\documentclass[12pt, twoside]{article}
\usepackage{amsfonts,mathrsfs,amsmath,amssymb,amsthm,upref,amscd}
\usepackage{xcolor}
\usepackage[all,cmtip]{xy}
\usepackage{setspace}
\usepackage{enumerate}
\usepackage[title,titletoc]{appendix}
\usepackage{times}
\usepackage{cite}
\usepackage{tikz}
\usepackage{pgfplots}
\usepackage{multirow}
\usepackage{array}
\usepackage{longtable}
\usepackage{needspace}
\usepackage{ulem}
\usepackage[colorinlistoftodos]{todonotes}
\usepackage{thmtools}
\usepackage{caption}
\usepackage{etoolbox}
\usetikzlibrary{arrows}
\numberwithin{equation}{section}

\makeatletter 
\@mparswitchfalse
\makeatother
\normalmarginpar %for right-handed notes and lines, or 
\def\titlerunning#1{\gdef\titrun{#1}}
\newcommand{\authoraddress}[2]{%
    \begin{minipage}[t]{0.48\linewidth}
        \footnotesize\raggedright
        {\scshape #1\par}
        \textit{Email address}: \texttt{#2}
    \end{minipage}}
\def\subjclass#1{{\renewcommand{\thefootnote}{}%
		\footnote{%\emph
			{Mathematics Subject Classification (2020):} #1}}}
\def\keywords#1{\par\medskip
	\noindent\textbf{Keywords.} #1}

\allowdisplaybreaks

\theoremstyle{plain}
\newtheorem{Thm}{Theorem}[section]
\newtheorem{Lem}[Thm]{Lemma}
\newtheorem{Cor}[Thm]{Corollary}
\newtheorem{Prop}[Thm]{Proposition}
\newtheorem*{Thm*}{Theorem}
\newtheorem*{claim*}{Claim}
\newtheorem*{Cor*}{Corollary}

\newtheorem*{Ques*}{Question}

\newtheorem*{Prob*}{Problem}
\newtheorem*{Fact*}{Fact}
\newtheorem{Fact}[Thm]{Fact}
\newtheorem*{OProb*}{Open Problem}
\theoremstyle{definition}
\newtheorem{Def}[Thm]{Definition}
\newtheorem*{Def*}{Definition}
\newtheorem{Rem}[Thm]{Remark}
\newtheorem{Ex}[Thm]{Example}

\AtEndEnvironment{Rem}{
	\unskip
	\hfill $\diamond$  
}

\newenvironment{revision}{\par\begingroup\emergencystretch=1em}{\par\endgroup}

\DeclareMathOperator{\tr}{tr}
\DeclareMathOperator{\real}{Re}

\DeclareMathOperator{\diver}{div}

\DeclareMathOperator{\sgn}{sgn}
\DeclareMathOperator{\vol}{vol}
\DeclareMathOperator{\dist}{dist}
\DeclareMathOperator{\supp}{supp}

\DeclareMathOperator{\Spin}{Spin}
\DeclareMathOperator{\SO}{SO}

\DeclareMathOperator{\diam}{diam}
\DeclareMathOperator{\spec}{Spec}

\newcommand{\equ}{equation}
\newcommand{\C}{\mathbb{C}}
\newcommand{\N}{\mathbb{N}}
\newcommand{\R}{\mathbb{R}}
\newcommand{\T}{\mathbb{T}}
\newcommand{\Z}{\mathbb{Z}}

\newcommand\ch{\mathcal{H}}

\newcommand\cj{\mathcal{J}}

\newcommand\cl{\mathcal{L}}
\newcommand\cm{\mathcal{M}}

\newcommand\co{\mathcal{O}}
\newcommand\cp{\mathcal{P}}

\newcommand\rr{\mathcal{R}}

\newcommand{\inp}[2]{\left\langle#1,#2\right\rangle}
\newcommand{\normm}[1]{{\left\vert\kern-0.25ex\left\vert\kern-0.25ex\left\vert #1 
		\right\vert\kern-0.25ex\right\vert\kern-0.25ex\right\vert}}

\def\mbs{\mathbb{S}}
\def\msa{\mathscr{A}}

\def\msd{\mathscr{D}}

\def\msf{\mathscr{F}}

\def\msi{\mathscr{I}}

\def\mfd{\mathfrak{D}}

\def\mft{\mathfrak{T}}
\def\mfs{\mathfrak{S}}
\def\mfj{\mathfrak{J}}
\def\mff{\mathfrak{F}}
\def\id{\text{Id}}
\def\ig{\textit{g}}
\def\ih{\textit{h}}

\def\pa {\partial}
\def\op{\oplus}
\def\ot{\otimes}

\def\De{\Delta}
\def\ka{\kappa}
\def\al{\alpha}
\def\bt{\beta}

\def\de{\delta}
\def\Ga{\Gamma}
\def\ga{\gamma}

\def\lm{\lambda}
\def\La{\Lambda}
\def\om{\omega}
\def\Om{\Omega}
\def\sa{\sigma}
\def\vr{\varepsilon}
\def\va{\varphi}

\begin{document}
	
	\titlerunning{Noncompactness and a sharp dimension threshold for the spinorial Yamabe equation}
	
	%\title{Examples of non-compactness for the spinorial Yamabe problem}
	
	\title{Noncompactness and a sharp dimension threshold for the spinorial Yamabe equation}
	
	\author{Juncheng Wei \qquad Tian Xu}%\thanks{Tian Xu is partially supported by the National Natural Science Foundation of China (NSFC 12371117) and the Natural Science Foundation of Zhejiang Province (LQK26A010001).}
	
	\date{}
	
	\maketitle
	
	\markboth{Juncheng Wei \qquad Tian Xu}{\titrun}
	
	
	\begin{abstract}
\begin{revision}
Motivated by the compactness conjecture for the classical Yamabe
equation, we study blow-up for the spinorial Yamabe equation on closed
spin manifolds. In every dimension $m\geq2$ we construct blow-up
sequences of nowhere-vanishing solutions whose maximum sets have any
prescribed codimension from $0$ to $m-1$.
Adapting the perturbation method of Isobe and Xu, we extend the
noncompactness construction of Isobe, Sire and Xu to smooth metrics.
For every $m\geq20$ we construct a smooth, non-locally conformally flat
metric on $S^m$ for which two single-bubble sequences concentrate at
the same point. Each solution in either sequence attains its maximum
length at a unique point, so its maximum set has codimension $m$.
Changing the unit spinor reverses the sign of the
leading energy correction, and the energies of the two sequences
approach the standard bubble energy from opposite sides.
For an admissible radial metric profile, differentiating the centered
reduced energy with respect to the bubble scale gives a quadratic form
in the radial coefficients. Under our hypotheses on the Weyl tensor
and the exterior terms, its definiteness rules out asymptotically
centered single-bubble concentration at the scale of the profile
through dimension $19$. We thus obtain a sharp threshold for this
class without invoking a Weyl vanishing theorem, a local Pohozaev
identity, or a positive mass theorem. Applied to the classical scalar
equations, the same scale calculation detects the known thresholds
of $25$ for the Yamabe and fourth-order $Q$-curvature equations and
$27$ for the sixth-order $Q$-curvature equation.
\end{revision}
{\keywords{Dirac operator, spinorial Yamabe equation, blow-up, dimension threshold.}}
\end{abstract}

\tableofcontents

\section{Introduction}

The Yamabe problem asks whether a closed Riemannian manifold of
dimension $m\geq3$ admits a metric of constant scalar curvature in a
prescribed conformal class. {The existence of a minimizing
metric follows from the work of Yamabe, Trudinger, Aubin and Schoen
(see \cite{Yamabe,Trudinger,Aubin,Schoen}).}
When the scalar curvature is fixed at a positive value, can a sequence of
positive conformal factors solving the Yamabe equation blow up on a fixed
manifold? The conformal transformations of the round sphere give such
sequences. In his 1988 Stanford lectures, Schoen conjectured that the full
set of positive solutions is compact unless the manifold is conformally
equivalent to the round sphere (see \cite[p.~144]{KMS}).
{For Schoen's proof in the locally conformally flat case,
see \cite{Schoen1991}. For general conformal classes the answer
depends on dimension. Khuri, Marques and Schoen proved compactness
through dimension $24$, assuming the positive mass theorem used in
their proof (see \cite{KMS}). Brendle constructed smooth metrics on
the sphere admitting blow-up in dimensions $m\geq52$ (see
\cite{Brendle}). For the extension to $25\leq m\leq51$, see his joint
work with Marques \cite{BrMa}.}

{Where does the dimension enter the compactness proof?
Near an isolated simple blow-up point, Khuri, Marques and Schoen use
a sharp approximation of the solution to expand a local Pohozaev
identity. The expansion yields a quadratic form in the Taylor
coefficients of the metric. Its coercivity in dimensions $6$ through
$24$ forces the Weyl tensor and its covariant derivatives up to a
dimension-dependent order to vanish at the blow-up point. This
vanishing gives sufficient control of the Green function expansion
to identify the mass in the limiting boundary term. The sign imposed
by blow-up is incompatible with the positive mass theorem, which
rules out concentration. We refer the reader to
\cite[Sections~6, 7 and~9]{KMS} for these steps. For the fourth- and
sixth-order constant $Q$-curvature
equations, Gong, Kim and Wei use the corresponding Pohozaev quadratic
forms and Green function masses to prove compactness through dimensions
$24$ and $26$, respectively, under the hypotheses of their theorems.
For the examples establishing sharpness, see the fourth-order
construction of Wei and Zhao \cite{WZ13} for $m\geq25$ and the
sixth-order construction of Gong, Kim and Wei \cite{GKW25} for $m\geq27$.}

On a closed Riemannian spin manifold $(M,\ig,\sa)$ of dimension $m\geq2$, where $\sa$ is a fixed spin structure,
let $D_{\ig}$ denote the Dirac operator for the metric $\ig$ and spin
structure $\sa$. We consider the spinorial Yamabe equation
\begin{\equ}\label{SY0}
 D_{\ig}\psi=|\psi|_{\ig}^{2/(m-1)}\psi.
\end{\equ}
Here $\psi$ is a section of the complex spinor bundle, and we call it
\textit{nowhere-vanishing} if its pointwise length satisfies
$|\psi(x)|_{\ig}>0$ for every $x\in M$. {Equation~\eqref{SY0}
is conformally covariant and admits a family of standard bubbles on
Euclidean space. For the scalar Yamabe and $Q$-curvature equations,
the center and scale determine the standard positive bubble. A
spinorial bubble also depends on a unit spinor. Varying this spinor
at fixed center and scale gives distinct bubbles with the same
pointwise length.}

A second feature of \eqref{SY0} is particularly clear on a flat torus. There one can
find a sequence of solutions that oscillate with increasing frequency,
although each has constant length. The maximum is attained everywhere
on the torus, while both the length and the energy tend to infinity.
To describe such blow-up, we record both the set where each solution
attains its maximum length and the energy it carries.

\begin{Def}\label{def:maximum-locus}
For a sequence of solutions $\{\psi_j\}$ to \eqref{SY0} on a fixed closed
spin manifold, let
\[
 \co_j=\{x\in M:|\psi_j(x)|_{\ig}=\max_M|\psi_j|_{\ig}\}.
\]
We call $\{\psi_j\}$ a \textit{blow-up sequence} if
$\max_M|\psi_j|_{\ig}\to\infty$ as $j\to\infty$.
When $\co_j$ is a submanifold, its codimension $\operatorname{codim}_M\co_j$ is taken in $M$.
\end{Def}
We record the codimension of $\co_j$ because it counts the independent directions normal to the maximum set in $M$. It is zero when the maximum is attained everywhere and equals $m$ when it is attained at a single point. In our examples, the maximum set may have several connected components, all of the same dimension. {Their number can vary with $j$ while their common codimension stays fixed. For instance, on the product manifold $S^1\times N$, the maximum set can consist of $2j$ disjoint copies of $N$, each of codimension one.}

On a fixed closed spin manifold $M$, elliptic estimates for \eqref{SY0}
give a subsequence converging in $C^1$ whenever
$\sup_j\max_M|\psi_j|_{\ig}<\infty$. If the maximum lengths have no
uniform bound, one can extract a subsequence along which they tend to
infinity. Thus a family of solutions is relatively compact in $C^1$ if
and only if the maximum lengths of its members are uniformly bounded.
{In this sense, noncompactness is equivalent to the existence
of a blow-up sequence in the family.} We also record the energy of the
solutions along each such sequence. At a solution $\psi$, the energy
equals $(2m)^{-1}\int_M|\psi|_{\ig}^{2m/(m-1)}\,d\vol_{\ig}$.

%\Needspace{18\baselineskip}
\begin{Thm}\label{main-thm}
For every $m\geq2$, the set of solutions to the spinorial Yamabe equation \eqref{SY0} on a closed spin $m$-manifold can fail to be compact. More precisely, the following three statements are valid:
\begin{itemize}
\item[$(1)$] For each $k\in\{0,\ldots,m-1\}$, there is a closed spin
$m$-manifold $(M,\ig,\sa)$ admitting a {blow-up sequence} $\{\psi_j\}$ of smooth, nowhere-vanishing solutions with
$\operatorname{codim}_M\co_j=k$ for every $j$.

\item[$(2)$] For every $m\geq20$, there is a smooth, non-locally conformally
flat metric $\ig$ on $S^m$ for which a {blow-up sequence} $\{\psi_j\}$ of smooth, nowhere-vanishing solutions
concentrates at one point as a single standard bubble. Each
$\co_j$ is a singleton, hence $\operatorname{codim}_{S^m}\co_j=m$, and
\[
 \frac1{2m}\int_{S^m}|\psi_j|_{\ig}^{2m/(m-1)}\,d\vol_{\ig}
 \to \frac1{2m}\Bigl(\frac m2\Bigr)^m
       \operatorname{Vol}(S^m,\ig_{S^m}), \quad \text{as } j\to\infty
\]
where $\ig_{S^m}$ is the unit round metric.

\item[$(3)$] The lower bound $m\geq20$ in \textup{(2)} is sharp in the sense
that asymptotically centered radial single-bubble concentration of this
type cannot occur in dimensions $m\leq19$.
\end{itemize}
\end{Thm}

The following Table~\ref{tab:noncompactness-examples} illustrates the blow-up sequences in Theorem~\ref{main-thm} with specific examples. In each row the metric and the spin structure are fixed. The
standard bubble energy is the limiting value displayed in Theorem~\ref{main-thm}~(2).
{In the last two rows the Weyl tensor of the metric at each
perturbation center is a nonzero multiple of the chosen algebraic tensor
$W$ in the local frame. In both rows we choose $W$ so that the positive
local maximum of the leading energy correction persists when the bubble
center and the unit spinor vary together. In the last row its negative
local minimum also persists.}

\begingroup
\small
\setlength{\tabcolsep}{4pt}
\renewcommand{\arraystretch}{1.5}
\setlength{\LTcapwidth}{\linewidth}
\captionsetup{labelfont={bf,color=black},textfont={bf,color=black}}
\color{black}
\newcommand{\exampletablehead}{%
  \textbf{Manifold and metric} & \textbf{Choice of data}
  & \textbf{Sequences and maximum sets $\co_j$} & \textbf{Energy of the solutions}}
\begin{longtable}{|>{\color{black}\centering\arraybackslash}m{3.30cm}|
                 >{\color{black}\centering\arraybackslash}m{4.10cm}|
                 >{\color{black}\centering\arraybackslash}m{3.90cm}|
                 >{\color{black}\centering\arraybackslash}m{3.45cm}|}
\caption{Some examples of noncompactness in Theorem~\ref{main-thm}}\label{tab:noncompactness-examples}\\
\hline
\exampletablehead\\
\hline
\endfirsthead
\multicolumn{4}{c}{\color{black}\bfseries Table~\thetable\ (continued)}\\
\hline
\exampletablehead\\
\hline
\endhead
$\T^m$, $m\geq2$, with a metric as in Corollary~\ref{cor:prescribed-maximum-set},
which can be chosen non-locally conformally flat when $m\geq3$.
  & Prescribe $k\in\{0,\ldots,m-1\}$.
  & One blow-up sequence~$\{\psi_j\}$ for each prescribed $k$. Every $\co_j$ equals the same coordinate subtorus of codimension $k$ (all of $\T^m$ if $k=0$).
  & Tends to infinity.\\
\hline
$S^1\times N$ with the product metric.
  & $N$ is odd-dimensional and carries a unit-length eigenspinor with nonzero eigenvalue of $D_{\ig_N}$. Examples include round spheres $S^{2\ell+1}$ and flat tori $\T^{2\ell+1}$, $\ell\in\N$.
  & Two blow-up sequences. For one, $\co_j=S^1\times N$. For the other, $\co_j$ consists of $2j$ disjoint fibers $\{t\}\times N$ of codimension one.
  & Tends to infinity along both sequences.\\
\hline
$S^k\times\T^{m-k}$ or $S^k\times S^3$, with a fixed conformal change of the product metric.
  & $k\geq2$. Use a unit parallel spinor on the flat torus ($m-k\geq1$), or a unit harmonic spinor on a Berger sphere $S^3$.
  & One blow-up sequence $\{\psi_j\}$ on each product. For every $j$, $\co_j$ consists of the two fibers over the poles of $S^k$, each of codimension $k$.
  & Tends to infinity.\\
\hline
$S^m$, $m\geq20$, with a smooth non-locally conformally flat metric.
  & {At each center, use the Weyl tensor of Example~\ref{rev-weyl-example} or a fixed tensor sufficiently close~to~it.}
  & One blow-up sequence~$\{\psi_j\}$ with $\co_j=\{x_j\}$, where $x_j$ converges to a concentration point.
  & Approaches the standard bubble energy from above.\\
\hline
$S^m$, $m\geq20$, with a smooth non-locally conformally flat metric.
  & {At each center, use the Weyl tensor of Example~\ref{Ex-2} or a fixed tensor sufficiently close~to~it.}
  & {Two blow-up sequences $\{\psi_j^+\}$ and $\{\psi_j^-\}$} on the same metric, with maximum sets $\co_j^+=\{x_j^+\}$ and $\co_j^-=\{x_j^-\}$. Both $x_j^+$ and $x_j^-$ converge to the same point.
  & Approaches the standard bubble energy from above along $\{\psi_j^+\}$ and from below along $\{\psi_j^-\}$.\\
\hline
\end{longtable}
\endgroup

In the first three rows, the minimum length tends to
infinity along each sequence. To see how the maximum set can be
prescribed in the first row, start with the flat torus and its trivial
spin structure, the choice we make on all torus factors in the table.
{Friedrich's description of the flat Dirac spectrum gives
unit-length eigenspinors $\varphi_j$ with positive eigenvalues
$\lambda_j\to\infty$ (see \cite{Friedrich84} or
Theorem~\ref{thm: Friedrich's result on torus} below).} Then $\lambda_j^{(m-1)/2}\varphi_j$ solves
\eqref{SY0} and has constant length $\lambda_j^{(m-1)/2}$.
More generally, whenever a fixed metric $\ig$ admits a blow-up sequence
with $|\psi_j|_\ig\equiv c_j$, a single conformal change can prescribe
the maximum set of every solution in the sequence. Given a nonempty
closed embedded submanifold $Z\subset M$, choose a smooth function
$u\geq0$ that vanishes precisely on $Z$, for instance by taking the
squared distance to $Z$ near $Z$ and extending it smoothly with positive
values away from $Z$. For the metric $e^{2u}\ig$, the transformed
solutions have lengths $c_j e^{-(m-1)u/2}$ and hence attain their
maxima precisely on $Z$.

In the second row, the same product metric on $S^1\times N$ supports
two blow-up sequences with different maximum sets. Take the trivial
spin structure on $S^1$, so that spinors are periodic in the circle
variable, and fix a unit eigenspinor on $N$ with nonzero eigenvalue.
Using Fourier modes on the circle gives the first sequence. Each
solution has constant length and attains its maximum on all of
$S^1\times N$. To obtain the second sequence, we keep the same
eigenspinor on $N$ and use the periodic orbits of the reduced system \eqref{eqs2} in Section \ref{sec:ansatz}.
Along these orbits the length
varies along the circle, and the $j$-th solution attains its maximum
on exactly $2j$ fibers. Thus the maximum set changes with $j$ even
though the metric stays fixed.
Both sequences occur on the flat product torus
$\T^m=S^1\times\T^{m-1}$ for even $m$. They also occur on
$S^1\times S^{2\ell+1}$ with the standard product metric for every
$\ell\geq1$, since the round sphere carries a unit Killing spinor
with nonzero Dirac eigenvalue.

In the third row all the solutions attain their maximum
on the same two fibers over the poles of $S^k$. To realize these
fibers as copies of $S^3$, we equip the second factor of
$S^k\times S^3$ with the Berger metric of parameter $4$ in
{Hitchin's normalization (see \cite{Hit74}).}
For either product in this row, the ratio of the lengths of any two
distinct members of the sequence varies over $S^k$, as we show in
Section~\ref{sec: product}. A common conformal change preserves this
ratio, so at most one of the transformed solutions can have constant
length. This distinguishes the third row from the first, where a
single conformal change recovers constant length for every member of
the sequence.

For the scalar Yamabe equation on a closed locally conformally flat
manifold of dimension at least $3$, {Schoen proved that the set
of positive solutions at a fixed positive scalar curvature is compact
unless the manifold is conformally equivalent to the round sphere
(see \cite{Schoen1991}).}
For the spinorial Yamabe equation, every codimension in
Theorem~\ref{main-thm}~(1) already occurs in the conformal class of the
flat torus $\T^m$, for every $m\geq2$. Indeed, a fixed conformal change
of the flat metric makes any prescribed coordinate subtorus the maximum
set of every spinor in the transformed sequence. When $m\geq3$,
Corollary~\ref{cor:prescribed-maximum-set} realizes the same range of
codimensions on $\T^m$ with non-locally conformally flat metrics as well.
The third row of Table~\ref{tab:noncompactness-examples} gives
another such example on $S^k\times\T^{m-k}$ when $k\geq2$ and $m-k\geq2$.
The product of the round and flat metrics has nonzero Weyl curvature on
every two-plane tangent to the torus factor, so its conformal class
is non-locally conformally flat. Thus, even away from the conformal
class of the round sphere, examples of noncompactness for the spinorial
Yamabe equation occur for both locally conformally flat and non-locally
conformally flat metrics.

A seven-dimensional factor for the second row is provided by any
closed simply connected $3$-Sasakian manifold $N$ with its standard
$3$-Sasakian metric. {Its canonical unit spinor has Dirac
eigenvalue $9/2$ in the normalization of Agricola and Friedrich
(see \cite[Theorem~4.1]{AF10}).} The circle construction therefore applies
on the eight-dimensional product $S^1\times N$.

For the construction in the third row with sphere factor $S^2$, it
suffices that $N$ carry a unit harmonic spinor.
{Agricola, Chiossi, Friedrich and H\"oll give a closed
six-dimensional example (see \cite[Example~3.8]{ACFH}). Their criterion
in \cite[Theorem~4.6]{ACFH} also applies, for instance, to any closed
seven-dimensional}
manifold $N$ with a torsion-free $G_2$-structure, whose defining unit
spinor is parallel. In either dimension this spinor reduces
\eqref{SY0} on $S^2\times N$ to an equation on $S^2$.
Our Theorem~\ref{main thm ODE2} supplies a sequence of radial solutions.
After the fixed conformal change described in Section~\ref{sec: product},
their maximum sets consist of the two copies of $N$ over the poles,
even though their lengths tend to infinity everywhere. {For
an example of non-isolated blow-up for a scalar curvature equation,
see the construction of Robert and V\'etois in \cite{RV14}, where
the potential varies.}

Of particular interest are the asymptotically centered radial
single-bubble sequences in the last two rows of the table.
For such sequences, Theorem~\ref{main-thm}~(3) gives the sharp lower
bound $m\geq20$.
{Isobe and Xu obtained multiple solutions by perturbing the
round metric (see \cite[Theorem~1.4]{Isobe-Xu23}). Using infinitely
many localized perturbations, Isobe, Sire and Xu constructed a metric
of finite regularity for which a blow-up sequence approaches the bubble
energy from below (see \cite[Theorem~1.2]{ISX24}). They asked in
\cite[Remark~1.3]{ISX24} whether a smooth metric could exhibit
noncompactness.}
Subsection~\ref{rev-gluing-subsection} gives such a metric for every
$m\geq20$.

{The smooth metric can also be chosen so that single-bubble
concentration occurs on both sides of the standard bubble energy.
At fixed center and scale the leading energy correction contains a
Clifford four-form term, quadratic in the Weyl tensor used to construct
the metric. For the tensor in Example~\ref{Ex-2}, the correction is
positive for one choice of unit spinor and negative for another.
Theorem~\ref{rev-noncompact-sphere} realizes both
signs on the same metric. The two sequences concentrate at the same
point, and their energies approach the value in
Theorem~\ref{main-thm}~(2) from above along one sequence and from
below along the other.}

{The sequence with energy below the standard bubble energy
has a consequence for the B\"ar-Hijazi-Lott invariant.}
Writing $\lambda_1^+(\widetilde\ig)$ for the first positive Dirac
eigenvalue of a metric $\widetilde\ig$ with spin structure $\sa$, set
\[
 \lambda_{\min}^+(M,[\ig],\sa)
 =\inf_{\widetilde\ig\in[\ig]}
   \lambda_1^+(\widetilde\ig)\operatorname{Vol}(M,\widetilde\ig)^{1/m}.
\]
Its value on the round sphere is
$\frac m2\operatorname{Vol}(S^m,\ig_{S^m})^{1/m}$. By the spinorial
Aubin inequality, this is an upper bound in every spin conformal
class {(see \cite{AGHM}).} Whether equality characterizes the round sphere up
to spin-conformal equivalence remains open in general.
{Applying the variational characterization of this invariant
(see, for instance, \cite[Remark~1.6]{Isobe-Xu23}) to any nonzero
solution $\psi$ of \eqref{SY0} gives}
\begin{\equ}\label{BHL-bound}
 \lambda_{\min}^+(M,[\ig],\sa)
 \leq\left(\int_M|\psi|_{\ig}^{2m/(m-1)}\,d\vol_{\ig}\right)^{1/m}.
\end{\equ}
For the metric constructed from Example~\ref{Ex-2}, a single solution
on the branch below the bubble energy therefore proves the strict
B\"ar-Hijazi-Lott inequality for its conformal class.

In the scalar Yamabe problem, Aubin showed that in dimension $m\geq6$
a test function centered at a point of nonzero Weyl curvature has
Yamabe quotient strictly below the spherical value {(see
\cite{Aubin} for the test function estimate).}
{In \cite[Theorem~1.1]{SX2021}, the second author claimed the corresponding
strict B\"ar-Hijazi-Lott inequality for every closed spin manifold of
dimension $m\geq6$ whose metric is not locally conformally flat.}
We acknowledge a gap in that proof.
The unit spinor chosen in Lemma~4.4 of that preprint cancels the
Clifford four-form term. And Lemma~4.6 of that preprint then requires the remaining metric contribution to lower the energy at fourth order, but the
calculation does not justify this sign.
%\footnote{In the normalization $\tau=0$ allowed in \cite{SX2021}, the metric matrix $G=(\ig_{ij})$ satisfies
%$\det G=1+O(|x|^5)$ in normal coordinates. Writing $h=\ln G$ gives
%$h=O(|x|^2)$ and $\operatorname{tr}h=O(|x|^5)$. For the standard
%bubble, contraction with the symmetric matrix $G^{-1/2}-I$ in the
%term denoted by $J_6$ leaves its trace times a positive radial factor.
%Since $\operatorname{tr}(G^{-1/2}-I)
%=\tfrac18\operatorname{tr}(h^2)+O(|x|^5)$, its degree-four term isnonnegative.} 
Here in this paper, our Example~\ref{rev-weyl-example} shows what can happen after this
cancellation. Its Weyl tensor is nonzero, and the Clifford four-form
term vanishes for every unit spinor. At the scale of our construction,
the remaining term raises the energy of every centered bubble. The
resulting solutions approach the standard bubble energy from above,
so the upper bounds they give in \eqref{BHL-bound} exceed the
spherical value. A different test spinor might still establish the
strict inequality for this metric, although we have not yet constructed
such a spinor.

{What changes when the dimension increases from $19$ to $20$?
In Theorem~\ref{main-thm}~(3), let $b_j$ be the length scale of the
$j$-th local metric profile and $\lambda_j$ the scale of its bubble.
The profile is determined by an algebraic Weyl tensor $W$ and a
polynomial in the squared distance from its center. We take this center
as the origin and measure distances in units of $b_j$. In these
coordinates the bubble center tends to the origin and the relative
bubble scale $\ell_j=\lambda_j/b_j$ tends to a finite limit
$\ell_*>0$. Under the operator and exterior estimates in
Theorem~\ref{centered-radial-bump-threshold}, testing the equation
along variations of the scale and unit spinor shows that the leading
correction to the centered bubble energy must be stationary at the
limiting scale and spinor.
For a unit spinor $\ga$, this correction is the product of a positive
function of the relative scale and the quadratic form
$\real(\msa_W\ga,\ga)$. Here $\msa_W$ is the Hermitian endomorphism
associated to $W$ in
Subsection~\ref{rev-cubic-subsection}. The spinor variation selects an
eigenvalue of $\msa_W$. Since the theorem assumes $\ker\msa_W=\{0\}$,
the scale-dependent factor must itself be stationary at $\ell_*$.}

{When the relative scale $\ell$ increases, radii smaller
and larger than $\ell$ contribute with opposite signs to the derivative
of this factor. Through dimension $19$, the outer contribution is
strictly larger in magnitude for every admissible radial polynomial
and every $\ell>0$. This excludes concentration of the type in
Theorem~\ref{main-thm}~(3). In dimension $20$ the cubic used in our
construction makes the two contributions balance at a strict maximum
of the radial factor. At that same scale the leading energy correction
has a maximum for a positive eigenvalue of $\msa_W$ and a minimum for
a negative one. To obtain solutions, we verify that for our choices of $W$
the extrema persist when the center and spinor vary and the cut-off
is restored. Theorem~\ref{rev-noncompact-sphere} then produces actual
single-bubble sequences in every dimension $m\geq20$. This proves sharpness within
the stated class. Single-bubble concentration on locally conformally
flat spin manifolds outside the spin conformal class of the round
sphere poses a separate question.}

{The spinorial result suggests a more general question for
conformally invariant equations that admit a variational formulation
and a standard bubble. ``If we fix the center of a standard bubble and vary
its size, in what dimension can the leading energy correction first
have a critical scale for some nonzero admissible radial metric
profile?'' In Subsection~\ref{subsec:centered-scalar-thresholds} we
examine this question for the scalar Yamabe equation and the fourth-
and sixth-order constant $Q$-curvature equations. In each case we
start from the dependence of the energy on the metric. Differentiating
the centered reduced energy with respect to the bubble scale gives a
quadratic form in the coefficients of the radial polynomial, without
using a local Pohozaev identity. This calculation can be carried out
before constructing any blow-up sequence or proving compactness in
lower dimensions. For the Yamabe equation, the centered energy formula
of Brendle and Marques again gives contributions of opposite signs from
inside and outside the bubble scale. For the centered energy and its
dependence on the scale, see \cite[Propositions~17 and~18]{BrMa}.
For the sixth-order constant $Q$-curvature equation, Gong, Kim and Wei
compute the matrix of the scale derivative in their noncompactness
construction (see \cite[Corollary~11.19]{GKW25}). We refer the reader
to \cite[Proposition~5.3]{GKW26} for their formula at general even
order. On the chosen radial family the scale
derivative agrees, up to sign and normalization, with the restriction
of the Pohozaev quadratic form. What surprised us was that in all three
scalar problems this family already contains a negative direction in
the first dimension where the full form loses positive definiteness.
For the Yamabe case, see \cite[Proposition~A.8]{KMS}. The comparison
across all three equations is given in \cite[Corollary~4.7]{GKWambient}.
This explains why the centered calculation detects dimension $25$ for
the Yamabe and fourth-order equations and $27$ for the sixth-order
equation. The compactness theorems and blow-up constructions recalled
above establish that these are the sharp thresholds for actual
solutions under the hypotheses of the respective theorems.}

The paper is organized as follows. Section~\ref{sec:preliminaries}
recalls the properties of the Dirac operator used below, including its
formula on a Riemannian product. In Section~\ref{sec:ansatz} we reduce
the spinorial Yamabe equation on products to ordinary differential
systems and write the equation on the sphere in stereographic
coordinates. In Section~\ref{sec: product} we construct blow-up sequences
on product manifolds and describe their maximum sets. In
Section~\ref{sec: perturbation} we construct smooth metrics on
$S^m$ admitting single-bubble blow-up for every $m\geq20$. In
Section~\ref{sec:radial-dimension-threshold} we use the scale equation
to prove that this dimension bound is sharp for asymptotically centered
radial single-bubble concentration and compare the calculation with its
counterparts for the scalar Yamabe and $Q$-curvature equations.
{Appendix~\ref{app:cubic-calculations} gives the polynomial
estimates and the derivative calculations used to locate the extrema in
Subsection~\ref{rev-cubic-subsection}.}
%\filbreak

% \todo{Can one fixed product metric carry both a sequence with uniformly diverging length and one concentrating at the standard bubble energy?}

\section{Preliminaries}\label{sec:preliminaries}
	
	\subsection{Notations and conventions}

Throughout the paper, $(M,\ig,\sa)$ denotes a closed, oriented Riemannian spin
manifold of dimension $m$, and $\mbs(M)$ its complex spinor bundle. We denote
the fiberwise Hermitian product and its norm by $(\cdot,\cdot)_\ig$ and
$|\cdot|_\ig$. We write $\Ga(\mbs(M))$ for the smooth sections, also denoted
by $C^\infty(M,\mbs(M))$. In the variational arguments we use the real
$L^2$ inner product
\[
 (\psi,\va)_2=\real\int_M(\psi,\va)_\ig\,d\vol_\ig.
\]
For $1\leq p<\infty$, the space $L^p(M,\mbs(M))$ consists of measurable
spinors endowed with the norm $\|\psi\|_p=(\int_M|\psi|_\ig^p\,d\vol_\ig)^{1/p}$. We denote the essential supremum of $|\psi|_\ig$ by $\|\psi\|_\infty$. For continuous spinors on $M$ it equals
$\max_M|\psi|_\ig$.

{Throughout the paper we use the Einstein summation convention
for an index that appears once as a superscript and once as a subscript, and we
write the summation sign explicitly when both occurrences are subscripts
or when this improves readability.}

\subsection{A brief review about the spinor bundle and the Dirac operator}\label{subsec:spinor bundle}
	
Let ${P_{\SO}(M)}$ be the principal bundle of positively oriented orthonormal
frames. A spin structure {$\sa$} is a principal $\Spin(m)$-bundle $P_{\Spin}(M)$
together with a map $\vartheta:P_{\Spin}(M)\to {P_{\SO}(M)}$ over the identity
on $M$ for which the diagram
\begin{displaymath}
		\xymatrix@R=0.3cm{
			{P_{\Spin}(M)}\times {\Spin(m)}  \ar[r] \ar[dd]^{\displaystyle\vartheta\times \Theta}&  {P_{\Spin}(M)}\ar[dd]^{\displaystyle\vartheta} \ar[dr] \\
			&&  M \\
			{P_{\SO}(M)}\times {\SO(m)} \ar[r] & {P_{\SO}(M)}  \ar[ur] }
	\end{displaymath}
commutes, where $\Theta:\Spin(m)\to {\SO(m)}$ is the spin double covering.
Such a structure exists precisely when the second Stiefel-Whitney class
$w_2(TM)\in H^2(M;\Z_2)$ vanishes, see for instance \cite{Friedrich00,LM} for more geometric materials.

The real Clifford algebra $Cl(\R^m)$ is generated by $\R^m$ with the relation
$x\cdot y+y\cdot x=-2(x,y)$, where $(\cdot,\cdot)$ is the Euclidean inner
product. Fix an irreducible module $\mbs_m$ for its complexification. The complex
dimension of $\mbs_m$ is $2^{\lfloor m/2\rfloor}$. For an oriented orthonormal basis
$(e_1,\dots,e_m)$, the complex volume element is
$\om_\C^{\R^m}=i^{\lfloor(m+1)/2\rfloor}e_1\cdots e_m$.
In even dimension the complex Clifford module is unique up to equivalence.
It splits under $\Spin(m)$ into the two eigenspaces $\mbs_m^\pm$ of
$\om_\C^{\R^m}$ for the eigenvalues $\pm1$. In odd dimension the two
inequivalent irreducible Clifford modules have equivalent restrictions to
$\Spin(m)$. They may be realized on the same vector space by reversing the
sign of Clifford multiplication by vectors.

Restricting Clifford multiplication on $\mbs_m$ to $\Spin(m)$ gives the
complex spin representation $\varrho:\Spin(m)\to\operatorname{GL}(\mbs_m)$.
Using this representation, we form the spinor bundle from the principal
bundle $P_{\Spin}(M)$ of the chosen spin structure $\sa$ by setting
\[
\mbs(M):=P_{\Spin}(M)\times_\varrho\mbs_m.
\]
At $x\in M$, choose an oriented orthonormal frame $\mathbf e=(e_1,\ldots,e_m)$
and a lift $\widetilde{\mathbf e}\in P_{\Spin}(M)$ of this frame.
For $a\in\Spin(m)$ and $v\in\mbs_m$, the pairs $(\widetilde{\mathbf e}a,v)$
and $(\widetilde{\mathbf e},\varrho(a)v)$ represent the same spinor in $\mbs_x(M)$.

The spinor bundle carries a Hermitian metric, a Clifford action and
a spin connection induced by the Levi-Civita connection.
We recall their compatibility relations:
\begin{itemize}
		\item[$(i)$] for any $x\in M$, $X,Y\in T_xM$ and $\psi\in\mbs_x(M)$
		\[
		X\cdot Y\cdot \psi + Y\cdot X\cdot\psi + 2\ig(X,Y)\psi = 0;
		\]
		\item[$(ii)$] for any $X\in T_xM$ and $\psi_1,\psi_2\in\mbs_x(M)$,
		\[
		(X\cdot \psi_1,\psi_2)_{\ig}=-(\psi_1,X\cdot\psi_2)_{\ig};
		\]
		%where $(\cdot,\cdot)_{\ig}$ is the hermitian metric on $\mbs(M)$ induced from $\ig$;
		\item[$(iii)$] for any $X,Y\in\Ga(TM)$ and $\psi\in\Ga(\mbs(M))$,
		\[
		\nabla_X^\mbs(Y\cdot\psi)=(\nabla_XY)\cdot\psi+Y\cdot\nabla_X^\mbs\psi,
		\]
		where $\nabla^\mbs$ stands for the metric connection on $\mbs(M)$ (if necessary, we will denote $\nabla^{\ig}$ to emphasize the background metric $\ig$).
	\end{itemize}
	The Dirac operator is then defined on the spinor bundle $\mbs(M)$ as the composition
	\begin{displaymath}
		\xymatrixcolsep{1.6pc}\xymatrix{
			D_\ig:\Ga(\mbs(M)) \ar[r]^-{\nabla^\mbs}
			& \Ga(T^*M\otimes \mbs(M)) \ar[r]
			& \Ga(TM\otimes \mbs(M)) \ar[r]^-{\mathfrak{m}}
			& \Ga(\mbs(M))}
	\end{displaymath}
where $\mathfrak m(X\otimes\psi)=X\cdot_\ig\psi$ denotes the Clifford multiplication and the middle arrow is the metric identification of covectors with vectors. In a local
orthonormal frame $\{e_1,\dots,e_m\}$ this gives
\[
D_\ig=\sum_{k=1}^m e_k\cdot_\ig\nabla^\mbs_{e_k}
	\]
as in the spin Dirac operator of Atiyah and Singer \cite{AS}.

Clifford multiplication depends on the metric. For $\tilde\ig=e^{2u}\ig$,
the map $X\mapsto e^{-u}X$ identifies the oriented orthonormal frame bundles.
Its spin lift induces a fiberwise isometry $\psi\mapsto\tilde\psi$ carrying
$X\cdot_\ig\psi$ to $e^{-u}X\cdot_{\tilde\ig}\tilde\psi$. In the sequel, we shall write
$D_\ig^M$ and $\cdot_\ig$ whenever we need to specify the manifold or metric.

\begin{Prop}[Ref. \cite{Hij86, Hit74}]\label{prop:conformal-invar}
Let $\ig=e^{2u}\ig_0$ be a conformal change of metric on an
$m$-dimensional spin manifold. The induced fiberwise isometry
$\iota:\mbs(M,\ig_0)\to\mbs(M,\ig)$ satisfies
\[
	e^{\frac{m+1}2 u}	D_\ig\big( \iota(\psi) \big)= \iota\big(D_{\ig_0}
		\big( e^{\frac{m-1}2 u}\psi \big)\big).
		\]
	\end{Prop}
\noindent
We point out particularly that if $\psi$ solves \eqref{SY0} for $\ig_0$ then
$e^{-(m-1)u/2}\iota(\psi)$ solves it for $e^{2u}\ig_0$.

\subsection{The spectrum of the Dirac operator}

On a closed spin manifold the self-adjoint realization of $D_\ig$ in
$L^2(M,\mbs(M))$ has domain $H^1(M,\mbs(M))$ and compact resolvent.
Its spectrum, denoted by $\spec(D_\ig)$, consists of real eigenvalues of
finite multiplicity, unbounded in both directions on the real line $\R$. The corresponding
eigenspaces admit a complete orthonormal eigenbasis \cite{Friedrich00,Ginoux,LM}.
The product constructions below use the following assumptions on the
Dirac operator $D_{\ig_N}^N$ of a closed spin factor $(N,\ig_N,\sa_N)$.
In each assumption, unit length is understood pointwise.
\begin{itemize}
 \item[(H1)] There are unit-length eigenspinors with eigenvalues
 $\lambda_j$ satisfying $|\lambda_j|\to\infty$, as $j\to\infty$.
 \item[(H2)] There is a unit-length harmonic spinor.
 \item[(H3)] There is a unit-length eigenspinor with nonzero eigenvalue.
\end{itemize}
Clearly, (H1) implies (H3). To construct constant-length
solutions on $N_1\times N_2$, we assume (H1) on $N_1$ and (H2) on $N_2$.
In this construction, either $N_1$ is even-dimensional or both factors
are odd-dimensional. The lengths of the solutions grow as $|\lambda_j|\to\infty$.
For the construction on $S^1\times N$ with $N$ odd-dimensional, (H3)
suffices on $N$. We keep one nonzero eigenvalue and its unit-length
eigenspinor fixed on $N$ and vary the least period of the solution on
the circle.

Round spheres already give examples of (H3). On the unit round
$n$-sphere, $n\geq2$, Killing spinors have constant length and Dirac
eigenvalues $\pm n/2$. Thus in the circle product $S^1\times N$ one
may take $N=S^{2\ell+1}$ for any $\ell\geq1$. The case $N=S^1$ also
satisfies (H3) for either spin structure, since the Fourier modes of
nonzero frequency have constant length.

The flat torus provides a standard model in which (H1) holds
for every spin structure and (H2) holds precisely for the trivial one.
Let $\Ga$ be a lattice in $\R^n$ with basis $\ga_1,\dots,\ga_n$. The
spin structures on the flat torus $\T^n=\R^n/\Ga$ correspond to characters
$\Ga\to\{\pm1\}$. Write $(-1)^{\de_j}$ for the value of the character on
$\ga_j$, with $\de_j\in\{0,1\}$. The following description of the
spectrum is due to Friedrich \cite{Friedrich84}.

\begin{Thm}\label{thm: Friedrich's result on torus}
For the spin structure determined by $(\de_1,\dots,\de_n)$, the spectrum
of the flat Dirac operator is
\[
	\spec(D_{\ig_\text{flat}}^{\T^n})=\Big\{ \pm2\pi\Big| \ga^*+\frac12\sum_{j=1}^n\de_j\ga_j^* \Big|:\, \ga^*\in \Ga^* \Big\},
	\]
	where $\Ga^*$ is the dual lattice and $\{\ga_1^*,\dots,\ga_n^*\}$ is the basis of $\Ga^*$ dual to $\{\ga_1,\dots,\ga_n\}$. Each eigenspace admits a basis of spinors of constant pointwise
length. The kernel consists of the parallel spinors for the trivial spin
structure and has complex dimension $2^{\lfloor n/2\rfloor}$ in that case.
For every other spin structure the kernel is trivial.
\end{Thm}

To see how the eigenspinors arise, fix a frequency $\nu$ in the shifted
dual lattice of Theorem~\ref{thm: Friedrich's result on torus} and a
constant spinor $\eta\in\mbs_n$. In the parallel trivialization on
$\R^n$, the function $x\mapsto e^{2\pi i(\nu,x)}\eta$ represents a
spinor on $\T^n$, since translation by $\ga_j$ multiplies it by the
prescribed sign $(-1)^{\de_j}$. In this trivialization the flat torus
Dirac operator acts by
\[
 D_{\ig_\text{flat}}^{\T^n}\bigl(e^{2\pi i(\nu,x)}\eta\bigr)
 =\sum_{j=1}^n e_j\cdot_{\ig_{\R^n}}
     \partial_{x^j}\bigl(e^{2\pi i(\nu,x)}\eta\bigr)
 =e^{2\pi i(\nu,x)}
     \bigl(2\pi i\,\nu\cdot_{\ig_{\R^n}}\eta\bigr).
\]
Clifford multiplication by $\nu$ is skew-adjoint and squares to
$-|\nu|^2\operatorname{Id}$. Thus $2\pi i\,\nu\cdot_{\ig_{\R^n}}$
is Hermitian, and each of its eigenvalues has absolute value
$2\pi|\nu|$. Choosing $\eta$ to be a unit eigenvector gives an
eigenspinor of pointwise length one, since $|e^{2\pi i(\nu,x)}|=1$.
Choosing an eigenbasis at each frequency gives the asserted basis
by Fourier expansion. Distinct frequencies are orthogonal, so their
multiplicities add whenever they yield the same eigenvalue. For every spin structure the allowed frequencies have arbitrarily
large norm, which justifies (H1). For the trivial spin structure one
may also take $\nu=0$. Any unit constant spinor then descends to a
unit-length parallel spinor on $\T^n$, justifying (H2).

Agricola, Chiossi, Friedrich and H\"oll also give geometric criteria for
(H2) in dimensions six and seven. To state their seven-dimensional
criterion, let $\psi$ be a unit real spinor on a Riemannian spin
seven-manifold $(N,\ig_N)$. It defines a $G_2$-structure whose
intrinsic torsion is described by the endomorphism $S$ of $TN$
determined by
\[
 \nabla_X^{\mbs}\psi=S(X)\cdot_{\ig_N}\psi
 \qquad\text{for every tangent vector }X.
\]
At each point, the infinitesimal rotations preserving the $G_2$-structure
form a $14$-dimensional space of skew-adjoint endomorphisms, denoted by
$\mathcal W_2$. The symmetric trace-free endomorphisms of the same
tangent space form a $27$-dimensional space, denoted by $\mathcal W_3$.

\begin{Thm}[{Agricola-Chiossi-Friedrich-H\"oll
\cite[Theorem~4.6]{ACFH}}]\label{thm:G2-harmonic-unit-spinor}
The spinor $\psi$ is harmonic if and only if, at every point,
the skew-adjoint part of $S$ lies in $\mathcal W_2$ and the symmetric
part of $S$ lies in $\mathcal W_3$.
\end{Thm}

In particular, a closed seven-manifold with a torsion-free
$G_2$-structure satisfies (H2) for the induced spin structure.
In dimension six, the corresponding $\mathrm{SU}(3)$ criterion is
given in \cite[Theorem~3.7]{ACFH}. The six-manifold used in the
introduction is their compact quotient of a Lie group. Its defining
unit spinor is harmonic and has nonzero covariant derivative
\cite[Example~3.8]{ACFH}.

For the Sobolev spaces, choose a complex orthonormal eigenbasis
$\{\psi_j\}$ with $D_\ig\psi_j=\lm_j\psi_j$. For $s>0$ the domain
of $|D_\ig|^s$ consists of the spinors $\psi=\sum_j\al_j\psi_j$ with
$\sum_j(1+|\lm_j|^{2s})|\al_j|^2<\infty$, and
$|D_\ig|^s\psi=\sum_j|\lm_j|^s\al_j\psi_j$. This domain is
$H^s(M,\mbs(M))$ (see \cite{Ammann}) and we equip it with the real inner product
\[
 \inp{\psi}{\va}_{s,2}
  =(|D_\ig|^s\psi,|D_\ig|^s\va)_2+(\psi,\va)_2.
\]
If $D_\ig$ is invertible, $\||D_\ig|^s\psi\|_2$ is an equivalent norm.
On $\R^m$ the homogeneous energy space $\msd^{1/2}(\R^m,\mbs(\R^m))$
is the completion of compactly supported smooth spinors with respect to the norm
$\||D_{\ig_{\R^m}}|^{1/2}\psi\|_2$.

\subsection{Spinors on product manifolds}\label{subsec: product}
	In this subsection, we follow the notation of \cite{Ba-Xu-CV22, SX2023}. Let $(M=N_1\times N_2,\ig=\ig_{N_1}\op\ig_{N_2})$ be a product of Riemannian spin $n_j$-manifolds $(N_j,\ig_{N_j},\sa_{N_j})$, $j=1,2$. We have
	\[
	{P_{\Spin}(M)=(P_{\Spin}(N_1)\times P_{\Spin}(N_2))\times_\zeta\Spin(n_1+n_2)}
	\]
	where $\zeta:\Spin(n_1)\times \Spin(n_2)\to\Spin(n_1+n_2)$ is the Lie group homomorphism lifting the standard embedding $\SO(n_1)\times \SO(n_2)\to\SO(n_1+n_2)$.
	
	The spinor bundle over $M$ can be identified with
	\[
	\mbs(M)=\left\{
	\aligned
	&(\mbs(N_1)\op\mbs(N_1))\ot \mbs(N_2) &\quad & \text{both } n_1 \text{ and } n_2 \text{ are odd}, \\[0.3em]
	&\qquad \mbs(N_1)\ot\mbs(N_2) &\quad & n_1 \text{ is even}.
	\endaligned \right.
	\]
	{When one factor is even-dimensional we place it first. In the
above identification the Clifford action is as follows. Take}
$X\in TN_1$, $Y\in TN_2$, $\va\in\Ga(\mbs(N_2))$ and
	\[
	\psi=\begin{cases}
		\psi_1\op\psi_2\in \Ga(\mbs(N_1)\op\mbs(N_1)) & \text{for both
			$n_1$ and $n_2$ odd} \\[0.3em]
		\qquad  \psi\in\Ga(\mbs(N_1)) & \text{for $n_1$ even}
	\end{cases}
	\]
	we have
	\begin{\equ}\label{product spinor representation}
		(X\op Y)\cdot_{\ig_M} (\psi\ot\va)=(X\cdot_{\ig_{N_1}}\psi)\ot\va+(\om_\C^{N_1}\cdot_{\ig_{N_1}}\psi)\ot(Y\cdot_{\ig_{N_2}}\va)
	\end{\equ}
	where in case $n_1$ and $n_2$ odd we set $X\cdot_{\ig_{N_1}}\psi=(X\cdot_{\ig_{N_1}}\psi_1)\op(-X\cdot_{\ig_{N_1}}\psi_2)$ and $\om_\C^{N_1}\cdot_{\ig_{N_1}}\psi=i(\psi_2\op-\psi_1)$. {Since $\om_\C^{N_1}$ interchanges the two summands and $X$ acts on them
with opposite signs, the Clifford actions of $X\oplus0$ and $0\oplus Y$
anticommute, as required for orthogonal vectors.
In tensor products over $M=N_1\times N_2$, each bundle on $N_j$ is pulled
back by the projection $M\to N_j$.}

Let $\nabla^{\mbs(N_1)}$ and $\nabla^{\mbs(N_2)}$ be the
	Levi-Civita connections on $\mbs(N_1)$ and $\mbs(N_2)$. By
	\[
	\nabla^{\mbs(N_1)\ot\mbs(N_2)}=\nabla^{\mbs(N_1)}\ot\id_{\mbs(N_2)}
	+\id_{\mbs(N_1)}\ot \nabla^{\mbs(N_2)}
	\]
	we mean the tensor product connection on $\mbs(N_1)\ot\mbs(N_2)$. Then, by \eqref{product spinor representation}, the Dirac operator on $M$ is given by
	\begin{\equ}\label{product-Dirac-operator}
		D_{\ig}^M = \tilde D_{\ig_{N_1}}^{N_1} \ot \id_{\mbs(N_2)}+ (\om_\C^{N_1}\cdot_{\ig_{N_1}}\id_{\mbs(N_1)})\ot D_{\ig_{N_2}}^{N_2}
	\end{\equ}
	where
	$\tilde D_{\ig_{N_1}}^{N_1}= D_{\ig_{N_1}}^{N_1}\op -D_{\ig_{N_1}}^{N_1}$ if both $n_1$ and $n_2$ are odd
	and $\tilde D_{\ig_{N_1}}^{N_1}={D_{\ig_{N_1}}^{N_1}}$ if $n_1$ is even.
	
%	For the case $n_1+n_2$ even, we have the decomposition
%	$\mbs(M)=\mbs(M)^+\op\mbs(M)^-$ and, moreover, when restrict $D_{\ig}^M$
%	on those half-spinor spaces we get
%	$D_{\ig}^M:\Ga(\mbs(M)^\pm)\to\Ga(\mbs(M)^\mp)$.

Suppose that $N_1$ satisfies (H1) and $N_2$ satisfies (H2).
Choose unit-length spinors $\psi_j$ on $N_1$ and $\va$ on $N_2$ such that
\[
 D_{\ig_{N_1}}^{N_1}\psi_j=\lm_j\psi_j,\qquad
 |\lm_j|\to\infty,\qquad D_{\ig_{N_2}}^{N_2}\va=0.
\]
After discarding finitely many terms we may assume $\lm_j\ne0$ for every
$j$. To realize $|\lm_j|$ as an eigenvalue of $D_\ig^M$, set
\[
 \phi_j=
 \begin{cases}
 \psi_j\ot\va,
     & n_1\text{ even},\ \lm_j>0,\\[0.2em]
 (\om_\C^{N_1}\cdot_{\ig_{N_1}}\psi_j)\ot\va,
     & n_1\text{ even},\ \lm_j<0,\\[0.2em]
 (\psi_j\op 0)\ot\va,
     & n_1,n_2\text{ odd},\ \lm_j>0,\\[0.2em]
 (0\op\psi_j)\ot\va,
     & n_1,n_2\text{ odd},\ \lm_j<0.
 \end{cases}
\]
Since $\va$ is harmonic, only the first term of
\eqref{product-Dirac-operator} contributes. When $n_1$ is even,
Clifford multiplication by $\om_\C^{N_1}$ reverses the sign of the
Dirac eigenvalue. When both dimensions are odd, the second summand
carries $-D_{\ig_{N_1}}^{N_1}$. Hence in all four cases
\[
 D_\ig^M\phi_j=|\lm_j|\phi_j,\qquad |\phi_j|_\ig\equiv1.
\]
For $m=n_1+n_2$, the spinor $|\lm_j|^{(m-1)/2}\phi_j$ therefore
solves \eqref{SY0} on the product $M=N_1\times N_2$. 
Its length is constant on $M$ and tends to infinity
as $j\to\infty$, so $\co_j=M$ for every $j$. In particular, when $N_1$
and $N_2$ are flat tori with their trivial spin structures, this
construction gives the example of $k=0$ in the
first row of Table~\ref{tab:noncompactness-examples}.

\section{Specific ansatz for the spinorial Yamabe equation}\label{sec:ansatz}

Having treated the case $\co_j=M$, we now seek examples on
different manifolds $M$ of blow-up sequences $\{\psi_j\}$ for \eqref{SY0}
such that each $|\psi_j|_\ig$ attains its maximum precisely on a proper subset
$\co_j$ of $M$. On product manifolds, our aim is to realize
the remaining codimensions $1\leq k\leq m-1$ in
Theorem~\ref{main-thm}~(1). For this purpose, we use a constant-length
spinor on one factor to reduce the equation to an ordinary differential
system on the other. On the sphere, stereographic projection gives a
critical equation on $\R^m$ whose coefficients record the choice of metric.

\subsection{For the Riemannian product spaces}

	\subsubsection{$M=S^1\times N$ with $\dim N$ odd}

	\begin{revision}
Write $S^1=\R/(2\pi\Z)$ with coordinate $t$ and take the trivial spin
structure. We choose Clifford multiplication by the positively oriented
unit tangent to act on $\C$ as multiplication by $i$. For odd
$n=\dim N$, formula \eqref{product-Dirac-operator} gives
\end{revision}

\[
	D_{\ig}^M\phi=i\Big( \frac{d}{dt}\psi_1\op -\frac{d}{dt}\psi_2\Big)\otimes\va
	+i(\psi_2\op -\psi_1)\otimes D_{\ig_N}^N\va
	\]
	for all $\phi=(\psi_1\op\psi_2)\otimes\va \in\Ga(\mbs(M))$
	where the spinor bundle of $M$ is given by
	\[
	\mbs(M)=\big( \mbs(S^1)\op \mbs(S^1)\big)\otimes \mbs(N).
	\]
	{Thus $\psi_1$ and $\psi_2$ are $2\pi$-periodic complex functions and
the equation on $S^1\times N$ is}
\begin{\equ}\label{SYS1N}
		D_{\ig}^M\phi=|\phi|_{\ig}^{\frac2n}\phi, \quad  n=\dim N
	\end{\equ}
	where $|\phi|_{\ig}=\sqrt{|\psi_1|^2+|\psi_2|^2}\,|\va|_{\ig_N}$.
	
	{Choose a unit eigenspinor $\va_\lm$ of $D^N_{\ig_N}$ with eigenvalue
$\lm$. Substitution of $\phi=(\psi_1\op\psi_2)\otimes\va_\lm$ into
\eqref{SYS1N} gives}
\begin{\equ}\label{system1}
		\left\{
		\aligned
		i\frac{d\psi_1}{dt}+i\lm\psi_2&=\big( |\psi_1|^2+|\psi_2|^2 \big)^{\frac1n}\psi_1 \\[0.3em]
		-i\frac{d\psi_2}{dt}-i\lm\psi_1&=\big( |\psi_1|^2+|\psi_2|^2 \big)^{\frac1n}\psi_2
		\endaligned
		\right.
	\end{\equ}
	on $\R/2\pi\Z$.

{The product formula also gives a sequence of constant-length
solutions for this fixed $\va_\lm$. For $j\geq1$ and
$(c_1,c_2)\in\C^2$, the Dirac operator on a spinor of the form
$e^{ijt}(c_1\oplus c_2)\ot\va_\lm$ acts on $(c_1,c_2)$ by the
Hermitian matrix
\[
 \begin{pmatrix}-j&i\lm\\-i\lm&j\end{pmatrix}.
\]
A unit eigenvector for its positive eigenvalue $\sqrt{j^2+\lm^2}$
gives a unit-length eigenspinor on $S^1\times N$. Multiplying this
spinor by $(j^2+\lm^2)^{n/4}$ gives a solution of \eqref{SYS1N}
whose constant length tends to infinity as $j\to\infty$.}

The system has an invariant real subspace given by
$\psi_2=\overline{\psi_1}$, where the bar denotes complex conjugation.
Indeed, on this subspace the second equation
is the complex conjugate of the first. Substituting
$\psi_1=u+iv$ and $\psi_2=u-iv$, with $u,v$ real, gives
\begin{\equ}\label{eqs2}
			\left\{
			\aligned
			u'+\lm u&=2^{\frac1n}(u^2+v^2)^{\frac1n}v \\
			-v'+\lm v&= 2^{\frac1n}(u^2+v^2)^{\frac1n}u
			\endaligned \right. 
		\end{\equ}
on $\R/(2\pi\Z)$. Since $|\va_\lm|_{\ig_N}\equiv1$, the squared
length of the product spinor is $2(u^2+v^2)$. If $u$ and $v$ vanished
at the same point, uniqueness for \eqref{eqs2} would give
$(u,v)\equiv(0,0)$, since the vector field is locally Lipschitz and
the origin is an equilibrium. Thus every nontrivial solution obtained
from this system is nowhere vanishing.
For $\lm=0$, the pairs
\[
 (u_j(t),v_j(t))=\frac{j^{n/2}}{\sqrt2}(\sin jt,\cos jt),
 \qquad j\in\N,
\]
give constant-length solutions. To obtain solutions of nonconstant
length, assume (H3) on $N$ and fix the corresponding eigenspinor with
eigenvalue $\lm\ne0$. The map $(u,v)\mapsto(v,-u)$ sends solutions
for $\lm$ to solutions for $-\lm$.
We shall therefore take $\lm>0$. Along solutions,
$(u^2+v^2)'=2\lm(v^2-u^2)$. Thus a nonzero solution with constant
$u^2+v^2$ has $u=\pm v$ and is stationary. Section~\ref{sec: product}
constructs nonstationary periodic solutions for a single fixed $\lm$.
Their norms depend on the circle variable alone. Odd-dimensional flat tori
and round spheres supply admissible choices of $N$.

The spinorial Weierstrass representation also gives examples on
$S^1\times N$. Suppose that the closed oriented manifold $N$ is
isometrically immersed as a hypersurface of constant mean curvature $H$
in a spin manifold $\widetilde N$. Equip $N$ with the induced spin structure and suppose that
$\widetilde N$ carries a nonzero Killing spinor with real Killing
constant $\al$. Normalize this spinor to have unit length and assume
$(H,\al)\ne(0,0)$. We use the restriction argument of
\cite[Section~5.3.2]{Ammann} in the odd-dimensional case considered here.
Writing the ambient spinor as $\Psi^++\Psi^-$ and choosing the unit
normal $\nu$, we identify the restrictions of $\Psi^+$ and
$\nu\cdot\Psi^-$ with two spinors on $N$ whose squared lengths sum to one.
Taking the Dirac trace of the hypersurface equations in
\cite[Theorem~3.1]{Conti-SegnanDalmasso} gives, in their conventions,
the constant symmetric matrix
$-n\left(\begin{smallmatrix}H/2&\al\\ \al&-H/2\end{smallmatrix}\right)$
on this pair, where $H$ is the average of the principal curvatures.
Its eigenvalues are $\pm\lm$, where $\lm=\frac n2\sqrt{H^2+4\al^2}$.
Choosing an orthonormal basis in which this matrix has off-diagonal
entries $\lm$ and zero diagonal gives a pair of spinors $\eta_1,\eta_2$
satisfying
\[
 D^N_{\ig_N}\eta_1=\lm\eta_2,\qquad
 D^N_{\ig_N}\eta_2=\lm\eta_1,\qquad
 |\eta_1|_{\ig_N}^2+|\eta_2|_{\ig_N}^2=1.
\]
For any solution $(u,v)$ of \eqref{eqs2} with this $\lm$, put $z=u+iv$
and take the product spinor
\[
 \phi=\sqrt{2}\bigl((z\oplus0)\ot\eta_1
                     +(0\oplus\overline z)\ot\eta_2\bigr).
\]
Then $|\phi|_\ig^2=2|z|^2$, and substitution in
\eqref{product-Dirac-operator} gives exactly \eqref{system1} with
$\psi_1=z$ and $\psi_2=\overline z$. Thus the same periodic solutions
give blow-up sequences on these products.

\subsubsection{$M=S^k\times N$ with $k>1$}\label{subsubsec:MSN}

Equip $S^k\times N$ with the product metric
$\ig=\ig_{S^k}\op\ig_N$, and put $n=\dim N$. Suppose first that $k$ is
even or that both $k$ and $n$ are odd. For spinors of product form the
equation is
\begin{\equ}\label{SYSkN}
		D_\ig^M\phi=|\phi|_{\ig}^{\frac{2}{k+n-1}}\phi, \quad \phi=\psi\ot\va\in\Ga(\mbs(M))
		\end{\equ}
		which can be expanded as
		\begin{\equ}\label{SYSkN2}
		(\tilde D_{\ig_{S^k}}^{S^k}\psi)\ot\va+\big( \om_\C^{S^k}\cdot_{\ig_{S^k}}\psi \big)\ot(D_{\ig_N}^N\va)=\big(|\psi|_{\ig_{S^k}}|\va|_{\ig_N}\big)^{\frac{2}{k+n-1}}\psi\ot\va.
		\end{\equ}
		
Suppose now that the Dirac operator $D_{\ig_N}^N$ on $N$
satisfies (H2), so that we may fix a harmonic spinor $\va_N$ with
$|\va_N|_{\ig_N}\equiv1$. For
$\phi=\psi\ot\va_N$, formula \eqref{SYSkN2} reduces to
\begin{\equ}\label{eqs2-2}
			\tilde D_{\ig_{S^k}}^{S^k}\psi=|\psi|_{\ig_{S^k}}^{\frac{2}{k+n-1}}\psi \quad \text{on } S^k.
		\end{\equ}
When $k$ is even, $\tilde D_{\ig_{S^k}}^{S^k}$
coincides with the Dirac operator $D_{\ig_{S^k}}^{S^k}$ acting on sections of $\mbs(S^k)$.
When both $k$ and $n$ are odd, $\tilde D_{\ig_{S^k}}^{S^k}$ acts on
sections of $\mbs(S^k)\oplus\mbs(S^k)$ by
\[
 \tilde D_{\ig_{S^k}}^{S^k}(\psi_1\oplus\psi_2)
 =(D_{\ig_{S^k}}^{S^k}\psi_1)\oplus(-D_{\ig_{S^k}}^{S^k}\psi_2).
\]

{In every dimension, we write $P_{North}$ for the north pole
of the standard sphere.} Suppose first that $k$ is even. Let
$\om:S^k\setminus\{P_{North}\}\to\R^k$ be stereographic projection from
the north pole, and put $\rho(x)=2/(1+|x|^2)$, so that
$(\om^{-1})^*\ig_{S^k}=\rho^2\ig_{\R^k}$.
{For $p\in S^k\setminus\{P_{North}\}$, put $x=\om(p)$
and let $\iota_p:\mbs_p(S^k)\to\mbs_x(\R^k)$ be the induced isometry
between the spinor fibers. By Proposition~\ref{prop:conformal-invar}
in dimension $k$, the spinor defined by
$\va(x)=\rho(x)^{(k-1)/2}\iota_p(\psi(p))$ satisfies}
\begin{\equ}\label{eqs2-2-1}
		D_{\ig_{\R^k}}\va=\rho^{\frac{n}{k+n-1}}|\va|_{\ig_{\R^k}}^{\frac{2}{k+n-1}}\va \quad \text{on } (\R^k,\ig_{\R^k}).
		\end{\equ}
Moreover, since $d\vol_{\ig_{S^k}}=\rho^k\,dx$ in these coordinates, we have
\[
 \int_{S^k}|\psi|_{\ig_{S^k}}^{\frac{2(k+n)}{k+n-1}}\,d\vol_{\ig_{S^k}}
 =\int_{\R^k}\rho^{\frac{n}{k+n-1}}
          |\va|_{\ig_{\R^k}}^{\frac{2(k+n)}{k+n-1}}\,dx .
\]
To construct solutions on $S^k$, we shall seek 
solutions $\va$ of \eqref{eqs2-2-1} for which the weighted integral
over $\R^k$ above is finite. We recover $\psi$ at each
$p\in S^k\setminus\{P_{North}\}$ by
\[
 \psi(p)=\rho(\om(p))^{-(k-1)/2}\,
          \iota_p^{-1}\bigl(\va(\om(p))\bigr).
\]
Since the nonlinearity has subcritical exponent (in the sense of Sobolev embedding) on $S^k$, Ammann's theorems on the removal of singularities and on regularity
\cite[Theorems~5.1 and~5.2]{Ammann2009} ensure that $\psi$
extends across $P_{North}$ to a $C^{1,\alpha}$ solution of
\eqref{eqs2-2} for some $\alpha>0$.

{The subcritical nonlinearity in \eqref{eqs2-2} allows
one to obtain infinitely many solutions on $S^k$ by variational methods, see for instance
\cite{Isobe-subcritical}. For the product equation \eqref{SYSkN}, however,
we need solutions that are nowhere vanishing. Indeed, if a solution $\psi$ of \eqref{eqs2-2} vanishes
at a point $p\in S^k$, then $\phi=\psi\ot\va_N$ vanishes on all of
$\{p\}\times N$. The variational argument alone gives no control over
these zeros. Spinor fibers have no natural order, so the usual maximum
principle argument for positive solutions of a scalar elliptic equation
is unavailable here.

To control the length pointwise, we seek solutions of
\eqref{eqs2-2-1} with two real radial coefficients. Fix a unit constant
spinor $\ga_0\in\mbs_k$ and set $r=|x|$. For real functions $f_1,f_2$
on $(0,\infty)$, consider the spinor}
\[
 \va(x)=f_1(r)\ga_0+f_2(r)\frac{x}{r}\cdot_{\ig_{\R^k}}\ga_0 .
\]
{Since Clifford multiplication by the unit vector $x/r$ is skew-adjoint
and preserves length, we have $|\va(x)|_{\ig_{\R^k}}^2=f_1(r)^2+f_2(r)^2$.
Thus a zero at radius $r>0$ would require both coefficients to vanish
there.} Direct differentiation gives
\[
 D_{\ig_{\R^k}}\va
 =-\left(f_2'(r)+\frac{k-1}{r}f_2(r)\right)\ga_0
    + f_1'(r)\frac{x}{r}\cdot_{\ig_{\R^k}}\ga_0 .
\]
{We now make the Emden-Fowler change of variables
\[
 r=e^{-t},\qquad
 f_1(r)=-e^{(k-1)t/2}u(t),\qquad
 f_2(r)=e^{(k-1)t/2}v(t),\qquad t\in\R.
\]
Since $r\rho(r)=(\cosh t)^{-1}$, \eqref{eqs2-2-1} reduces to}
\begin{\equ}\label{eqs2-2-1-1}
			\left\{
			\aligned
			u'+\frac{k-1}2 u &= \cosh(t)^{-\frac{n}{k+n-1}}(u^2+v^2)^{\frac1{k+n-1}}v \\[0.3em]
			-v' +\frac{k-1}2v &= \cosh(t)^{-\frac{n}{k+n-1}}(u^2+v^2)^{\frac1{k+n-1}}u
			\endaligned
			\right.
		\end{\equ}
on $\R$. Conversely, every solution $(u,v)$ of \eqref{eqs2-2-1-1} on $\R$
gives a solution of \eqref{eqs2-2-1} on $\R^k\setminus\{0\}$ by setting
\[
 \va(x)=r^{-(k-1)/2}
 \left[-u(t)\ga_0+v(t)\frac{x}{r}\cdot_{\ig_{\R^k}}\ga_0\right],
 \qquad r=|x|,\quad t=-\ln r.
\]
The corresponding product spinor $\phi=\psi\ot\va_N$ has squared length
\[
 |\phi|_{\ig}^2=(\cosh t)^{k-1}\big(u(t)^2+v(t)^2\big).
\]
If this expression is bounded as $t\to+\infty$, the singularity of
$\va$ at $x=0$ is removable. If it is also bounded as $t\to-\infty$,
the spherical spinor $\psi$ extends across $P_{North}$ by the same
removability theorem, and $\phi$ solves \eqref{SYSkN} on all of $S^k\times N$.

When $k$ and $n$ are both odd, $\psi$ is a section of
$\mbs(S^k)\oplus\mbs(S^k)$. Taking $\psi=\psi_1\oplus 0$ gives
$\tilde D_{\ig_{S^k}}^{S^k}\psi=(D_{\ig_{S^k}}^{S^k}\psi_1)\oplus 0$,
so \eqref{eqs2-2} reduces to the same equation for $\psi_1$ as in the
even-dimensional sphere case. Hence the
same radial system and length formula apply.
{If $k$ is odd and $n$ is even, we place $N$ first in
\eqref{product-Dirac-operator}. In this case assume that $N$ carries
a unit harmonic spinor $\va_N$ of positive chirality, meaning that
$\om_\C^N\cdot_{\ig_N}\va_N=\va_N$. Since $D_{\ig_N}^N\va_N=0$,
the product formula gives
\[
 D_\ig^M(\va_N\ot\psi)=\va_N\ot D_{\ig_{S^k}}^{S^k}\psi.
\]
Thus the same equation on $S^k$ applies in this case as well.
An even-dimensional flat torus with the trivial spin structure carries
such a spinor.}

The Hamiltonian for \eqref{eqs2-2-1-1} is
\[
 H(t,u,v)=-\frac{k-1}{2}uv+\frac{k+n-1}{2(k+n)}
   (\cosh t)^{-\frac{n}{k+n-1}}(u^2+v^2)^{\frac{k+n}{k+n-1}} .
\]
Along each solution $dH/dt=\partial_t H$, which is negative for $t>0$
along any nonzero solution and positive for $t<0$. Section~\ref{sec: product}
uses the monotonicity on the positive half-line and the symmetry under
$t\mapsto-t$, $u\leftrightarrow v$ to construct the required global
solutions.

\subsection{For the perturbed metric on $M=S^m$}\label{subsec: perturbed-bg}
	
For the construction on $S^m$ in dimension $m\ge20$, we seek
nontrivial solutions of
\begin{\equ}\label{SY-m-sphere}
		D_{\ig_\vr}\psi=|\psi|_{\ig_\vr}^{\frac2{m-1}}\psi \quad \text{on } S^m,
	\end{\equ}
where $\ig_\vr$ is a smooth metric close to the round metric.
Let $\om:S^m\setminus\{P_{North}\}\to\R^m$ be stereographic projection from
the north pole. We write the metric in Euclidean coordinates as
\begin{\equ}\label{the metric g}
		\tilde\ig_\vr=\rho^{-2}\cdot(\om^{-1})^*\ig_\vr, \quad \rho(x)=\frac2{1+|x|^2}, \ x\in\R^m.
	\end{\equ}

The choice $\tilde\ig_\vr=\ig_{\R^m}$ corresponds to the round
sphere. To construct $\ig_\vr$, choose points $x_j\in\R^m$ tending to
infinity and smooth symmetric matrix fields $h^{(j)}$ whose translates
$h^{(j)}(\,\cdot-x_j)$ have pairwise disjoint supports. We set
\begin{\equ}\label{h-expansion}
 \tilde\ig_\vr(x)=\exp\big(\vr\tilde h(x)\big),\qquad
 \tilde h(x)=\sum_{j=1}^\infty h^{(j)}(x-x_j).
\end{\equ}
The sizes and locations of these supports will be chosen together with
the amplitudes of the tensors
so that the compactified metric is smooth at $P_{North}$.

For $p\in S^m\setminus\{P_{North}\}$ and $x=\om(p)$, let
$\iota_p:\mbs_p(S^m,\ig_\vr)\to\mbs_x(\R^m,\tilde\ig_\vr)$
be the spinor isometry induced by \eqref{the metric g}.
By conformal covariance, $\va(x)=\rho(x)^{(m-1)/2}\iota_p(\psi(p))$
satisfies
\begin{\equ}\label{reduced Dirac problem 1}
 D_{\tilde\ig_\vr}\va=|\va|_{\tilde\ig_\vr}^{2/(m-1)}\va
 \qquad\text{on }\R^m .
\end{\equ}
This identification preserves the $L^{\frac{2m}{m-1}}$ norm of the spinor. For a
smooth compactification, a finite-energy solution on $\R^m$ extends
weakly across the north pole by the cutoff argument for the Dirac
equation \cite[Theorem~5.1]{Ammann2009}.

For a fixed small $\vr\ne0$, Section~\ref{sec: perturbation} constructs
solutions approaching standard Euclidean bubbles in
$\msd^{1/2}(\R^m,\mbs(\R^m))$. Write $\xi_j$ and $\lm_j$ for the
center and scale of the $j$-th bubble. Since $|\xi_j|\to\infty$ and
$\lm_j\to0$, the sets $U_j=\om^{-1}(B_{\lm_j}(\xi_j))$ have spherical
volumes tending to zero. After rescaling each ball $B_{\lm_j}(\xi_j)$
to the unit ball, the standard bubbles all have the same positive
$L^{\frac{2m}{m-1}}$-norm there. Convergence in this norm therefore
gives a constant $c>0$, independent of $j$, such that the corresponding
spherical solutions $\psi_j$ satisfy
\[
 c\leq \int_{U_j}|\psi_j|_{\ig_\vr}^{\frac{2m}{m-1}}\,d\vol_{\ig_\vr}
 \leq \operatorname{Vol}(U_j,\ig_\vr)
      \left(\max_{S^m}|\psi_j|_{\ig_\vr}\right)^{\frac{2m}{m-1}}
\]
for all sufficiently large $j$. Hence $\max_{S^m}|\psi_j|_{\ig_\vr}\to\infty$.

	\section{Noncompactness on product manifolds}\label{sec: product}

The reductions \eqref{eqs2} and \eqref{eqs2-2-1-1} lead to different boundary conditions. On the circle we choose the period of a closed orbit. On the sphere the two poles correspond to the ends of the real line, and we seek solutions approaching the origin at both ends. In either case, if $z=(u,v)$ and $H$ is the corresponding Hamiltonian, the system takes the form
\begin{\equ}\label{Hamiltionian system}
 z'=J\nabla_zH(t,z),\qquad J=\begin{pmatrix}0&1\\-1&0\end{pmatrix}.
\end{\equ}

	\subsection{The autonomous case}

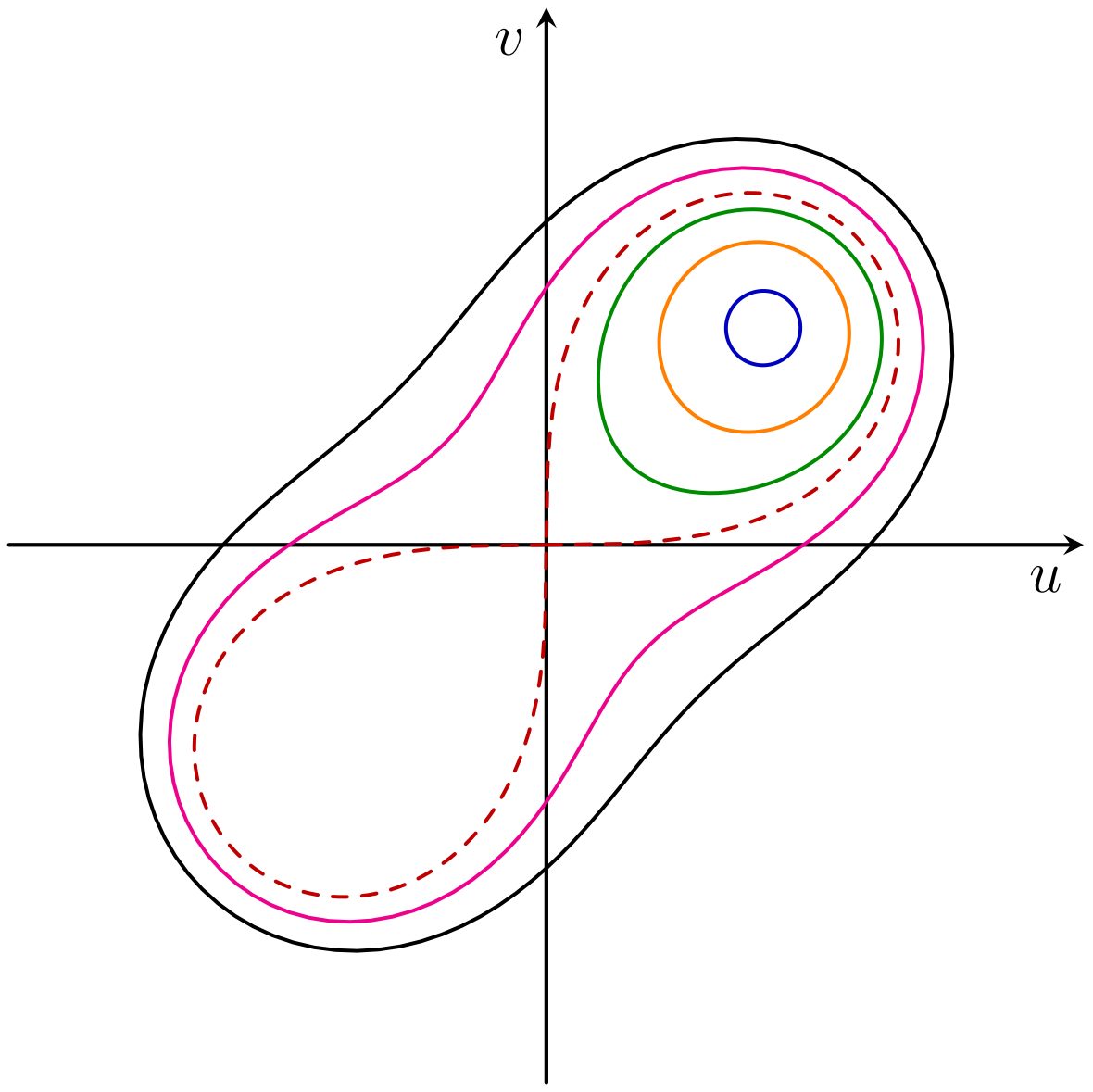
\begin{figure}[h!]
	\centering
	\begin{minipage}{0.5\linewidth}
		\centering
		\begin{tikzpicture}[
			x=.39\linewidth, y=.39\linewidth, >=stealth,
			line width=.75pt, line join=round, line cap=round,
			declare function={
				outerradius(\anglevar,\energy)=
					sqrt((sin(2*\anglevar)+sqrt((sin(2*\anglevar))^2+8*\energy))/2);
				sumuv(\anglevar,\amplitude)=
					sqrt(1+2*\amplitude*cos(\anglevar)-(\amplitude*sin(\anglevar))^2);
			}]
			% For n=lambda=1, H=(u^2+v^2)^2/2-uv.
			\draw[->, thick] (-1.23,0) -- (1.23,0) node[below left] {$u$};
			\draw[->, thick] (0,-1.23) -- (0,1.23) node[below left] {$v$};
			% Positive levels: solve H=E for the polar radius.
			\foreach \energy/\orbitcolor in {.15/black,.06/magenta}{
				\draw[\orbitcolor] plot[domain=0:360,samples=145,variable=\t]
					({outerradius(\t,\energy)*cos(\t)},{outerradius(\t,\energy)*sin(\t)})
					-- cycle;
			}
			% At H=0, r^2=sin(2t) in the first and third quadrants.
			\foreach \rotation in {0,180}{
				\draw[rotate=\rotation,red!75!black,dashed]
					plot[domain=0:90,samples=101,variable=\t]
					({sqrt(max(0,sin(2*\t)))*cos(\t)},
					 {sqrt(max(0,sin(2*\t)))*sin(\t)});
			}
			% On a negative level, (u^2+v^2-1/2)^2+(u-v)^2=1/4+2E.
			% Parametrize this circle with radius \amplitude and recover u,v.
			\foreach \amplitude/\orbitcolor in {.43/green!55!black,.30/orange,.12/blue!75!black}{
				\draw[\orbitcolor] plot[domain=0:360,samples=101,variable=\t]
					({(sumuv(\t,\amplitude)+\amplitude*sin(\t))/2},
					 {(sumuv(\t,\amplitude)-\amplitude*sin(\t))/2}) -- cycle;
			}
		\end{tikzpicture}
		\caption{The positive and negative levels of
				\eqref{H1} in the $(u,v)$-plane.  The dashed curves represent the zero level and the solid curves represent periodic orbits. The periodic orbits in the third quadrant are omitted.}
			\label{fig:eqs2}
	\end{minipage}
\end{figure}

Let $n=\dim N$. Fix a unit eigenspinor on the odd-dimensional factor $N$ with eigenvalue $\lm\ne0$, and take the trivial spin structure on $S^1=\R/2\pi\Z$. A rotation of the $(u,v)$-plane by $\pi/2$ exchanges $\lm$ and $-\lm$ in \eqref{eqs2}, so we may take $\lm>0$. The Hamiltonian
\begin{\equ}\label{H1}
 H(u,v)=-\lm uv+\frac{n}{4(n+1)}
            \bigl[2(u^2+v^2)\bigr]^{\frac{n+1}{n}}
\end{\equ}
is constant along each orbit, since $J$ is skew-symmetric. The zero level consists of the origin and two homoclinic orbits. Up to translation, one of them is
\[
 u_0(t)=\frac{\bigl(\frac{n+1}{n}\lm\bigr)^{n/2}e^{\lm t/n}}{2\cosh(2\lm t/n)^{(n+1)/2}},
 \qquad v_0(t)=u_0(-t),
\]
and its negative gives the other. They enclose the two bounded components of $\{H<0\}$ shown in Figure~\ref{fig:eqs2}. The two-dimensional torus calculation in \cite[Section~7]{SX2023} concerns orbits in the positive quadrant. We also consider orbits on the positive energy levels, which cross the coordinate axes.

\begin{Prop}\label{prop:tauk}
For each $E>0$, the level $\{(u,v)\,:\, H(u,v)=E\}$ corresponds to a closed orbit of \eqref{eqs2}. Its least period $T(E)$ is continuous and strictly decreasing, with $T(E)\to\infty$ as $E\searrow0$ and $T(E)\to0$ as $E\nearrow+\infty$.
\end{Prop}
\begin{proof}
In polar coordinates write $(u,v)=\sqrt{s}(\cos\theta,\sin\theta)$, where $s=u^2+v^2$.
For each fixed angle $\theta$, the energy equation $H(u,v)=E$ becomes
\[
 \frac{n}{4(n+1)}(2s)^{\frac{n+1}{n}}-\frac{\lm}{2}s\sin(2\theta)=E.
\]
For $E>0$, this equation has exactly one positive root $s=s_E(\theta)$.
At this root, dividing the energy equation by $s$ gives
$\frac{\lm}{2}\sin(2\theta)=\frac{n}{2(n+1)}(2s)^{1/n}-E/s$.
Thus the derivative of the left-hand side of the energy equation with
respect to $s$ at this root is
\begin{\equ}\label{ds}
 \frac{(2s)^{1/n}}{2}-\frac{\lm}{2}\sin(2\theta)
 =\frac{(2s)^{1/n}}{2(n+1)}+\frac{E}{s}>0.
\end{\equ}
The implicit function theorem now gives smooth dependence of $s_E(\theta)$
on $(E,\theta)$. Differentiating the energy equation with respect to $E$
at fixed $\theta$ gives
\[
 \partial_E s_E(\theta)=
 \left(\frac{\bigl(2s_E(\theta)\bigr)^{1/n}}{2(n+1)}
       +\frac{E}{s_E(\theta)}\right)^{-1}>0.
\]
To compute the period, observe that $uv'-vu'=s\theta'$ in polar
coordinates. With $H$ as in \eqref{H1}, system \eqref{eqs2} takes
the Hamiltonian form $u'=\partial_v H$, $v'=-\partial_u H$. Hence
\[
 -\theta'(t)=\frac{u\,\partial_u H+v\,\partial_v H}{s}
 =2\frac{\partial}{\partial s}H(\sqrt{s}\cos\theta,\sqrt{s}\sin\theta).
\]
Here $\theta'(t)$ is the time derivative of the angle along the orbit.
On the right, $\partial/\partial s$ differentiates the Hamiltonian
along a ray, varying $s$ while holding the angle $\theta$ fixed.
Along a solution $(u(t),v(t))$ of energy $E>0$, we have
$s(t)=u(t)^2+v(t)^2=s_E(\theta(t))$. Substituting this identity
into the formula for the angular velocity gives
\[
 -\theta'(t)=\bigl(2s_E(\theta(t))\bigr)^{1/n}-\lm\sin(2\theta(t))>0,
\]
where the last inequality follows from \eqref{ds}.
Thus the orbit runs clockwise. Since the energy level meets each ray
at exactly one point, the orbit first returns to its starting point
when $\theta$ has decreased by $2\pi$. Its least period is therefore
\begin{\equ}\label{period-T}
 T(E)=\int_{2\pi}^{0}\frac{d\theta}{\theta'}
     =\int_0^{2\pi}
       \frac{d\theta}{\bigl(2s_E(\theta)\bigr)^{1/n}-\lm\sin(2\theta)}.
\end{\equ}
To compare two energies $0<E_1<E_2$, fix an angle $\theta$. Since
$\partial_E s_E(\theta)>0$, we have $s_{E_2}(\theta)>s_{E_1}(\theta)$.
The term $\lm\sin(2\theta)$ is the same for both orbits, so the orbit
of energy $E_2$ crosses this ray with greater angular speed. 
Integrating over one period gives $T(E_2)<T(E_1)$.
The period also depends continuously on $E>0$, since the integrand is
jointly continuous in $(E,\theta)$ and $\theta$ ranges over $[0,2\pi]$.
As $E\to\infty$, $s_E(\theta)\to\infty$ uniformly, which gives $T(E)\to0$.
For the limit as $E\searrow0$, consider the time spent in the second
quadrant of the $(u,v)$-plane. For each fixed $\theta\in(\pi/2,\pi)$, we have $\sin(2\theta)<0$,
so both terms on the left-hand side of the energy equation are positive.
In particular,
\[
 0<s_E(\theta)\leq\frac{2E}{\lm|\sin(2\theta)|}.
\]
Thus $s_E(\theta)$ decreases to zero as $E\searrow0$, and the
integrand in the period formula \eqref{period-T} increases to $1/(\lm|\sin(2\theta)|)$.
Since $|\sin(2\theta)|\leq 2(\theta-\pi/2)$ for $\theta\in(\pi/2,\pi)$, the integral of this limit
diverges at $\theta=\pi/2$. By the monotone convergence theorem, the time spent in the second quadrant already tends to infinity, and hence so does $T(E)$.
\end{proof}

For each $j\geq1$ there is therefore a unique energy $E_j$ with $T(E_j)=2\pi/j$. We have $E_j\to\infty$, and the energy equation gives
\begin{\equ}\label{autonomous-orbit-length-asymptotic}
 2s_{E_j}(\theta)=\left(\frac{4(n+1)}{n}E_j\right)^{\frac{n}{n+1}}
 \bigl(1+o(1)\bigr)
 \qquad\text{as }j\to\infty,
\end{\equ}
uniformly in $\theta$. Let $(u_j(t),v_j(t))$ be a solution on the orbit
of energy $E_j$, and let $\theta_j(t)$ be its polar angle. Since the
eigenspinor $\va_\lm$ on $N$ has unit length, the corresponding product
spinor $\phi_j$ in the reduction \eqref{system1} satisfies
\[
 |\phi_j(t,y)|_{\ig}
 =\sqrt{2\bigl(u_j(t)^2+v_j(t)^2\bigr)}
 =\sqrt{2s_{E_j}(\theta_j(t))},\qquad (t,y)\in S^1\times N.
\]
Taking square roots in \eqref{autonomous-orbit-length-asymptotic} thus
gives the same leading term
$\bigl(4(n+1)E_j/n\bigr)^{\frac{n}{2(n+1)}}$ for the length at every
point, with a relative error tending to zero uniformly on $S^1\times N$.
In particular, the minimum and maximum lengths are both asymptotic to
this quantity. Since $E_j\to\infty$, the minimum length tends to infinity,
and the ratio of the maximum length to the minimum length tends to one.
For each fixed $j$, the pointwise length $|\phi_j(t,y)|_{\ig}$ is
independent of $y\in N$ and varies with $t\in S^1$. To locate its
maximum, write the energy equation on the orbit of energy $E_j$ as
\[
 \frac{\lm}{2}\sin(2\theta)
 =\frac{n}{2(n+1)}(2s)^{1/n}-\frac{E_j}{s},
 \qquad s=s_{E_j}(\theta).
\]
Since $E_j>0$, both terms on the right are strictly increasing in $s$.
As $\lm>0$ and $|\phi_j|_{\ig}^2=2s$, the length is therefore largest
precisely when $\sin(2\theta)=1$. In one period this occurs at
$\theta=\pi/4$ and $\theta=5\pi/4$ modulo $2\pi$. Since $T(E_j)=2\pi/j$, the interval $0\leq t<2\pi$ contains precisely $j$ periods. Thus $ |\phi_j(t,y)|_{\ig}$ attains its maximum exactly $2j$ times on this interval. Since its value is independent of $y$, the maximum set in $S^1\times N$ consists of $2j$ copies of $N$.

\medskip

On the flat torus $\ell S^1\times S^1$, let $2\pi\ell$ be the length of
the first circle and let $\lm>0$ be a Dirac eigenvalue on the second.
The period calculation underlying \cite[Theorem~7.3(1)]{SX2023} concerns
the negative levels of the Hamiltonian \eqref{H1}. Indeed, for $n=1$
and $u,v>0$, the integration constant $K$ in formula~(7.7) of
\cite{SX2023} satisfies
\[
 K=4uv-\frac{2}{\lm}(u^2+v^2)^2=-\frac{4}{\lm}H(u,v).
\]
Thus $0<K<\lm/2$, the range in \cite[Lemma~7.1]{SX2023}, corresponds
to $-\lm^2/8<H<0$. Every orbit in this range has least period
greater than $\pi/\lm$. Hence, when $\ell\leq1/(2\lm)$, the negative
levels give only the constant-length spinor at $H=-\lm^2/8$, up to
a constant phase.

Positive levels of the same Hamiltonian give additional solutions even
for the same circle lengths and spin structures. For example, take both
circles of length $2\pi$, with periodic spinors on the first and
antiperiodic spinors on the second. If $\tau$ denotes the angle on
the second circle, then $e^{-i\tau/2}$ has eigenvalue $\lm=1/2$,
so $\ell=1=1/(2\lm)$. For each $j\geq1$, the orbit of
Proposition~\ref{prop:tauk} with $H=E_j>0$ and least period $2\pi/j$
yields the solution
\[
 \phi_j(t,\tau)=
 \bigl((u_j(t)+iv_j(t))\oplus(u_j(t)-iv_j(t))\bigr)e^{-i\tau/2}.
\]
Its length is nonconstant, and its orbit crosses both coordinate axes.

\subsection{The non-autonomous case}

Fix $k\geq2$ and $n\geq1$. We state the shooting result for the system \eqref{eqs2-2-1-1} before returning to spinors on $S^k\times N$.

\begin{Thm}\label{main thm ODE2}
For every integer $j\geq0$ there is a solution $(u_j,v_j)$ of \eqref{eqs2-2-1-1} such that $u_j(0)=v_j(0)=\mu_j>0$, $u_j(-t)=v_j(t)$, and $v_j$ has exactly $j$ zeros on $(0,\infty)$. The initial values can be chosen so that $\mu_j$ increases to infinity. For each $j$ there is a finite $C_j>2\mu_j^2$ such that
\begin{\equ}\label{decay-c1-c2-j}
 2\mu_j^2\leq e^{(k-1)|t|}\bigl(u_j(t)^2+v_j(t)^2\bigr)\leq C_j
 \quad \text{for } t\in\R.
\end{\equ}
Moreover,
\[
 u_j(t)^2+v_j(t)^2
 =\bigl(C_j+o(1)\bigr)e^{-(k-1)|t|}
 \qquad\text{as } |t|\to\infty.
\]
\end{Thm}

For the ease of notations, in the shooting argument, we set $a=(k-1)/2$, $d=k+n-1$ and $b(t)=\cosh(t)^{-n/d}$. For the solution with $u_\mu(0)=v_\mu(0)=\mu>0$ write $R_\mu=u_\mu^2+v_\mu^2$. Then the Hamiltonian along this solution is
\[
 H_\mu(t)=-a u_\mu(t)v_\mu(t)+\frac{d}{2(d+1)} b(t)R_\mu(t)^{1+1/d}.
\]
Hamilton's equations \eqref{eqs2-2-1-1} and the product rule give
\begin{\equ}\label{product-shooting-identities}
\left\{
\aligned
&R_\mu'=-2a(u_\mu^2-v_\mu^2),\\[0.3em]
&(u_\mu v_\mu)'=b R_\mu^{1/d}(v_\mu^2-u_\mu^2),\\[0.3em]
&H_\mu'=-\frac{n}{2(d+1)}\tanh(t)b R_\mu^{1+1/d}.
\endaligned
\right.
\end{\equ}
In particular $H_\mu$ decreases strictly for $t>0$.

\begin{Lem}\label{ODE2-lemma1}
For each $\mu>0$, system \eqref{eqs2-2-1-1} has a unique solution
$(u_\mu,v_\mu)$ on $\R$ with $u_\mu(0)=v_\mu(0)=\mu$. For every $t\in\R$,
\[
 u_\mu(t)^2+v_\mu(t)^2>0,\qquad u_\mu(-t)=v_\mu(t).
\]
Each component has only simple zeros and hence finitely many zeros on
each compact interval. As $\nu\to\mu$, the solutions $(u_\nu,v_\nu)$
and their first derivatives with respect to $t$ converge uniformly on
each compact interval to $(u_\mu,v_\mu)$ and $(u_\mu',v_\mu')$, respectively.
\end{Lem}
\begin{proof}
Since $\mu>0$, the vector field in \eqref{eqs2-2-1-1} is smooth in a
neighborhood of the initial point $(\mu,\mu)$. Hence the initial value
problem has a unique local solution.
The first identity in \eqref{product-shooting-identities} gives
$-2aR_\mu\leq R_\mu'\leq 2aR_\mu$. Since $R_\mu(0)=2\mu^2$,
integration in either direction from $t=0$ yields
\begin{\equ}\label{R-bdd}
 2\mu^2 e^{-2a|t|}\leq R_\mu(t)\leq 2\mu^2 e^{2a|t|}.
\end{\equ}
In particular, the two components never vanish simultaneously. On each
bounded time interval these bounds keep $(u_\mu,v_\mu)$ in a compact
subset of $\R^2\setminus\{(0,0)\}$, where the vector field is smooth.
The solution therefore extends to all of $\R$.

Since $b$ is even, $(v_\mu(-t),u_\mu(-t))$ solves the same system
with the same initial data. Uniqueness gives the stated symmetry.
If $v_\mu(t_0)=0$, then $u_\mu(t_0)\ne0$ and the second equation gives
$v_\mu'(t_0)=-b(t_0)R_\mu(t_0)^{1/d}u_\mu(t_0)\ne0$. Thus every
zero of $v_\mu$ is simple, and the symmetry gives the same conclusion
for $u_\mu$. Consequently, each component has only finitely many zeros
on each compact interval. Indeed, if a component had infinitely many
zeros in such an interval, a sequence of distinct zeros would converge
to a point $t_*$. By continuity the component would
vanish at $t_*$. The difference quotients along this sequence would
then all be zero, contradicting the nonzero derivative at $t_*$.

Finally, continuous dependence on the initial data for ordinary
differential equations gives uniform convergence of $(u_\nu,v_\nu)$
to $(u_\mu,v_\mu)$ on each compact interval as $\nu\to\mu$.
Substitution into the two equations gives uniform convergence of
the derivatives on the same interval.
\end{proof}

By Lemma~\ref{ODE2-lemma1}, reflection across $t=0$
interchanges $u_\mu$ and $v_\mu$. Their values on $[0,+\infty)$ therefore
determine the solution on the whole $\R$, so we restrict the shooting
argument to this half-line.

\begin{Lem}\label{ODE2-lemma-exponential-decay}
Suppose that $H_\mu(t)\geq0$ for all $t\geq0$ and that $v_\mu$ has
only finitely many zeros on $[0,\infty)$. Then $e^{2at}R_\mu(t)$
is strictly increasing on $[0,\infty)$ and has a finite positive limit
as $t\to\infty$. Moreover, $u_\mu(t)v_\mu(t)>0$ for all sufficiently
large $t$, and $H_\mu(t)\to0$ as $t\to\infty$.
\end{Lem}
\begin{proof}
Choose $T\geq0$ so that $v_\mu$ has no zeros on $[T,\infty)$.
Changing the signs of both components if necessary, we may assume that
$v_\mu(t)>0$ for $t\geq T$. At any zero of $u_\mu$ in this interval,
the first equation in \eqref{eqs2-2-1-1} gives
$u_\mu'=bR_\mu^{1/d}v_\mu>0$. Hence once $u_\mu$ becomes positive,
it stays positive. Suppose instead that $u_\mu$ stays negative on
$[T,\infty)$. The second equation then gives $v_\mu'\geq a v_\mu$, so
\[
 v_\mu(t)\geq v_\mu(T)e^{a(t-T)}\qquad(t\geq T).
\]
Since $-a u_\mu v_\mu>0$, $R_\mu\geq v_\mu^2$ and
$b(t)=(\cosh t)^{-n/d}\geq e^{-nt/d}$ for $t\geq0$, it follows that
\[
 H_\mu(t)\geq \frac{d}{2(d+1)}b(t)v_\mu(t)^{2+2/d}
 \geq C\exp\!\left(\frac{2a(d+1)-n}{d}\,t\right)
 \qquad(t\geq T)
\]
for some $C>0$. Here $2a=k-1\geq1$ and $d+1=k+n$, so
$2a(d+1)-n\geq k>0$. This would force $H_\mu(t)\to\infty$,
whereas \eqref{product-shooting-identities} gives $H_\mu(t)\leq H_\mu(T)$
for all $t\geq T$. Thus we can choose $T_1\geq T$ so that
$u_\mu(t)v_\mu(t)>0$ for every $t\geq T_1$.

On $[T_1,\infty)$, the two equations in \eqref{eqs2-2-1-1} and the
definition of $H_\mu$ give
\[
 \frac{d}{dt}\ln\frac{u_\mu}{v_\mu}
 =-2a+\frac{bR_\mu^{1+1/d}}{u_\mu v_\mu}
 =\frac{2(d+1)}{d}\frac{H_\mu}{u_\mu v_\mu}+\frac{2a}{d}
 \geq\frac{2a}{d}.
\]
Integrating from $T_1$ to $t$ gives
\[
 \frac{u_\mu(t)}{v_\mu(t)}
 \geq \frac{u_\mu(T_1)}{v_\mu(T_1)}e^{2a(t-T_1)/d}.
\]
Both sides are positive, so taking reciprocals gives an upper bound
for $v_\mu(t)/u_\mu(t)$. Squaring this bound and using
$R_\mu(t)\geq u_\mu(t)^2$, we obtain
\begin{\equ}\label{v2R}
 0\leq\frac{v_\mu(t)^2}{R_\mu(t)}
 \leq\left(\frac{v_\mu(t)}{u_\mu(t)}\right)^2
 \leq\left(\frac{v_\mu(T_1)}{u_\mu(T_1)}\right)^2 e^{-4a(t-T_1)/d}.
\end{\equ}
Thus $v_\mu^2/R_\mu$ is integrable on $[T_1,\infty)$ and is bounded
by one on $[0,T_1]$. By the first identity in
\eqref{product-shooting-identities}, we have that
$(\ln(e^{2at}R_\mu))'=4a v_\mu^2/R_\mu$. And since $R_\mu(0)=2\mu^2$, it follows that
\[
 e^{2at}R_\mu(t)=2\mu^2
 \exp\!\left(4a\int_0^t\frac{v_\mu(s)^2}{R_\mu(s)}\,ds\right).
\]
The integral increases strictly because $v_\mu$ has only finitely many
zeros, and it converges as $t\to\infty$ by \eqref{v2R}.
This proves the asserted monotonicity and the existence of a finite
positive limit. In particular $R_\mu(t)\to0$. Since
$|u_\mu v_\mu|\leq R_\mu/2$ and $0<b(t)\leq1$, both terms in the
definition of $H_\mu(t)$ tend to zero as well.
\end{proof}

To compare the number of zeros of $v_\mu$ for nearby initial values
$\mu$, we need to know how many zeros can occur after a given time
$T\geq1$. At a zero of either component, the term $-a u_\mu v_\mu$
vanishes and the remaining term in $H_\mu$ is positive since $R_\mu>0$.
Thus, if $H_\mu(T)\leq0$, the monotonicity of $H_\mu$ rules out any
further zeros. The next lemma shows that if $H_\mu(T)$ is positive
but sufficiently small and the two components have the same sign at
$T$, then $v_\mu$ has at most one further zero.

\begin{Lem}\label{ODE2-lemma7}
There exists $\de>0$, depending only on $k$ and $n$, with the following property: If $T\geq1$, $u_\mu(T)v_\mu(T)>0$ and $H_\mu(T)\leq\de$, then $v_\mu$ has at most one zero on $(T,\infty)$.
\end{Lem}
\begin{proof}
Changing the signs of both components if necessary, we may assume that
$u_\mu(T)>0$ and $v_\mu(T)>0$. Suppose, for a contradiction, that
$v_\mu$ has at least two zeros after $T$, and denote the first two by
$t_1<t_2$. Both zeros are simple by Lemma~\ref{ODE2-lemma1}. Since
$v_\mu(T)>0$, it follows that $v_\mu'(t_1)<0$, that $v_\mu<0$ on
$(t_1,t_2)$, and that $v_\mu'(t_2)>0$. At a zero of $v_\mu$, the
second equation in \eqref{eqs2-2-1-1} gives
$v_\mu'=-bR_\mu^{1/d}u_\mu$. Thus $u_\mu(t_1)>0$ and $u_\mu(t_2)<0$.

Put $f_\mu=(u_\mu-v_\mu)/\sqrt2$ and $g_\mu=(u_\mu+v_\mu)/\sqrt2$.
In these variables
\[
 f_\mu'=(bR_\mu^{1/d}-a)g_\mu,\qquad g_\mu'=-(bR_\mu^{1/d}+a)f_\mu,\qquad
 H_\mu=\frac a2(f_\mu^2-g_\mu^2)+\frac{d}{2(d+1)}bR_\mu^{1+1/d}.
\]
Since $v_\mu$ vanishes at $t_1$ and $t_2$, we have
$f_\mu(t_1)=u_\mu(t_1)/\sqrt2>0$ and $f_\mu(t_2)=u_\mu(t_2)/\sqrt2<0$.
Let $\widehat T$ be the first zero of $f_\mu$ in $(t_1,t_2)$.
On $[t_1,\widehat T)$ we have $f_\mu>0$, so $g_\mu$ decreases strictly
from $g_\mu(t_1)=f_\mu(t_1)>0$. At $\widehat T$ the two components
$u_\mu$ and $v_\mu$ are equal. Since $v_\mu<0$ throughout $(t_1,t_2)$,
we obtain $g_\mu(\widehat T)=\sqrt2\,v_\mu(\widehat T)<0$.

At $t_2$ the product $u_\mu v_\mu$ vanishes and $R_\mu(t_2)>0$.
Since $H_\mu$ decreases and $\widehat T<t_2$, we have
\[
 H_\mu(\widehat T)\geq H_\mu(t_2)
 =\frac{d}{2(d+1)}b(t_2)R_\mu(t_2)^{1+1/d}>0.
\]
To use this energy inequality to locate $g_\mu(\widehat T)$, we choose
an auxiliary lower bound for $H_\mu$ that is exact at $\widehat T$
and whose negative root can be computed explicitly. One such choice
is obtained by fixing the time at $\widehat T$ in the Hamiltonian
and evaluating it at the points $(0,s)$ in the $(f_\mu,g_\mu)$-plane.
Here $s\in\R$ represents the $g_\mu$-coordinate. Since
$R_\mu=f_\mu^2+g_\mu^2$, the resulting function is
\[
 F_{\widehat T}(s)=-\frac a2 s^2+
       \frac{d}{2(d+1)}b(\widehat T)|s|^{2+2/d}.
\]
For $t_1\leq t\leq\widehat T$, we have $b(t)\geq b(\widehat T)$ and
$R_\mu(t)\geq g_\mu(t)^2$. Comparing the two energy expressions and
using $f_\mu(\widehat T)=0$ at the endpoint gives
\[
 F_{\widehat T}(g_\mu(t))\leq H_\mu(t),\qquad
 F_{\widehat T}(g_\mu(\widehat T))=H_\mu(\widehat T)>0.
\]
For $s<0$, it follows that $F_{\widehat T}(s)>0$
precisely when $s<-G$, where
\[
 G=\left(\frac{a(d+1)}{d\,b(\widehat T)}\right)^{d/2}.
\]
Since $g_\mu(\widehat T)<0$, the endpoint identity yields
$g_\mu(\widehat T)<-G$. The minimum of $F_{\widehat T}$ on the negative half-line occurs at $-\alpha G$, where $\alpha=(d/(d+1))^{d/2}$. Figure~\ref{fig-Fs} shows the portion of the graph between its minimum and its negative root.

\begin{figure}[h!]
	\centering
	\captionsetup{singlelinecheck=false, justification=centering}
	\begin{minipage}{0.55\linewidth}
\begin{tikzpicture}[scale=1.05]
	\begin{axis}[
		xlabel={$s$},
		ylabel={$F_{\widehat T}(s)$},
		xmin=-10, xmax=1.5,
		ymin=-3, ymax=5,
		grid=both,
		major grid style={gray!50},
		minor grid style={gray!25},
		axis lines=middle,
		samples=200,
		xtick=\empty, 
		ytick=\empty,
		xlabel style={at={(ticklabel cs:1)}, anchor=north}, 
		ylabel style={at={(ticklabel cs:1)}, anchor=east}, 
		clip=false,
		axis line style={black, thick}, 
		]
		\node at (axis cs:0.4,0.4) {$O$};
		\addplot[blue, thick, domain=-10:0] {-0.03*x^3 - 0.25*x^2};
		\addplot[black, dashed, thick] coordinates {(-9,0) (-9,1.62)};
        \addplot[only marks, mark=*, black] coordinates {(-9,0) (-8.3333,0)};
        \node[anchor=east,blue] at (axis cs:-9.1,0.8) {$g_\mu(\widehat T)$};
        \node[anchor=south west,blue] at (axis cs:-8.35,0.05) {$-G$};
        \addplot[red, dashed, thick] coordinates {(-5.5556,0) (-5.5556,-2.572)};
        \addplot[only marks, mark=*, red] coordinates {(-5.5556,0) (-5.5556,-2.572)};
        \node[anchor=south,blue] at (axis cs:-5.5556,0.05) {$-\alpha G$};
    \end{axis}
\end{tikzpicture}
\caption{The graph of $(s,F_{\widehat T}(s))$ for $s\leq0$. The energy estimate uses the interval $[-G,-\alpha G]$.}\label{fig-Fs}
\end{minipage}
\end{figure}
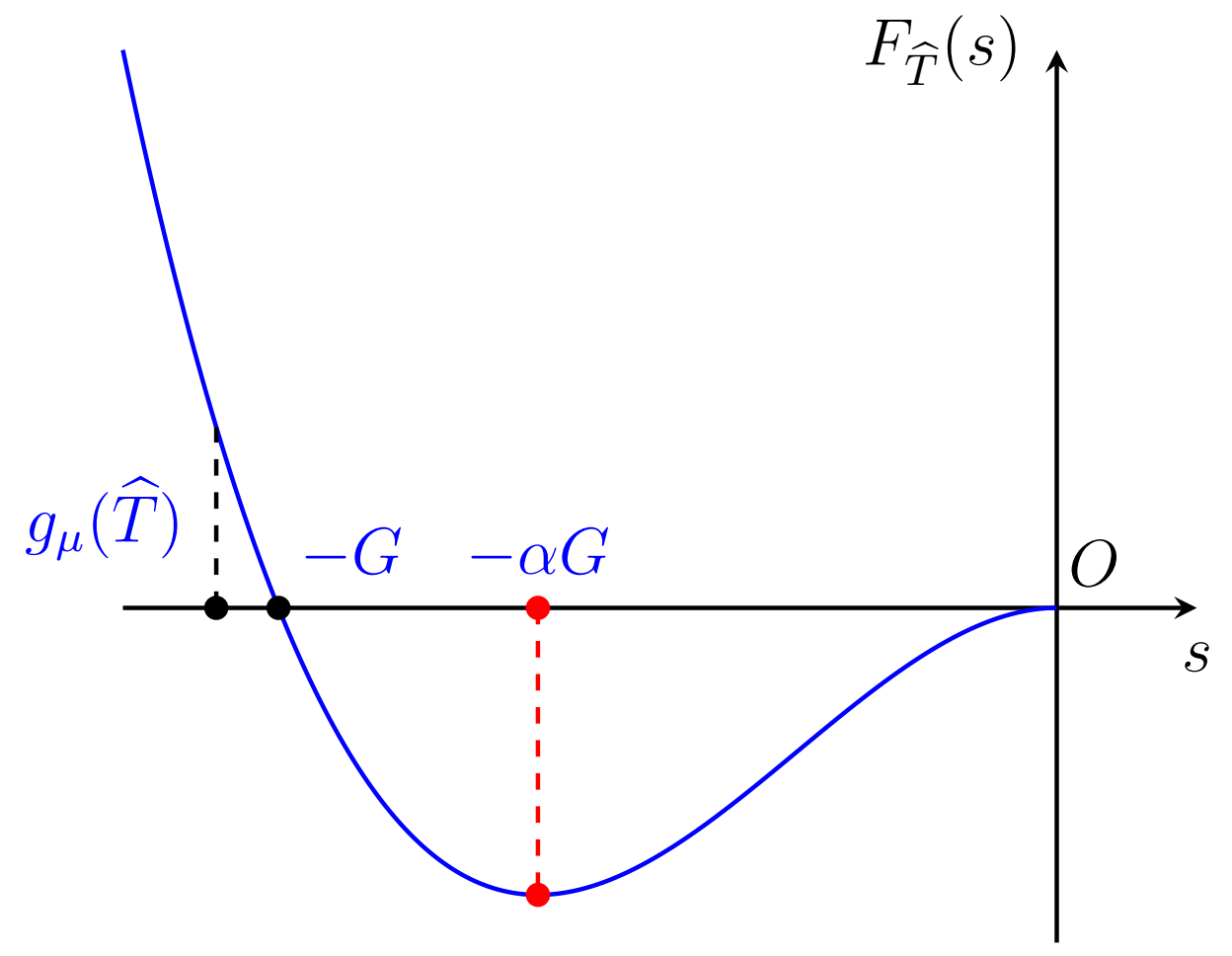

On the interval where $g_\mu$ decreases from $-\alpha G$ to $-G$, the
bound $H_\mu\leq\de\leq1$ gives $f_\mu^2\leq g_\mu^2+2/a$. Since
$G\geq(a(d+1)/d)^{d/2}$ and $|g_\mu|\leq G$, it follows that
$R_\mu=f_\mu^2+g_\mu^2$ is bounded by a constant multiple of $G^2$
depending only on $k$ and $n$. Moreover,
$b(t)G^{2/d}\geq b(\widehat T)G^{2/d}=a(d+1)/d$.
Using $|f_\mu|\leq R_\mu^{1/2}$ in the equation for $g_\mu'$, we may
therefore fix $C=C(k,n)>0$ such that
\[
 |g_\mu'|\leq C b(t)G^{1+2/d}.
\]
To estimate $-H_\mu'$, use $\tanh(t)\geq\tanh(1)$ and
$R_\mu\geq g_\mu^2\geq\alpha^2G^2$ in \eqref{product-shooting-identities}.
With $c=\frac{n}{2(d+1)}\tanh(1)\alpha^{2+2/d}$ and $c'=c/C$, this gives
\[
 -H_\mu'\geq c\,b(t)G^{2+2/d}\geq c'G|g_\mu'|.
\]
Both $c$ and $c'$ are positive and depend only on $k$ and $n$.
Integrating over this interval gives an energy loss of at least
$c'(1-\alpha)G^2$. Since $H_\mu$ decreases on $[T,\widehat T]$, we obtain
\[
 H_\mu(T)-H_\mu(\widehat T)
 \geq c'(1-\alpha)\left(\frac{a(d+1)}{d}\right)^d.
\]
We can thus choose $\de\leq1$ smaller than this positive lower bound,
independently of $\mu$, $T$ and $\widehat T$. This contradicts
$0<H_\mu(\widehat T)\leq H_\mu(T)\leq\de$ and proves the lemma.
\end{proof}

\begin{Lem}\label{ODE2-lemma8}
For every $\mu>0$ the number of zeros of $v_\mu$ on $(0,\infty)$ is finite. If this number is $j$, all sufficiently close initial values yield at most $j+1$ zeros. When $H_\mu$ takes a negative value, the number of zeros is locally constant.
\end{Lem}
\begin{proof}
Fix $\mu>0$. Suppose that $v_\mu$ has infinitely many zeros on
$(0,\infty)$. By Lemma~\ref{ODE2-lemma1}, only finitely many can lie
in any compact interval, so we may list them as $t_1<t_2<\cdots$ with
$t_i\to\infty$. At each $t_i$ the product $u_\mu v_\mu$ vanishes and
$R_\mu(t_i)>0$, so $H_\mu(t_i)>0$. Since $H_\mu$ is decreasing,
it stays positive and has a limit $H_\infty\geq0$ as $t\to\infty$.

We claim that $H_\infty=0$. Integrating the identity for $H_\mu'$ in
\eqref{product-shooting-identities} and using $\tanh(t)\geq\tanh(1)$
for $t\geq1$, we obtain
\begin{\equ}\label{bR-integrable}
 \int_1^\infty b(t)R_\mu(t)^{1+1/d}\,dt
 \leq \frac{2(d+1)}{n\tanh(1)}\bigl(H_\mu(1)-H_\infty\bigr)<\infty.
\end{\equ}
Also, \eqref{product-shooting-identities} gives that $|(u_\mu v_\mu)'|\leq bR_\mu^{1+1/d}$.
Thus $u_\mu v_\mu$ has a finite limit at infinity, and that limit is
zero because $v_\mu(t_i)=0$ for every $i$. The definition of $H_\mu$
now gives $bR_\mu^{1+1/d}\to 2(d+1)H_\infty/d$. A positive limit would
contradict \eqref{bR-integrable}, so $H_\mu(t)\to0$ as $t\to\infty$.

Let $\de>0$ be the constant in Lemma~\ref{ODE2-lemma7}. At each zero
$t_i$, we have
\[
 (u_\mu v_\mu)'(t_i)=-b(t_i)R_\mu(t_i)^{1/d}u_\mu(t_i)^2<0.
\]
Hence $u_\mu v_\mu$ is positive immediately to the left of $t_i$.
For $i$ sufficiently large we can therefore choose $T\geq1$ just before
$t_i$ with $u_\mu(T)v_\mu(T)>0$ and $H_\mu(T)<\de$.
Lemma~\ref{ODE2-lemma7} then allows at most one zero after $T$,
although both $t_i$ and $t_{i+1}$ lie after $T$. This proves finiteness.

Let $j$ be the number of zeros of $v_\mu$ on $(0,\infty)$. For any
$T\geq1$ with all $j$ zeros in $(0,T)$ and $v_\mu(T)\ne0$, the number
of zeros in $[0,T]$ remains $j$ for all initial values $\nu$
sufficiently close to $\mu$. Indeed, each zero is simple, so we can
enclose the zeros in disjoint closed intervals inside $(0,T)$ on each
of which $v_\mu'$ is bounded away from zero and $v_\mu$ has opposite
signs at the endpoints. By Lemma~\ref{ODE2-lemma1}, $v_\nu$ and $v_\nu'$
converge uniformly on $[0,T]$ to $v_\mu$ and $v_\mu'$, respectively.
Hence $v_\nu$ has exactly one zero in each of these intervals.
On the compact set outside their interiors, $|v_\mu|$ has a positive
minimum, so $v_\nu$ has no other zeros in $[0,T]$.

If $H_\mu$ takes a negative value, choose $T$ as above large enough that
$H_\mu(T)<0$. For $\nu$ close to $\mu$, continuous dependence gives
$H_\nu(T)<0$ as well. The energy then remains negative for $t\geq T$,
whereas it is positive at any zero of $v_\nu$. Thus $v_\nu$ has
exactly $j$ zeros on $(0,\infty)$.

If $H_\mu(t)\geq0$ for every $t\geq0$, we can now apply
Lemma~\ref{ODE2-lemma-exponential-decay}, since finiteness has already
been proved. It gives $H_\mu(t)\to0$ and $u_\mu(t)v_\mu(t)>0$ for all
sufficiently large $t$. Choose $T$ beyond all the zeros of $v_\mu$ so
that $T\geq1$, $H_\mu(T)<\de$ and $u_\mu(T)v_\mu(T)>0$. These
strict inequalities persist for $\nu$ sufficiently close to $\mu$.
Lemma~\ref{ODE2-lemma7} gives at most one zero of $v_\nu$ on
$(T,\infty)$ in addition to its $j$ zeros on $[0,T]$, as claimed.
\end{proof}

\begin{proof}[Proof of Theorem~\ref{main thm ODE2}]
Fix an integer $j\geq0$ and let $S_j$ be the set of initial values
$\mu>0$ for which $v_\mu$ has at most $j$ zeros on $(0,\infty)$. Evidently, we have $S_0\subset S_1\subset S_2\subset\cdots$.
For all sufficiently small $\mu>0$, the initial energy satisfies
\[
 H_\mu(0)=-a\mu^2+\frac{d}{2(d+1)}(2\mu^2)^{1+1/d}<0.
\]
Since $H_\mu$ decreases and is positive at every zero of $v_\mu$,
the corresponding $v_\mu$ has no positive zeros. Thus $S_0$ is nonempty and so is $S_j$, $j\geq1$.
Moreover, $S_j$ is closed in $(0,\infty)$. Indeed, if $\mu\notin S_j$,
choose disjoint intervals around the first $j+1$ zeros of $v_\mu$.
By Lemma~\ref{ODE2-lemma1} the zeros are simple, so the intervals can
be chosen with opposite signs of $v_\mu$ at their endpoints. These signs
persist for all nearby initial values, which therefore also lie outside $S_j$.

To show that $S_j$ is bounded above, it suffices to prove that for
every sufficiently large $\mu$, the function $v_\mu$ has at least
$j+1$ positive zeros. And, of course, it is clear that the existence of $j+1$ zeros in any fixed interval
$(0,T)$ would immediately exclude $\mu$ from $S_j$. To proceed, let us take $T=1$. 
On $[0,1]$ we have
$b(t)\geq b(1)>0$ and, by the estimate in \eqref{R-bdd}, $R_\mu(t)\geq 2\mu^2e^{-2a}$.
To count the zeros, write the solution in polar coordinates with a
continuous angle $\theta_\mu$ chosen so that
\[
 (u_\mu,v_\mu)=\sqrt{R_\mu}\,(\cos\theta_\mu,\sin\theta_\mu) \quad \text{and} \quad
 \theta_\mu(0)=\pi/4.
\]
By \eqref{eqs2-2-1-1}, we deduce that
\[
 \theta_\mu'
 =\frac{u_\mu v_\mu'-v_\mu u_\mu'}{R_\mu}
 =a\sin(2\theta_\mu)-b(t)R_\mu^{1/d}.
\]
The two lower bounds on $[0,1]$ give
$\theta_\mu'\leq a-b(1)(2\mu^2e^{-2a})^{1/d}$. For all sufficiently
large $\mu$ this upper bound is less than $-(j+1)\pi$. Starting from $t=0$, the angle $\theta_\mu(t)$ then passes through $0,-\pi,\dots,-j\pi$ before
$t=1$. Each of these $j+1$ angles gives a distinct zero of $v_\mu$,
so every sufficiently large $\mu$ has to lie outside $S_j$.

Since $S_j$ is nonempty, bounded above and closed in $(0,\infty)$,
its supremum is positive and belongs to $S_j$. Set $\mu_j=\max S_j$
and write $(u_j,v_j)=(u_{\mu_j},v_{\mu_j})$. If $v_j$ had $q<j$
zeros, Lemma~\ref{ODE2-lemma8} would give at most $q+1\leq j$ zeros
for all initial values sufficiently close to $\mu_j$, including some
larger than $\mu_j$. This contradicts the choice of $\mu_j$, so $v_j$
has exactly $j$ zeros. If $H_{\mu_j}$ took a negative value on
$[0,\infty)$, the same lemma would give exactly $j$ zeros for all
nearby initial values, again contradicting maximality. Thus
$H_{\mu_j}(t)\geq0$ for every $t\geq0$.

Lemma~\ref{ODE2-lemma8} now gives at most $j+1$ zeros for initial
values close to $\mu_j$. In particular, $S_{j+1}$ contains values
larger than $\mu_j$, so $\mu_{j+1}>\mu_j$. To see that $\mu_j\to\infty$,
fix any $\mu>0$. By the finiteness assertion of Lemma~\ref{ODE2-lemma8},
$v_\mu$ has some finite number $q$ of zeros. Hence $\mu\in S_j$
and $\mu_j\geq\mu$ for every $j\geq q$. Since $\mu$ was arbitrary,
the assertion follows.

Finally, Lemma~\ref{ODE2-lemma-exponential-decay} applies to
$(u_j,v_j)$. It shows that $e^{2at}R_{\mu_j}(t)$ increases strictly
from $R_{\mu_j}(0)=2\mu_j^2$ to a finite limit $C_j>2\mu_j^2$ as
$t\to\infty$. By Lemma~\ref{ODE2-lemma1}, $u_j(-t)=v_j(t)$, so
$R_{\mu_j}(-t)=R_{\mu_j}(t)$. Since $2a=k-1$, this proves
\eqref{decay-c1-c2-j} on $\R$ and gives
\[
 \lim_{|t|\to\infty}e^{(k-1)|t|}\bigl(u_j(t)^2+v_j(t)^2\bigr)=C_j,
\]
which completes the whole proof.
\end{proof}

Starting from the solutions $(u_j,v_j)$ of Theorem~\ref{main thm ODE2},
we follow Section~\ref{subsubsec:MSN} to construct solutions of the
spinorial Yamabe equation on $S^k\times N$ and determine their maximum
sets after a fixed conformal change of the product metric.
We use the product metric $\ig=\ig_{S^k}\op\ig_N$ and the unit harmonic
spinor $\va_N$ from that construction. In particular, we may take $N$
to be a flat torus with the trivial spin structure in any dimension $n\geq1$.
Reversing the Emden-Fowler change of variables and stereographic
projection gives a spinor $\psi_j$ on $S^k$ away from the poles.
Tensoring $\psi_j$ with $\va_N$ as in the product construction gives
a spinor $\phi_j$ of squared length
\[
 |\phi_j|_\ig^2=(\cosh t)^{k-1}\bigl(u_j(t)^2+v_j(t)^2\bigr),
 \qquad r=e^{-t}.
\]
The upper bound in \eqref{decay-c1-c2-j} keeps this length bounded
near both poles, so the removability argument in Section~\ref{subsubsec:MSN}
gives a solution of \eqref{SYSkN} on all of $S^k\times N$.
The lower bound in \eqref{decay-c1-c2-j} and
$\cosh t\geq e^{|t|}/2$ give, by continuity at the poles,
\[
 \min_{S^k\times N}|\phi_j|_\ig^2\geq 2^{2-k}\mu_j^2.
\]
Thus each $\phi_j$ is nowhere vanishing, and its minimum length tends
to infinity as $j\to\infty$ because $\mu_j\to\infty$.

To make the length attain its maximum precisely over the two poles,
let $\theta$ be the polar angle measured from the north pole of $S^k$.
For stereographic projection from this pole, $r=\cot(\theta/2)$, so
$\sin\theta=2r/(1+r^2)=1/\cosh t$. We now equip $S^k\times N$ with
the fixed metric
\[
 \widehat\ig=(1+\sin^2\theta)^{2(k-1)/(k+n-1)}\ig.
\]
By Proposition~\ref{prop:conformal-invar} in dimension $k+n$, the
corresponding solution $\widehat\phi_j$ has length
$|\widehat\phi_j|_{\widehat\ig}=(1+\sin^2\theta)^{-(k-1)/2}|\phi_j|_\ig$.
Since the factor is positive everywhere, $\widehat\phi_j$ has the same
zero set as $\phi_j$, which is empty. Substituting the formulas for $\sin\theta$ and
$|\phi_j|_\ig$ gives
\[
 |\widehat\phi_j|_{\widehat\ig}^2
 =2^{1-k}e^{(k-1)|t|}\bigl(u_j(t)^2+v_j(t)^2\bigr)
 \left(\frac{1+e^{-2|t|}}{1+(\cosh t)^{-2}}\right)^{k-1}.
\]
The fraction is at most one because $\cosh t\leq e^{|t|}$, and it tends
to one as $|t|\to\infty$. As shown in the proof of
Theorem~\ref{main thm ODE2}, $e^{(k-1)|t|}(u_j(t)^2+v_j(t)^2)$
increases strictly to $C_j$ as $|t|$ increases. Hence
$|\widehat\phi_j|_{\widehat\ig}^2$ is strictly less than $2^{1-k}C_j$
away from the poles and equals $2^{1-k}C_j$ over each pole.
For every $j$, the maximum set $\co_j$ therefore consists of the two
fibers over the poles, each a copy of $N$ of codimension $k$ in $S^k\times N$.

Although the maximum set is the same for every $j$, the pointwise
lengths of two distinct members of this sequence are never proportional.
Indeed, their squared lengths have ratio $R_{\mu_i}/R_{\mu_j}$, since
the metric factors in the length formula cancel. To see what this
ratio records, use \eqref{product-shooting-identities} and $2a=k-1$
to write
\[
 \frac{d}{dt}\bigl(e^{(k-1)t}R_{\mu_j}(t)\bigr)
 =2(k-1)e^{(k-1)t}v_j(t)^2.
\]
By Theorem~\ref{main thm ODE2}, the derivative on the left vanishes
at exactly $j$ points of $(0,\infty)$. If $R_{\mu_i}$ were a constant
multiple of $R_{\mu_j}$, the corresponding derivatives would vanish
at the same points, forcing $i=j$. By
Proposition~\ref{prop:conformal-invar}, a further conformal change
multiplies all the lengths by the same positive function. It therefore
preserves their ratios and cannot make two distinct members of the
sequence simultaneously of constant length.

To obtain the example on $S^k\times S^3$ in the third row of
Table~\ref{tab:noncompactness-examples}, it remains to exhibit a unit
harmonic spinor on the second factor. Equip $S^3$ with the Berger
metric obtained from the unit round metric by multiplying lengths
along the Hopf fibers by $4$ and keeping the horizontal metric fixed.
For this metric the Dirac kernel has complex dimension two, see
\cite[Proposition~3.2 and p.~37]{Hit74}. Let $\eta$ be an
$L^2$-normalized harmonic spinor. The quaternionic structure $J$ on
spinors in dimension three commutes with the Dirac operator, so
$\eta,J\eta$ form an orthonormal basis of the kernel. The pointwise
sum $|\eta|_{\ig_N}^2+|J\eta|_{\ig_N}^2$ is independent of the
choice of orthonormal basis. It is therefore invariant under left
translations, which preserve the metric and act transitively on $S^3$.
Since $J$ preserves pointwise length, this sum equals
$2|\eta|_{\ig_N}^2$. Thus $\eta$ has positive constant length and
can be rescaled to unit length. The construction above now applies
for every $k\geq2$ and gives the two copies of $S^3$ as the common
maximum set after the fixed conformal change.

On the torus $\T^m$, constant-length eigenspinors give another way to
prescribe the maximum sets in Theorem~\ref{main-thm}~(1).

\begin{Cor}\label{cor:prescribed-maximum-set}
Let $m\geq2$ and $0\leq k\leq m-1$. On the torus $\T^m$ there is a
smooth metric and a sequence of smooth solutions to the spinorial Yamabe equation \eqref{SY0} whose
minimum lengths tend to infinity and whose maximum sets are all the
same coordinate subtorus of codimension $k$. One may choose the metric
locally conformally flat in every dimension, or non-locally conformally
flat if $m\geq3$.
\end{Cor}
\begin{proof}
Use angular coordinates $x^1,\ldots,x^m$ of period $2\pi$ and the
trivial spin structure. On the flat torus, we put $t=x^1$ for simplicity. 
If $\chi$ is a constant unit spinor with $i e_1\cdot\chi=\chi$, the spinor
$\phi_j=e^{ijt}\chi$ has unit length and Dirac eigenvalue $j$ for
every positive integer $j$.

For $m\geq3$, we now construct a non-locally conformally flat metric
on $\T^m$ for which the spinors represented by $e^{ijt}\chi$ in the
orthonormal frame below still have Dirac eigenvalue $j$.
Set $a(t)=\sin^3 t$ and consider the metric
\[
 \ig_a=dt^2+e^{2a(t)}(dx^2)^2+e^{-2a(t)}(dx^3)^2+
       \sum_{\nu=4}^m(dx^\nu)^2 .
\]
Take the orthonormal frame $e_1=\partial_t$, $e_2=e^{-a(t)}\partial_{x^2}$,
$e_3=e^{a(t)}\partial_{x^3}$, and $e_\nu=\partial_{x^\nu}$ for $\nu\geq4$.
The Levi-Civita connection of $\ig_a$ satisfies
$\nabla_{e_2}e_1=a'(t)e_2$ and $\nabla_{e_3}e_1=-a'(t)e_3$.
For a spinor $\phi$ whose components in the lifted spin frame depend
only on $t$, the local formula for the spin connection
\cite[Chapter~II, Theorem~4.14]{LM} therefore gives
\[
 e_2\cdot_{\ig_a}\nabla^{\mbs}_{e_2}\phi
   =\frac12 a'(t)e_1\cdot_{\ig_a}\phi,\qquad
 e_3\cdot_{\ig_a}\nabla^{\mbs}_{e_3}\phi
   =-\frac12 a'(t)e_1\cdot_{\ig_a}\phi.
\]
These two terms cancel in $D_{\ig_a}\phi$. The frame is parallel in
the $e_1$ direction and in each $e_\nu$ direction for $\nu\geq4$,
so $D_{\ig_a}\phi=e_1\cdot_{\ig_a}\partial_t\phi$. In particular,
$D_{\ig_a}\phi_j=j\phi_j$ for the spinor represented by $e^{ijt}\chi$.
Write $K_{\ig_a}(e_i,e_j)$ for the sectional curvature of the plane
spanned by $e_i$ and $e_j$. In the orthonormal frame above,
\[
 K_{\ig_a}(e_1,e_2)=-a''-(a')^2,\qquad
 K_{\ig_a}(e_1,e_3)=a''-(a')^2,\qquad
 K_{\ig_a}(e_2,e_3)=(a')^2.
\]
All other planes spanned by pairs of frame vectors have zero sectional
curvature. Since $a(t)=\sin^3 t$, we have $a'(\pi/2)=0$ and
$a''(\pi/2)=-3$. Hence, for $m\geq4$, the Weyl tensor of $\ig_a$
satisfies $|W_{\ig_a}(e_1,e_2,e_2,e_1)|=3(m-3)/(m-2)>0$
at every point with $t=\pi/2$.
For $m=3$, the Schouten tensor of $\ig_a$ is
$P=\operatorname{Ric}_{\ig_a}-\frac14\operatorname{Scal}_{\ig_a}\,\ig_a$,
where $\operatorname{Ric}_{\ig_a}$ and $\operatorname{Scal}_{\ig_a}$
denote the Ricci tensor and scalar curvature of $\ig_a$, respectively.
At $t=0$ we have $a'(0)=a''(0)=0$ and
$(\nabla_1P)_{22}-(\nabla_2P)_{12}=-a'''(0)\ne0$, so the Cotton
tensor is nonzero.

Let $\ig$ denote either the flat metric on $\T^m$ or, when $m\geq3$,
the metric $\ig_a$ constructed above. In both cases, $\phi_j$ has unit
length and Dirac eigenvalue $j$, so $\psi_j=j^{(m-1)/2}\phi_j$
solves \eqref{SY0}.
Put $u(x)=\sum_{\nu=1}^k(1-\cos x^\nu)$, with $u=0$ if $k=0$,
and equip the torus with the fixed metric $\widehat \ig=e^{2u}\ig$.
By Proposition~\ref{prop:conformal-invar}, the corresponding solutions
have length $j^{(m-1)/2}e^{-(m-1)u(x)/2}$. Their minimum tends to
infinity, and their maximum is attained precisely where $u=0$, namely
on $\{x^1=\cdots=x^k=0\}$. The resulting metric $\widehat\ig$ is
locally conformally flat when $\ig$ is flat and non-locally conformally
flat when $\ig=\ig_a$ with $m\geq3$.
\end{proof}

\section{Noncompactness on $S^m$ with a non-flat metric}\label{sec: perturbation}

For each $m\geq20$, we construct a smooth metric on $S^m$ for which the
spinorial Yamabe equation \eqref{SY-m-sphere} admits a sequence of
{smooth, nowhere-vanishing solutions} whose supremum norms tend to infinity.
The metric is obtained from localized Weyl-type perturbations of the Euclidean metric and conformal compactification.
We adapt the abstract perturbation framework of Isobe and Xu~\cite{Isobe-Xu23},
who obtained multiple solutions by placing finitely many perturbations
far apart and seeking a solution near a standard bubble at each center.
Here we use infinitely many centers and prescribe bubble scales tending
to zero. To obtain all these solutions for a single metric, we estimate,
at each scale, the combined effect of the other perturbations on the
bubble near the chosen center. The same construction applies at a
positive local maximum or a negative local minimum of the reduced
energy. In Example~\ref{Ex-2} both occur for one choice of metric.

	The following setting is due to Brendle \cite{Brendle}. We fix $W=\{W_{ijkl}\}: (\R^m)^4\to\R$ be a multi-linear form such that its tensor norm 
	\[
	|W|=\Big( \sum_{i,j,k,l=1}^m (W_{ijkl}+W_{ilkj})^2 \Big)^{1/2}
	\]
	is strictly positive and it satisfies all algebraic properties of the Weyl tensor: $W_{ijkl}=-W_{jikl}=-W_{ijlk}=W_{klij}$ (symmetry and anti-symmetry), $W_{ijkl}+W_{iklj}+W_{iljk}=0$ (the Bianchi identity) and  $\sum_{i=1}^mW_{ijil}=0$ (trace free). For abbreviation, we set (under the Einstein’s summation convention)
   \[
   H_{ij}(x)= W_{i\al j\bt}x^\al x^\bt
   \]
   for $x=(x^1,x^2,\dots,x^m)\in\R^m$. Given $\vr,a,b>0$ (where these values are small), we consider a Riemannian metric of the form 
   \begin{\equ}\label{perturb-g}
   \tilde \ig(x)=\tilde\ig_{\vr,a,b}(x)=\exp(\vr h(x))=I +\vr\ih(x)+\frac{\vr^2\ih^2(x)}{2!}+\vr^3\ih(x)\rr(\vr,\ih(x))
   \end{\equ}
   defined on $\R^m$. Here $h(x)$ is a trace-free symmetric two-tensor given by
   \begin{\equ}\label{hij}
   h_{ij}(x)=a\eta(r^{-1}|x|)f(b^{-2}|x|^2)H_{ij}(x)
   \end{\equ}
where we choose the cubic polynomial
\begin{\equ}\label{rev-cubic-profile}
f(t)=1-\frac9{14}t+\frac9{91}t^2-\frac2{637}t^3.
\end{\equ}
The cut-off function $\eta\in C^\infty(\R,[0,1])$ satisfies
   \[
   \eta(t)=\begin{cases}
   	1 & \text{if } t\in(-\frac12,\frac12), \\[0.3em]
   	0 & \text{if } t\in\R\setminus(-1,1),
   \end{cases}
   \]
   and $\rr(\vr,A)$ is defined for an arbitrary symmetric two-tensor $A$ as 
   \[
   \rr(\vr,A)=\frac{A^2}{3!}+\frac{\vr A^3}{4!}+\frac{\vr^2 A^4}{5!}+\cdots=\sum_{k=1}^\infty \frac{\vr^{k-1}A^{k+1}}{(k+2)!}.
   \] 
   The parameters $\vr,a,b,r$ will be determined later. And
we impose $0<a<b^6<b^2<b<r\leq1$.
The cut-off radius $r$ will remain fixed. For every fixed $k$,
\[
\|h\|_{C^k(\R^m)}\leq C_k a b^{-6},
\]
Here $C_k$ depends only on $k,r,\eta,W$ and the polynomial $f$ in \eqref{rev-cubic-profile}.
The Weyl identities give $\tr h=0$, $\diver h=0$ and $h(x)x=0$.

\begin{Rem}\label{rem:cubic-coefficient-range}
To explain the choice of $f$ in \eqref{rev-cubic-profile},
consider a bubble of scale $\lm$ centered at the origin. We keep the
cut-off radius $r$ fixed and let $b\to0$. In the coordinates $y=x/b$,
the cut-off radius $r/b$ tends to infinity, while the bubble scale
$\lm/b$ must tend to a finite positive limit. We choose $f$ so that
the leading energy correction, as a function of the relative scale,
has a strict local extremum of nonzero value in this limit. It appears natural to
try a linear or quadratic $f$ first. Such an extremum requires
$m\geq37$ in the linear case and $m\geq23$ in the quadratic case.
The cubic in \eqref{rev-cubic-profile} gives such an extremum in dimension $20$, and
increasing the degree cannot lower this dimension. We explain the
choice of coefficients in Subsection~\ref{rev-cubic-subsection} and
compare the degrees in Subsection~\ref{rev-other-choices}.

Roughly speaking, the three coefficients in
\eqref{rev-cubic-profile} multiplying $t$, $t^2$, and $t^3$ can be
varied independently while the constant term stays equal to one.
For every $m\geq20$, the strict inequalities used in the construction
continue to hold for any polynomial
$1+\sigma_1t+\sigma_2t^2+\sigma_3t^3$ satisfying
\[
 \max\left\{
 \left|\sigma_1+\frac9{14}\right|,
 \left|\sigma_2-\frac9{91}\right|,
 \left|\sigma_3+\frac2{637}\right|
 \right\}\leq\frac1{7000}.
\]
For example, rounding the three coefficients
to three decimal places gives the admissible polynomial
$1-0.643t+0.099t^2-0.003t^3$.
Appendix~\ref{app:cubic-perturbations} shows why the bound $1/7000$
suffices to preserve the required strict inequalities uniformly for $m\geq20$.
\end{Rem}

   \subsection{Some basic materials}\label{subsec:bg-identification}
   
   Here, we will briefly review the Bourguignon-Gauduchon construction of a suitable local trivialization for spinor bundles, which is particularly useful for our purposes. For more details, we refer readers to the original paper \cite{BG} by Bourguignon and Gauduchon, as well as \cite{ADHH, AGHM, Isobe11, Ba-Xu-JFA21, Ba-Xu-CV22} for its various applications. 
   
   Let $\ig_1$ and $\ig_2$ be Riemannian metrics on the spin $m$-manifold $M$. There exists a unique bundle endomorphism $\mft:(TM,\ig_2)\to(TM,\ig_1)$ such that $\ig_1(\mft X,Y)=\ig_2(X,Y)$ for all $X,Y\in TM$. It is $\ig_1$-symmetric and positive definite. Set $B:=(\mft)^{-1/2}$, the inverse of $(\mft)^{1/2}$, where $(\mft)^{1/2}$ is the unique positive definite square root of $\mft$. It follows that $B$ maps $\ig_1$-orthonormal frames to $\ig_2$-orthonormal frames. This, in other words, gives rise to an ${\SO(m)}$-equivariant bundle morphism $B: {P_{\SO}(M,\ig_1)}\to {P_{\SO}(M,\ig_2)}$ of the principal bundles of orthonormal frames. The map $B$ lifts to a ${\Spin(m)}$-equivariant bundle morphism $\mfs: {P_{\Spin}(M,\ig_1)}\to {P_{\Spin}(M,\ig_2)}$ of the corresponding spin structures. And this induces an isometry of vector bundles
   \[
   \mfj: \mbs(M,\ig_1)\to \mbs(M,\ig_2)
   \]
   which is a fiberwise isometry. In particular, any spinor $\va\in\Ga(\mbs(M,\ig_2))$ can be expressed in terms of a spinor $\psi\in\Ga(\mbs(M,\ig_1))$ so that $\va=\mfj(\psi)$, and $B(X)\cdot_{\ig_2} \mfj(\psi)=\mfj(X\cdot_{\ig_1}\psi)$ for all $X\in TM$ and $\psi\in \Ga(\mbs(M,\ig_1))$, where ``$\,\cdot_{\ig_1}$'' and ``$\,\cdot_{\ig_2}$" are Clifford multiplication with respect to $\ig_1$ and $\ig_2$.
   
   We remark that the spin connections $\nabla^{\ig_1}$ and $\nabla^{\ig_2}$ on $\mbs(M,\ig_1)$ and $\mbs(M,\ig_2)$ are different, as they depend on the respective metrics (see e.g. \cite[Theorem 4.14]{LM}). As a result, the map $\mfj$ does not preserve the Dirac operators. The core idea of the Bourguignon-Gauduchon trivialization is to transport the $\ig_2$-Dirac operator onto the $\ig_1$-spinor bundle via $\mfj$, so that the discrepancy between the two operators can be analyzed on a single vector bundle. More precisely, we can define the transported Dirac operator $\widetilde D_{\ig_2}$ acting on sections of $\mbs(M,\ig_1)$ by 
   \[
   \widetilde D_{\ig_2}= \mfj^{-1}\circ D_{\ig_2}\circ\mfj.
   \]
   Hence, an explicit expression for $\widetilde D_{\ig_2}$ can be formulated locally in terms of a $\ig_1$-orthonormal frame $\{e_k\}_{k=1}^m$ as
   \[
   \widetilde D_{\ig_2}\psi=\mfj^{-1}\Big( \sum_{k=1}^m B(e_k) \cdot_{\ig_2} \nabla^{\ig_2}_{B(e_k)}\mfj(\psi)\Big).
   \]
   
   In our settings, let us consider $M=\R^m$ with $\ig_1=\ig_{\R^m}$ being the Euclidean metric and $\ig_2=\tilde\ig$ the perturbed metric defined by \eqref{perturb-g}. Let $x=(x^1,\dots,x^m)$ be the $\ig_1$-normal coordinate chart on $\R^m$, then the metrics $\tilde\ig$ can be written as $\tilde\ig=\tilde\ig_{ij}dx^i dx^j$, where
   \[
   \tilde\ig_{ij}(x)=\big(\exp(\vr\ih(x))\big)_{ij}=\de_{ij}+\vr\ih_{ij}(x)+\frac{\vr^2}{2!}(\ih^2(x))_{ij}+\vr^3\big( \ih(x)\rr(\vr,\ih(x)) \big)_{ij}.
   \]
   Applying this to the Bourguignon-Gauduchon trivialization, we can see that the induced map $B$ is pointwisely given by  $B_x=(b_{ij}(x))=\exp(-\vr\ih(x)/2)$, where 
   \begin{\equ}\label{B-bij}
   \aligned
   b_{ij}(x)&=\de_{ij}-\frac{\vr\ih_{ij}(x)}{2}+\frac1{2!}\Big(\frac{\vr^2\ih^2(x)}4\Big)_{ij}-\frac{\vr^3}2\big( \ih(x)\rr(\vr,-\ih(x)/2) \big)_{ij} \\[0.3em]
   &=\de_{ij}-\frac{\vr\ih_{ij}(x)}{2}+\frac{\vr^2}8\sum_{d=1}^m \ih_{id}(x)\ih_{dj}(x)-\frac{\vr^3}2\big( \ih(x)\rr(\vr,-\ih(x)/2) \big)_{ij}.
   \endaligned
   \end{\equ}
   Denoting, for the sake of notational convenience, the entries of the inverse matrix of $B_x$ by $\bar b_{ij}(x)$, we have
   \begin{\equ}\label{B-inverse-ij}
   \bar b_{ij}(x)=\de_{ij}+\frac{\vr\ih_{ij}(x)}{2}+\frac{\vr^2}8\sum_{d=1}^m \ih_{id}(x)\ih_{dj}(x)+\frac{\vr^3}2\big( \ih(x)\rr(\vr,\ih(x)/2) \big)_{ij}.
   \end{\equ}
   Then, it has been shown in \cite[Proposition 3.2]{AGHM} that
   \begin{\equ}\label{Dirac-expansion}
   \widetilde D_{\tilde\ig} \psi = D_{\ig_{\R^m}}\psi+\Theta\cdot_{\ig_{\R^m}}\psi+\La\cdot_{\ig_{\R^m}}\psi+\sum_{i,j}(b_{ij}-\de_{ij}) e_i\cdot_{\ig_{\R^m}}\nabla^{\ig_{\R^m}}_{e_j}\psi
   \end{\equ}
   for any $\psi\in\Ga(\mbs(\R^m,\ig_{\R^m}))$, where $\big\{e_i=\frac{\pa }{\pa x^i}\big\}_{i=1}^m$ stands for the Euclidean orthonormal frame,
   \[
   \Theta=\frac14\sum_{\substack{i,j,k \\ i\neq j\neq k\neq i}} \sum_{\al,\bt}b_{i\al}(\pa_\al b_{j\bt})\bar b_{\bt k}\, e_i\cdot_{\ig_{\R^m}} e_j\cdot_{\ig_{\R^m}} e_k
   \]
   and
   \[
   \La=\frac12\sum_{i,k}\tilde\Ga_{ik}^i e_k
   \]
   with $\tilde\Ga_{ij}^k(x)=\tilde\ig_x\big(\nabla_{B_x(e_i)}^{\tilde\ig}B_x(e_j),B_x(e_k)\big)$ being the Christoffel symbols of the second kind.
   
   \begin{Rem}
   		On a spin manifold, since the tangent bundle is embedded in the bundle of Clifford algebra, vector fields act on spinors in two distinct ways: through Clifford multiplication and through covariant derivatives. To distinguish these two actions on a spinor $\psi$, we adopt the following notations. For $e_i :=\frac{\pa }{\pa x^i}= \pa_i \in T\mathbb{R}^m$:
   		\begin{itemize}
   			\item $e_i \cdot_{\ig_{\mathbb{R}^m}} \psi$ or $\pa_i \cdot_{\ig_{\mathbb{R}^m}} \psi$ denote the Clifford multiplication of $e_i$ with $\psi$,
   			
   			\item $\nabla_{e_i} \psi$ or $\nabla_{\pa_i} \psi$ denote the covariant derivative of $\psi$ under the metric $\ig_{\mathbb{R}^m}$.
   		\end{itemize}
   		For scalar functions, we will simplify the notation by writing $\pa_i u$ to represent the $i$-th partial derivative of $u$.
   \end{Rem}
	
	\subsection{Expansion of the perturbed functional}
	
	To provide readers with a precise expansion of the perturbed functional, we define the space
	\[
	\msd^{1/2}(\R^m,\mbs(\R^m))=\Big\{ \psi\in L^{\frac{2m}{m-1}}(\R^m,\mbs(\R^m))\,:\, \int_{\R^m} \big||D_{\ig_{\R^m}}|^{1/2}\psi\big|_{\ig_{\R^m}}^2d\vol_{\ig_{\R^m}}<\infty \Big\}
	\]
	over the Euclidean space $(\R^m,\ig_{\R^m})$, where the operator $|D_{\ig_{\R^m}}|^{1/2}$ is defined via the Fourier transformation. More specifically,  we have $\msf\big( |D_{\ig_{\R^m}}|^{1/2}\psi \big)(\xi)=|\xi|^{1/2}\msf(\psi)(\xi)$, where $\msf$ denotes the Fourier transform on $\R^m$. By the Sobolev's inequality, there exists a positive constant $C_m$, depending only on $m$, such that
	\[
	\Big( \int_{\R^m}|\psi|_{\ig_{\R^m}}^{\frac{2m}{m-1}}d\vol_{\ig_{\R^m}} \Big)^{\frac{m-1}{m}}\leq C_m \int_{\R^m} \big||D_{\ig_{\R^m}}|^{1/2}\psi\big|_{\ig_{\R^m}}^2d\vol_{\ig_{\R^m}}
	\]
	for all $\psi\in\msd^{1/2}(\R^m,\mbs(\R^m))$. We equip  $\msd^{1/2}(\R^m,\mbs(\R^m))$ with the norm 
	\[
	\|\psi\|_\msd=\Big(\int_{\R^m} \big||D_{\ig_{\R^m}}|^{1/2}\psi\big|_{\ig_{\R^m}}^2d\vol_{\ig_{\R^m}}\Big)^{1/2}.
	\]
	It is straightforward to verify that $\msd^{1/2}(\R^m,\mbs(\R^m))$, equipped with this norm, is complete. By introducing the inner product
	\[
	\inp{\psi}{\va}_\msd=\real\int_{\R^m}(|D_{\ig_{\R^m}}|^{1/2}\psi,|D_{\ig_{\R^m}}|^{1/2}\va)_{\ig_{\R^m}}d\vol_{\ig_{\R^m}},
	\]
	we have $\msd^{1/2}(\R^m,\mbs(\R^m))$ is a real Hilbert space.
    
For a Fr\'echet differentiable functional $J\colon\msd^{1/2}(\R^m,\mbs(\R^m))\to\R$,
we write $J'(\psi)$ for its derivative at $\psi$ and $\nabla J(\psi)$
for the unique element of $\msd^{1/2}(\R^m,\mbs(\R^m))$ satisfying
\[
J'(\psi)[v]=\inp{\nabla J(\psi)}{v}_{\msd}
\qquad\text{for every }v\in\msd^{1/2}(\R^m,\mbs(\R^m)).
\]
 
	For the Riemannian metric $\tilde\ig$ on $\R^m$ of the form \eqref{perturb-g}, let us consider the existence of weak solutions of the spinorial Yamabe equation
	\begin{\equ}\label{SY-tilde-g}
		D_{\tilde\ig}\va = |\va|_{\tilde\ig}^{\frac2{m-1}}\va \quad \text{on } (\R^m,\tilde\ig).
	\end{\equ}
	Here, a weak solution is a spinor $\va\in L^{\frac{2m}{m-1}}(\R^m,\mbs(\R^m))$ satisfying
	\[
	\int_{\R^m}(\va, D_{\tilde\ig}\phi)_{\tilde\ig}\, d\vol_{\tilde\ig}=\int_{\R^m}|\va|_{\tilde\ig}^{\frac{2}{m-1}}(\va,\phi)_{\tilde\ig}\,d\vol_{\tilde\ig}
	\]
	for all $\phi\in C_c^\infty(\R^m,\mbs(\R^m))$. From a variational perspective, a solution of \eqref{SY-tilde-g} corresponds to a critical point of the following energy functional 
	\[
	\cj(\va):=\frac12\int_{\R^m}(D_{\tilde\ig}\va,\va)_{\tilde\ig}\,d\vol_{\tilde\ig}-\frac{m-1}{2m}\int_{\R^m}|\va|_{\tilde\ig}^{\frac{2m}{m-1}}d\vol_{\tilde\ig}
	\] 
	define for spinors in $\mbs(\R^m,\tilde\ig)$. Using the Bourguignon-Gauduchon trivialization, we can apply the induced bundle isometry $\mfj:\mbs(\R^m,\ig_{\R^m})\to\mbs(\R^m,\tilde\ig)$. This allows us to define an equivalent functional $\cl$ on $\msd^{1/2}(\R^m,\mbs(\R^m))$ in the way that
	\begin{\equ}\label{the-cl}
	\cl(\psi)=\cj\circ\mfj(\psi)=\frac12\int_{\R^m}(\widetilde D_{\tilde\ig}\psi,\psi)_{\ig_{\R^m}}\,d\vol_{\ig_{\R^m}}-\frac{m-1}{2m}\int_{\R^m}|\psi|_{\ig_{\R^m}}^{\frac{2m}{m-1}}d\vol_{\ig_{\R^m}}
	\end{\equ}
	where we have used the fact $d\vol_{\ig_{\R^m}}=d\vol_{\tilde\ig}$. This construction enables the variational framework for studying weak solutions of \eqref{SY-tilde-g} in the setting of the Sobolev space $\msd^{1/2}(\R^m,\mbs(\R^m))$.
	
	\medskip

We next expand the quadratic part of the functional using \eqref{Dirac-expansion}. The exponential parametrization is important here: since $\tr h=0$, the volume element is unchanged, and only the Dirac quadratic form depends on the metric perturbation.

	\begin{Lem}\label{lemma: cl-expansion}
		Let $\tilde\ig$ be given in \eqref{perturb-g}, then
		\begin{\equ}\label{cl-expanded}
			\cl(\psi)=\cl_0(\psi)+\vr\Ga(\psi)+\vr^2\Phi(\psi)+o(\vr^2)
		\end{\equ}
		where 
		\[
		\cl_0(\psi)=\frac12\int_{\R^m}( D_{\ig_{\R^m}}\psi,\psi)_{\ig_{\R^m}}\,d\vol_{\ig_{\R^m}}-\frac{m-1}{2m}\int_{\R^m}|\psi|_{\ig_{\R^m}}^{\frac{2m}{m-1}}d\vol_{\ig_{\R^m}},
		\]
		\[
		\Ga(\psi)=-\frac{1}4\sum_{i,j}\int_{\R^m}\ih_{ij}\real(e_i\cdot_{\ig_{\R^m}}\nabla_{e_j}\psi,\psi)_{\ig_{\R^m}}d\vol_{\ig_{\R^m}}
		\]
		and
		\[
		\aligned
		\Phi(\psi)&=-\frac{1}{32}\sum_{\substack{i,j,k \\ i\neq j\neq k\neq i}}\sum_d\int_{\R^m}\pa_i\ih_{jd}\ih_{dk}\real (e_i\cdot_{\ig_{\R^m}} e_j\cdot_{\ig_{\R^m}}e_k\cdot_{\ig_{\R^m}}\psi,\psi)_{\ig_{\R^m}} d\vol_{\ig_{\R^m}} \\[0.3em]
		&\qquad + \frac{1}{16}\sum_{i,j}\sum_{d}\int_{\R^m}\ih_{id}\ih_{dj}\real(e_i\cdot_{\ig_{\R^m}}\nabla_{e_j}\psi,\psi)_{\ig_{\R^m}}d\vol_{\ig_{\R^m}}
		\endaligned
		\]
		for $\psi\in\msd^{1/2}(\R^m,\mbs(\R^m))$.
	\end{Lem}
	\begin{proof}
	We only need to expand the quadratic part of $\cl$ in terms of $\vr$. Using \eqref{Dirac-expansion} and the fact $\La\in T\R^m$, we get $\real(\La\cdot_{\ig_{\R^m}}\psi,\psi)_{\ig_{\R^m}}\equiv0$ and
	\[
	\aligned
	\real\big(\widetilde D_{\tilde\ig}\psi,\psi\big)_{\ig_{\R^m}}&=\real( D_{\ig_{\R^m}}\psi,\psi)_{\ig_{\R^m}}+\real(\Theta\cdot_{\ig_{\R^m}}\psi,\psi)_{\ig_{\R^m}} \\[0.3em]
	&\qquad +\sum_{i,j}(b_{ij}-\de_{ij})\real(e_i\cdot_{\ig_{\R^m}}\nabla_{e_j}\psi,\psi)_{\ig_{\R^m}},
	\endaligned
	\]
	where
	\[
	\real(\Theta\cdot_{\ig_{\R^m}}\psi,\psi)_{\ig_{\R^m}}=\frac14\sum_{\substack{i,j,k \\ i\neq j\neq k\neq i}} \Big(\sum_{\al,\bt}b_{i\al}(\pa_\al b_{j\bt})\bar b_{\bt k}\Big)\real (e_i\cdot_{\ig_{\R^m}} e_j\cdot_{\ig_{\R^m}}e_k\cdot_{\ig_{\R^m}}\psi,\psi)_{\ig_{\R^m}}
	\]
	Plainly, $\real(\Theta\cdot_{\ig_{\R^m}}\psi,\psi)_{\ig_{\R^m}}\equiv0$ for $m=2$. For $m\geq3$, by \eqref{B-bij} and \eqref{B-inverse-ij}, we have that
	\[
	\left\{
	\aligned
	&b_{i\al}=\de_{i\al}-\frac\vr2\ih_{i\al}+\frac{\vr^2}8\sum_{d}\ih_{id}\big[ \ih_{d\al}-4\vr \rr(\vr,-\ih/2)_{d\al}\big], \\[0.3em]
	&\pa_\al b_{j\bt}=-\frac\vr2\pa_\al \ih_{j\bt}+\frac{\vr^2}8\sum_{d}\Big( \pa_\al\ih_{jd}\ih_{d\bt}+\ih_{jd}\pa_\al\ih_{d\bt} \Big)-\frac{\vr^3}2\sum_d\pa_\al\big[ \ih_{jd}\rr(\vr,-\ih/2)_{d\bt} \big], \\[0.3em]
	&\bar b_{\bt k}=\de_{\bt k}+\frac\vr2\ih_{\bt k}+\frac{\vr^2}8\sum_{d}\ih_{\bt d}\big[\ih_{d k}+4\vr\rr(\vr,\ih/2)_{dk}\big].
	\endaligned
	\right.
	\]
	Hence
\begin{\equ}\label{bbb}
	\aligned
	b_{i\al}(\pa_\al b_{j\bt})\bar b_{\bt k}&=-\frac\vr2\de_{i\al}\de_{\bt k}\pa_\al\ih_{j\bt}+\frac{\vr^2}4\pa_\al\ih_{j\bt}\big( \de_{\bt k}\ih_{i\al}-\de_{i\al}\ih_{\bt k} \big)  \\[0.3em]
	&\qquad +\frac{\vr^2\de_{i\al}\de_{\bt k}}{8}\sum_{d}\Big( \pa_\al\ih_{jd}\ih_{d\bt}+\ih_{jd}\pa_\al\ih_{d\bt} \Big) + o(\vr^2),
	\endaligned
\end{\equ}
For $|\vr|\|h\|_\infty$ sufficiently small, the remainder $o(\vr^2)$
in the preceding expansion is bounded pointwise by $C|\vr|^3|h|^2|\nabla h|$.
An exact quadratic expression for the remainder is given in Remark \ref{rev-exact-remainder}.

	Since we have $\pa_i\ih_{jk}=\pa_i\ih_{kj}$ and $ e_j\cdot_{\ig_{\R^m}} e_k=- e_k\cdot_{\ig_{\R^m}} e_j$, it follows that
	\[
	\sum_{\substack{i,j,k \\ i\neq j\neq k\neq i}}\pa_i\ih_{jk} \, e_i\cdot_{\ig_{\R^m}} e_j\cdot_{\ig_{\R^m}} e_k=0
	\]
	Therefore, $\real(\Theta\cdot_{\ig_{\R^m}}\psi,\psi)_{\ig_{\R^m}}$ has no term of order $1$ of $\vr$. To compute the terms of order $2$, we first note that
	\[
	\sum_{\substack{i,j,k \\ i\neq j\neq k\neq i}}\pa_\al\ih_{jk}\ih_{i\al} \, e_i\cdot_{\ig_{\R^m}} e_j\cdot_{\ig_{\R^m}} e_k=0
	\]
	and
	\[
	\sum_{\substack{i,j,k \\ i\neq j\neq k\neq i}}\sum_d\Big( \pa_i\ih_{jd}\ih_{dk}+\ih_{jd}\pa_i\ih_{dk} \Big)e_i\cdot_{\ig_{\R^m}} e_j\cdot_{\ig_{\R^m}} e_k=0.
	\]
	These are because, in the above summations, we iterate through the cases where $j$ and $k$ are swapped.
	Thus, in the expression of $b_{i\al}(\pa_\al b_{j\bt})\bar b_{\bt k}$, we only need to focus on the terms arising from $-\frac{\vr^2}4\de_{i\al}\pa_\al\ih_{j\bt}\ih_{\bt k}$. As a result, we obtain
	\[
	\real(\Theta\cdot_{\ig_{\R^m}}\psi,\psi)_{\ig_{\R^m}}=-\frac{\vr^2}{16}\sum_{\substack{i,j,k \\ i\neq j\neq k\neq i}}\sum_d\pa_i\ih_{jd}\ih_{dk}\real (e_i\cdot_{\ig_{\R^m}} e_j\cdot_{\ig_{\R^m}}e_k\cdot_{\ig_{\R^m}}\psi,\psi)_{\ig_{\R^m}}+o(\vr^2).
	\]
	And it follows that
	\[
	\aligned
	\cl(\psi)&=\cl_0(\psi)-\frac{\vr}4\sum_{i,j}\int_{\R^m}\ih_{ij}\real(e_i\cdot_{\ig_{\R^m}}\nabla_{e_j}\psi,\psi)_{\ig_{\R^m}}d\vol_{\ig_{\R^m}}\\[0.3em]
	&\qquad -\frac{\vr^2}{32}\sum_{\substack{i,j,k \\ i\neq j\neq k\neq i}}\sum_d\int_{\R^m}\pa_i\ih_{jd}\ih_{dk}\real (e_i\cdot_{\ig_{\R^m}} e_j\cdot_{\ig_{\R^m}}e_k\cdot_{\ig_{\R^m}}\psi,\psi)_{\ig_{\R^m}} d\vol_{\ig_{\R^m}} \\[0.3em]
	&\qquad + \frac{\vr^2}{16}\sum_{i,j}\sum_{d}\int_{\R^m}\ih_{id}\ih_{dj}\real(e_i\cdot_{\ig_{\R^m}}\nabla_{e_j}\psi,\psi)_{\ig_{\R^m}}d\vol_{\ig_{\R^m}}+o(\vr^2),
	\endaligned
	\]
	which completes the proof.
	\end{proof}
	
\begin{Rem}\label{rev-exact-remainder}
For later uniform estimates it is useful to retain an exact expression
for the remainder rather than an unweighted $o(\vr^2)$ in \eqref{bbb}.
Write
\[
A_{ijk}(\vr)=\sum_{\al,\bt}b_{i\al}
(\partial_\al b_{j\bt})\bar b_{\bt k},
\qquad
C_{ijk}=\sum_d(\partial_i h_{jd})h_{dk}.
\]
The remainder in \eqref{cl-expanded} is the sum $R_{\rm pr}+R_{\rm Cl}$.
Here the subscript ``$\mathrm{pr}$'' refers to the principal part of the Dirac operator,
which differentiates the spinor, and ``$\mathrm{Cl}$" refers to Clifford
multiplication by the three-form $\Theta$ in \eqref{Dirac-expansion}.
Explicitly,
\[
R_{\rm pr}(\psi)=\frac12\sum_{i,j}\int_{\R^m}
\left(B-I+\frac{\vr}{2}h-\frac{\vr^2}{8}h^2\right)_{ij}
\real(e_i\cdot\partial_j\psi,\psi)\,dx
\]
and
\[
R_{\rm Cl}(\psi)=\frac18
\sum_{\substack{i,j,k \\ i\neq j\neq k\neq i}}\int_{\R^m}
\left(A_{ijk}(\vr)+\frac{\vr^2}{4}C_{ijk}\right)
\real(e_i\cdot e_j\cdot e_k\cdot\psi,\psi)\,dx.
\]
In $R_{\rm Cl}$, the terms linear in $\vr$ cancel in pairs.
Indeed, interchanging $j$ and $k$ leaves $\partial_i h_{jk}$ unchanged
and reverses the sign of the Clifford product.
After summation, the term $+\frac{\vr^2}{4}C_{ijk}$ cancels the
remaining quadratic contribution from $A_{ijk}(\vr)$.
{When $|\vr|\|h\|_\infty$ is small, the coefficients in
$R_{\rm pr}$ are bounded by $C|\vr|^3|h|^3$. In $R_{\rm Cl}$, summing
the six permutations of each triple $i<j<k$ gives a coefficient
bounded by $C|\vr|^3|h|^2|\nabla h|$ for the ordered Clifford product
$e_i\cdot e_j\cdot e_k$.}
All metric perturbation terms are exactly quadratic in $\psi$.
\end{Rem}

\medskip

	Recall that the spinor module $\mbs_m$ has complex dimension $2^{[\frac m2]}$ and equipped with a hermitian metric $|\cdot|_{\ig_{\R^m}}$. Let us consider $S(\mbs_m)$, the unit sphere in $\mbs_m$, which has real dimension $2^{[\frac m2]+1}-1$. Given any triple $(\lm,\xi,\ga)\in(0,+\infty)\times\R^m\times S(\mbs_m)$, we define a spinor $\psi_{\lm,\xi,\ga}:\R^m\to\mbs_m$ by
	\begin{\equ}\label{bubble}
	\psi_{\lm,\xi,\ga}(x)=\frac{m^{\frac{m-1}2}\lm^{\frac{m-1}2}}{(\lm^2+|x-\xi|^2)^{\frac m2}}(\lm-(x-\xi))\cdot_{\ig_{\R^m}}\ga
	\end{\equ}
	and set 
	\begin{\equ}\label{critical-mfd}
	\cm=\big\{\psi_{\lm,\xi,\ga}:\, (\lm,\xi,\ga)\in(0,+\infty)\times\R^m\times S(\mbs_m)\big\}.
	\end{\equ}
	The spinor $\psi_{\lm,\xi,\ga}$ satisfies the PDE
	\[
	D_{\ig_{\R^m}}\psi_{\lm,\xi,\ga}=|\psi_{\lm,\xi,\ga}|_{\ig_{\R^m}}^{\frac2{m-1}}\psi_{\lm,\xi,\ga} \quad \text{on } (\R^m,\ig_{\R^m}).
	\]
	It is also well known that
	\[
	\int_{\R^m}|\psi_{\lm,\xi,\ga}|^{\frac{2m}{m-1}}_{\ig_{\R^m}}d\vol_{\ig_{\R^m}}=\Big(\frac m2\Big)^m\om_m
	\]
	for all $(\lm,\xi,\ga)\in(0,+\infty)\times\R^m\times S(\mbs_m)$, where $\om_m$ stands for the volume of the $m$-sphere $(S^m,\ig_{S^m})$. 
	
%	We also remark that
%	\[
%	\aligned
%	\pa_i\ih_{jd}\ih_{dk}&=a^2f(b^{-2}|x|^2)f'(b^{-2}|x|^2)\frac{2x^i}{b^2}H_{jd}H_{dk}+a^2f(b^{-2}|x|^2)^2(W_{ji\bt d}+W_{j\bt id})W_{kpqd}x^px^qx^\bt \\
%	&=a^2f(b^{-2}|x|^2)f'(b^{-2}|x|^2)\frac{2x^i}{b^2}H_{jd}H_{dk}+a^2f(b^{-2}|x|^2)^2\frac{W_{j\bt id}W_{kpqd}x^px^qx^\bt x^i}{x^i} \\
%	&=a^2f(b^{-2}|x|^2)f'(b^{-2}|x|^2)\frac{2x^i}{b^2}H_{jd}H_{dk}+a^2f(b^{-2}|x|^2)^2\frac{W_{jpqd}W_{k\bt id}x^px^qx^\bt x^i}{x^i} \\
%	&=a^2f(b^{-2}|x|^2)f'(b^{-2}|x|^2)\frac{2x^i}{b^2}H_{jd}H_{dk}+a^2f(b^{-2}|x|^2)^2(W_{ki\bt d}+W_{k\bt id})W_{jpqd}x^px^qx^\bt
%	\endaligned
%	\]
% The cancellation for Table~\ref{tab:weyl-vanishing-spinor-term} is proved after \eqref{rev-table1-structure}. For Table~\ref{tab:weyl-indefinite-spinor-term} the Clifford term survives.
	
	\begin{Lem}\label{lemma: critical-manifold}
	Let $\cl_0$, $\Ga$ and $\Phi$ be as in Lemma \ref{lemma: cl-expansion} and $\cm$ be as above, then
	\begin{itemize}
		\item[$(1)$]  $\cm\subset \msd^{1/2}(\R^m,\mbs(\R^m))$ is a finite dimensional $C^2$-submanifold such that $\nabla\cl_0(\psi)=0$ for all $\psi\in\cm$;
		
		\item[$(2)$] $\cm$ is non-degenerate with respect to $\cl_0$ in the sense that $T_{\psi}\cm=\ker\nabla^2\cl_0(\psi)$ for all $\psi\in\cm$;
		
		\item[$(3)$]  $\nabla^2\cl_0(\psi):\msd^{1/2}(\R^m,\mbs(\R^m))\to\msd^{1/2}(\R^m,\mbs(\R^m))$ is a Fredholm operator with index zero for all $\psi\in\cm$;
		
		\item[$(4)$] $\Ga|_{\cm}\equiv0$, i.e., $\Ga(\psi)=0$ for all $\psi\in\cm$, moreover, $\nabla \Ga(\psi_{\lm,0,\ga})=0$ for all $(\lm,\ga)\in(0,+\infty)\times S(\mbs_m)$;
		
		\item[$(5)$] $\Phi(\psi_{\lm,0,\ga})$ has the expression
		\[
		\aligned
		\Phi(\psi_{\lm,0,\ga})&=\frac{m^{m-1}\lm^m}{16}\bigg[\int_{\R^m}\frac{\tr(\ih^2)(x)}{(\lm^2+|x|^2)^m}d\vol_{\ig_{\R^m}} \\[0.5em]
		&\qquad +P(W,\ga) \int_{\R^m}\frac{a^2\eta(r^{-1}|x|)^2f(b^{-2}|x|^2)^2|x^1|^2 |x^2|^2}{(\lm^2+|x|^2)^m} d\vol_{\ig_{\R^m}}  \bigg]
		\endaligned
		\]
		for all $(\lm,\ga)\in(0,+\infty)\times S(\mbs_m)$, where
		\[
		P(W,\ga)=9\sum_{\substack{i,j,k \\ i\neq j\neq k\neq i}}\sum_{d,p,q}
		W_{jpid}W_{kpqd}S_{ijkq}(\ga).
		\]
		and 
		\[
		S_{ijkq}(\ga)=\real(e_i\cdot_{\ig_{\R^m}} e_j\cdot_{\ig_{\R^m}}e_k\cdot_{\ig_{\R^m}}e_q\cdot_{\ig_{\R^m}}\ga,\ga)_{\ig_{\R^m}}.
		\]
	\end{itemize}
	\end{Lem}
	\begin{proof}
	It is obvious to verify $(1)$ since we already have the explicit expression of $\cl_0$ and $\cm$. The proof for $(2)$ and $(3)$ can be found in \cite[Sections 5 and 6]{Isobe13}, thus we only need to check $(4)$ and $(5)$.
	
	To see $(4)$, for an arbitrary $\psi_{\lm,\xi,\ga}\in\cm$, we first compute
	\begin{\equ}\label{pa-j-psi}
	\aligned
	\nabla_{e_j}\psi_{\lm,\xi,\ga}(x)&=-\frac{m^{\frac{m+1}2}\lm^{\frac{m-1}2}}{(\lm^2+|x-\xi|^2)^{\frac m2+1}}(x^j-\xi^j)(\lm-(x-\xi))\cdot_{\ig_{\R^m}}\ga\\[0.3em]
	&\qquad -\frac{m^{\frac{m-1}2}\lm^{\frac{m-1}2}}{(\lm^2+|x-\xi|^2)^{\frac m2}}e_j\cdot_{\ig_{\R^m}}\ga.
	\endaligned
	\end{\equ}
	Then, by using $\real(e_i\cdot_{\ig_{\R^m}}\ga,\ga)_{\ig_{\R^m}}=0$ and $\real(e_i\cdot_{\ig_{\R^m}}e_j\cdot_{\ig_{\R^m}}\ga,\ga)_{\ig_{\R^m}}=0$ for $i\neq j$, we obtain
	\begin{\equ}\label{e-d-psi}
	\aligned
	&\real(e_i\cdot_{\ig_{\R^m}}\nabla_{e_j}\psi_{\lm,\xi,\ga},\psi_{\lm,\xi,\ga})_{\ig_{\R^m}}= \\[0.5em]
	&\qquad \qquad
	\left\{
	\aligned
	&\displaystyle \frac{m^{m-1}\lm^m}{(\lm^2+|x-\xi|^2)^m}, && \text{if } i=j, \\[0.5em]
	&\displaystyle \frac{m^{m-1}\lm^{m-1}}{(\lm^2+|x-\xi|^2)^m}\real(e_i\cdot_{\ig_{\R^m}}e_j\cdot_{\ig_{\R^m}}\ga,(x-\xi)\cdot_{\ig_{\R^m}}\ga)_{\ig_{\R^m}}, && \text{if } i\neq j.
	\endaligned
	\right.
	\endaligned
	\end{\equ}
	Since $\ih(x)=(\ih_{ij}(x))$ is a trace-free symmetric two-tensor and $e_i\cdot_{\ig_{\R^m}}e_j=-e_j\cdot_{\ig_{\R^m}} e_i$ for $i\neq j$, it follows immediately that $\Ga(\psi_{\lm,\xi,\ga})\equiv0$ for all $(\lm,\xi,\ga)\in(0,+\infty)\times\R^m\times S(\mbs_m)$.

To compute $\nabla\Ga(\psi_{\lm,0,\ga})$, let $z\in\msd^{1/2}(\R^m,\mbs(\R^m))$ be smooth
and let $w$ be a smooth compactly supported test spinor. Differentiating
the expression for $\Ga$ in Lemma~\ref{lemma: cl-expansion} and
integrating by parts gives
\[
\begin{aligned}
\Ga'(z)[w]
&=-\frac14\sum_{i,j}\int_{\R^m}h_{ij}
 \{\real(e_i\cdot\partial_jz,w)+\real(e_i\cdot\partial_jw,z)\}\,dx\\
&=-\frac14\sum_{i,j}\int_{\R^m}(\partial_jh_{ij})
 \real(e_i\cdot z,w)\,dx\\
&\qquad-\frac12\sum_{i,j}\int_{\R^m}h_{ij}
 \real(e_i\cdot\partial_jz,w)\,dx.
\end{aligned}
\]
Thus $\displaystyle \Ga'(z)[w]=\real\int_{\R^m}(F_h(z),w)_{\ig_{\R^m}}\,dx$, where
\begin{\equ}\label{rev-source}
F_h(z)=-\frac14\sum_{i,j}(\pa_jh_{ij})e_i\cdot_{\ig_{\R^m}}z
-\frac12\sum_{i,j}h_{ij}e_i\cdot_{\ig_{\R^m}}\nabla_{e_j}z.
\end{\equ}
By density, this identifies $\Ga'(z)$ with $F_h(z)$ in the dual of
$\msd^{1/2}(\R^m,\mbs(\R^m))$. Its Hilbert gradient is therefore
$\nabla\Ga(z)=|D_{\ig_{\R^m}}|^{-1}F_h(z)$, so it suffices to show
that $F_h(\psi_{\lm,0,\ga})=0$.

Set $\xi=0$ in \eqref{pa-j-psi} and substitute the resulting expression
for $\nabla_{e_j}\psi_{\lm,0,\ga}$ into \eqref{rev-source}.
The term containing $x^j$ vanishes upon summation since
$\sum_j h_{ij}(x)x^j=0$ for every $i$. In the other term, the
off-diagonal summands cancel in pairs by the symmetry of $h$ and the
Clifford relations, so
\[
\sum_{i,j}h_{ij}\,e_i\cdot_{\ig_{\R^m}}e_j\cdot_{\ig_{\R^m}}\ga
=-(\tr h)\ga=0.
\]
Since $\diver h=0$ as well, \eqref{rev-source} gives
$F_h(\psi_{\lm,0,\ga})=0$ and hence $\nabla\Ga(\psi_{\lm,0,\ga})=0$.

	\medskip
	
	Finally, to derive the specific form of $\Phi(\psi_{\lambda,0,\ga})$, a more precise calculation is required by substituting $\psi_{\lambda,\xi,\ga}$ into the expression of $\Phi$ and then take $\xi=0$. First, note that by using the anti-commutativity property, where $e_i \cdot_{\ig_{\R^m}} e_j = -e_j \cdot_{\ig_{\R^m}} e_i$ for $i \neq j$, the second part of $\Phi(\psi_{\lm,\xi,\ga})$ (see its expression in Lemma \ref{lemma: cl-expansion}) can be calculated as
	\begin{\equ}\label{Phi-on-M}
	\aligned
	&\sum_{i,j}\sum_{d}\int_{\R^m}\ih_{id}\ih_{dj}\real(e_i\cdot_{\ig_{\R^m}}\nabla_{e_j}\psi_{\lm,\xi,\ga},\psi_{\lm,\xi,\ga})_{\ig_{\R^m}}d\vol_{\ig_{\R^m}} \\
	&\qquad =\sum_{i,d}\int_{\R^m}\ih_{id}\ih_{di}\real(e_i\cdot_{\ig_{\R^m}}\nabla_{e_i}\psi_{\lm,\xi,\ga},\psi_{\lm,\xi,\ga})_{\ig_{\R^m}}d\vol_{\ig_{\R^m}} \\
	&\qquad =m^{m-1}\lm^m \int_{\R^m}\frac{\tr(\ih^2)}{(\lm^2+|x-\xi|^2)^m}d\vol_{\ig_{\R^m}}.
	\endaligned
	\end{\equ}
	Meanwhile, to analyze the first part of $\Phi(\psi_{\lm,\xi,\ga})$ with $\pa_i\ih_{jd}\ih_{dk}$ in detail, we compute
	\[
	\aligned
	\pa_i\ih_{jd}(x)\ih_{dk}(x)&=\frac{2a^2}{b^2}\eta(r^{-1}|x|)^2f(b^{-2}|x|^2)f'(b^{-2}|x|^2)x^iH_{jd}(x)H_{dk}(x)\\[0.3em]
	&\qquad + \frac{a^2}{r|x|}\eta(r^{-1}|x|)\eta'(r^{-1}|x|)f(b^{-2}|x|^2)^2x^iH_{jd}(x)H_{dk}(x)\\[0.3em]
	&\qquad +a^2\eta(r^{-1}|x|)^2f(b^{-2}|x|^2)^2\pa_i H_{jd}(x)H_{dk}(x) \\[0.3em]
	&=a^2\Big[\frac{2}{b^2}\eta(r^{-1}|x|)^2f(b^{-2}|x|^2)f'(b^{-2}|x|^2)\\[0.3em]
	&\qquad +\frac{1}{r|x|}\eta(r^{-1}|x|)\eta'(r^{-1}|x|)f(b^{-2}|x|^2)^2
	\Big]x^iH_{jd}(x)H_{dk}(x)\\[0.3em]
	&\qquad +a^2\eta(r^{-1}|x|)^2f(b^{-2}|x|^2)^2(W_{ji\al d}+W_{j\al id})W_{kpqd}x^{\al}x^px^q.
	\endaligned
	\]
	In the sequel, for the sake of simplicity, let us denote
	\[
	\aligned
	A_{ijdk}(x)&=a^2\Big[\frac{2}{b^2}\eta(r^{-1}|x|)^2f(b^{-2}|x|^2)f'(b^{-2}|x|^2)\\[0.3em]
	&\qquad +\frac{1}{r|x|}\eta(r^{-1}|x|)\eta'(r^{-1}|x|)f(b^{-2}|x|^2)^2\Big]x^iH_{jd}(x)H_{dk}(x)
	\endaligned
	\]
	and
    \begin{\equ}\label{B-ijdk}
    B_{ijdk}(x)=a^2\eta(r^{-1}|x|)^2f(b^{-2}|x|^2)^2(W_{ji\al d}+W_{j\al id})W_{kpqd}x^{\al}x^px^q
    \end{\equ}
    so that $\pa_i\ih_{jd}(x)\ih_{dk}(x)=A_{ijdk}(x)+B_{ijdk}(x)$ (with no summation over $i,j,d,k$ in this identity).
	Since $e_j\cdot_{\ig_{\R^m}}e_k=-e_k\cdot_{\ig_{\R^m}} e_j$ for $j\neq k$ and $A_{ijdk}(x)=A_{ikdj}(x)$, the contributions involving $A_{ijdk}$ to
	$\Phi(\psi_{\lm,\xi,\ga})$ cancel in pairs upon interchanging
	$j$ and $k$. Consequently,
	\begin{\equ}\label{Phi-on-M2}
	\aligned
	&\sum_{\substack{i,j,k \\ i\neq j\neq k\neq i}}\sum_d\int_{\R^m}\pa_i\ih_{jd}\ih_{dk}\real (e_i\cdot_{\ig_{\R^m}} e_j\cdot_{\ig_{\R^m}}e_k\cdot_{\ig_{\R^m}}\psi_{\lm,\xi,\ga},\psi_{\lm,\xi,\ga})_{\ig_{\R^m}} d\vol_{\ig_{\R^m}} \\[0.3em]
	&\qquad =\sum_{\substack{i,j,k \\ i\neq j\neq k\neq i}}\sum_d\int_{\R^m}B_{ijdk}(x)\real (e_i\cdot_{\ig_{\R^m}} e_j\cdot_{\ig_{\R^m}}e_k\cdot_{\ig_{\R^m}}\psi_{\lm,\xi,\ga},\psi_{\lm,\xi,\ga})_{\ig_{\R^m}} d\vol_{\ig_{\R^m}}.
	\endaligned
	\end{\equ}
	We now set $\xi=0$ and substitute the explicit expression
	for $\psi_{\lm,0,\ga}$ from \eqref{bubble}.
	For pairwise distinct $i,j,k$, the operator
	$e_i\cdot_{\ig_{\R^m}}
	e_j\cdot_{\ig_{\R^m}}
	e_k\cdot_{\ig_{\R^m}}$ is self-adjoint. Thus, the two cross terms in the spinor
	pairing combine, and each integral on the right-hand side
	of \eqref{Phi-on-M2} is the sum of the following three terms:
	\begin{\equ}\label{3order}
    m^{m-1}\lm^{m+1}\int_{\R^m}\frac{B_{ijdk}(x)}{(\lm^2+|x|^2)^m}\real(e_i\cdot_{\ig_{\R^m}} e_j\cdot_{\ig_{\R^m}}e_k\cdot_{\ig_{\R^m}}\ga,\ga)_{\ig_{\R^m}} d\vol_{\ig_{\R^m}},
	\end{\equ}
	\begin{\equ}\label{4order}
	-2m^{m-1}\lm^{m}\int_{\R^m}\frac{B_{ijdk}(x)}{(\lm^2+|x|^2)^m}\real(e_i\cdot_{\ig_{\R^m}} e_j\cdot_{\ig_{\R^m}}e_k\cdot_{\ig_{\R^m}}x\cdot_{\ig_{\R^m}}\ga,\ga)_{\ig_{\R^m}}  d\vol_{\ig_{\R^m}},
	\end{\equ}
	and
	\begin{\equ}\label{5order}
	m^{m-1}\lm^{m-1}\int_{\R^m}\frac{B_{ijdk}(x)}{(\lm^2+|x|^2)^m}\real(e_i\cdot_{\ig_{\R^m}} e_j\cdot_{\ig_{\R^m}}e_k\cdot_{\ig_{\R^m}}x\cdot_{\ig_{\R^m}}\ga,x\cdot_{\ig_{\R^m}}\ga)_{\ig_{\R^m}}  d\vol_{\ig_{\R^m}}.
	\end{\equ}
	By \eqref{B-ijdk}, $B_{ijdk}(x)$ is the product of a radial
	function and a homogeneous polynomial of degree three.
	In particular, it is odd under $x\mapsto -x$.
	The spinor pairings in \eqref{3order} and \eqref{5order}
	are, respectively, independent of $x$ and quadratic in $x$,
	and hence are both even.
	The integrands in these two expressions are therefore odd,
	so both integrals vanish. And it remains to evaluate the term in \eqref{4order}.
	
	Expanding $x \cdot_{\ig_{\R^m}} \ga = x^l e_l \cdot_{\ig_{\R^m}} \ga$, we observe that any term in \eqref{4order} containing an odd power of at least one coordinate integrates to zero by reflection symmetry. Thus, only terms in which every coordinate occurs to an even power can contribute. For convenience, write $S_{ijkl}(\ga)=\real(e_i\cdot_{\ig_{\R^m}} e_j\cdot_{\ig_{\R^m}}e_k\cdot_{\ig_{\R^m}}e_l\cdot_{\ig_{\R^m}}\ga,\ga)_{\ig_{\R^m}}$, for $\ga\in S(\mbs_m)$. For fixed, pairwise distinct $i,j,k$ and fixed $d$, the potentially nonzero contributions are therefore of the following three types, with the common prefactor $-2m^{m-1}\lm^m$ omitted:
\begin{itemize}
		\item[(i)] when $\al=p$, $q=l$ and $p\neq q$
	\[
(W_{ji p d}+W_{jp id})W_{kp qd}S_{ijkq}(\ga)	\int_{\R^m}\frac{a^2\eta(r^{-1}|x|)^2f(b^{-2}|x|^2)^2|x^p|^2 |x^q|^2}{(\lm^2+|x|^2)^m} d\vol_{\ig_{\R^m}};
	\]

\item[(ii)] when $\al=q$, $p=l$ and $p\neq q$
	\[
(W_{jiqd}+W_{jq id})W_{kpq d}	S_{ijkp}(\ga) \int_{\R^m}\frac{a^2\eta(r^{-1}|x|)^2f(b^{-2}|x|^2)^2|x^p|^2 |x^q|^2}{(\lm^2+|x|^2)^m} d\vol_{\ig_{\R^m}};
	\]

\item[(iii)] when $\al=l$ and $p=q$
	\[
(W_{ji\al d}+W_{j\al id})W_{kpp d} S_{ijk\al}(\ga)	\int_{\R^m}\frac{a^2\eta(r^{-1}|x|)^2f(b^{-2}|x|^2)^2|x^{\al}|^2 |x^p|^2}{(\lm^2+|x|^2)^m} d\vol_{\ig_{\R^m}}.
	\]
\end{itemize}
For $1\leq i,j,k,l\leq m$, set
\begin{\equ}\label{I12}
 I_{ijkl}=m^{m-1}\lm^m\int_{\R^m}
 \frac{a^2\eta(r^{-1}|x|)^2 f(b^{-2}|x|^2)^2\,x^i x^j x^k x^l}
 {(\lm^2+|x|^2)^m}\,d\vol_{\ig_{\R^m}}.
\end{\equ}
Write $I_{12}=I_{1122}$ and $I_{11}=I_{1111}$. By radial symmetry,
$I_{ijkl}$ equals $I_{12}$ when two distinct indices each occur twice,
$3I_{12}$ when all four coincide, and zero otherwise.
In particular, $I_{11}=3I_{12}$.
Substitution in \eqref{4order} gives
\begin{eqnarray*}
\eqref{4order}&=&-2\bigg[I_{12}\sum_{p\neq q}(W_{ji p d}+W_{jp id})W_{kp qd}S_{ijkq}(\ga) +I_{12}\sum_{p\neq q}(W_{ji q d}+W_{jq id})W_{kp qd}S_{ijkp}(\ga) \\[0.3em]
& &+I_{12}\sum_{\al\neq p}(W_{ji\al d}+W_{j\al id})W_{kpp d}S_{ijk\al}(\ga) +I_{11}\sum_{p}(W_{jipd}+W_{jp id})W_{kpp d}S_{ijkp}(\ga) \bigg] \\[0.3em]
&=&-2\bigg[I_{12}\sum_{p,q}(W_{ji p d}+W_{jp id})(W_{kp qd}+W_{kqpd})S_{ijkq}(\ga) \\[0.3em]
& &-2I_{12}\sum_p(W_{jipd}+W_{jpid})W_{kppd}S_{ijkp}(\ga) 
+I_{12}\sum_{\al\neq p}(W_{ji\al d}+W_{j\al id})W_{kpp d}S_{ijk\al}(\ga)
\\[0.3em]
&&+ I_{11}\sum_{p}(W_{jipd}+W_{jp id})W_{kpp d}S_{ijkp}(\ga) \bigg] \\[0.3em]
&=&-2I_{12}\sum_{p,q}\Big[(W_{ji p d}+W_{jp id})(W_{kp qd}+W_{kqpd})+(W_{jiq d}+W_{jq id})W_{kpp d}\Big]S_{ijkq}(\ga)
\end{eqnarray*}
where, from the first equality to the second, we merely extend the first two sums to all $p,q$ and subtract the resulting diagonal contributions; while for the third equality, we have used precisely the relation $I_{11}=3I_{12}$.

Since $i,j,k$ are pairwise distinct, we have $S_{ijkq}(\ga)=0$ whenever $q\in\{i,j,k\}$. Thus, only terms with four pairwise distinct indices contribute to the sum in \eqref{4order}. On these terms,  interchanging any two indices changes the sign of $S_{ijkq}(\ga)$.

Note that, the trace-free property of $W$ gives
$\sum_p W_{kppd}=0$. And consequently,
\[
\sum_{\substack{i,j,k \\ i\neq j\neq k\neq i}}
\sum_{d,p,q}
(W_{jiqd}+W_{jqid})W_{kppd}S_{ijkq}(\ga)=0.
\]
For the remaining terms, the first Bianchi identity implies
\[
W_{jpid}-W_{ipjd}=W_{jipd},
\qquad
W_{kpqd}-W_{qpkd}=W_{kqpd}.
\]
Antisymmetrizing in $i,j$ and in $k,q$, we obtain
\[
\aligned
	&\sum_{\substack{i,j,k \\ i\neq j\neq k\neq i}}\sum_{d,p,q}
	(W_{jipd}+W_{jpid})(W_{kpqd}+W_{kqpd})
	S_{ijkq}(\ga)
	\\[0.3em]
	&\qquad=
	3
	\sum_{\substack{i,j,k \\ i\neq j\neq k\neq i}}\sum_{d,p,q}
	W_{jpid}(W_{kpqd}+W_{kqpd})S_{ijkq}(\ga)
	\\[0.3em]
	&\qquad=
	9\sum_{\substack{i,j,k \\ i\neq j\neq k\neq i}}\sum_{d,p,q}
	W_{jpid}W_{kpqd}S_{ijkq}(\ga).
\endaligned
\]
Substituting these identities, together with the preceding computations, into the expansion in Lemma \ref{lemma: cl-expansion} yields $(5)$ and completes the proof.
\end{proof}

Building on the abstract framework in \cite[Section 2]{Isobe-Xu23} and combining Lemmas \ref{lemma: cl-expansion} and \ref{lemma: critical-manifold}, we know that finding solutions to the spinorial Yamabe equation \eqref{SY-tilde-g} reduces to finding critical points of the reduced functional
\begin{\equ}\label{cl-red}
\aligned
\cl^{red}(\lm,\xi,\ga)&= \cl(\psi_{\lm,\xi,\ga}+w_{\lm,\xi,\ga})  \\[0.3em]
&=\cl_0(\psi_{\lm,\xi,\ga})+\vr^2\Big( \Phi(\psi_{\lm,\xi,\ga})-\frac12\inp{K_{\lm,\xi,\ga}(\nabla\Ga(\psi_{\lm,\xi,\ga}))}{\nabla\Ga(\psi_{\lm,\xi,\ga})}_{\msd}\Big)+o(\vr^2) \\[0.3em]
&=\frac1{2m}\Big(\frac m2\Big)^m\om_m+\vr^2\,\widehat\Phi(\lm,\xi,\ga)+o(\vr^2)
\endaligned
\end{\equ}
defined for $(\lm,\xi,\ga)\in(0,+\infty)\times\R^m\times S(\mbs_m)$, where we use the symbol $K_{\lm,\xi,\ga}$ to denote the inverse of $\nabla^2\cl_0(\psi_{\lm,\xi,\ga})$ restricted to $T_{\psi_{\lm,\xi,\ga}}\cm^\bot$, the orthogonal complement of $T_{\psi_{\lm,\xi,\ga}}\cm$ in $\msd^{1/2}(\R^m,\mbs(\R^m))$,
\[
w_{\lm,\xi,\ga}=-\vr K_{\lm,\xi,\ga}(\nabla\Ga(\psi_{\lm,\xi,\ga}))+o(\vr)
\] 
solves the equation $P_{\lm,\xi,\ga}\nabla\cl(\psi_{\lm,\xi,\ga}+w)=0$ with $P_{\lm,\xi,\ga}:\msd^{1/2}(\R^m,\mbs(\R^m))\to T_{\psi_{\lm,\xi,\ga}}\cm^\bot$ being the orthogonal projection 
and 
\begin{\equ}\label{widehat-Phi}
\widehat\Phi(\lm,\xi,\ga)=\Phi(\psi_{\lm,\xi,\ga})-\frac12\inp{K_{\lm,\xi,\ga}(\nabla\Ga(\psi_{\lm,\xi,\ga}))}{\nabla\Ga(\psi_{\lm,\xi,\ga})}_{\msd}.
\end{\equ}
For a fixed $h$, this second-order expansion holds on compact sets of the
bubble parameters $(\lm,\xi,\ga)$. Later, in Lemma~\ref{rev-local-estimate}, we will give the corresponding uniform estimates as $a,b\to0$.

\subsection{The role of the Weyl tensor and fixed-metric limits}\label{subsec:weyl}

In scalar equations such as the Yamabe equation and the constant
$Q$-curvature equations, the center and scale determine the standard
positive bubble. A standard spinorial bubble also depends on a unit
spinor $\ga$, and this parameter can affect even the sign of the
second-order term $\widehat\Phi(\lm,\xi,\ga)$ in \eqref{cl-red}. Formula~\eqref{widehat-Phi} applies to the perturbation
\eqref{hij} for any algebraic Weyl tensor $W$, including those in
Examples~\ref{rev-weyl-example} and~\ref{Ex-2} below. In particular,
Example~\ref{Ex-2} shows that for a single choice of $W$ and each
fixed $\lm>0$, the function $\ga\mapsto\widehat\Phi(\lm,0,\ga)$
attains both a positive maximum and a negative minimum on the unit
spinor sphere.

The following limits concern a fixed compactly supported tensor $h$.
For the polynomial in \eqref{rev-cubic-profile}, Proposition~\ref{rev-cubic-max}
will also locate critical points at a scale comparable to $b$ as
$b\to0$. Their type under translations depends on the Weyl tensor.

\begin{Lem}\label{lemma:widehat-Phi lm 0}
For a fixed smooth compactly supported symmetric trace-free tensor $h$ and $m\geq4$,
\[
\widehat\Phi(\lm,\xi,\ga)\to0\quad\hbox{as }\lm\to0,
\]
uniformly in $(\xi,\ga)\in\R^m\times S(\mbs_m)$.
\end{Lem}
\begin{proof}
Consider the two integrals defining $\Phi$
in Lemma~\ref{lemma: cl-expansion}. The first integral involves
Clifford multiplication by a three-form. Since each unit vector $e_i$
acts isometrically on spinors, the spinor pairing in each summand has
absolute value at most $|\psi_{\lm,\xi,\ga}(x)|_{\ig_{\R^m}}^2$.
Under the change of variables
$x=\xi+\lm y$, the bubble scales as
$\psi_{\lm,\xi,\ga}(\xi+\lm y)=\lm^{-(m-1)/2}\psi_{1,0,\ga}(y)$, so that
\[
 |\psi_{\lm,\xi,\ga}(x)|^2\,dx=\lm\,|\psi_{1,0,\ga}(y)|^2\,dy.
\]
Since the coefficients of this integral are products of the form
$(\pa_i h_{jd})h_{dk}$, its absolute value is bounded by
\[
 C_m\lm\,\sup_{\R^m}|\nabla h|\,\sup_{\R^m}|h|\,
 \|\psi_{1,0,\ga}\|_{L^2(\R^m,\mbs(\R^m))}^2.
\]
For every unit $\ga$, the bubble formula gives
$|\psi_{1,0,\ga}(y)|^2=m^{m-1}(1+|y|^2)^{1-m}$. This is integrable for $m>2$,
so the first integral is $O(\lm)$ uniformly in $\xi$ and $\ga$.

By \eqref{Phi-on-M}, the second integral, including its factor
$1/16$, becomes
\[
 \frac{m^{m-1}}{16}\int_{\R^m}
 \frac{\tr(h(\xi+\lm y)^2)}{(1+|y|^2)^m}\,dy.
\]
Since $h$ is smooth and compactly supported, $\tr(h^2)$ has bounded
gradient. Replacing $\tr(h(\xi+\lm y)^2)$ by $\tr(h(\xi)^2)$ thus
changes the numerator by at most $C\lm |y|$, uniformly in $\xi$.
As $|y|(1+|y|^2)^{-m}$ is integrable, this replacement contributes
an error of order $\lm$ after integration. Setting
$I_0=\int_{\R^m}(1+|y|^2)^{-m}\,dy$ and combining the two estimates
gives, uniformly in $(\xi,\ga)\in\R^m\times S(\mbs_m)$,
\begin{\equ}\label{convergence1}
 \Phi(\psi_{\lm,\xi,\ga})\to
 \frac{m^{m-1}}{16}\tr(h(\xi)^2)I_0
 \qquad\text{as }\lm\to0.
\end{\equ}

Linearizing the Euclidean spinorial Yamabe equation
$D_{\ig_{\R^m}}\psi=|\psi|_{\ig_{\R^m}}^{2/(m-1)}\psi$ at the bubble
$\psi_{1,0,\ga}$ gives the real-linear operator
\[
 L_{\psi_{1,0,\ga}}w
 =D_{\ig_{\R^m}}w-|\psi_{1,0,\ga}|_{\ig_{\R^m}}^{2/(m-1)}w
 -\frac{2}{m-1}|\psi_{1,0,\ga}|_{\ig_{\R^m}}^{(4-2m)/(m-1)}
   \real(\psi_{1,0,\ga},w)_{\ig_{\R^m}}\psi_{1,0,\ga}.
\]
This operator takes values in the dual of $\msd^{1/2}(\R^m,\mbs(\R^m))$.
Its Hilbert space representative is the Hessian
$\nabla^2\cl_0(\psi_{1,0,\ga})=|D_{\ig_{\R^m}}|^{-1}L_{\psi_{1,0,\ga}}$.
We now turn to the second term in \eqref{widehat-Phi}. Set
\[
 \bar w_{\lm,\xi,\ga}
 =-K_{\lm,\xi,\ga}\nabla\Ga(\psi_{\lm,\xi,\ga}).
\]
Since $\Ga$ vanishes on $\cm$ by Lemma~\ref{lemma: critical-manifold}(4),
its derivative vanishes on $T_{\psi_{\lm,\xi,\ga}}\cm$. Hence
$\nabla\Ga(\psi_{\lm,\xi,\ga})$ belongs to the orthogonal complement
of this tangent space. By the definition of $K_{\lm,\xi,\ga}$,
$\bar w_{\lm,\xi,\ga}$ is the unique spinor in that complement satisfying
\[
 \nabla^2\cl_0(\psi_{\lm,\xi,\ga})[\bar w_{\lm,\xi,\ga}]
 =-\nabla\Ga(\psi_{\lm,\xi,\ga}).
\]
In differential form, this equation reads
$L_{\psi_{\lm,\xi,\ga}}\bar w_{\lm,\xi,\ga}=-F_h(\psi_{\lm,\xi,\ga})$,
where $L_{\psi_{\lm,\xi,\ga}}$ is the linearization at
$\psi_{\lm,\xi,\ga}$ and $F_h$ is given by \eqref{rev-source}.

To compare these equations at different scales, write $x=\xi+\lm y$ again
and define the rescaled spinor by
\[
 \widetilde w_{\lm,\xi,\ga}(y)
 =\lm^{(m-1)/2}\bar w_{\lm,\xi,\ga}(\xi+\lm y).
\]
This rescaling preserves the bubble manifold and the
$\msd^{1/2}(\R^m,\mbs(\R^m))$ inner product, and takes
$\psi_{\lm,\xi,\ga}$ to $\psi_{1,0,\ga}$.
Thus $\widetilde w_{\lm,\xi,\ga}\in T_{\psi_{1,0,\ga}}\cm^\perp$.
Substituting \eqref{rev-source} and using
$\partial_{x^j}=\lm^{-1}\partial_{y^j}$ gives the exact equation
\[
 L_{\psi_{1,0,\ga}}\widetilde w_{\lm,\xi,\ga}(y)
 =\frac{\lm}{4}\sum_{i,j}(\partial_jh_{ij})(\xi+\lm y)
       e_i\cdot_{\ig_{\R^m}}\psi_{1,0,\ga}(y)
 +\frac12\sum_{i,j}h_{ij}(\xi+\lm y)
       e_i\cdot_{\ig_{\R^m}}\partial_j\psi_{1,0,\ga}(y).
\]
The first sum tends to zero because of the factor $\lm$. In the
second sum we may replace $h_{ij}(\xi+\lm y)$ by $h_{ij}(\xi)$
with an error bounded by $\lm |y|\sup_{\R^m}|\nabla h|$ in each
coefficient. Indeed, by the bubble formula and \eqref{pa-j-psi},
both $|\psi_{1,0,\ga}(y)|_{\ig_{\R^m}}$ and
$|y|\,|\partial_j\psi_{1,0,\ga}(y)|_{\ig_{\R^m}}$ are bounded by
$C_m(1+|y|)^{1-m}$ for every unit $\ga$. Since $m>3$, this bound
belongs to $L^{2m/(m+1)}(\R^m)$. The right-hand side therefore
converges in $L^{2m/(m+1)}(\R^m,\mbs(\R^m))$, and hence in the
dual of $\msd^{1/2}(\R^m,\mbs(\R^m))$, uniformly in $\xi$ and $\ga$.

Recall that $K_{1,0,\ga}$ is the inverse of the Hessian
$\nabla^2\cl_0(\psi_{1,0,\ga})$ restricted to
$T_{\psi_{1,0,\ga}}\cm^\perp$. Its operator norm depends continuously
on $\ga$ and is therefore bounded on the compact unit spinor sphere.
Since the operator $L_{\psi_{1,0,\ga}}$ in the rescaled equation
is independent of $\lm$, this uniform bound and the convergence of
the right-hand side in the dual norm imply that
$\widetilde w_{\lm,\xi,\ga}\to w_{0,\xi,\ga}$ in
$\msd^{1/2}(\R^m,\mbs(\R^m))$ as $\lm\to0$, uniformly in $\xi$ and
$\ga$, where $w_{0,\xi,\ga}$ is the unique solution in
$T_{\psi_{1,0,\ga}}\cm^\perp$ of
\begin{\equ}\label{equ-w0}
L_{\psi_{1,0,\ga}}w_{0,\xi,\ga}
=\frac12\sum_{i,j}h_{ij}(\xi)e_i\cdot_{\ig_{\R^m}}\partial_j\psi_{1,0,\ga}.
\end{\equ}

To evaluate the limiting pairing, observe that
\[
 y\mapsto\frac12\sum_{i,j}h_{ij}(\xi)y^i\partial_j\psi_{1,0,\ga}(y)
\]
is a particular solution of \eqref{equ-w0}. Indeed, the translation directions satisfy
$L_{\psi_{1,0,\ga}}(\partial_j\psi_{1,0,\ga})=0$ and the product rule gives
\[
 L_{\psi_{1,0,\ga}}(y^i\partial_j\psi_{1,0,\ga})
 =e_i\cdot_{\ig_{\R^m}}\partial_j\psi_{1,0,\ga}.
\]
This particular solution differs from $w_{0,\xi,\ga}$ by an element
of $\ker L_{\psi_{1,0,\ga}}$. For every $\eta$ in this kernel, testing
\eqref{equ-w0} with $\eta$ gives
\[
 \frac12\sum_{i,j}h_{ij}(\xi)\int_{\R^m}
 \real\bigl(e_i\cdot_{\ig_{\R^m}}\partial_j\psi_{1,0,\ga},\eta
 \bigr)_{\ig_{\R^m}}\,dy=0,
\]
by the symmetry of $L_{\psi_{1,0,\ga}}$ and the identity
$L_{\psi_{1,0,\ga}}\eta=0$. Thus either solution gives the same
value when integrated against the right-hand side of \eqref{equ-w0}.
By definition, the pairing in \eqref{widehat-Phi} equals
$-\Ga'(\psi_{\lm,\xi,\ga})[\bar w_{\lm,\xi,\ga}]$.
After rescaling, its limit is the pairing of the right-hand side of
\eqref{equ-w0} with $w_{0,\xi,\ga}$. Substituting the particular
solution and using \eqref{pa-j-psi} reduces this limit to a radial integral,
\begin{\equ}\label{limitK}
\lim_{\lm\to0}
\inp{K_{\lm,\xi,\ga}\nabla\Ga(\psi_{\lm,\xi,\ga})}
{\nabla\Ga(\psi_{\lm,\xi,\ga})}_{\msd}
=\frac{m^{m-1}}4\tr(h(\xi)^2)
 \int_{\R^m}\frac{|y|^2}{(1+|y|^2)^{m+1}}\,dy.
\end{\equ}
To evaluate the integral in terms of $I_0$, write $r=|y|$ and set
$u=r^2/(1+r^2)$ in both radial integrals. This gives
\begin{\equ}\label{BB-Euler}
\int_{\R^m}\frac{|y|^2}{(1+|y|^2)^{m+1}}\,dy
 =\frac{B(m/2+1,m/2)}{B(m/2,m/2)}I_0,
\end{\equ}
where $B(s,t)$ is the beta function, defined by Euler's integral
of the first kind,
\[
 B(s,t)=\int_0^1 u^{s-1}(1-u)^{t-1}\,du,\qquad s,t>0.
\]
The property 
\[
B(s+1,t)=\frac{s}{s+t}B(s,t)
\] 
gives a ratio of $1/2$ in \eqref{BB-Euler}. Thus the limit in \eqref{limitK} equals $\displaystyle \frac{m^{m-1}}8\tr(h(\xi)^2)I_0$.
Substituting this limit and \eqref{convergence1} into
\eqref{widehat-Phi} proves the assertion.
\end{proof}

Lemma~\ref{lemma:widehat-Phi lm 0} shows that
$\widehat\Phi(\lm,\xi,\ga)\to0$ as $\lm\to0$. When
$\lm+|\xi|\to\infty$, the bubble and its first derivatives decay
uniformly on $\supp h$. Combining this decay with the uniform bound
on $K_{1,0,\ga}$ used in the proof of the Lemma~\ref{lemma:widehat-Phi lm 0} gives the following corollary.

\begin{Cor}\label{corollary: widehat-Phi-0}
Let $m\geq4$ and let $h$ be a fixed smooth, compactly supported,
symmetric trace-free tensor on $\R^m$. Then
$\widehat\Phi(\lm,\xi,\ga)\to0$ as $\lm+|\xi|\to\infty$,
uniformly for $\ga\in S(\mbs_m)$.
\end{Cor}
\begin{proof}
Choose $R\geq1$ such that $\supp h\subset B_R(0)$. Both terms in
the formula for $\Phi(\psi_{\lm,\xi,\ga})$ in
Lemma~\ref{lemma: cl-expansion} are integrals over this fixed ball.
We first estimate the bubble and its first derivatives there.
By \eqref{bubble} and \eqref{pa-j-psi}, for $k=0,1$ and
$\ga\in S(\mbs_m)$,
\[
 |\nabla^k\psi_{\lm,\xi,\ga}(x)|_{\ig_{\R^m}}
 \leq C_m
 \frac{\lm^{(m-1)/2}}{(\lm^2+|x-\xi|^2)^{(m-1+k)/2}}
 \leq \frac{C_m}{(\lm+|x-\xi|)^{(m-1)/2+k}}.
\]
For $x\in B_R(0)$ and $\lm+|\xi|\geq2R$, the triangle inequality
gives $\lm+|x-\xi|\geq\lm+|\xi|-R\geq(\lm+|\xi|)/2$.
Hence these bounds tend to zero throughout $B_R(0)$ as
$\lm+|\xi|\to\infty$, uniformly in the unit spinor $\ga$.

In the Clifford three-form term of $\Phi(\psi_{\lm,\xi,\ga})$, each
spinor pairing is bounded by $|\psi_{\lm,\xi,\ga}|_{\ig_{\R^m}}^2$.
In the term containing $h^2$, it is bounded by
$|\nabla\psi_{\lm,\xi,\ga}|_{\ig_{\R^m}}\,|\psi_{\lm,\xi,\ga}|_{\ig_{\R^m}}$.
Since the coefficients are bounded and the ball has finite volume,
the preceding estimates give
\[
 |\Phi(\psi_{\lm,\xi,\ga})|\leq C_{m,h,R}(\lm+|\xi|)^{1-m}
 \qquad\text{for }\lm+|\xi|\geq2R.
\]

It remains to estimate the term involving $\nabla\Ga$ in
\eqref{widehat-Phi}. The source $F_h(\psi_{\lm,\xi,\ga})$ in
\eqref{rev-source} is supported in $B_R(0)$. Its first sum contains
$(\partial_j h_{ij})\psi_{\lm,\xi,\ga}$ and its second contains
$h_{ij}\nabla_{e_j}\psi_{\lm,\xi,\ga}$. Clifford multiplication by
$e_i$ preserves the norm of each summand. The same pointwise bounds, followed by
integration on $B_R(0)$, therefore give
\[
 \|F_h(\psi_{\lm,\xi,\ga})\|_{L^{2m/(m+1)}(\R^m,\mbs(\R^m))}
 \leq C_{m,h,R}(\lm+|\xi|)^{-(m-1)/2}.
\]
By H\"older's inequality and the Sobolev embedding into
$L^{2m/(m-1)}(\R^m,\mbs(\R^m))$, this bounds the dual norm of
$\Ga'(\psi_{\lm,\xi,\ga})$, up to a constant depending only on $m$.
That dual norm equals
$\|\nabla\Ga(\psi_{\lm,\xi,\ga})\|_{\msd^{1/2}(\R^m,\mbs(\R^m))}$.

Finally, translation and dilation are isometries of
$\msd^{1/2}(\R^m,\mbs(\R^m))$ and identify $K_{\lm,\xi,\ga}$ with
$K_{1,0,\ga}$. The latter has uniformly bounded operator norm on the
unit spinor sphere, as in the proof of Lemma~\ref{lemma:widehat-Phi lm 0}.
Thus the absolute value of the second term in \eqref{widehat-Phi}
is bounded by a constant times
$\|\nabla\Ga(\psi_{\lm,\xi,\ga})\|_{\msd^{1/2}(\R^m,\mbs(\R^m))}^2$.
Combining the two estimates gives
\[
 \sup_{\ga\in S(\mbs_m)}|\widehat\Phi(\lm,\xi,\ga)|
 \leq C_{m,h,R}(\lm+|\xi|)^{1-m}
 \qquad\text{for }\lm+|\xi|\geq2R,
\]
which proves the assertion.
\end{proof}

\begin{Rem}\label{remark:unified-domain-lm-xi-ga}
In Lemma~\ref{lemma:widehat-Phi lm 0} and
Corollary~\ref{corollary: widehat-Phi-0}, the tensor $h$ was fixed.
We now return to the family \eqref{hij}, keeping the Weyl tensor $W$,
the polynomial $f$, and the cut-off function $\eta$ fixed while varying
$a,b,r$. Write $\Ga_{a,b,r}$ and $\Phi_{a,b,r}$ for the functionals
in Lemma~\ref{lemma: cl-expansion} evaluated with this choice of $h$.
We write $\widehat\Phi_{a,b,r}$ for the second-order reduced functional
defined by \eqref{widehat-Phi} with $\Ga_{a,b,r}$ and $\Phi_{a,b,r}$.
Since $\Ga_{a,b,r}$ vanishes on the bubble manifold $\cm$, the
expansion \eqref{cl-red} gives
\[
 \cl^{red}(\lm,\xi,\ga)-\cl_0(\psi_{\lm,\xi,\ga})
 =\vr^2\widehat\Phi_{a,b,r}(\lm,\xi,\ga)+o(\vr^2).
\]
To express $\widehat\Phi_{a,b,r}$ on the length scale $b$, set $x=by$.
Since $H_{ij}$ is homogeneous of degree two, \eqref{hij} gives
\[
 h_{ij}(by)=ab^2\,\eta\!\left(\frac{|y|}{r/b}\right)f(|y|^2)H_{ij}(y).
\]
Thus the rescaled tensor is $ab^2$ times the tensor obtained from
\eqref{hij} with parameters $(1,1,r/b)$. For $b>0$, define the
dilation operator $U_b$ on $\msd^{1/2}(\R^m,\mbs(\R^m))$ by
\[
 (U_b\varphi)(y)=b^{(m-1)/2}\varphi(by),\qquad y\in\R^m.
\]
The bubble formula \eqref{bubble} then gives
$U_b\psi_{\lm,\xi,\ga}=\psi_{\lm/b,\xi/b,\ga}$.
Substituting $x=by$ in the integrals defining $\Ga$ and $\Phi$ gives
\[
 \Ga_{a,b,r}(\varphi)=ab^2\Ga_{1,1,r/b}(U_b\varphi),\qquad
 \Phi_{a,b,r}(\varphi)=a^2b^4\Phi_{1,1,r/b}(U_b\varphi).
\]
Indeed, the two spinor factors contribute $b^{-(m-1)}$, the single
derivative contributes $b^{-1}$, and $dx=b^m\,dy$. The remaining
powers of $ab^2$ follow from the linear dependence of $\Ga$ and the
quadratic dependence of $\Phi$ on $h$. The dilation operator $U_b$ is an isometry of $\msd^{1/2}(\R^m,\mbs(\R^m))$.
It carries $\nabla\Ga_{a,b,r}(\psi_{\lm,\xi,\ga})$ to
$ab^2\nabla\Ga_{1,1,r/b}(\psi_{\lm/b,\xi/b,\ga})$ and identifies
the inverse Hessian $K_{\lm,\xi,\ga}$ with $K_{\lm/b,\xi/b,\ga}$.
Since the gradient occurs twice in the second term of \eqref{widehat-Phi},
that term also acquires the factor $a^2b^4$. We obtain
\begin{\equ}\label{ab-a1}
\widehat\Phi_{a,b,r}(\lm,\xi,\ga)
=a^2b^4\widehat\Phi_{1,1,r/b}(\lm/b,\xi/b,\ga).
\end{\equ}
In the rescaled tensor, $r/b$ appears only as the cut-off radius.
We will keep $r$ fixed and let $b\to0$, so this radius tends to infinity.
For fixed $a,b,r$, Lemma~\ref{lemma:widehat-Phi lm 0} gives
$\widehat\Phi_{a,b,r}(\lm,\xi,\ga)\to0$ as $\lm\to0$, and
Corollary~\ref{corollary: widehat-Phi-0} gives the same limit as
$\lm+|\xi|\to\infty$. The small-scale estimates use the suprema of
$|h|$ and $|\nabla h|$. In the decay estimate at the end of the
corollary, the constant $C_{m,h,R}$ also depends on the radius $R$
of a ball containing $\supp h$. For the rescaled tensor in
\eqref{ab-a1}, we may take $R=r/b$. Thus the two results give
convergence for each fixed $r/b$, with no uniform rate established
as $r/b\to\infty$. For
the cubic polynomial in \eqref{rev-cubic-profile} and $m\geq20$,
we will study the limit on compact neighborhoods of the centered
critical points in the rescaled variables.
\end{Rem}

We now return to the notation $\widehat\Phi$ of \eqref{widehat-Phi}.
To find a positive maximum of $\widehat\Phi$, it suffices to exhibit
one bubble at which $\widehat\Phi$ is positive. Indeed, for fixed
$a,b,r$, the decay
in Lemma~\ref{lemma:widehat-Phi lm 0} and
Corollary~\ref{corollary: widehat-Phi-0}, together with the compactness
of $S(\mbs_m)$, ensures that the maximum is attained. The same
argument gives a negative minimum whenever $\widehat\Phi$ takes a
negative value.

We first look for such values among bubbles centered at the origin.
By Lemma~\ref{lemma: critical-manifold}(4),
$\nabla\Ga(\psi_{\lm,0,\ga})=0$, so \eqref{widehat-Phi} reduces to
$\widehat\Phi(\lm,0,\ga)=\Phi(\psi_{\lm,0,\ga})$. Part~(5) of the
same lemma then shows that the dependence on the unit spinor $\ga$
enters entirely through $P(W,\ga)$. The two examples below illustrate
the role of this term. In the first, $P(W,\ga)$ vanishes for every
$\ga$. For the nonzero perturbation in the second example, varying
$\ga$ alone gives both positive and negative values of
$\widehat\Phi(\lm,0,\ga)$ at each fixed scale $\lm$.
These algebraic calculations apply to any radial factor $f(|x|^2)$
in \eqref{hij}. For our first example, consider the following two
conditions on $W$.
\begin{itemize}
\item[(W1)] $W_{ijil}=0$ for all $i,j,l$;
\item[(W2)] $W_{ijkl}=0$ if any index is at least $6$.
\end{itemize}

\begin{Ex}\label{rev-weyl-example}
	A rather simple and clear example is as follows. Consider the Weyl-type tensor $W=\{W_{ijkl}\}$ whose non-zero components only arise when the indices $i, j, k, l \in \{1, 2, 3, 4\}$. Under condition (W1), to guarantee the nontriviality of the tensor $W$, it suffices to assign values to $W_{1234}$, $W_{1342}$, and $W_{1423}$ --- with the proviso that such assignments must satisfy the Bianchi identity $W_{1234}+W_{1342}+W_{1423}=0$. See the following {Table~\ref{tab:weyl-vanishing-spinor-term}} for instance.
	\begin{center}
		\refstepcounter{table}\label{tab:weyl-vanishing-spinor-term}
		\renewcommand{\arraystretch}{1.3}
		\begin{tabular}{p{2.2cm} p{2.2cm} p{2.2cm}p{2.2cm} p{2.2cm} p{2.2cm}}
			\hline
			\multicolumn{6}{c}{{\tablename~\thetable}: A specific example of $W_{ijkl}$ satisfying (W1) and (W2)} \\
			\hline 
			$W_{1212}=0$ & $W_{1313}=0$ & $W_{1414}=0$ & $W_{2323}=0$ & $W_{2424}=0$ & $W_{3434}=0$ \\
			$W_{1213}=0$ & $W_{1314}=0$ & $W_{1423}=-2$ & $W_{2324}=0$ & $W_{2434}=0$ &  \\
			$W_{1214}=0$ & $W_{1323}=0$ & $W_{1424}=0$ & $W_{2334}=0$ & & \\
			$W_{1223}=0$ & $W_{1324}=-1$ & $W_{1434}=0$ & & & \\
			$W_{1224}=0$ & $W_{1334}=0$ &  & & & \\
			$W_{1234}=1$ &  &  & & & \\
			\hline
		\end{tabular}
	\end{center}
\end{Ex}

\medskip

More generally, beyond the case listed in {Table~\ref{tab:weyl-vanishing-spinor-term}}, the algebraic symmetries of $W$, together with (W1), imply that any component of $W$ with a repeated index vanishes. Moreover, the Clifford relations imply that $S_{ijkq}(\ga)=0$ whenever $i,j,k$ are pairwise distinct and $q\in\{i,j,k\}$. Thus, a nonzero summand in $P(W,\ga)$ would require $i,j,k,q$ to be pairwise distinct. The condition $W_{jpid}W_{kpqd}\neq0$ would further force $p$ and $d$ to be distinct from each other and from all four indices $i,j,k,q$. Hence all six indices $i,j,k,q,p,d$ would have to be pairwise distinct, contradicting the support condition (W2), which requires them to lie in $\{1,2,3,4,5\}$. Consequently, every summand in $P(W,\ga)$ vanishes, and therefore $\Phi(\psi_{\lm,0,\ga})$ can be expanded as
\begin{\equ}\label{W1W2-Phi}
\Phi(\psi_{\lm,0,\ga})
=\frac{m^{m-1}\lm^m}{16}
 \int_{\R^m}\frac{\tr(h^2)(x)}{(\lm^2+|x|^2)^m}\,dx
=\frac{|W|^2}{32}I_{12}>0.
\end{\equ}
For the last equality, expand $H_{ij}(x)=\sum_{p,q}W_{ipjq}x^p x^q$
and use the values of $I_{ijkl}$ in \eqref{I12} together with
$\sum_p W_{ipjp}=0$. For {Table~\ref{tab:weyl-vanishing-spinor-term}}, $|W|^2=144$, so the value is
$9I_{12}/2$. This is strictly positive and independent of $\ga$. At the center, $\widehat\Phi=\Phi$ by Lemma \ref{lemma: critical-manifold}.
For fixed $a,b,r$, Lemma~\ref{lemma:widehat-Phi lm 0} and
Corollary~\ref{corollary: widehat-Phi-0} give
$\widehat\Phi(\lm,0,\ga)\to0$ both as $\lm\to0$ and as $\lm\to\infty$.
Thus the positive integral in \eqref{W1W2-Phi} attains a maximum at a
finite scale, even for a linear $f$ with $f(0)=1$ in dimension $20$.
Since the weight $\tr(h^2)$ has fixed compact support and $\lm$ enters
only through $\lm^m/(\lm^2+|x|^2)^m$, the integral is real analytic in
$\lm>0$. It is nonconstant, so each maximizing scale is a strict local maximum.
For the blow-up construction we also need to control the location of
the extremum as $b\to0$, which is why Remark~\ref{rem:cubic-coefficient-range}
considers the relative scale $\lm/b$.

\medskip

%It is worth noting that, when the terms such as $W_{ijil}$ are non-zero (violating the condition (W1)), we are able to construct nontrivial examples of $C_{ijkl}$ creating more critical points for $\Phi(\psi_{\lm,\xi,\ga})$.

In the next example some components $W_{ijil}$ are nonzero. For this
tensor the sign of the centered energy changes as the unit spinor varies.
\begin{Ex}\label{Ex-2}
Let us consider the Weyl-type tensor $W=\{W_{ijkl}\}$, for which the non-zero components appear only with indices $i, j, k, l \in \{1, 2, 3, 4\}$ (namely, we assume $W_{ijkl} = 0$ if any of its indices is at least $5$). In this case, the tensor has a total of $21$ components to specify, although only $10$ of these are independent due to the trace-free property and the Bianchi identity.
\begin{center}
	\refstepcounter{table}\label{tab:weyl-indefinite-spinor-term}
	\renewcommand{\arraystretch}{1.3}
\begin{tabular}{p{2.2cm} p{2.2cm} p{2.2cm}p{2.2cm} p{2.2cm} p{2.2cm}}
	\hline
	\multicolumn{6}{c}{{\tablename~\thetable}: A specific example of $W_{ijkl}$ that generates nontrivial $P(W,\ga)$} \\
	\hline 
	$W_{1212}=0$ & $W_{1313}=0$ & $W_{1414}=0$ & $W_{2323}=0$ & $W_{2424}=0$ & $W_{3434}=0$ \\
	$W_{1213}=0$ & $W_{1314}=1$ & $W_{1423}=0$ & $W_{2324}=-1$ & $W_{2434}=0$ &  \\
	$W_{1214}=1$ & $W_{1323}=1$ & $W_{1424}=-1$ & $W_{2334}=1$ & & \\
	$W_{1223}=1$ & $W_{1324}=0$ & $W_{1434}=1$ & & & \\
	$W_{1224}=1$ & $W_{1334}=-1$ &  & & & \\
	$W_{1234}=0$ &  &  & & & \\
	\hline
\end{tabular}
\end{center}
\end{Ex}

\medskip

Using the datum specified above, we compute $P(W,\ga)$ as follows. Since $S_{ijkq}(\ga)=0$ whenever $i,j,k$ are pairwise distinct and $q\in\{i,j,k\}$, only permutations of $\{1,2,3,4\}$ contribute. Direct substitution of the prescribed components gives $P(W,\ga)=-288S_{1234}(\ga)$, $\,|W|^2=240$ and
\begin{eqnarray}\label{nontrivial-Phi}
\Phi(\psi_{\lm,0,\ga})&=&\frac{m^{m-1}\lm^m}{16}\int_{\R^m}\frac{\tr(\ih^2)(x)}{(\lm^2+|x|^2)^m}d\vol_{\ig_{\R^m}}  \nonumber  \\[0.3em]
& & -18m^{m-1}\lm^m\int_{\R^m}\frac{a^2\eta(r^{-1}|x|)^2f(b^{-2}|x|^2)^2|x^1|^2 |x^2|^2}{(\lm^2+|x|^2)^m} d\vol_{\ig_{\R^m}} S_{1234}(\ga)    \nonumber\\[0.3em]
&=&\frac{m^{m-1}\lm^m}{16}\sum_{i,d}W_{i\al d\bt}W_{ikdl}\int_{\R^m}\frac{a^2\eta(r^{-1}|x|)^2f(b^{-2}|x|^2)^2x^{\al}x^{\bt} x^k x^l}{(\lm^2+|x|^2)^m}d\vol_{\ig_{\R^m}} \nonumber\\[0.5em]   
& & -18I_{12} S_{1234}(\ga)    \nonumber\\[0.5em]
&=&\frac{1}{16}\sum_{i,d} \bigg[I_{12}\sum_{k>l} (W_{i k d l}+ W_{ildk})^2 +I_{12}\sum_{k\neq l}W_{ikdk}W_{ildl} + I_{11}\sum_kW_{ikdk}^2   \bigg] \nonumber\\[0.5em]
& & -18I_{12} S_{1234}(\ga)    \nonumber\\[0.5em]
&=&\frac{1}{16}\sum_{i,d} \bigg[ I_{12}\sum_{k>l} (W_{i k d l}+ W_{ildk})^2   + (I_{11}-I_{12})\sum_kW_{ikdk}^2  \bigg] -18I_{12} S_{1234}(\ga)    \nonumber\\[0.3em]
&=&\frac{I_{12}}{16}\sum_{i,d} \bigg[ \sum_{k>l} (W_{i k d l}+ W_{ildk})^2   + 2\sum_kW_{ikdk}^2  \bigg] -18I_{12} S_{1234}(\ga)    \nonumber\\[0.3em]
&=&\frac{I_{12}}{32}|W|^2-18I_{12} S_{1234}(\ga) \nonumber\\[0.5em]
%&=&m^{m-1}\lm^m\Big[\frac{15}2-18S_{1234}(\ga)\Big]\int_{\R^m}\frac{a^2\eta(r^{-1}|x|)^2f(b^{-2}|x|^2)^2|x^1|^2|x^{2}|^2}{(\lm^2+|x|^2)^m}d\vol_{\ig_{\R^m}}.
&=&\Big[\frac{15}2-18S_{1234}(\ga)\Big] I_{12},
\end{eqnarray}
where $I_{12}$ stands for the integral given in \eqref{I12}.
These constants follow by finite contraction of the prescribed tensor entries, after imposing the curvature symmetries. The resulting identities are algebraic and do not depend on the choice of the radial polynomial.

For {Table~\ref{tab:weyl-vanishing-spinor-term}}, \eqref{W1W2-Phi} is independent of the unit spinor.
For {Table~\ref{tab:weyl-indefinite-spinor-term}}, the dependence of
$\widehat\Phi(\lm,0,\ga)=\Phi(\psi_{\lm,0,\ga})$ on $\ga$ enters
through the single factor $S_{1234}(\ga)$ in \eqref{nontrivial-Phi}.
Its range determines both extreme values at each fixed scale.
\begin{Fact}\label{fact}
For $m\geq4$, the function $\ga\mapsto S_{1234}(\ga)$, which is defined for a unit spinor $\ga\in S(\mbs_m)$, has a range of $[-1,1]$.
\end{Fact}
\begin{proof}
Let $A=e_1\cdot e_2\cdot e_3\cdot e_4$. Then $A^*=A$ and $A^2=1$.
Since $e_1A=-Ae_1$, the invertible map $e_1$ exchanges its two eigenspaces, so both are nonzero.
Choose unit eigenspinors $\ga_\pm$ with $A\ga_\pm=\pm\ga_\pm$.
Since they are orthogonal,
\[
 S_{1234}(\cos\theta\,\ga_++\sin\theta\,\ga_-)
 =\cos^2\theta-\sin^2\theta=\cos(2\theta).
\]
As $\theta$ runs from $0$ to $\pi/2$, this gives every value in $[-1,1]$.
The Cauchy-Schwarz inequality gives $|S_{1234}(\ga)|\leq1$ for every
unit spinor $\ga$, since $A$ is an isometry.
\end{proof}

For the fixed perturbation in Example~\ref{Ex-2}, formula~\eqref{nontrivial-Phi}
and Fact~\ref{fact} show that, at each fixed scale $\lm>0$,
$\widehat\Phi(\lm,0,\ga)$ takes both positive and negative values as
$\ga$ varies over the unit spinor sphere. By Lemma~\ref{lemma:widehat-Phi lm 0}
and Corollary~\ref{corollary: widehat-Phi-0}, $\widehat\Phi(\lm,\xi,\ga)$ tends uniformly
to zero as $\lm\to0$ or $\lm+|\xi|\to\infty$. Continuity of $\widehat\Phi$
and compactness of the unit spinor sphere therefore imply that
$\widehat\Phi$ attains a positive maximum and a negative minimum on
the full parameter space. We next
examine how to choose $f$ so that there are critical points for which
the ratio $\lm/b$ stays between two fixed positive constants as $b\to0$.

\subsection{The role of the cubic polynomial in locating critical points}
\label{rev-cubic-subsection}

\begin{revision}
We locate the critical scale by first letting the rescaled cut-off
radius tend to infinity. For a bubble centered at the origin, the
limiting energy is a product of two factors. The cubic polynomial $f$
determines the dependence on the scale, and the Weyl tensor determines
a quadratic form on the unit spinors.
\end{revision}

{Fix a nonzero algebraic Weyl tensor $W=(W_{ijkl})$ on the
Euclidean vector space $\R^m$. Its components in the standard basis
serve as the constant coefficients of the quadratic polynomials
$H_{ij}(x)=W_{i\al j\bt}x^\al x^\bt$
in \eqref{hij}. For a concrete choice of $W$, one may use either
of the examples in Tables~\ref{tab:weyl-vanishing-spinor-term}
and~\ref{tab:weyl-indefinite-spinor-term}.}
{Recall from Remark~\ref{remark:unified-domain-lm-xi-ga}
that the cut-off radius in the rescaled coordinates $x/b$ is $R=r/b$.
Setting $a=b=1$ and taking $R$ as the cut-off radius in \eqref{hij}
gives the second-order reduced functional $\widehat\Phi_{1,1,R}$
defined by \eqref{ab-a1} with
$h_{ij}(x)=\eta(|x|/R)f(|x|^2)H_{ij}(x)$. The cut-off factor $\eta$ equals
$1$ for $|x|\le R/2$ and vanishes for $|x|\ge R$. For $m\ge18$, we can
remove the cut-off by considering}
\begin{\equ}\label{Phi-infty}
 \widehat\Phi_\infty(\lm,\xi,\ga)
 =\lim_{R\to\infty}\widehat\Phi_{1,1,R}(\lm,\xi,\ga).
\end{\equ}

\begin{revision}
For the uncut tensor $h(x)=f(|x|^2)H(x)$, the radial integrands in
\eqref{Phi-on-M} and \eqref{Phi-on-M2} are respectively $O(\rho^{15-m})$
and $O(\rho^{16-m})$ as $\rho=|x|\to\infty$. These bounds and the
corresponding estimate for the source in \eqref{rev-source} give
convergence uniformly on compact sets of bubble parameters when
$m\ge18$. The estimates are given in Appendix~\ref{app:cubic-translations}.
\end{revision}

{For $m\ge18$, let us set}
\begin{\equ}\label{rev-moments}
 A_m=\int_0^\infty\frac{\rho^{m+3}}{(1+\rho^2)^m}\,d\rho,
 \quad M_0=1, \quad
 M_k=\prod_{j=1}^k\frac{m+2j+2}{m-2j-4}\quad \text{for }1\le k\le6,
\end{\equ}
and
\[
 \kappa_m=\frac{m^{m-1}}{16}
 \frac{\om_{m-1}}{m(m+2)}A_m>0.
\]
For the cubic polynomial $f$ in \eqref{rev-cubic-profile}, write
\[
 f(s)^2=q_0+q_1s+q_2s^2+q_3s^3+q_4s^4+q_5s^5+q_6s^6.
\]
Expanding the square gives
\begin{\equ}\label{rev-cubic-square}
 (q_0,\ldots,q_6)=
 \left(1,-\frac97,\frac{1557}{2548},-\frac{85}{637},
       \frac{801}{57967},-\frac{36}{57967},\frac4{405769}\right).
\end{\equ}
For a bubble centered at the origin, Lemma~\ref{lemma: critical-manifold}~(4)
gives $\nabla\Ga(\psi_{\lm,0,\ga})=0$ for every cut-off radius.
Thus the inverse-Hessian term in \eqref{widehat-Phi} vanishes in the
limit as well, and to evaluate $\widehat\Phi_\infty$ it suffices to evaluate the two integrals
in Lemma~\ref{lemma: critical-manifold}~(5) with the cut-off removed. Write
$x=\rho\theta$ in polar coordinates, where $\rho=|x|$ and
$\theta=(\theta^1,\ldots,\theta^m)\in S^{m-1}$ is the unit vector in
the direction of $x$. We write $d\theta$ for surface measure on $S^{m-1}$.
The angular integrals satisfy the same index relations as in \eqref{I12}.
Integrating $(\sum_i(\theta^i)^2)^2=1$ fixes their common normalization.
Contracting with $W$ as in \eqref{W1W2-Phi}, using its trace-free property,
then gives
\[
 \int_{S^{m-1}}\tr(H(\theta)^2)\,d\theta
 =\frac{\om_{m-1}}{2m(m+2)}|W|^2,\qquad
 \int_{S^{m-1}}(\theta^1\theta^2)^2\,d\theta
 =\frac{\om_{m-1}}{m(m+2)} .
\]
Since $H(\rho\theta)=\rho^2H(\theta)$, both $\tr(H(x)^2)$ and
$(x^1x^2)^2$ contribute a factor of $\rho^4$. Together with the volume
element $\rho^{m-1}\,d\rho\,d\theta$, this gives
\begin{\equ}\label{rev-central-profile0}
 \widehat\Phi_\infty(\lm,0,\ga)
 =\frac{\kappa_m}{A_m}\Big(\frac12|W|^2+P(W,\ga)\Big)
 \lm^m\int_0^\infty
 \frac{\rho^{m+3}f(\rho^2)^2}{(\lm^2+\rho^2)^m}\,d\rho .
\end{\equ}
Expand $f(\rho^2)^2=\sum_{k=0}^6q_k\rho^{2k}$ and set $\rho=\lm\tau$.
For each $0\le k\le6$ the corresponding radial integral becomes
\begin{\equ}\label{rev-radial-integral-k}
 \lm^m\int_0^\infty\frac{\rho^{m+3+2k}}{(\lm^2+\rho^2)^m}\,d\rho
 =\lm^{2k+4}\int_0^\infty\frac{\tau^{m+3+2k}}{(1+\tau^2)^m}\,d\tau
 =\lm^{2k+4}A_mM_k .
\end{\equ}
Indeed, $s=\tau^2$ gives $\frac12 B(m/2+2+k,m/2-2-k)$ for the integral
in $\tau$. Its value at $k=0$ is $A_m$, and the ratio of consecutive
values is $(m+2k+2)/(m-2k-4)$ for $1\le k\le6$, as in \eqref{rev-moments}.
Since $\lm^{2k+4}=(\lm^2)^{k+2}$, setting
$P_m(t)=\sum_{k=0}^6q_kM_kt^{k+2}$ gives
\begin{\equ}\label{rev-central-profile}
 \widehat\Phi_\infty(\lm,0,\ga)
 =\kappa_m\Big(\frac12|W|^2+P(W,\ga)\Big)P_m(\lm^2).
\end{\equ}
In this formula, $\tfrac12|W|^2+P(W,\ga)$ depends only on the Weyl
tensor $W$ and the unit spinor $\ga$. The dependence on the bubble
scale $\lm$ is entirely through $P_m(\lm^2)$.

For fixed $\lm>0$, the integrand on the left of \eqref{rev-radial-integral-k}
is asymptotic to $\rho^{3+2k-m}$ as $\rho\to\infty$.
For the choice of $f$ in \eqref{rev-cubic-profile}, the coefficient $q_6$
is positive, so the term with $k=6$ gives the leading power
$\rho^{15-m}$ in the full radial integrand. {Thus all seven
integrals in \eqref{rev-radial-integral-k} converge precisely when
$m>16$, in particular throughout the range $m\ge18$. We also remark that, by comparing \eqref{rev-central-profile0} with \eqref{rev-central-profile}, $P_m(\lm^2)$ is strictly positive for $\lm>0$.}
For fixed $W$ and $\ga$ with $\tfrac12|W|^2+P(W,\ga)\ne0$,
\eqref{rev-central-profile} shows that $\lm>0$ is a critical point of
$\lm\mapsto\widehat\Phi_\infty(\lm,0,\ga)$ if and only if $P_m'(\lm^2)=0$.
{Since $P_m'(t)$ is divisible by $t$ and we seek only
positive roots $t=\lm^2$, it suffices to consider the polynomial}
\begin{\equ}\label{rev-radial-polynomials}
 Q_m(t)=\frac{P_m'(t)}t
       =\sum_{k=0}^6(k+2)q_kM_kt^k .
\end{\equ}
Thus each positive zero $t$ of $Q_m$ gives a critical scale $\lm=\sqrt t$.
At such a zero, $P_m''(t)=tQ_m'(t)$, so $Q_m'(t)<0$ gives a strict
local maximum of $P_m$.

At this point the center is still fixed. For the construction we must
restore the cut-off and allow both the center and the spinor to vary.
We therefore look for a neighborhood on whose boundary the limiting
energy is strictly below its value at the centered maximum. For a
minimum the inequality is reversed. Uniform convergence then gives
an interior extremum for every sufficiently large finite $R$.
This is why we compute second derivatives. The scale derivative
locates a candidate, and the second variation transverse to the set
of extremizing parameters gives the strict boundary inequality.
In particular, although $\nabla\Ga(\psi_{\lm,0,\ga})=0$, the term
involving $K_{\lm,\xi,\ga}$ in \eqref{widehat-Phi} can contribute
to the second derivative in $\xi$.

The two example tensors are supported on
$E=\operatorname{span}\{e_1,e_2,e_3,e_4\}$. We introduce the radial
polynomials that arise when the bubble is translated. For translations
perpendicular to $E$, set
\begin{\equ}\label{rev-transverse-translation-polynomial}
 T_m(t)=\sum_{k=1}^6k(m-4-2k)q_kM_kt^{k-1},\qquad
 J_m(t)=\frac{2t}{m+4}T_m(t).
\end{\equ}

For the tensor of Example~\ref{rev-weyl-example}, let $\ch_m(t)$ be
the second derivative of $\widehat\Phi_\infty$ in the $\xi^1$ direction
at $(\lm,0,\ga)$, divided by $72\kappa_m\lm^2$, where $t=\lm^2$.
It is a polynomial in $t$, independent of the unit spinor $\ga$.
The four coordinate directions in $E$ have the same second derivative,
and the full translation Hessian is
\begin{\equ}\label{rev-full-translation-hessian}
 \partial_{\xi^\ell\xi^k}^2\widehat\Phi_\infty(\lambda,0,\ga)
 =\begin{cases}
 72\kappa_m\lambda^2\ch_m(t),&1\le\ell=k\le4,\\[0.3em]
 72\kappa_m\lambda^2 J_m(t),&5\le\ell=k\le m,\\[0.3em]
 0,&\ell\ne k.
 \end{cases}
\end{\equ}

Appendix~\ref{app:cubic-translations} derives this formula and gives
the coefficients of $\ch_m$. Thus a zero of $Q_m$ at which
$Q_m'<0$, $\ch_m<0$ and $T_m<0$ gives a strict maximum in scale
and translation for Example~\ref{rev-weyl-example}. In Example~\ref{Ex-2}, the Clifford three-form term depends on the
unit spinor. The stronger inequality involving $J_m+P_m/t^2$ below controls
the mixed variation of the center and the spinor in that example.

\begin{Lem}\label{rev-cubic-signs}
For the polynomial $f$ in \eqref{rev-cubic-profile} and every integer
$m\ge20$, we have
\[
 Q_m(17/20)>0,\qquad Q_m(1)<0,
\]
and, for $17/20\le t\le1$,
\[
 Q_m'(t)<0,\qquad T_m(t)<0,\qquad\ch_m(t)<0,\qquad
 J_m(t)+P_m(t)/t^2<0.
\]
\end{Lem}

In the range $m\ge18$ where the cut-off can be removed, the endpoint
value $Q_m(1)$ is positive for $m=18,19$ and becomes negative at
$m=20$. Appendix~\ref{app:cubic-signs} proves all six inequalities.
The following remark explains how the endpoint condition at $m=20$
leads to the cubic in \eqref{rev-cubic-profile}.

\begin{Rem}\label{rem:choice-cubic-at-twenty}
First allow $f$ to be any polynomial for which the radial integral in
\eqref{rev-central-profile0} converges. Multiplying $f$ by a
nonzero constant multiplies $P_m$ and $Q_m=P_m'/t$ by its square.
We choose $f(0)\ne0$ so that the prescribed Weyl tensor contributes to
the quadratic Taylor term $a f(0)H(x)$ at the center in \eqref{hij}.
We may then normalize $f(0)=1$ without changing the critical scales.

We also want the critical bubble scale to stay in a fixed compact
interval after rescaling by $b$. In terms of $t=(\lm/b)^2$, replacing
$f(s)$ by $f(cs)$ for $c>0$ replaces $Q_m(t)$ by $Q_m(ct)$ in
\eqref{rev-radial-polynomials}. We can thus take $t=1$ as the upper
endpoint at which to test the sign. Since $Q_m(0)=2$ under our
normalization, a first condition on the coefficients is $Q_m(1)<0$.

To choose a cubic profile at $m=20$, write $f(s)=\sum_{j=0}^3 c_js^j$.
Expanding $f^2$ in \eqref{rev-radial-polynomials} gives
\[
 Q_{20}(1)=\sum_{i,j=0}^3(i+j+2)M_{i+j}c_ic_j,
\]
where the $M_j$ in \eqref{rev-moments} are evaluated at $m=20$.
By homogeneity we may temporarily fix $c_3=1$ and minimize this
expression over $c_0,c_1,c_2$. Its quadratic part in these three
variables is positive definite, since the leading principal minors of
its matrix are $2$, $160/49$, and $7904/343$. In particular, a profile
of degree at most $2$ cannot give $Q_{20}(1)<0$. For $c_3=1$, the
unique minimizer satisfies
\[
 \sum_{j=0}^2(i+j+2)M_{i+j}c_j=-(i+5)M_{i+3},
 \qquad i=0,1,2 .
\]
Its solution is $(c_0,c_1,c_2)=(-637/2,819/4,-63/2)$. Dividing the
polynomial by $c_0$ gives precisely \eqref{rev-cubic-profile}, for
which $Q_{20}(1)=-114/31213$.

This calculation selects the coefficients using only the endpoint
at $m=20$. Lemma~\ref{rev-cubic-signs} shows that the same choice gives
$Q_m(17/20)>0$, $Q_m(1)<0$, and $Q_m'<0$ on $[17/20,1]$ for every
$m\ge20$. {We use the unique zero of $Q_m$ in this interval.
The corresponding centered critical parameters lie in the compact set}
$[\sqrt{17/20},1]\times\{0\}\times S(\mbs_m)$ in the rescaled coordinates
$(\lm/b,\xi/b,\ga)$, independently of $b$. The remaining
inequalities in the lemma control the second variations involving
translations for the two Weyl tensors of
Examples~\ref{rev-weyl-example} and~\ref{Ex-2}.
{For the cubic in \eqref{rev-cubic-profile} at $m=20$,
$Q_{20}(1)=-114/31213<0$ and $Q_{20}(2)=40682/31213>0$.
Thus $Q_{20}$ has another zero in $(1,2)$, giving a further critical
scale of the radial factor with $1<\lm/b<\sqrt{2}$.}
\end{Rem}

The inequalities of Lemma~\ref{rev-cubic-signs} persist when the
coefficients vary in the range stated in Remark~\ref{rem:cubic-coefficient-range},
as proved in Appendix~\ref{app:cubic-perturbations}. In particular, each
such profile has a unique critical scale in the indicated interval. To identify the spinors
giving the extrema at that scale, we return to the other factor in
\eqref{rev-central-profile}.

Extend $P(W,\ga)$ to
arbitrary constant spinors by the formula in
Lemma~\ref{lemma: critical-manifold}\,(5). Its Clifford pairings define a
Hermitian endomorphism $\msa_W$ by polarization,
\begin{\equ}\label{quadratic-weyl-spinor-form}
 9\real(\msa_W\ga,\ga)_{\ig_{\R^m}}
 =\frac12|W|^2|\ga|_{\ig_{\R^m}}^2+P(W,\ga),\qquad \ga\in\mbs_m .
\end{\equ}
Clearly, the constant coefficient $9$ is rather harmless. 
For a unit spinor the centered energy is therefore
$9\kappa_m P_m(\lm^2)\real(\msa_W\ga,\ga)_{\ig_{\R^m}}$.
Since $P_m(\lm^2)>0$, the largest and smallest values at each scale
are attained on the corresponding eigenspaces of $\msa_W$.

\begin{Prop}\label{rev-cubic-max}
Let $m\ge20$ and let $f$ be the cubic profile in \eqref{rev-cubic-profile},
with coefficients allowed to vary as in Remark~\ref{rem:cubic-coefficient-range}.
Take $P_m$ as in \eqref{rev-central-profile}, using the coefficients
of $f^2$. Then $P_m'$ has a unique zero $t_m$ in $(17/20,1)$,
and $P_m''(t_m)<0$. The corresponding bubble scale $\sqrt{t_m}$
depends only on $m$ and $f$.

For an algebraic Weyl tensor $W$, let $\ga_+$ and $\ga_-$ be unit
eigenspinors for the largest and smallest eigenvalues of $\msa_W$,
respectively. With the center fixed at $\xi=0$,
every $W\ne0$ gives a positive local maximum of $\widehat\Phi_\infty$
at $(\sqrt{t_m},0,\ga_+)$. If the smallest eigenvalue is negative, there is also a negative
local minimum at
$(\sqrt{t_m},0,\ga_-)$. These extrema are strict in the scale variable.
Both occur, for example, for the tensor in Example~\ref{Ex-2} and for
all algebraic Weyl tensors sufficiently close to it. If
$\tfrac12|W|^2+P(W,\ga)=0$, then $\widehat\Phi_\infty(\lm,0,\ga)=0$
at every scale for that spinor.

At each of these centered extrema, the first derivatives with respect
to the translation variables also vanish. Moreover, for $R=r/b$
sufficiently large, there is a scale $\ell_R$ such that
$\ell_R\to\sqrt{t_m}$ as $R\to\infty$ and $(b\ell_R,0,\ga)$ is
a critical point of $\widehat\Phi_{a,b,r}$ for each of the corresponding
extremizing spinors $\ga$.
\end{Prop}

Appendix~\ref{app:cubic-centered} proves Proposition~\ref{rev-cubic-max},
including stationarity with respect to translations and persistence of
the centered critical scale for a finite cut-off. To obtain an extremum
in the full parameter space, we must also check what happens when the
center and the spinor vary together.

Put $\ell_m=\sqrt{t_m}$ and $G_m=9\kappa_mP_m(t_m)>0$.
Let $\bt_*$ be the largest eigenvalue of $\msa_W$, which is positive
for $W\ne0$. If the smallest eigenvalue is negative, we may instead
take it as $\bt_*$ to study the minimum. Write $E_{\bt_*}$ for its
eigenspace and $S(E_{\bt_*})$ for its unit sphere. The energy is
$G_m\bt_*$ on $\mathcal C_{\bt_*}=\{\ell_m\}\times\{0\}\times S(E_{\bt_*})$.
For $\ga\in S(E_{\bt_*})$ and $v\in\R^m$, define
\begin{\equ}\label{rev-joint-translation-test}
 \ch_{W,\ga}(v)
 =D^2_{\xi\xi}\widehat\Phi_\infty(\ell_m,0,\ga)[v,v]
 -\frac{G_m}{2\ell_m^2}
 \real\bigl((\msa_W-\bt_*\operatorname{Id})(v\cdot_{\ig_{\R^m}}\ga),
             v\cdot_{\ig_{\R^m}}\ga\bigr)_{\ig_{\R^m}}.
\end{\equ}

The first term measures a translation with the spinor fixed. The second
accounts for the mixed derivative in translation and spinor directions.
Completing the square in the spinor variation, as in
Appendix~\ref{app:cubic-joint}, gives the following sufficient condition
for a strict extremum transverse to $\mathcal C_{\bt_*}$,
\begin{\equ}\label{rev-weyl-sign-condition}
 \sgn(\bt_*)\,\ch_{W,\ga}(v)<0
 \quad\text{for }\ga\in S(E_{\bt_*}),\ v\in\R^m\setminus\{0\}.
\end{\equ}

The spinor Hessian already has the required sign on $E_{\bt_*}^\perp$,
and the sign in the scale direction follows from $P_m''(t_m)<0$.
By compactness of $S(E_{\bt_*})$, there is then a relatively compact open
neighborhood $\mathcal U$ of $\mathcal C_{\bt_*}$ in the full parameter
space such that
\begin{\equ}\label{rev-uniform-local-gap}
 \sgn(\bt_*)\widehat\Phi_\infty<G_m|\bt_*|
 \quad\text{on }\overline{\mathcal U}\setminus\mathcal C_{\bt_*},
 \qquad
 \widehat\Phi_\infty=G_m\bt_*\quad\text{on }\mathcal C_{\bt_*} .
\end{\equ}
In particular, the values on $\partial\mathcal U$ are separated from
the extremal value by a positive constant.

In Example~\ref{rev-weyl-example} we have $\msa_W=8\operatorname{Id}$.
The mixed correction in \eqref{rev-joint-translation-test} vanishes, so
\eqref{rev-full-translation-hessian} and Lemma~\ref{rev-cubic-signs}
give \eqref{rev-uniform-local-gap} for the positive maximum.

In Example~\ref{Ex-2}, let $A$ denote Clifford multiplication by the
oriented volume form of $E$. By \eqref{nontrivial-Phi},
$9\msa_W=120\operatorname{Id}-288A$. For a unit spinor with
$A\ga=\sigma\ga$, where $\sigma\in\{-1,1\}$, put $p=P(W,\ga)=-288\sigma$.
Write $v_E$ and $v_{E^\perp}$ for the orthogonal components of $v$.
The translation Hessian with $\ga$ fixed is
\begin{\equ}\label{table-two-translation-hessian}
 \frac{D^2_{\xi\xi}\widehat\Phi_\infty(\lm,0,\ga)[v,v]}{\ka_m\lm^2}
 =\Big[\Big(120+\frac{p}{2}\Big)\ch_m(t)+pJ_m(t)\Big]|v_E|^2
+\Big(120+p\Big)J_m(t)|v_{E^\perp}|^2 .
\end{\equ}

When $\ga$ also varies, the coefficient of $|v_E|^2$ in this
normalization acquires the additional term $pP_m(t_m)/t_m^2$ at
$t=t_m$. This explains the last inequality in Lemma~\ref{rev-cubic-signs}.
For $p=288$ the two coefficients in the joint test are negative,
giving the positive maximum. For $p=-288$ both are positive, giving the
negative minimum. Appendix~\ref{app:cubic-joint} verifies
\eqref{table-two-translation-hessian} and the joint test in coordinates.
The strict inequalities persist for tensors sufficiently close to
the respective examples.

In Subsection~\ref{rev-gluing-subsection} the strict boundary
inequality \eqref{rev-uniform-local-gap} gives a blow-up sequence for
each of the extrema just found. In Section~\ref{sec:radial-dimension-threshold}
the spinor equation selects an eigenvalue of $\msa_W$ at a limiting
centered bubble. When this eigenvalue is nonzero, the scale equation
reduces to a condition on the radial polynomial alone.

\subsection{A blow-up sequence for a smooth metric}\label{rev-gluing-subsection}

Fix $m\ge20$ and use the cubic profile in \eqref{rev-cubic-profile}, or any
cubic polynomial in the coefficient range of Remark~\ref{rem:cubic-coefficient-range}.
{For a fixed nonzero Weyl tensor $W$, we construct a single smooth
metric on $S^m$. For each centered critical set $\mathcal C_{\bt_*}$
satisfying \eqref{rev-uniform-local-gap}, we obtain a blow-up sequence
for this metric. Subsection~\ref{rev-cubic-subsection} already verifies this condition for the maximum in
Example~\ref{rev-weyl-example} and for both extrema in Example~\ref{Ex-2}.}
The cut-off radius
$r\in(0,1)$ is fixed throughout.

For $\mu\in\N$, let
\begin{\equ}\label{rev-h-mu}
h^{(\mu)}_{ij}(x)=a_\mu\eta(|x|/r)
f(|x|^2/b_\mu^2)H_{ij}(x),
\qquad
\tilde\ig^{(\mu)}(x)=\exp\big(\vr h^{(\mu)}(x-x_\mu)\big).
\end{\equ}
We shall choose $a_\mu,b_\mu\to0$ with $a_\mu/b_\mu^6\to0$. For ease of notation, let us write
\begin{\equ}\label{rev-small-quantities}
c_\mu=\frac{a_\mu}{b_\mu^6}\quad \text{and} \quad
\de_{\vr,\mu}=|\vr|a_\mu b_\mu^2.
\end{\equ}
Throughout this subsection, constants in the estimates may depend on
the fixed data $m,r,\eta,f,W$ and, when the bubble parameters vary,
on the compact set in which the rescaled parameters range.
If a constant depends on the index $\mu$, we indicate this by a
subscript $\mu$. For each integer $k\geq0$, the compact support and
the fixed cut-off radius $r$ give
\begin{\equ}\label{rev-coefficient-bound}
\|h^{(\mu)}\|_{C^k(\R^m)}\leq C_kc_\mu,
\end{\equ}
where $C_k$ may also depend on $k$ and is uniform in $\mu$.

Let $\cl^{(\mu)}$ be the functional of the metric in \eqref{rev-h-mu}. Its difference from $\cl_0$ is exactly quadratic. Thus there is a bounded self-adjoint operator $A_\mu$ on the real Hilbert space $\msd^{1/2}(\R^m,\mbs(\R^m))$ such that
\begin{\equ}\label{rev-exact-operator}
\cl^{(\mu)}(u)=\cl_0(u)+\frac12\inp{A_\mu u}{u}_{\msd}.
\end{\equ}
We use Hilbert gradients in this formula. In particular, if $A_{1,\mu}$ denotes the first coefficient of $A_\mu$, then
$\Ga^{(\mu)}(u)=\frac12\inp{A_{1,\mu}u}{u}_{\msd}$ and
$\nabla\Ga^{(\mu)}(u)=A_{1,\mu}u$.
We denote the operator norm on this Hilbert space by
$\|\cdot\|_{\mathrm{op}}$, so that
$\|A_\mu\|_{\mathrm{op}}=\sup_{\|u\|_{\msd}=1}\|A_\mu u\|_{\msd}$.

\begin{Lem}\label{rev-local-estimate}
Fix a domain $\mfd\Subset(0,+\infty)\times\R^m$ for the bubble scale
and center in the rescaled coordinates, and put
\[
\Om_\mu=\{(\lm,\xi):\ (\lm/b_\mu,(\xi-x_\mu)/b_\mu)\in\mfd\}.
\]
For sufficiently small $|\vr|$, the small normal correction
$w^{(\mu)}_z\in T_z\cm^\bot$ is well defined for
$z=\psi_{\lm,\xi,\ga}$ with $(\lm,\xi)\in\Om_\mu$ and satisfies, uniformly in $\ga\in S(\mbs_m)$,
\begin{\equ}\label{rev-local-w-bound}
\|w^{(\mu)}_z\|_{\msd}\leq C\de_{\vr,\mu},
\qquad
\|w^{(\mu)}_z+\vr K_z\nabla\Ga^{(\mu)}(z)\|_{\msd}
\leq C\big(|\vr|c_\mu\de_{\vr,\mu}+
\de_{\vr,\mu}^{(m+1)/(m-1)}\big).
\end{\equ}
Moreover,
\begin{\equ}\label{rev-local-reduced-expansion}
\cl^{(\mu),red}(\lm,\xi,\ga)
=E_0+\vr^2\widehat\Phi^{(\mu)}(\lm,\xi,\ga)
+R_\mu(\lm,\xi,\ga),
\end{\equ}
with $E_0=\frac1{2m}\big(\frac m2\big)^m\om_m$ and
\begin{\equ}\label{rev-relative-remainder}
|R_\mu|\leq
C\big(|\vr|c_\mu+\de_{\vr,\mu}^{2/(m-1)}\big)\de_{\vr,\mu}^2.
\end{\equ}
\end{Lem}
\begin{proof}
We first establish the estimates
\begin{\equ}\label{rev-soft-source-bound}
\|\nabla\Ga^{(\mu)}(z)\|_{\msd}\leq Ca_\mu b_\mu^2,
\quad
\|A_\mu\|_{\mathrm{op}}\leq C|\vr|c_\mu,
\quad
\|(A_\mu-\vr A_{1,\mu})z\|_{\msd}
\leq C|\vr|c_\mu\de_{\vr,\mu}.
\end{\equ}
For the first estimate, write $h(x)=h^{(\mu)}(x-x_\mu)$ and set
$x=x_\mu+b_\mu y$. By
\eqref{rev-h-mu},
\[
 h(x_\mu+b_\mu y)=a_\mu b_\mu^2\eta(b_\mu|y|/r)f(|y|^2)H(y).
\]
The polynomial $f(|y|^2)H(y)$ has degree $8$, and its first derivatives
have degree at most $7$. The derivative of the cut-off satisfies the same
growth estimate after multiplication by the polynomial, since on the
support of $\eta'(b_\mu|y|/r)$ we have $b_\mu/r\leq 1/|y|$.
On the other hand, \eqref{bubble} and \eqref{pa-j-psi} show that the
rescaled bubble $b_\mu^{(m-1)/2}z(x_\mu+b_\mu y)$ and its first
derivatives in $y$ are bounded by $C(1+|y|)^{1-m}$ and
$C(1+|y|)^{-m}$, respectively, uniformly for the rescaled parameters in $\mfd$.
Thus both terms in the bound for $F_h(z)$ in \eqref{rev-source} satisfy
\[
 b_\mu^{(m+1)/2}\bigl(|h|\,|\nabla z|+|\nabla h|\,|z|\bigr)(x_\mu+b_\mu y)
 \leq C a_\mu b_\mu^2(1+|y|)^{8-m}.
\]
For $q=2m/(m+1)$ we have $m/q=(m+1)/2$, so the factor on the left
is exactly the one that preserves the $L^q$ norm under this change of
variables. The weight $(1+|y|)^{8-m}$ belongs to $L^q(\R^m)$ precisely
when $\int_1^\infty \rho^{m-1-q(m-8)}\,d\rho$ converges, that is,
when $q(m-8)>m$, or equivalently $m>17$. Testing $F_h(z)$ against a
spinor and applying H\"older's inequality and the Sobolev embedding
into $L^{2m/(m-1)}(\R^m,\mbs(\R^m))$ now gives
$\|\nabla\Ga^{(\mu)}(z)\|_{\msd}\leq Ca_\mu b_\mu^2$.

For the operator norm, we use \eqref{rev-coefficient-bound} in the
original coordinates. The coefficients of the Dirac difference in
\eqref{Dirac-expansion} are supported in $B_r(x_\mu)$ and bounded by
$C|\vr|c_\mu$. Its bilinear form on smooth compactly supported spinors
$u,v$ is therefore bounded by $C|\vr|c_\mu\|\nabla u\|_{L^2}\|v\|_{L^2}$,
using $\|u\|_{L^2(B_r(x_\mu))}\leq C_r\|\nabla u\|_{L^2(\R^m)}$ for the
term without derivatives of $u$. By symmetry the same bound holds
with $u$ and $v$ interchanged. Interpolation gives
$|\inp{A_\mu u}{v}_{\msd}|\leq C|\vr|c_\mu\|u\|_{\msd}\|v\|_{\msd}$.

Finally, after subtraction of the term linear in $\vr$, the coefficients
of $\nabla z$ and $z$ in the Dirac difference are bounded by
$C|\vr|^2|h|^2$ and $C|\vr|^2|h|\,|\nabla h|$, respectively,
by \eqref{B-bij}, \eqref{bbb} and Remark~\ref{rev-exact-remainder}.
Taking one factor $|h|$ in $L^\infty$ leaves precisely the two terms
estimated above in $L^q$. The same Sobolev estimate therefore bounds
$\|(A_\mu-\vr A_{1,\mu})z\|_{\msd}$ by
$C|\vr|^2c_\mu a_\mu b_\mu^2=C|\vr|c_\mu\de_{\vr,\mu}$.

Write $P_z$ for the orthogonal projection onto $T_z\cm^\bot$. As in the perturbation framework in \cite{Isobe-Xu23},
we seek $w\in T_z\cm^\bot$ such that $P_z\nabla\cl^{(\mu)}(z+w)=0$.
Since $\nabla\cl_0(z)=0$, write
$\mathcal R_z(w)=\nabla\cl_0(z+w)-\nabla^2\cl_0(z)w$ for the Taylor
remainder. By \eqref{rev-exact-operator}, applying the inverse $K_z$ of
the Hessian on $T_z\cm^\bot$ gives the equivalent equation
\begin{\equ}\label{rev-local-normal-equation}
 w=-K_zP_z\bigl(A_\mu z+A_\mu w+\mathcal R_z(w)\bigr).
\end{\equ}
Translation and dilation give $\|K_z\|_{\mathrm{op}}=\|K_{1,0,\ga}\|_{\mathrm{op}}$.
The latter norm depends continuously on $\ga$ and is therefore bounded on
the compact unit spinor sphere.
The derivative of $\phi\mapsto|\phi|_{\ig_{\R^m}}^{2/(m-1)}\phi$ is
H\"older continuous with exponent $2/(m-1)$. Hence, by the Sobolev
embedding, $\mathcal R_z$ has Lipschitz constant at most $C s^{2/(m-1)}$
on $\|w\|_{\msd}\leq s$, uniformly in $z$.
By \eqref{rev-soft-source-bound}, the right-hand side at $w=0$ has
norm at most $C\de_{\vr,\mu}$, and its Lipschitz constant on this
ball is at most $C(|\vr|c_\mu+s^{2/(m-1)})$. Choosing $s$ small
and then $|\vr|$ small makes it a contraction of the ball into itself
with Lipschitz constant at most $1/2$. Its unique fixed point thus
satisfies $\|w_z^{(\mu)}\|_{\msd}\leq 2\|K_zP_z A_\mu z\|_{\msd}
\leq C\de_{\vr,\mu}$.

For the second estimate in \eqref{rev-local-w-bound}, recall that
$A_{1,\mu}z=\nabla\Ga^{(\mu)}(z)$ lies in $T_z\cm^\bot$ because
$\Ga^{(\mu)}$ vanishes on $\cm$ by Lemma~\ref{lemma: critical-manifold}~(4).
Moving its contribution to the left of the fixed point equation gives,
with $w=w_z^{(\mu)}$,
\[
 w+\vr K_z\nabla\Ga^{(\mu)}(z)
 =-K_zP_z\bigl((A_\mu-\vr A_{1,\mu})z+A_\mu w+\mathcal R_z(w)\bigr).
\]
By \eqref{rev-soft-source-bound} and the first estimate for $w$, each
of the first two terms inside the parentheses is bounded in norm by
$C|\vr|c_\mu\de_{\vr,\mu}$. Since $\mathcal R_z(0)=0$, its Lipschitz
estimate with $s=\|w\|_{\msd}$ bounds the last term by
$C\|w\|_{\msd}^{(m+1)/(m-1)}\leq C\de_{\vr,\mu}^{(m+1)/(m-1)}$.
This proves \eqref{rev-local-w-bound}.

To estimate the energy at $z$, expand the exact quadratic perturbation to second order. The omitted metric terms are bounded in absolute value by a constant times
\[
|\vr|^3c_\mu
\int_{\R^m}
\big(|h|^2|\nabla z||z|
+|h||\nabla h||z|^2\big)\,dx.
\]
The integral is $O(a_\mu^2b_\mu^4)$: after rescaling its worst radial exponent is $16-m$, again integrable for $m>17$. Thus this energy error is $O(|\vr|c_\mu\de_{\vr,\mu}^2)$. Finally, the expansion of $\cl_0(z+w)$ has remainder $O(\|w\|_{\msd}^{2m/(m-1)})$. Substitution of the first approximation $-\vr K_z\nabla\Ga^{(\mu)}(z)$, and completion of the normal quadratic form, yield the correction
$-\frac12\vr^2\inp{K_z\nabla\Ga^{(\mu)}(z)}{\nabla\Ga^{(\mu)}(z)}_{\msd}$.
All remaining terms are bounded by the right-hand side of \eqref{rev-relative-remainder}. This proves the lemma.
\end{proof}

Assume now that the balls $B_r(x_\nu)$ are pairwise disjoint and that $\sum_\nu c_\nu<\infty$, where $c_\nu$ is given in \eqref{rev-small-quantities} for $\nu\in\N$. Define
\begin{\equ}\label{rev-h-star}
h^*(x)=\sum_{\nu=1}^\infty h^{(\nu)}(x-x_\nu),
\qquad \tilde\ig^*(x)=\exp(\vr h^*(x)).
\end{\equ}
The associated functional has the exact form
\[
\cl^*(u)=\cl_0(u)+\frac12\inp{A^*u}{u}_{\msd},
\qquad A^*=\sum_{\nu=1}^\infty A_\nu.
\]
The local quadratic forms add because the supports are disjoint.
By \eqref{rev-soft-source-bound} and $\sum_\nu c_\nu<\infty$, the series
defining $A^*$ converges in operator norm, with
\[
 \|A^*\|_{\mathrm{op}}\leq\sum_\nu\|A_\nu\|_{\mathrm{op}}
 \leq C|\vr|\sum_\nu c_\nu .
\]
For sufficiently small $|\vr|$, the same contraction argument gives
a unique small solution $w_z^*\in T_z\cm^\bot$ of the global normal equation
\begin{\equ}\label{rev-global-normal-equation}
 w=-K_zP_z\bigl(A^*z+A^*w+\mathcal R_z(w)\bigr).
\end{\equ}
We write $\cl^{*,red}(z)=\cl^*(z+w_z^*)$.

For the comparison with the other centers, choose $\mfd$ in
Lemma~\ref{rev-local-estimate} to contain the centered parameter $(\ell_m,0)$,
where $\ell_m=\sqrt{t_m}$ is the critical bubble scale identified in
Proposition~\ref{rev-cubic-max} and $t_m$ is the unique zero of $P_m'$ in $(17/20,1)$.
We use the Euclidean distance $\dist$
on $(0,+\infty)\times\R^m$ and write
\[
 \diam(\mfd)=\sup_{p,p'\in\mfd}\dist(p,p')
            =\max_{p,p'\in\overline{\mfd}}\dist(p,p').
\]
\begin{Lem}\label{rev-gluing-comparison}
Fix $\mu\in\N$ and suppose that
\[
 |x_\nu-x_\mu|\geq 2\bigl(r+b_\mu\diam(\mfd)\bigr)\qquad\text{for every }\nu\ne\mu.
\]
For $z=\psi_{\lm,\xi,\ga}$ with $(\lm,\xi)\in\Om_\mu$, set
\begin{\equ}\label{rev-positive-tail}
\mathcal T_\mu=|\vr|b_\mu^{(m-1)/2}
\sum_{\nu\ne\mu}\frac{c_\nu}{|x_\nu-x_\mu|^{m-1}}.
\end{\equ}
Then
\begin{\equ}\label{rev-comparison-w}
\|w_z^*-w_z^{(\mu)}\|_{\msd}
\leq C\big(\mathcal T_\mu+|\vr|\de_{\vr,\mu}\big),
\end{\equ}
and
\begin{\equ}\label{rev-comparison-energy}
|\cl^{*,red}(z)-\cl^{(\mu),red}(z)|
\leq C\big(\mathcal T_\mu^2+
\mathcal T_\mu\de_{\vr,\mu}+|\vr|\de_{\vr,\mu}^2\big).
\end{\equ}
\end{Lem}
\begin{proof}
Let $A_{*,\mu}$ denote the contribution from all centers $x_\nu$ with $\nu\ne\mu$, so that
\[
 A_{*,\mu}=A^*-A_\mu=\sum_{\nu\ne\mu}A_\nu.
\]
Applying \eqref{rev-soft-source-bound} to each $A_\nu$ gives
\[
 \|A_{*,\mu}\|_{\mathrm{op}}
 \leq \sum_{\nu\ne\mu}\|A_\nu\|_{\mathrm{op}}
 \leq C|\vr|\sum_{\nu\ne\mu}c_\nu
 \leq C'|\vr|,
\]
where $C'=C\sum_{\nu=1}^\infty c_\nu$ is finite and independent of $\mu$.
For $(\lm,\xi)\in\Om_\mu$, the rescaled parameters
$(\ell,\zeta)=(\lm/b_\mu,(\xi-x_\mu)/b_\mu)$ belong to $\mfd$.
Since $(\ell_m,0)\in\mfd$, we have
$|\zeta|\leq\dist\bigl((\ell,\zeta),(\ell_m,0)\bigr)\leq\diam(\mfd)$.
Hence, for $x\in B_r(x_\nu)$ with $\nu\ne\mu$,
\[
 |x-\xi|\geq |x_\nu-x_\mu|-r-b_\mu\diam(\mfd)\geq \tfrac12|x_\nu-x_\mu|.
\]
The bubble formula \eqref{bubble} and its derivative \eqref{pa-j-psi} then give
\[
|z|\leq Cb_\mu^{(m-1)/2}|x_\nu-x_\mu|^{-(m-1)},
\quad
|\nabla z|\leq Cb_\mu^{(m-1)/2}|x_\nu-x_\mu|^{-m}.
\]
The local coefficient and derivative bounds, followed by the Sobolev dual estimate, imply
\begin{\equ}\label{rev-remote-source}
\|A_{*,\mu}z\|_{\msd}\leq C\mathcal T_\mu.
\end{\equ}
For the direct bubble energy a stronger estimate holds. The identity
$\Ga^{(\nu)}|_{\cm}=0$ removes the first-order energy on each remote patch. Expanding the exact metric coefficients only for this estimate gives
\begin{\equ}\label{rev-direct-tail-square}
\frac12|\inp{A_{*,\mu}z}{z}_{\msd}|
\leq C\vr^2b_\mu^{m-1}
\sum_{\nu\ne\mu}
\frac{c_\nu^2}{|x_\nu-x_\mu|^{2m-2}}
\leq C\mathcal T_\mu^2.
\end{\equ}
Higher powers of the coefficients are absorbed since $|\vr|\sup c_\nu$ is small.
The last inequality follows because all cross terms in the expansion of
$\mathcal T_\mu^2$ from \eqref{rev-positive-tail} are nonnegative.

To estimate $w_z^*-w_z^{(\mu)}$, fix $z$ and compare the local equation
\eqref{rev-local-normal-equation} with the global equation
\eqref{rev-global-normal-equation}. At $w=w_z^{(\mu)}$, the right-hand
side of the global equation differs from that of the local equation by
$-K_zP_z A_{*,\mu}(z+w_z^{(\mu)})$. Since the global map has Lipschitz
constant at most $1/2$ on the ball containing both corrections, its fixed
point satisfies
\[
\|w_z^*-w_z^{(\mu)}\|_{\msd}
\leq C\|A_{*,\mu}(z+w_z^{(\mu)})\|_{\msd}
\leq C(\mathcal T_\mu+|\vr|\de_{\vr,\mu}),
\]
where the last inequality uses \eqref{rev-remote-source} and
\eqref{rev-local-w-bound}. This proves \eqref{rev-comparison-w}.

The local correction $w_z^{(\mu)}$ from \eqref{rev-local-normal-equation}
and the global correction $w_z^*$ from \eqref{rev-global-normal-equation}
both lie in $T_z\cm^\bot$. Their difference is therefore a normal
direction at $z$. Since $P_z\nabla\cl^*(z+w_z^*)=0$, the linear term
in the Taylor expansion of $\cl^*$ at $z+w_z^*$ in this direction
vanishes, giving
\[
|\cl^*(z+w_z^*)-\cl^*(z+w_z^{(\mu)})|
\leq C\|w_z^*-w_z^{(\mu)}\|_{\msd}^2.
\]
It remains to compare the two functionals at the same spinor $z+w_z^{(\mu)}$.
Since $A_{*,\mu}=A^*-A_\mu$ is self-adjoint, expanding the quadratic form gives
\[
 \cl^*(z+w_z^{(\mu)})-\cl^{(\mu)}(z+w_z^{(\mu)})
 =\frac12\inp{A_{*,\mu}z}{z}_{\msd}
  +\inp{A_{*,\mu}z}{w_z^{(\mu)}}_{\msd}
  +\frac12\inp{A_{*,\mu}w_z^{(\mu)}}{w_z^{(\mu)}}_{\msd}.
\]
The first term is bounded in absolute value by $C\mathcal T_\mu^2$
by \eqref{rev-direct-tail-square}. For the mixed term, Cauchy-Schwarz
and the bounds on $A_{*,\mu}z$ in \eqref{rev-remote-source} and
on $w_z^{(\mu)}$ in \eqref{rev-local-w-bound} give
$C\mathcal T_\mu\de_{\vr,\mu}$. For the last term we use
$\|A_{*,\mu}\|_{\mathrm{op}}\leq C|\vr|$ and obtain
$C|\vr|\de_{\vr,\mu}^2$. Combining these bounds with the estimate for
the change between the two corrections proves \eqref{rev-comparison-energy}.
\end{proof}

The estimates \eqref{rev-local-w-bound} and \eqref{rev-comparison-w}
control the correction $w_z^*$ to the standard bubble $z$ in the energy
norm and hence, by Sobolev embedding, in $L^{2m/(m-1)}(\R^m,\mbs(\R^m))$.
To show that the solutions we seek have no zeros, we also need to
control their length pointwise, even far from the center of concentration.
On the sphere, we can normalize the center and scale by a conformal
map that takes the standard bubble to a Killing spinor of constant
positive length. Comparing the Green representations of the transformed
solution and this Killing spinor will then give a positive lower bound
for the length of the transformed solution throughout $S^m$.

Let $\ig_0$ be the unit round metric on $S^m$. We measure the $C^2$
distance between a smooth metric $\ig$ and $\ig_0$ on the whole sphere,
using $\ig_0$ to take the supremum norms of $\ig-\ig_0$ and of its first
two covariant derivatives. As in Subsection~\ref{subsec:bg-identification},
the Bourguignon-Gauduchon isometry
$\mfj:\mbs(S^m,\ig_0)\to\mbs(S^m,\ig)$ allows us to work on the fixed
spinor bundle $\mbs(S^m,\ig_0)$ of the round metric. On this bundle
the transported Dirac operator is $\widetilde D_\ig=\mfj^{-1}D_\ig\mfj$.
For $\ig$ sufficiently close to $\ig_0$, it is invertible. Define its
Green kernel $G_\ig$ using the same volume form for every $\ig$, so that
\[
 (\widetilde D_\ig^{-1}\varphi)(x)
 =\int_{S^m}G_\ig(x,y)\varphi(y)\,d\vol_{\ig_0}(y)
\]
for every smooth section $\varphi$ of $\mbs(S^m,\ig_0)$.
Thus $G_\ig(x,y)$ maps the spinor fiber at $y$ to the fiber at $x$,
both for $\ig_0$. With the operator norm on these maps and the chordal
distance $d(x,y)=|x-y|$ for $x,y\in S^m\subset\R^{m+1}$, we claim that
\begin{\equ}\label{rev-round-kernel-stability}
 \sup_{x\ne y}d(x,y)^{m-1}
      \|G_\ig(x,y)-G_{\ig_0}(x,y)\|_{\mathrm{op}}\to0
 \quad\text{as }\|\ig-\ig_0\|_{C^2(S^m,\ig_0)}\to0 .
\end{\equ}
The weight $d(x,y)^{m-1}$ compensates for the singularity at $x=y$.
Using the chordal distance will also preserve this bound when we
normalize the center and scale of the bubble by a conformal map.

To prove \eqref{rev-round-kernel-stability}, fix $s>m$ and take the
Sobolev norms of sections of $\mbs(S^m,\ig_0)$. As $\ig\to\ig_0$ in
$C^2(S^m,\ig_0)$, $\widetilde D_\ig\to D_{\ig_0}$ in operator norm
from $W^{1,s}$ to $L^s$. Invertibility of $D_{\ig_0}$ therefore gives
invertibility of $\widetilde D_\ig$ and continuity of its inverse from
$L^s$ to $W^{1,s}$. To describe the kernel near $x=y$, write the
principal part of $\widetilde D_\ig$ in a local spin frame and evaluate
its coefficients at $y$. In linear coordinates orthonormal for $\ig(y)$,
the resulting constant-coefficient operator is the Euclidean Dirac
operator. Its fundamental solution is homogeneous of degree $1-m$ in
the coordinate difference. Since the coordinate distance is locally
comparable to $d(x,y)$, the kernel is $O(d(x,y)^{1-m})$ and its first
derivatives in $x$ are $O(d(x,y)^{-m})$.

These local fundamental solutions give an approximate inverse on $S^m$
as follows. Let $\{\chi_j\}$ be a partition of unity subordinate to a
finite coordinate cover. For a smooth spinor $\varphi$, we apply the
local kernel in the $j$-th chart to $\chi_j\varphi$. To write this
integral with respect to the round volume $d\vol_{\ig_0}(y)=J_0(y)\,dy$,
we divide the coordinate kernel by $J_0(y)$. A smooth cut-off supported
in the chart and equal to one near $\supp\chi_j$ allows us to extend the
resulting spinor by zero. Summing these contributions over $j$ defines
$K_\ig\varphi$ on the sphere.

When $\widetilde D_\ig$ acts on $K_\ig\varphi$ in the $x$ variable,
the point mass at $x=y$ in the $j$-th chart gives $\chi_j\varphi$,
since the cut-off is one on $\supp\chi_j$. These terms sum to $\varphi$,
which accounts for the identity in $\widetilde D_\ig K_\ig=I+R_\ig$.
In the remainder $R_\ig$, a derivative of the frozen kernel is
multiplied by the difference between the principal coefficients at $x$
and $y$. The derivative is $O(d(x,y)^{-m})$ and the difference is
$O(d(x,y))$, so their product is $O(d(x,y)^{1-m})$. The lower-order terms already
have this bound, and derivatives of the cut-off contribute only away
from the diagonal. Thus the kernels of $K_\ig$ and $R_\ig$ satisfy the
same bound. The frozen fundamental solution depends smoothly on the
principal coefficients at $y$. Since these coefficients vary continuously
with $\ig$ in $C^1$ and the lower-order coefficients in $C^0$, the same
estimates applied to differences give continuity of both kernels in the
weighted norm in \eqref{rev-round-kernel-stability}.

We can now improve the approximate inverse by successive compositions
with $R_\ig$. Integrating in local coordinates shows that a product of $j$ such
integral operators has kernel bounded by $C d^{j-m}$ for $j<m$, by
$C(1+|\ln d|)$ for $j=m$, and by $C$ for $j=m+1$. The same bounds
apply to differences of products, with a factor tending to zero,
by taking the difference in one kernel at a time. Now write
\[
 \widetilde D_\ig^{-1}=K_\ig\sum_{j=0}^{m}(-R_\ig)^j
                    +\widetilde D_\ig^{-1}(-R_\ig)^{m+1}.
\]
The finite sum is continuous in the weighted kernel norm. In the last
term, the kernel of $(-R_\ig)^{m+1}$ converges uniformly to the one
for $\ig_0$. Applying $\widetilde D_\ig^{-1}$ in the $x$ variable and
using $W^{1,s}\subset C^0$ gives uniform convergence of the last
kernel as well. This proves \eqref{rev-round-kernel-stability}.

\begin{Thm}\label{rev-noncompact-sphere}
Let $m\ge20$, let $W$ be a nonzero algebraic Weyl tensor on $\R^m$, and
let $f$ be the cubic polynomial in \eqref{rev-cubic-profile}, with its
coefficients allowed to vary as in Remark~\ref{rem:cubic-coefficient-range}.
Suppose that at least one centered critical set $\mathcal C_{\bt_*}$ of
$\widehat\Phi_\infty$ satisfies \eqref{rev-uniform-local-gap} in the full
space of bubble parameters. There is a smooth, non-locally conformally
flat metric $\ig^*$ on $S^m$ such that, for each such $\mathcal C_{\bt_*}$,
the equation
\[
 D_{\ig^*}\psi=|\psi|_{\ig^*}^{2/(m-1)}\psi
\]
admits a sequence $\{\psi_\mu\}$ of distinct smooth solutions such that
$|\psi_\mu(x)|_{\ig^*}>0$ for every $x\in S^m$.
{For all sufficiently large $\mu$, the function
$|\psi_\mu|_{\ig^*}$ attains its maximum at exactly one point.}
In the coordinates of the metric profiles their bubble parameters satisfy
\[
 \frac{\lm_\mu}{b_\mu}\to\ell_m,\qquad
 \frac{\xi_\mu-x_\mu}{b_\mu}\to0,\qquad
 \operatorname{dist}(\ga_\mu,S(E_{\bt_*}))\to0 .
\]
Moreover, it holds that $\|\psi_\mu\|_{L^\infty}\to\infty$ and
\begin{\equ}\label{rev-sphere-blow-up}
 \frac1{2m}\int_{S^m}|\psi_\mu|_{\ig^*}^{2m/(m-1)}\,d\vol_{\ig^*}
 =E_0+\de_{\vr,\mu}^2 G_m\bt_*+o(\de_{\vr,\mu}^2).
\end{\equ}
In Example~\ref{rev-weyl-example} the energies approach $E_0$ from above.
In Example~\ref{Ex-2} the same metric admits two such sequences, one
with energies above $E_0$ and the other with energies below $E_0$.
\end{Thm}
\begin{proof}
For a clear view, we divide the proof into five steps: 1. Choice of the metric parameters; 2. The limiting reduced energy and the construction of solutions; 3. Extension of the metric and the solutions to the sphere; 4. Blow-up and uniqueness of the maximum point; 5. Nonvanishing on the whole sphere.

\medskip

\noindent
\textit{Step 1. Choice of the metric parameters.}\quad
Fix a neighborhood $\mathcal U$ of $\mathcal C_{\bt_*}$ as in
\eqref{rev-uniform-local-gap}. The estimates in
Lemmas~\ref{rev-local-estimate} and~\ref{rev-gluing-comparison} apply on
$\overline{\mathcal U}$ by taking $\mfd$ to contain its projection onto the scale and center parameters.

Take $r=3/4$ and choose
\begin{\equ}\label{rev-final-parameters}
b_\mu=2^{-\mu/2},\quad a_\mu=b_\mu^7=2^{-7\mu/2},\quad
x_\mu=L\mu e_1,\quad L=|\vr|^{-1/(m-1)}.
\end{\equ}
Since $|x_\nu-x_\mu|=L|\nu-\mu|\geq L$ for $\nu\ne\mu$ and $b_\mu\leq1$,
the separation in Lemma~\ref{rev-gluing-comparison} holds for every $\mu$
as soon as $L\geq2\bigl(r+\diam(\mfd)\bigr)$. Here $c_\mu=b_\mu$, $\de_{\vr,\mu}=|\vr|b_\mu^9$ and, hence, $\sum_\mu c_\mu<\infty$.
Since $m-1>1$, the elementary convergent series
$\sum_{j\ne0}|j|^{-(m-1)}$ gives
\[
\sum_{\nu\ne\mu}\frac{c_\nu}{|x_\nu-x_\mu|^{m-1}}
\leq CL^{-(m-1)}.
\]
Consequently,
\begin{\equ}\label{rev-tail-ratio}
\frac{\mathcal T_\mu}{\de_{\vr,\mu}}
\leq Cb_\mu^{(m-19)/2}L^{-(m-1)}
\to0\quad \text{for }m\geq20.
\end{\equ}
Equivalently, the choice $a_\mu=2^{-\alpha\mu}$ only requires $\alpha>3$ for coefficient decay and $m>4\alpha+5$ for this tail comparison; \eqref{rev-final-parameters} uses $\alpha=7/2$.

\smallskip\noindent
\textit{Step 2. The limiting reduced energy and the construction of solutions.}\quad
In the coordinates $(\ell,\zeta)=(\lm/b_\mu,(\xi-x_\mu)/b_\mu)$,
the scaling law \eqref{ab-a1} and the cut-off limit give
\begin{\equ}\label{rev-local-limit}
 \frac{\widehat\Phi^{(\mu)}(b_\mu\ell,x_\mu+b_\mu\zeta,\ga)}
      {a_\mu^2b_\mu^4}
 \to\widehat\Phi_\infty(\ell,\zeta,\ga)
 \qquad\text{uniformly on }\overline{\mathcal U}.
\end{\equ}
We claim that the same limit holds for the full reduced energy,
\begin{\equ}\label{rev-glued-limit}
 \frac{\cl^{*,red}(b_\mu\ell,x_\mu+b_\mu\zeta,\ga)-E_0}{\de_{\vr,\mu}^2}
 \to\widehat\Phi_\infty(\ell,\zeta,\ga)
 \qquad\text{uniformly on }\overline{\mathcal U}.
\end{\equ}
By \eqref{rev-relative-remainder}, it suffices to show that the remote
profiles change the reduced energy by $o(\de_{\vr,\mu}^2)$ for fixed
$\vr\ne0$. After rescaling, the normalized local corrections
$w_z^{(\mu)}/\de_{\vr,\mu}$ form a relatively compact set in
$\msd^{1/2}(\R^m,\mbs(\R^m))$. Indeed, \eqref{rev-local-w-bound}
approximates them uniformly by the normalized linear responses up to
$O(|\vr|c_\mu+\de_{\vr,\mu}^{2/(m-1)})$. The normalized cut-off sources
agree with their uncut limits on each fixed ball for large $\mu$.
Their pointwise norms are bounded uniformly by
$C(1+|y|)^{8-m}\in L^{2m/(m+1)}(\R^m)$, so convergence in the
dual norm is uniform on the compact parameter set.
The normal inverses depend continuously on the bubble parameters.
Thus the limiting linear responses range over a compact set.

In the same coordinates the remote operator $A_{*,\mu}=A^*-A_\mu$
is uniformly bounded and eventually annihilates each fixed smooth
compactly supported spinor, since the remote supports leave every fixed
ball. Density gives strong convergence to zero, uniform on compact sets.
Hence $\|A_{*,\mu}w_z^{(\mu)}\|_{\msd}=o(\de_{\vr,\mu})$. Since
$\|A_{*,\mu}z\|_{\msd}=O(\mathcal T_\mu)=o(\de_{\vr,\mu})$, comparison
of the normal equations as in Lemma~\ref{rev-gluing-comparison} gives
$\|w_z^*-w_z^{(\mu)}\|_{\msd}=o(\de_{\vr,\mu})$. The resulting change
in energy between the two normal corrections is quadratic in this
difference. The direct remote energy is $O(\mathcal T_\mu^2)$ by
\eqref{rev-direct-tail-square}, its cross term with $w_z^{(\mu)}$
is $O(\mathcal T_\mu\de_{\vr,\mu})$, and its value on $w_z^{(\mu)}$
alone is $o(\de_{\vr,\mu}^2)$. This proves \eqref{rev-glued-limit}.

Maximize $\sgn(\bt_*)(\cl^{*,red}-E_0)$ in the rescaled
parameters on $\overline{\mathcal U}$. By \eqref{rev-uniform-local-gap}
and \eqref{rev-glued-limit}, all maximizers are interior for large $\mu$,
and every limit point belongs to $\mathcal C_{\bt_*}$.
This gives the three parameter limits and the expansion of the reduced
energy at the selected points.

Let $z$ be the bubble at one of these interior critical points.
The global normal equation \eqref{rev-global-normal-equation} gives
$P_z\nabla\cl^*(z+w_z^*)=0$, so the full gradient lies in $T_z\cm$.
To see that it also vanishes in the tangent directions, recall that
$\cl^{*,red}(z)=\cl^*(z+w_z^*)$. The implicit function theorem applied
to the normal equation gives $C^1$ dependence of $w_z^*$ on $z$,
since its linearization in the normal variable is invertible by the
same contraction estimate. For $Z\in T_z\cm$, differentiating
$P_z w_z^*=w_z^*$ gives
\[
 (I-P_z)\,\mathrm{d}w_z^*[Z]=(\mathrm{d}P_z[Z])w_z^* .
\]
The chain rule and the tangency of the full gradient now give
\[
 0=\bigl(\cl^{*,red}\bigr)'(z)[Z]
   =\inp{\nabla\cl^*(z+w_z^*)}{Z+(\mathrm{d}P_z[Z])w_z^*}_{\msd}.
\]
As $Z$ varies, the vectors in the second slot span $T_z\cm$ for all
sufficiently large $\mu$. Indeed, on the compact set of rescaled
parameters $\overline{\mathcal U}$ the tangent spaces of the bubble
manifold vary smoothly. Since translation and dilation are Hilbert
space isometries, $\|\mathrm{d}P_z[Z]\|_{\mathrm{op}}\le C\|Z\|_{\msd}$
with $C$ independent of $\mu$. Together with
$\|w_z^*\|_{\msd}=O(\de_{\vr,\mu})\to0$, this makes the map
$Z\mapsto Z+(\mathrm{d}P_z[Z])w_z^*$ an invertible perturbation of
the identity on $T_z\cm$. Thus $\nabla\cl^*(z+w_z^*)=0$, and the
spinor corresponding to $z+w_z^*$ under the Bourguignon-Gauduchon
isometry is a weak solution of \eqref{SY-tilde-g} for the metric
$\tilde\ig^*$ in \eqref{rev-h-star}.

For Example~\ref{Ex-2} we carry out this selection on the two disjoint
neighborhoods in the same parameter space. The profiles and their
locations are identical for both selections, and one choice of $\vr$
works on the union of their closures. The two energy coefficients are
$408\kappa_m P_m(t_m)$ and $-168\kappa_m P_m(t_m)$.

In each sequence $\|w_z^*\|_{\msd}=O(\de_{\vr,\mu})\to0$. Bubbles on distinct profiles
have mutual Hilbert inner products tending to zero as their scales
shrink. Passing to a subsequence gives pairwise distinct solutions.

\smallskip\noindent
\textit{Step 3. Extension of the metric and the solutions to the sphere.}\quad
Let
\[
\ka(y)=\frac{y}{|y|^2},\qquad
Q(y)=I-2\frac{y\otimes y}{|y|^2}.
\]
For $y\ne0$, the Jacobian matrix of the inversion is
$\nabla\ka(y)=(\pa\ka^i/\pa y^j)_{i,j=1}^m=|y|^{-2}Q(y)$.
Pulling back the metric by $\ka$ and factoring out the round conformal
factor $4/(1+|y|^2)^2$ gives the metric coefficient matrix in the chart
at the north pole,
\begin{\equ}\label{rev-correct-inversion}
G^\sharp(y)=Q(y)^T\tilde\ig^*(\ka(y))Q(y).
\end{\equ}
On the inverted support of the $\mu$-th patch, $|\ka(y)-x_\mu|\leq r$,
so $|x_\mu|-r\leq |y|^{-1}\leq |x_\mu|+r$. The chain rule, including derivatives of $Q$, gives for every fixed $k\geq0$
\begin{\equ}\label{rev-smooth-north-estimate}
\|Q^T(\tilde\ig^{(\mu)}\circ\ka-I)Q\|_{C^k}
\leq C_k|\vr|\frac{a_\mu}{b_\mu^6}|x_\mu|^{2k}
\leq C_k|\vr|^{1-2k/(m-1)}2^{-\mu/2}\mu^{2k}.
\end{\equ}
For fixed $\vr\ne0$, the last expression is summable in $\mu$ for every $k$. Each finite partial sum equals $I$ near $y=0$, and the partial sums converge in every $C^k$ norm. Hence $G^\sharp$ extends smoothly across $0$, and the compactified metric $\ig^*$ is smooth on $S^m$. At the centre of any patch, $f(0)=1$ and the curvature formula in normal coordinates gives a nonzero multiple of $\vr a_\mu W$ for its Weyl tensor. Thus $\ig^*$ is not locally conformally flat.

The solutions have finite $L^{2m/(m-1)}$ norm, which is preserved under conformal compactification. They extend weakly across the north pole. Indeed, for a cut-off supported in a ball of radius $2s$, with derivative bounded by $C/s$, its derivative has $L^{2m/(m+1)}$ norm bounded by $Cs^{(m-1)/2}\to0$. H\"older's inequality therefore removes the extra cut-off term in the weak equation. The regularity results in \cite[Appendix]{Isobe11} and \cite[Chapter~3]{Ammann} then give $C^1$ solutions on the whole sphere.

At the critical point $u=z+w_z^*$, conformal compactification gives
\[
 2m\cl^*(u)=\int_{\R^m}|u|^{2m/(m-1)}\,dx
 =\int_{S^m}|\psi_\mu|_{\ig^*}^{2m/(m-1)}\,d\vol_{\ig^*}.
\]
The reduced energy expansion from Step~2 now gives \eqref{rev-sphere-blow-up}.

{
\smallskip\noindent
\textit{Step 4. Blow-up and uniqueness of the maximum point.}\quad
Let $p_\mu$ be the point with stereographic coordinate $x_\mu$, put
$\rho(x)=2/(1+|x|^2)$, and set $s_\mu=b_\mu\rho(x_\mu)$.
This is the length scale on the sphere corresponding to $b_\mu$ in
the chart. In the coordinates $x=x_\mu+b_\mu y$, the metric
$s_\mu^{-2}\ig^*$ tends smoothly to the Euclidean metric on compact
sets. Indeed, $\rho(x_\mu+b_\mu y)/\rho(x_\mu)\to1$ smoothly there,
and the rescaled log metric tends smoothly to zero by \eqref{rev-h-mu}.
Since $\|w_z^*\|_{L^{2m/(m-1)}(\R^m)}\to0$ by the Sobolev inequality,
the rescaled spinors converge in $L^{2m/(m-1)}$ on each fixed ball to
$\psi_{\ell_m,0,\ga}$ along any subsequence on which $\ga_\mu\to\ga$.
The global $L^{2m/(m-1)}$ convergence of the correction to zero and the integrability
of $|\psi_{\ell_m,0,\ga}|^{2m/(m-1)}$ give
\begin{\equ}\label{rev-sphere-bubble-tail}
 \lim_{R\to\infty}\limsup_{\mu\to\infty}
 \int_{S^m\setminus B_{R s_\mu}(p_\mu)}
 |\psi_\mu|_{\ig^*}^{2m/(m-1)}\,d\vol_{\ig^*}=0 .
\end{\equ}

Write $M_\mu=\max_{S^m}|\psi_\mu|_{\ig^*}$ and
$d_\mu=M_\mu^{-2/(m-1)}$. For some fixed $R>0$, local convergence to
the nonzero bubble gives
$\int_{B_{R s_\mu}(p_\mu)}|\psi_\mu|_{\ig^*}^{2m/(m-1)}\,d\vol_{\ig^*}
\ge c>0$ for all large $\mu$. Hence
$M_\mu\ge c s_\mu^{-(m-1)/2}\to\infty$ and $d_\mu\le C s_\mu$.
At any maximum point, rescale by $d_\mu$ and divide the spinor by
$M_\mu$. Its norm is then at most one and equals one at the origin.
Interior $W^{1,p}$ estimates for the Dirac equation with $p>m$
give a uniform H\"older bound on the unit ball in these coordinates.
Consequently, at every maximum point $q_\mu$,
$\int_{B_{d_\mu}(q_\mu)}|\psi_\mu|_{\ig^*}^{2m/(m-1)}\,d\vol_{\ig^*}$
has a lower bound $c_0>0$ independent of $\mu$ and $q_\mu$.
By \eqref{rev-sphere-bubble-tail}, we can fix $R$ large enough that
the integral outside $B_{R s_\mu}(p_\mu)$ is less than $c_0$ for all
sufficiently large $\mu$. Thus $B_{d_\mu}(q_\mu)$ must meet
$B_{R s_\mu}(p_\mu)$ for every maximum point $q_\mu$.
Since $d_\mu\le C s_\mu$, the triangle inequality places every such
point within $(R+C)s_\mu$ of $p_\mu$. Moreover,
$d_\mu/s_\mu$ is bounded below: otherwise these balls would shrink
to points in the $s_\mu$ coordinates, contradicting the strong
local $L^{2m/(m-1)}$ convergence just proved.

The spinors at scale $s_\mu$ are thus uniformly bounded.
Interior elliptic estimates first give local uniform convergence
to the limiting bubble. Since its norm is positive on every compact
set, the nonlinearity is smooth on the range there for large $\mu$,
and the convergence improves to $C^2$ on compact sets.
The squared norm of the limiting bubble is a positive constant
times $(\ell_m^2+|y|^2)^{1-m}$. It has a unique maximum at $0$ and
a negative definite Hessian there. All rescaled maximum points
therefore tend to $0$, and the squared norms are strictly concave
on a fixed small Euclidean ball about $0$ for large $\mu$.
This proves uniqueness of the maximum. The argument applies to every
convergent subsequence of unit spinor parameters.}

\smallskip\noindent
\textit{Step 5. Nonvanishing on the whole sphere.}\quad
By
\eqref{rev-smooth-north-estimate} with $k=0,1,2$ and the bounds on
each Euclidean patch, we have that $\ig^*\to\ig_0$ in $C^2(S^m,\ig_0)$ as $\vr\to0$.
Indeed, $1-4/(m-1)>0$ and $\sum_\mu 2^{-\mu/2}\mu^4<\infty$.
Write $e(\vr)$ for the supremum in \eqref{rev-round-kernel-stability}
with $\ig=\ig^*$, so $e(\vr)\to0$.

In stereographic coordinates the map $y\mapsto\xi_\mu+\lm_\mu y$
extends to a conformal diffeomorphism $\phi_\mu$ of the round sphere.
Apply its conformal spinor pullback to the compactified standard bubble
with parameters $(\lm_\mu,\xi_\mu,\ga_\mu)$ and to the section
$\mfj^{-1}\psi_\mu$ of $\mbs(S^m,\ig_0)$ representing the exact
solution for $\ig^*$. Denote the resulting spinors by $Z_\mu$
and $U_\mu$, respectively. The bubble formula gives a Killing spinor
of constant length $|Z_\mu|_{\ig_0}\equiv c_m=(m/2)^{(m-1)/2}$.
The $L^{2m/(m-1)}$ norm is conformally invariant, so the Sobolev estimate
on $w_z^*$ gives
\[
 \|U_\mu-Z_\mu\|_{L^{2m/(m-1)}(S^m)}\to0
 \quad  \text{as }\mu\to\infty.
\]
All norms in the rest of the proof use the round metric.

Let $\upsilon_\mu$ be the positive conformal factor on $S^m$ defined by
$\phi_\mu^*\ig_0=\upsilon_\mu^2\ig_0$.
In the Green representation of $U_\mu$ the kernel, denoted by
$G_\mu$, is $G_{\ig^*}(\phi_\mu x,\phi_\mu y)$ multiplied by
$[\upsilon_\mu(x)\upsilon_\mu(y)]^{(m-1)/2}$ and composed at both ends
with the spin isometries. This follows by changing variables in the
integral, using the conformal spinor weight $(m-1)/2$ and the
homogeneity $(m+1)/(m-1)$ of $u\mapsto |u|^{2/(m-1)} u$. The same
transformation preserves $G_{\ig_0}$ by conformal covariance.
Since stereographic projection gives
\[
 d(\phi_\mu x,\phi_\mu y)
 =\sqrt{\upsilon_\mu(x)\upsilon_\mu(y)}\,d(x,y),
\]
we obtain $|G_\mu(x,y)-G_{\ig_0}(x,y)|\le e(\vr)d(x,y)^{1-m}$,
uniformly in the bubble parameters. In particular
$|G_\mu(x,y)|\le C_m d(x,y)^{1-m}$ for all sufficiently small $|\vr|$.

To pass from $L^{2m/(m-1)}$ convergence to a pointwise estimate, observe that
$0<2/(m-1)<1$ and hence
\[
 \big\||U_\mu|^{2/(m-1)}-c_m^{2/(m-1)}\big\|_{L^m}
 \le \|U_\mu-Z_\mu\|_{L^{2m/(m-1)}}^{2/(m-1)}=o(1).
\]
The Hardy-Littlewood-Sobolev inequality in round charts gives
a uniform bound from $L^{2m/3}$ to $L^{2m}$ for the integral operator
with kernel $G_\mu$. Applying it to
$U_\mu=G_\mu(|U_\mu|^{2/(m-1)} U_\mu)$ and using H\"older's inequality
therefore gives
\[
 \|U_\mu\|_{L^{2m}}\le C_m\|U_\mu\|_{L^{2m/3}}
                          +o(1)\|U_\mu\|_{L^{2m}} .
\]
Each $U_\mu$ is continuous. Absorbing the last term and interpolating
$L^{2m/3}$ between $L^{2m/(m-1)}$ and $L^{2m}$, using the uniform
$L^{2m/(m-1)}$ bound on $U_\mu$, gives
$\sup_\mu\|U_\mu\|_{L^{2m}}<\infty$. Since
$(m+1)^2/(m-1)<2m$, interpolation now gives
$|U_\mu|^{2/(m-1)} U_\mu-|Z_\mu|^{2/(m-1)} Z_\mu\to0$ in $L^{m+1}$.
The bound by $C_m d^{1-m}$ also gives a uniform operator bound from
$L^{m+1}$ to $C^0$ on continuous inputs. Subtracting
$Z_\mu=G_{\ig_0}(|Z_\mu|^{2/(m-1)} Z_\mu)$ from the equation for
$U_\mu$, we conclude that
\begin{\equ}\label{rev-global-bubble-length}
 \limsup_{\mu\to\infty}\|U_\mu-Z_\mu\|_{C^0(S^m)}
 \le C_m e(\vr).
\end{\equ}
Here the kernel difference acts on $|Z_\mu|^{2/(m-1)} Z_\mu$, whose
length is the fixed constant $c_m^{(m+1)/(m-1)}$. Choose $|\vr|>0$
small enough that $C_m e(\vr)<c_m/2$. Then $|U_\mu|\ge c_m/2$
everywhere for all sufficiently large $\mu$. The inverse conformal
spinor pullback preserves the zero set, so $\psi_\mu$ is nowhere
vanishing, including at the north pole. Elliptic regularity now gives
$\psi_\mu\in C^\infty$, since the nonlinearity is smooth away from
the zero spinor. Discarding finitely many terms gives the stated
sequences on the same metric.
\end{proof}

\begin{Rem}\label{rem:weyl-data-noncompactness}
{The construction in Theorem~\ref{rev-noncompact-sphere}
applies to any nonzero algebraic Weyl tensor $W$ for which a centered
extremum satisfies \eqref{rev-uniform-local-gap} in the full space of
bubble parameters. To establish noncompactness in a given dimension,
it suffices to find one such $W$. Example~\ref{rev-weyl-example} does
so in every dimension $m\ge20$. In Example~\ref{Ex-2}, the same smooth
metric admits two blow-up sequences whose energies approach the
standard bubble energy from opposite sides. The strict inequalities
verified in Subsection~\ref{rev-cubic-subsection} also give the
respective conclusions for all algebraic Weyl tensors sufficiently
close to those of the two examples.}
\end{Rem}

\subsection{Other choices of $f$}\label{rev-other-choices}

The cubic polynomial in \eqref{rev-cubic-profile} can be understood through a finite-dimensional optimization problem. Increasing the degree gives more freedom in the radial profile, but also increases the growth of the perturbation at infinity. These two effects determine the lowest dimension obtainable by the centered construction used above.

In this subsection we allow
\begin{\equ}\label{rev-choice-polynomial}
f(s)=\sum_{j=0}^{d}v_js^j,\qquad v_0=1,\qquad v_d\ne0.
\end{\equ}
Recall that in \eqref{hij} the parameter $b$ sets the length scale of the
radial factor $f(|x|^2/b^2)$ in $h_{ij}(x)$. For any algebraic Weyl
tensor $W$, we keep the coefficients $v_0,\ldots,v_d$ and the cut-off
radius $r$ fixed as $b\to0$, with $\lm/b$ in a compact subset of
$(0,\infty)$. The untruncated centered profile contains the positive integral
\begin{\equ}\label{rev-choice-integral}
I_m(f)=\int_0^\infty
\frac{\rho^{m+3}f(\rho^2)^2}{(1+\rho^2)^m}\,d\rho.
\end{\equ}
Since its integrand is asymptotic to $v_d^2\rho^{4d+3-m}$, one has
\begin{\equ}\label{rev-choice-integrability}
I_m(f)<\infty\quad\Longleftrightarrow\quad m>4d+4.
\end{\equ}
In particular, for a quartic polynomial with nonzero leading coefficient
and $m=20$, the radial integral $I_m(f)$ diverges logarithmically, since
its integrand is asymptotic to $v_4^2/\rho$.

Suppose henceforth that \eqref{rev-choice-integrability} holds.
For $0\le k\le 2d$, define the constants $M_k(m)$ by
\begin{\equ}\label{rev-choice-moments}
M_k(m)=
\frac{\displaystyle\int_0^\infty
s^{m/2+1+k}(1+s)^{-m}\,ds}
{\displaystyle\int_0^\infty
s^{m/2+1}(1+s)^{-m}\,ds}
=\prod_{j=1}^{k}\frac{m+2j+2}{m-2j-4},
\qquad M_0(m)=1.
\end{\equ}
For these polynomials the same radial integration as in
\eqref{rev-central-profile} gives the centered coefficient of the
uncut model, divided by $a^2b^4$, as the product of
$\kappa_m(\tfrac12|W|^2+P(W,\ga))$ and the radial polynomial
\begin{\equ}\label{rev-choice-radial-profile}
P_{m,f}(t)=t^2\sum_{i,j=0}^{d}v_iv_jM_{i+j}(m)t^{i+j},
\qquad t=(\lambda/b)^2.
\end{\equ}
Define the two matrices
\[
\mathsf G_d(m)=\big(M_{i+j}(m)\big)_{i,j=0}^{d},
\qquad
\mathsf S_d(m)=\big((i+j+2)M_{i+j}(m)\big)_{i,j=0}^{d}.
\]
Collecting the terms $v_jt^j$ in a column vector, we can write the
derivative of \eqref{rev-choice-radial-profile} in matrix form as
\begin{\equ}\label{rev-choice-radial-matrix}
V(t)=(v_0,v_1t,\ldots,v_dt^d)^T,\qquad
\frac{P_{m,f}'(t)}{t}=V(t)^T\mathsf S_d(m)V(t).
\end{\equ}
By \eqref{rev-choice-moments}, $\mathsf G_d(m)$ is the Gram matrix of
$1,s,\ldots,s^d$ for the weight $s^{m/2+1}(1+s)^{-m}$ normalized to
have integral one on $(0,\infty)$. It is therefore positive definite.
Thus a natural coefficient selection is obtained by minimizing
\begin{\equ}\label{rev-choice-rayleigh}
\frac{v^T\mathsf S_d(m)v}{v^T\mathsf G_d(m)v},
\qquad v\ne0,
\end{\equ}
or, equivalently, by solving the generalized eigenvalue problem
$\mathsf S_d(m)v=\kappa\mathsf G_d(m)v$.
If the minimum of the Rayleigh quotient \eqref{rev-choice-rayleigh}
is nonnegative, \eqref{rev-choice-radial-matrix} gives $P_{m,f}'(t)\ge0$
for every choice of coefficients and every $t>0$.
If the generalized eigenvalue problem admits an eigenvalue $\kappa<0$,
take a corresponding eigenvector $v$ as the coefficient vector of $f$.
Since $V(1)=v$, \eqref{rev-choice-radial-matrix} gives
$P_{m,f}'(1)=v^T\mathsf S_d(m)v=\kappa\,v^T\mathsf G_d(m)v<0$.
This strict inequality persists under small changes of the coefficients,
so we may arrange that $v_0\ne0$ and then normalize it to one.
Since $P_{m,f}'(t)>0$ for sufficiently small $t$, this provides a radial local maximum after an arbitrarily small generic perturbation if necessary. At this scale the centered energy has a local maximum if
$\tfrac12|W|^2+P(W,\ga)>0$ and a local minimum if
$\tfrac12|W|^2+P(W,\ga)<0$. The translation directions and the normal
correction must also be checked for a spatial extremum, as in
\eqref{rev-full-translation-hessian}.

\begin{Prop}\label{rev-choice-optimal-degree}
For fixed polynomial coefficients and a finite untruncated radial
profile, a strict scale extremum of nonzero centered energy requires
$m\ge20$ for either sign of $\tfrac12|W|^2+P(W,\ga)$. The cubic polynomial
\eqref{rev-cubic-profile} attains dimension $20$ for both
Examples~\ref{rev-weyl-example} and~\ref{Ex-2}. The latter realizes
both signs of the extremal energy.
\end{Prop}
\begin{proof}
Fix the degree $d$ of $f$ in \eqref{rev-choice-polynomial}. By
\eqref{rev-choice-radial-matrix}, it suffices for the lower bound to
show that $\mathsf S_d(m)$ is positive definite whenever $m\le19$ and
the radial integral is finite. The integrability condition $m>4d+4$
then gives $d\le3$, so the matrix has size at most $4$.

To apply Sylvester's criterion, let $N$ denote the order of a leading
principal minor of $\mathsf S_d(m)$. Thus $1\le N\le d+1$. The first
$N$ rows and columns are indexed by $0,\ldots,N-1$ and form the matrix
$\mathsf S_{N-1}(m)$. In particular, $N=d+1$ gives the determinant of
the full matrix. We compare each minor with the corresponding Gram
determinant by writing
\[
\De_N(m)=\det\mathsf S_{N-1}(m),\qquad
\mathcal D_N(m)=\det\mathsf G_{N-1}(m)>0.
\]
The largest index of $M_k(m)$ in either determinant is $k=2N-2$,
so the following formula for each order $N$ applies when $m>4N$.
Substituting \eqref{rev-choice-moments} into the determinants gives
\begin{\equ}\label{rev-choice-determinant-ratios}
\begin{aligned}
\frac{\De_1}{\mathcal D_1}&=2,\\
\frac{\De_2}{\mathcal D_2}&=
\frac{m(36-m)}{4(m-1)},\\
\frac{\De_3}{\mathcal D_3}&=
\frac{3m^2(22-m)}{(m-1)(m-2)},\\
\frac{\De_4}{\mathcal D_4}&=
\frac{3m^3(m^2-202m+3600)}{16(m-1)(m-2)(m-3)}.
\end{aligned}
\end{\equ}
Since $\mathcal D_N(m)>0$, each ratio has the sign of the corresponding
leading minor. For $d=0$, the only ratio is $2$. For $d=1$, the second
ratio is positive for $9\le m\le35$, vanishes at $m=36$, and becomes
negative at $m=37$. For $d=2$, the first two ratios are positive
throughout $13\le m\le23$, while the third vanishes at $m=22$ and
first becomes negative at $m=23$. In both cases the smaller leading
minors are positive when the full determinant vanishes, so the matrix
is positive semidefinite at that endpoint. Thus the first dimensions
admitting a negative direction for degrees $1$ and $2$ are $37$ and
$23$, respectively.

For $d=3$, integrability requires $m\ge17$. The first three ratios
are positive for $17\le m\le20$, so in this range it remains to check
the fourth. Its sign is that of $m^2-202m+3600$, whose smaller root
$101-\sqrt{6601}$ lies strictly between $19$ and $20$. Hence all four
leading minors are positive at $m=17,18,19$, whereas $\De_4(20)<0$.
Taken together, these sign computations give positive definiteness of
$\mathsf S_d(m)$ for every admissible degree when $m\le19$.
The vector $V(t)$ in \eqref{rev-choice-radial-matrix} is nonzero for
every $t>0$ because its first entry is $v_0=1$. Positive definiteness
therefore gives $P_{m,f}'(t)>0$ for every $t>0$.
For fixed $W$ and $\ga$ with $\tfrac12|W|^2+P(W,\ga)\ne0$, the centered
energy is therefore strictly monotone in the scale.

It remains to see that the lower bound is attained. For the specific
cubic polynomial \eqref{rev-cubic-profile},
Proposition~\ref{rev-cubic-max} gives a strict local maximum of the
radial polynomial for every $m\ge20$, in particular for $m=20$.
The Weyl data of Examples~\ref{rev-weyl-example} and~\ref{Ex-2}
satisfy \eqref{rev-weyl-sign-condition}, so Theorem~\ref{rev-noncompact-sphere}
gives the corresponding sequences. In Example~\ref{Ex-2} the two
spinor choices give opposite signs of the extremal energy correction.
\end{proof}

For completeness, the determinant formulas also have an orthogonal-polynomial interpretation. Put
\[
w_m(s)=s^{m/2+1}(1+s)^{-m}.
\]
Integration by parts, whose boundary term vanishes by \eqref{rev-choice-integrability}, yields
\begin{\equ}\label{rev-choice-integration-by-parts}
\int_0^\infty\big(2f^2+s(f^2)'\big)w_m\,ds
=\frac m2\int_0^\infty\frac{s-1}{s+1}f(s)^2w_m(s)\,ds.
\end{\equ}
Set $x=(s-1)/(s+1)$. In the new variable the polynomial in
\eqref{rev-choice-polynomial} satisfies
\[
(1-x)^d f\!\left(\frac{1+x}{1-x}\right)
=\sum_{j=0}^d v_j(1+x)^j(1-x)^{d-j}.
\]
The right-hand side has degree at most $d$. Absorbing $(1-x)^{-2d}$
from the square of $f$ into the weight, we see from
\eqref{rev-choice-integration-by-parts} that the quotient
\eqref{rev-choice-rayleigh} is $m/2$ times the Rayleigh quotient of
multiplication by $x$ on polynomials of degree at most $d$ with Jacobi weight
\[
(1-x)^{\alpha_J}(1+x)^{\beta_J},\qquad
\alpha_J=\frac m2-2d-3,\quad
\beta_J=\frac m2+1.
\]
Both parameters exceed $-1$. If $J_N^{(\alpha_J,\beta_J)}$ denotes the classical Jacobi polynomial, normalized by
$J_N^{(\alpha_J,\beta_J)}(1)=\binom{N+\alpha_J}{N}$, then the eigenvalues of the compressed multiplication operator are its zeros for $N=d+1$. This follows from the three-term relation for the Jacobi polynomials. The leading coefficient of $J_N^{(\alpha_J,\beta_J)}$ is
$2^{-N}\binom{m}{N}$, so
\begin{\equ}\label{rev-choice-jacobi-determinant}
\frac{\De_N(m)}{\mathcal D_N(m)}
=\frac{m^N(-1)^N}{\binom{m}{N}}
J_N^{(\alpha_J,\beta_J)}(0),\qquad d=N-1.
\end{\equ}
The finite binomial formula
\[
J_N^{(\alpha_J,\beta_J)}(0)
=2^{-N}\sum_{k=0}^{N}(-1)^{N-k}
\binom{N+\alpha_J}{k}
\binom{N+\beta_J}{N-k}
\]
gives, for $N=1,2,3,4$, respectively,
\[
-2,\qquad \frac{36-m}{8},\qquad
\frac{m-22}{2},\qquad
\frac{m^2-202m+3600}{128},
\]
which proves \eqref{rev-choice-determinant-ratios} as well. In particular, the minimum of \eqref{rev-choice-rayleigh} equals $m/2$ times the smallest zero of the corresponding Jacobi polynomial. This gives a systematic way to optimize the radial coefficients before checking the remaining directions.

For $d=3$ and $m=20$, the minimization in
Remark~\ref{rem:choice-cubic-at-twenty} amounts to fixing $v_3=1$
in $v^T\mathsf S_3(20)v$ and completing the square in $v_0,v_1,v_2$.

\begin{Rem}\label{rev-choice-scope}
The preceding optimality statement concerns extrema of the centered second-order term, with polynomial coefficients independent of the concentrating scale. It is not a lower-dimension obstruction for the spinorial Yamabe equation itself. It does not exclude critical points with $\xi\ne0$, a different relation between the cut-off and concentration scales, coefficients depending on $b$, non-polynomial choices of $f$, or a construction based on higher-order terms. In the present route, however, the obstruction for higher polynomial degrees is intrinsic: their positive leading square forces stronger growth at infinity before any possible gain from additional coefficients can be used.
\end{Rem}

\section{A sharp dimension threshold}\label{sec:radial-dimension-threshold}

Section~\ref{sec: perturbation} constructs a concentrating solution near each
of a sequence of localized metric profiles. Proposition~\ref{rev-choice-optimal-degree}
determines the first dimension in which the centered reduced energy can
have a strict scale extremum of nonzero value. We now ask the
converse question for solutions. We consider a single bubble whose
scale is comparable to that of the profile and whose rescaled center
tends to its origin. Testing the equation in the scale and spinor
directions will give a necessary condition for such concentration.

\subsection{The metric and the concentrating sequence}
\label{subsec:centered-radial-bump-threshold}

Fix an algebraic Weyl tensor $W$. We use the tensor and the norm of
Section~\ref{sec: perturbation}, that is,
\[
 H_{ij}(y)=W_{ikjl}y^ky^l \quad \text{and} \quad
 |W|^2=\sum_{i,j,k,l=1}^m(W_{ijkl}+W_{ilkj})^2.
\]

To compare $H$ with the quadratic term in the expansion of a metric,
define, for each algebraic Weyl tensor $W$,
\begin{\equ}\label{radial-quadratic-weyl-tensor}
 (J_W)_{ij}(y)=-\frac13 W_{ikjl}y^ky^l .
\end{\equ}
Let $G(y)=(\ig_{ij}(y))$ be the matrix of a metric $\ig$ in conformal
normal coordinates centered at $a$, with $\operatorname{Ric}_{\ig}(a)=0$.
At $a$ the Riemann tensor then equals the Weyl tensor $W_{\ig}(a)$.
With our curvature convention, the normal coordinate expansion
in \cite[formula~(5.4)]{LeeParker} gives
\begin{\equ}\label{radial-log-metric-expansion}
 \ln G(y)=J_{W_{\ig}(a)}(y)+O(|y|^3).
\end{\equ}
Comparing \eqref{radial-quadratic-weyl-tensor} with the definition of $H$
gives $J_W=-H/3$. Thus, when $W=W_{\ig}(a)$, the quadratic Taylor term
of $\ln G$ in \eqref{radial-log-metric-expansion} is $-H/3$.
For a general metric, the higher Taylor coefficients require
further curvature data. In our local model we prescribe them as in
\eqref{hij} by multiplying $H(y)$ by a polynomial in $|y|^2$.
For instance, near the center, where the cut-off equals one, the terms
of degrees $4$ and $6$ in the log metric are scalar multiples of
$|y|^2H(y)$ and $|y|^4H(y)$, respectively. This is where the radial
assumption enters.

The endomorphism $\msa_W$ was defined in
\eqref{quadratic-weyl-spinor-form}. The factor $9$ accounts for the
normalization $H=-3J_W$ in the quadratic energy. For
Example~\ref{rev-weyl-example}, $P(W,\ga)=0$ and
$|W|^2=144$, so $\msa_W=8\operatorname{Id}$.
For Example~\ref{Ex-2}, let $A$ be Clifford multiplication by the
oriented volume form of $\operatorname{span}\{e_1,e_2,e_3,e_4\}$.
Then \eqref{nontrivial-Phi} gives $9\msa_W=120\operatorname{Id}-288A$.
By Fact~\ref{fact} and its proof, the eigenvalues of $\msa_W$
are $408/9$ and $-168/9$. Thus the trivial-kernel hypothesis below
allows both signs of the centered energy, even for the same tensor.

Let $f(s)=\sum_{j=0}^d v_js^j$ be a fixed real polynomial with
$v_d\ne0$ and $m>4d+4$. This is precisely the integrability condition
in \eqref{rev-choice-integrability}. Let $x_n$ be the center of a
metric profile and $b_n\to0$ its length scale. In coordinates
$y=(x-x_n)/b_n$, let $R_n\to\infty$ be its cut-off radius and
$\de_n>0$ its amplitude. Using the cut-off $\eta$ from \eqref{hij},
write the local metric as
\begin{\equ}\label{block-local-metric}
 \tilde\ig_n^{\mathrm{loc}}=\exp(h_n),\qquad
 h_n(y)=\de_n\eta(|y|/R_n) f(|y|^2)H(y),\qquad
 \de_n R_n^{2d+2}\to0.
\end{\equ}
For example, if $f(0)=0$, the quadratic jet of this local metric
vanishes at its center even though its higher homogeneous terms can be
nonzero. In the construction \eqref{rev-h-mu}, $R_n=r/b_n$ and
$\de_n=|\vr|a_nb_n^2$ by the homogeneity of $H$. A fixed sign
of $\vr$ can be absorbed into $W$. Let $\tilde\ig_n$ denote the
full rescaled metric, also of determinant one in these coordinates.
The difference between the full and local operators accounts, for
example, for the other profiles in \eqref{rev-h-star}.

As in \eqref{rev-exact-operator}, let $A_n$ and $A_{*,n}$ be
the real self-adjoint operators on $\msd^{1/2}(\R^m,\mbs(\R^m))$
representing the two differences
\[
 \widetilde D_{\tilde\ig_n^{\mathrm{loc}}}-D_{\ig_{\R^m}},
 \qquad
 \widetilde D_{\tilde\ig_n}-\widetilde D_{\tilde\ig_n^{\mathrm{loc}}}.
\]
Thus pairing either operator with a test spinor in the real
$\msd^{1/2}(\R^m,\mbs(\R^m))$ inner product gives the weak action of the corresponding
Dirac difference. We use $\|\cdot\|_{\msd}$ for the Hilbert norm and
$\|\cdot\|_{\mathrm{op}}$ for its operator norm. Let $K_z$ be the
inverse of the Euclidean Hessian on $T_z\cm^\perp$, as in
\eqref{widehat-Phi}. By Lemma~\ref{lemma: critical-manifold} and the
symmetries of \eqref{bubble}, $\sup_{z\in\cm}\|K_z\|_{\mathrm{op}}$
is finite. We assume that for all sufficiently large $n$,
\begin{\equ}\label{threshold-small-operator}
 \|A_n+A_{*,n}\|_{\mathrm{op}}\,
 \sup_{z\in\cm}\|K_z\|_{\mathrm{op}}\leq\frac12 .
\end{\equ}

Consider weak solutions of
$\widetilde D_{\tilde\ig_n}u_n=|u_n|^{2/(m-1)}u_n$ in the rescaled
coordinates. Write $\ell_n$ and $\zeta_n$ for the scale and center of
the bubble in these coordinates. Thus a bubble of original scale
$\lambda_n$ and center $\xi_n$ has
$\ell_n=\lambda_n/b_n$ and $\zeta_n=(\xi_n-x_n)/b_n$.
We assume that, for unit spinors $\ga_n$,
\begin{\equ}\label{block-centered-solution}
 \begin{gathered}
 u_n=z_n+w_n,\qquad z_n=\psi_{\ell_n,\zeta_n,\ga_n},\qquad
 w_n\perp_{\msd}T_{z_n}\cm,\qquad \|w_n\|_{\msd}\to0,\\[0.3em]
 \ell_n\to\ell_*\in(0,\infty),\qquad \zeta_n\to0, \qquad \text{as } n\to\infty.
 \end{gathered}
\end{\equ}
The last two limits express concentration at the scale of the profile
and asymptotic centering. By the critical Sobolev embedding, the energy
$(2m)^{-1}\int_{\R^m}|u_n|^{2m/(m-1)}\,dy$ then tends to the
standard bubble energy $E_0$ in \eqref{rev-local-reduced-expansion}.

Let $\mathcal V_n$ be the real span of $z_n$, its scale derivative,
and its derivatives along the unit spinor sphere. We require the
exterior part of the operator to satisfy
\begin{\equ}\label{block-tail-conditions}
 \sup_{\substack{Z\in\mathcal V_n\\\|Z\|_{\msd}=1}}
      \|A_{*,n}Z\|_{\msd}=o(\de_n) \quad \text{and} \quad
 \sup_{\substack{Z\in\mathcal V_n\\\|Z\|_{\msd}=1}}
      |\inp{A_{*,n}z_n}{Z}_{\msd}|=o(\de_n^2).
\end{\equ}
The first estimate controls the forcing of the normal correction.
The second makes the exterior contribution negligible in the scale
and spinor equations at order $\de_n^2$.

\begin{Thm}\label{centered-radial-bump-threshold}
Suppose $\ker\msa_W=\{0\}$. A sequence satisfying
\eqref{block-local-metric}-\eqref{block-tail-conditions} can exist
only if $m\geq20$. This holds for every fixed polynomial of degree
$d$ with $m>4d+4$. Conversely, in every dimension $m\geq20$,
each extremum satisfying \eqref{rev-uniform-local-gap} for the cubic
in \eqref{rev-cubic-profile} gives a sequence in this class on the
sphere of Theorem~\ref{rev-noncompact-sphere}. In particular, both
Examples~\ref{rev-weyl-example} and~\ref{Ex-2} attain the threshold,
and the latter gives two sequences for a single metric.
\end{Thm}

{Theorem~\ref{centered-radial-bump-threshold} makes precise
the phrase ``of this type'' in Theorem~\ref{main-thm}~(3). In the coordinates
of the polynomial radial metric profile \eqref{block-local-metric}, the
center of the single bubble tends to the origin and its scale tends to a
finite limit $\ell_*>0$, as in \eqref{block-centered-solution}.
We give the proof in Subsection~\ref{subsec:proof-centered-radial-threshold}.
The operator bound \eqref{threshold-small-operator} and the exterior
estimates \eqref{block-tail-conditions} allow us to pass to the limit
when we test the full equation in the spinor and scale directions.
The limiting unit spinor must be an eigenspinor of $\msa_W$. Since
$\ker\msa_W=\{0\}$, its eigenvalue is nonzero, and the scale equation
therefore gives $P_{m,f}'(\ell_*^2)=0$. The calculation in the proof of
Proposition~\ref{rev-choice-optimal-degree} gives $P_{m,f}'(t)>0$ for
every $t>0$ when $4d+4<m\le19$, ruling out concentration in this range.
For $m\ge20$, we verify that the sequences constructed in
Theorem~\ref{rev-noncompact-sphere} satisfy the same hypotheses.}

\subsection{The quadratic energy at a centered bubble}

We use the two leading coefficients of the Dirac operator in
\eqref{Dirac-expansion}. For a symmetric trace-free tensor $h$, write
$D_1[h]$ and $D_2[h]$ for the coefficients of $t$ and $t^2$ in
$\widetilde D_{\exp(th)}$. For a tensor $B$ we use the convention
$(\diver B)^i=\sum_j\partial_j B_{ij}$. The two coefficients are
\begin{\equ}\label{next-D2}
\begin{aligned}
D_1[h]u={}&-\frac12\sum_{i,j}h_{ij}e_i\cdot_{\ig_{\R^m}}\partial_j u
            -\frac14(\diver h)\cdot_{\ig_{\R^m}}u,\\
D_2[h]u={}&\frac18\sum_{i,j}(h^2)_{ij}e_i\cdot_{\ig_{\R^m}}\partial_j u
             +\frac1{16}\diver(h^2)\cdot_{\ig_{\R^m}}u\\
&-\frac1{16}\sum_{\substack{i,j,k \\ i\neq j\neq k\neq i}}\sum_d
 (\partial_i h_{jd})h_{dk}\,
 e_i\cdot_{\ig_{\R^m}}e_j\cdot_{\ig_{\R^m}}e_k\cdot_{\ig_{\R^m}}u .
\end{aligned}
\end{\equ}
At $t=1$, for $|h|$
sufficiently small, the remaining terms are bounded pointwise by
$C(|h|^3|\nabla u|+|h|^2|\nabla h||u|)$.
The quadratic coefficient of the energy is
$\Phi(u)=\tfrac12\int_{\R^m}\real(D_2[h]u,u)\,dx$.

\begin{Lem}\label{quadratic-weyl-evaluation}
Take $H$ as in Section~\ref{sec: perturbation}, and put
$h(x)=F(|x|^2)H(x)$ for a smooth real function $F$.
Let $\ell>0$, let $\ga$ be a unit spinor, and set $z=\psi_{\ell,0,\ga}$.
On the Euclidean ball $B_R(0)$ of any radius $R>0$,
\begin{\equ}\label{radial-centered-pairing}
\frac12\int_{B_R(0)}\real(D_2[h]z,z)\,dx
=\frac{m^{m-1}\om_{m-1}}{16m(m+2)}
 \left(\frac{|W|^2}{2}+P(W,\ga)\right)
\cdot\ell^m\int_0^R
 \frac{\rho^{m+3}F(\rho^2)^2}{(\ell^2+\rho^2)^m}\,d\rho .
\end{\equ}
\end{Lem}

\begin{proof}
The angular calculation in Lemma~\ref{lemma: critical-manifold}\,(5)
applies on each centered sphere. Indeed, in the alternating term of
\eqref{next-D2} every derivative of $F$ multiplies $(H^2)_{jk}$,
which is symmetric in $j,k$. The terms of odd parity in the remaining
Clifford product integrate to zero on the sphere. Thus the same
contraction $P(W,\ga)$ multiplies the integral of $(x^1)^2(x^2)^2$
over the sphere. On the unit sphere this integral equals
$\om_{m-1}/[m(m+2)]$. For the scalar term, the same calculation
as in \eqref{W1W2-Phi} gives
\[
 \int_{S^{m-1}}\tr(H(\theta)^2)\,d\theta
 =\frac{\om_{m-1}}{2m(m+2)}|W|^2,
\]
with the norm of Section~\ref{sec: perturbation}. Both angular
contributions are homogeneous of degree four in $x$. Integrating in
$\rho$ yields \eqref{radial-centered-pairing}.
\end{proof}

The spinor dependence in this formula can also be written intrinsically.
For $1\le \al,\bt\le m$ let
$\omega_{\al\bt}=\sum_{i<j}W_{\al\bt ij}\,dx^i\wedge dx^j$ and identify
exterior forms with their Euclidean Clifford actions.
In the formula for $P(W,\ga)$, only the completely alternating part
in the four Clifford indices survives the real pairing. Antisymmetrizing
in the first and last pairs and using the algebraic Bianchi identity gives
\[
 P(W,\ga)=-9\sum_{\al,\bt}
 \real\bigl((\omega_{\al\bt}\wedge\omega_{\al\bt})
       \cdot_{\ig_{\R^m}}\ga,\ga\bigr).
\]
Thus \eqref{quadratic-weyl-spinor-form} retains both the scalar
and the four-form contributions of the Weyl tensor.

For $F=f$ of degree $d$ and $m>4d+4$, the radial integral in
\eqref{radial-centered-pairing} converges at infinity. In particular
the formula passes to the limit from radial cut-offs even at
$m=4d+5$, since the odd terms have already vanished on each sphere.
Writing $A_m=\int_0^\infty s^{m+3}(1+s^2)^{-m}\,ds$ as in
\eqref{rev-moments}, the substitution $\rho=\ell s$ and
\eqref{rev-choice-radial-profile} give
\begin{\equ}\label{threshold-centered-energy}
 \lim_{R\to\infty}\frac12\int_{B_R(0)}\real(D_2[f(|x|^2)H]z,z)\,dx
 =\frac{9m^{m-1}\om_{m-1}A_m}{16m(m+2)}
 \real(\msa_W\ga,\ga) P_{m,f}(\ell^2).
\end{\equ}
For the tensor in Example~\ref{rev-weyl-example}, the coefficient
is $72\kappa_m$, in agreement with \eqref{rev-central-profile}.

\subsection{The scale derivative at the center}

In \eqref{threshold-centered-energy}, fixing $W$ and the unit spinor
$\ga$ fixes the value of $\real(\msa_W\ga,\ga)$. The scale
dependence is then entirely through $P_{m,f}(\ell^2)$.
For the scale test we need the following consequence of the matrix
calculation in Section~\ref{rev-other-choices}. It also covers
polynomials whose constant coefficient vanishes.

\begin{Cor}\label{radial-scale-positivity}
Let $f$ be a nonzero real polynomial of degree $d$ and assume
$m>4d+4$. Then $P_{m,f}(t)>0$ for every $t>0$. If $m\leq19$,
then $P_{m,f}'(t)>0$ for every $t>0$ as well.
\end{Cor}
\begin{proof}
For $V(t)=(v_0,v_1t,\ldots,v_dt^d)^T$, the identities in
\eqref{rev-choice-radial-profile} and \eqref{rev-choice-radial-matrix} give
\[
 P_{m,f}(t)=t^2 V(t)^T\mathsf G_d(m)V(t),\qquad
 P_{m,f}'(t)=t V(t)^T\mathsf S_d(m)V(t).
\]
The first matrix is the Gram matrix of $1,s,\ldots,s^d$ for the
positive weight in \eqref{rev-choice-moments}. If $4d+4<m\leq19$,
then $d\leq3$ and every leading principal minor of $\mathsf S_d(m)$ is
positive by \eqref{rev-choice-determinant-ratios}. Since $V(t)\ne0$ for
$t>0$, the two claims follow.
\end{proof}

\subsection{Proof of Theorem~\ref{centered-radial-bump-threshold}}
\label{subsec:proof-centered-radial-threshold}

\begin{proof}[Step 1. Exclusion for $m\leq19$]
We work in the rescaled coordinates under the metric \eqref{block-local-metric}.
The identities $H(y)y=0$, $\sum_j\partial_j H_{ij}=0$ and
$\operatorname{tr}H=0$ also hold for $h_n$, since its scalar factor
is radial. For $z=\psi_{\ell,\zeta,\ga}$ the first variation of
the Dirac operator at the Euclidean metric consequently reduces to
\[
 -\frac12\sum_{i,j}(h_n)_{ij}e_i\cdot_{\ig_{\R^m}}\partial_j z
 =-\frac{m}{2(\ell^2+|y-\zeta|^2)}
       (h_n\zeta)\cdot_{\ig_{\R^m}}z .
\]
Its $L^{2m/(m+1)}$-norm is $O(\de_n|\zeta|)$, uniformly for
$\ell$ in a compact subset of $(0,\infty)$ and $\zeta$ near zero.
Indeed, the pointwise bound is
$C\de_n|\zeta|(1+|y|)^{2d+1-m}$.
The same bound holds after differentiating in scale or in the
unit spinor. The terms of order two and higher are bounded by
$C\de_n^2(1+|y|)^{4d+4-m}$ on the cut-off support.
Their $L^{2m/(m+1)}$-norm is at most
$C\de_n^2(1+R_n^{2d+2})=o(\de_n)$,
also when $m=4d+5$.

It follows from $\zeta_n\to0$ and \eqref{block-tail-conditions} that
$\|(A_n+A_{*,n})z_n\|_{\msd}=o(\de_n)$.
Projecting the equation onto $T_{z_n}\cm^\bot$ and using the
uniform bound on the inverse of the normal Hessian gives
$\|w_n\|_{\msd}=o(\de_n)$. Here the small operator bound
absorbs the linear metric term in $w_n$, and differentiability of
$u\mapsto |u|^{2/(m-1)}u$ from $\msd^{1/2}(\R^m,\mbs(\R^m))$ to its dual gives a
remainder $o(\|w_n\|_{\msd})$, uniformly over the bubble parameters.

Now test the full equation with a scale or unit-spinor tangent $Z$
at $z_n$ of bounded norm. The Euclidean Hessian annihilates $Z$.
Since $|Z|\leq C|z_n|$, the absolute value of the integrand in the
pairing of the nonlinear Taylor remainder with $Z$ is bounded by a
constant multiple of $|z_n|^{2/(m-1)}|w_n|^2$, uniformly in $n$.
This follows by Taylor expansion where $|w_n|\leq |z_n|/2$ and by
$2/(m-1)\leq1$ on the complement. All standard bubbles have the same
$L^{2m/(m-1)}$ norm, so H\"older's inequality and the Sobolev embedding give
\[
 \int_{\R^m}|z_n|^{2/(m-1)}|w_n|^2\,dy
 =O\bigl(\|w_n\|_{\msd}^2\bigr)=o(\de_n^2).
\]
Self-adjointness and the source estimates above give
$\big|\inp{w_n}{(A_n+A_{*,n})Z}_{\msd}\big|=o(\de_n^2)$.
The second bound in \eqref{block-tail-conditions} handles the direct
exterior pairing. The first-order local energy vanishes on $\cm$
by Lemma~\ref{lemma: critical-manifold}, so its derivative in the
direction $Z$ vanishes as well.

Thus only the derivatives of the quadratic local energy on $\cm$
can contribute at order $\de_n^2$. We give the estimates that allow
them to be evaluated at the center. The absolute contribution of the
higher terms in the metric expansion to the pairing with $Z$ is
bounded by a constant multiple of the radial integral below, uniformly
in $n$. Since $m\geq4d+5$, we have $6d+6-m\leq2d+1$, so
\[
 \de_n^3\int_0^{R_n}(1+r)^{6d+6-m}\,dr
 =O\bigl(\de_n^3(1+R_n^{2d+2})\bigr)=o(\de_n^2),
\]
where the last step follows from $R_n\to\infty$ and
$\de_n R_n^{2d+2}\to0$ in \eqref{block-local-metric}.
In the quadratic term a derivative in $\zeta$ gains one power of
decay from the bubble or tangent. Including the polar volume factor, the
absolute bound is $C\de_n^2(1+r)^{4d+3-m}\,dr$.
Replacing $\zeta_n$ by zero therefore costs $O(\de_n^2|\zeta_n|)$.
At the center, \eqref{radial-centered-pairing} applies to the radial
factor $\eta(r/R_n)f(r^2)$. The odd part vanishes on each sphere,
including in the borderline dimension $m=4d+5$. Since $m>4d+4$,
the remaining radial integrals and their scale derivatives converge
as $R_n\to\infty$.
In particular, after a scale change the term differentiated on
the cut-off annulus is $O(R_n^{4d+4-m})$.

Pass to a subsequence with $\ga_n\to\ga_*$. By
\eqref{radial-centered-pairing} and \eqref{rev-choice-radial-profile},
the limiting scale and spinor equations are precisely stationarity of
\[
 (\ell,\ga)\longmapsto P_{m,f}(\ell^2)
                    \real(\msa_W\ga,\ga)
\]
at $(\ell_*,\ga_*)$, after dividing by a positive dimensional constant.
Since $P_{m,f}(\ell_*^2)>0$, the spinor equation makes $\ga_*$ an
eigenspinor of $\msa_W$. The scale equation gives
$P'_{m,f}(\ell_*^2)\real(\msa_W\ga_*,\ga_*)=0$.
For $m\leq19$ and $m>4d+4$, \eqref{rev-choice-radial-matrix} and
Corollary~\ref{radial-scale-positivity} give $P'_{m,f}(\ell_*^2)>0$.
The selected eigenvalue must therefore vanish, contrary to the
assumption on $\msa_W$.
\end{proof}

\begin{proof}[Step 2. Attainment for every $m\geq20$]
Take a tensor $W$ with $\ker\msa_W=\{0\}$ and a centered
extremum as in \eqref{rev-uniform-local-gap}. Theorem~\ref{rev-noncompact-sphere}
gives the required scale and center limits. Its normal corrections lie
in $T_{z_n}\cm^\perp$ and tend to zero in norm by
\eqref{rev-local-w-bound} and the comparison proving \eqref{rev-glued-limit}.
This proves \eqref{block-centered-solution}. With $\vr>0$ sufficiently
small and the parameters of \eqref{rev-final-parameters}, the
$n$-th block in the coordinates $y=(x-x_n)/b_n$ has amplitude
$\de_n=\vr b_n^9$ and cut-off radius $R_n=r/b_n$. Hence
$\de_nR_n^8\to0$. The bound on the full metric operator follows
from \eqref{rev-coefficient-bound} and $\sum_n b_n<\infty$.

The exterior estimate \eqref{rev-remote-source} also holds with the
bubble replaced by its rescaled scale or spinor tangents. These have
the same pointwise bounds as the bubble and its spatial derivative
on each remote patch. The corresponding direct pairings with the
bubble are bounded by $C\mathcal T_n^2$, with $\mathcal T_n$ given
by \eqref{rev-positive-tail}. Indeed, differentiating
$\Ga^{(j)}|_{\cm}=0$ removes the linear term on every remote patch.
The remaining terms satisfy the product estimate used in
\eqref{rev-direct-tail-square}, with one bubble factor replaced by
its tangent. The estimates are uniform on a compact set of rescaled
parameters and for unit vectors in their span. Now
\eqref{rev-tail-ratio} gives $\mathcal T_n=o(\de_n)$, proving
both parts of \eqref{block-tail-conditions}.

For Example~\ref{rev-weyl-example}, the only eigenvalue is $8$.
For Example~\ref{Ex-2}, they are $408/9$ and $-168/9$.
Thus both tensors meet the trivial-kernel hypothesis, and the same
estimates apply to each of their centered extrema. Since the eigenvalues
of $\msa_W$ depend continuously on $W$, the requirement
$\ker\msa_W=\{0\}$ is an open condition on the algebraic Weyl tensor $W$.
Thus the nearby tensors considered in Subsection~\ref{rev-cubic-subsection}
also satisfy this hypothesis when the neighborhoods are chosen sufficiently small.
\end{proof}

\subsection{Dimension thresholds from the scale equations}
\label{subsec:centered-scalar-thresholds}

\begin{revision}
All Euclidean standard bubbles have the same energy. Once the metric
varies, the energy of a bubble can depend on its size even if its center
is fixed. In which dimension can the leading correction first have a
critical scale for some choice of the radial polynomial? We
consider this question for the spinorial Yamabe equation, the scalar
Yamabe equation and the fourth- and sixth-order constant $Q$-curvature
equations, using the same family of metric variations in each case.

Fix a nonzero algebraic Weyl tensor $W$ and let
$H_{ij}(y)=W_{i\al j\bt}y^\al y^\bt$ as in Section~\ref{sec: perturbation}.
For a nonzero real polynomial $f(s)=\sum_{j=0}^d v_js^j$, take the log metric
to be $\vr f(|y|^2)H(y)$ and place the bubble at the origin. We work
in the range where the integrals below converge as the cut-off is
removed. Write $F(\ell)$ for the coefficient of $\vr^2$ in the reduced
energy of the bubble at relative scale $\ell>0$, keeping the unit spinor
fixed in the spinorial case. To compute $F(1)$ one expands the energy
in the given metric and includes the first correction to the bubble
obtained from the linearized equation, as in \eqref{widehat-Phi}.
For the radial profiles considered here, the first variation of the
operator applied to the centered bubble vanishes. In the spinorial
case this follows from $H(y)y=0$, $\diver H=0$ and $\tr H=0$, as in Step~1 of
the proof of Theorem~\ref{centered-radial-bump-threshold}. The normal
correction therefore makes no contribution to $F(1)$.

Keeping the coefficients $v_j$ free, write the resulting quadratic form
at unit scale as $F(1)=\sum_{i,j=0}^d B_{ij}v_iv_j$, with $B_{ij}=B_{ji}$.
Thus $B_{ij}$ can be calculated by pairing the two metric terms
$|y|^{2i}H(y)$ and $|y|^{2j}H(y)$ in the quadratic energy. The term
$|y|^{2j}H(y)$ is homogeneous of degree $k_j=2j+2$, so the change of
variables $y=\ell x$ replaces $v_j$ by $c_j(\ell)=v_j\ell^{k_j}$.
Consequently,
\[
 F(\ell)=\sum_{i,j=0}^d B_{ij}c_i(\ell)c_j(\ell),\qquad
 \ell F'(\ell)=\sum_{i,j=0}^d(k_i+k_j)B_{ij}c_i(\ell)c_j(\ell).
\]
The dependence of the operator on the metric determines $B_{ij}$
at one scale, and homogeneity gives the matrix of the scale derivative.
For a prescribed $f$, substituting $t=\ell^2$ gives the polynomial
$F(\sqrt t)=\sum_{i,j=0}^d B_{ij}v_iv_j t^{i+j+2}$.

For the spinorial Yamabe equation, the angular integration has already
been carried out in \eqref{threshold-centered-energy}. With
$t=\ell^2$, the result is
\[
 F(\ell)=9\kappa_m\,\real(\msa_W\ga,\ga)\,P_{m,f}(t).
\]
The endomorphism $\msa_W$ acts on spinors. By
Corollary~\ref{radial-scale-positivity}, $P_{m,f}(t)>0$ for $m>4d+4$.
Thus stationarity with respect to $\ga$ on the unit spinor sphere
selects an eigenspinor of $\msa_W$. If the corresponding eigenvalue
is $\bt_*\ne0$, stationarity in $\ell$ requires $P_{m,f}'(t)=0$.
The hypothesis $\ker\msa_W=\{0\}$ in
Theorem~\ref{centered-radial-bump-threshold} ensures this for every
limiting eigenspinor. The matrix that tests for a critical scale,
$\mathsf S_d(m)$ in \eqref{rev-choice-radial-matrix}, acts on the
coefficients of $f$. In fact, for
$V(t)=(v_0,v_1t,\ldots,v_dt^d)^T$ we have
\[
 \frac{P_{m,f}'(t)}{t}=V(t)^T\mathsf S_d(m)V(t).
\]
The same matrix $\mathsf S_d(m)$ therefore applies to every choice of
$W$ and every spinor branch with nonzero eigenvalue. Its definiteness
tests all polynomials of the allowed degree at once.

There is a spatial interpretation of this last formula. With $r=|y|$,
\[
 \ell\,\partial_\ell\ln\!\left(\frac{\ell}{\ell^2+r^2}\right)
 =\frac{r^2-\ell^2}{r^2+\ell^2}.
\]
Differentiating the integral in \eqref{radial-centered-pairing} shows
that $P_{m,f}'(\ell^2)=0$ precisely when
\begin{\equ}\label{threshold-inner-outer-balance}
 \int_0^\ell \frac{(\ell^2-r^2)r^{m+3}f(r^2)^2}{(\ell^2+r^2)^{m+1}}\,dr
 =
 \int_\ell^\infty \frac{(r^2-\ell^2)r^{m+3}f(r^2)^2}{(\ell^2+r^2)^{m+1}}\,dr .
\end{\equ}
Both sides are nonnegative. For $4d+4<m\leq19$ the contribution from
$r>\ell$ is strictly larger for every nonzero $f$ and every $\ell>0$,
by Corollary~\ref{radial-scale-positivity}. For $d=3$ the determinant
in \eqref{rev-choice-determinant-ratios} changes sign at
$101-\sqrt{6601}$, between $19$ and $20$. At $m=20$ the cubic
in \eqref{rev-cubic-profile} realizes equality at a strict local
maximum of $P_{m,f}$. For the two eigenspaces in Example~\ref{Ex-2},
the same balance gives a scale maximum on the positive branch and a
scale minimum on the negative branch.

We now perform the comparison for the scalar equations. Write $n$ for
the dimension and $2\sigma\in\{2,4,6\}$ for the order of the operator.
The standard positive bubble at the origin is a constant multiple of
$(\ell/(\ell^2+|y|^2))^{(n-2\sigma)/2}$. For these equations the
first-order term vanishes on this bubble for the tensor $f(|y|^2)H(y)$.
After integration over the sphere the dependence on $W$ reduces to a
positive multiple of its squared norm. For these two calculations, see
\cite[Propositions~17 and~18]{BrMa} for the Yamabe equation,
\cite[Lemma~9.15 and Proposition~9.16]{WZ13} for the fourth-order
equation, and \cite[Lemma~11.16 and Corollary~11.19]{GKW25} for the
sixth-order equation.
Thus any nonzero $W$ gives the same scale matrix up to a positive
factor, and there is no spinor direction to select. We take
$4d+4<n-2\sigma$, so the quadratic integrals converge without a
cut-off. At equality the highest degree requires a logarithmic estimate.

To use the sign convention of the scalar compactness calculations,
let $M_{ij}=-(k_i+k_j)B_{ij}$ and write
$c(\ell)=(c_0(\ell),\ldots,c_d(\ell))^T$. Then
$-\ell F'(\ell)=c(\ell)^TMc(\ell)$. Suppose that a sequence of
metric profiles has amplitudes $\vr_\nu\to0$, and that in rescaled
coordinates the bubble centers tend to zero and the scales tend to
$\ell_*>0$. If the errors due to the normal correction, displacement
of the center, cut-off and exterior metric are $o(\vr_\nu^2)$ in
the scale equation, testing with the scale derivative of the bubble gives
\begin{\equ}\label{scalar-centered-scale-test}
 0=\vr_\nu^2\,c(\ell_*)^T M c(\ell_*)+o(\vr_\nu^2).
\end{\equ}
A definite matrix excludes such concentration for every nonzero $f$.
To exhibit a critical scale in a given dimension, it suffices to find
one choice of coefficients for which $F'(\ell)$ vanishes.

For the scalar Yamabe equation, Brendle and Marques show in
\cite[Proposition~17]{BrMa} that $F(\ell)$ is a negative constant times
\[
 \int_{\R^n}\left|\nabla\bigl(f(|y|^2)H(y)\bigr)\right|^2
 \left(\frac{\ell}{\ell^2+|y|^2}\right)^{n-2}\,dy .
\]
Here the squared gradient supplies a
nonnegative density, so differentiation gives the same division at
$r=\ell$ as in \eqref{threshold-inner-outer-balance}. The entries of
$H$ are harmonic homogeneous polynomials of degree $2$, with spherical
Laplacian eigenvalue $-2n$. Integration in polar coordinates gives,
up to a common positive factor,
\[
 M_{ij}=(k_i+k_j)(k_i k_j+2n)I_{k_i+k_j},\qquad
 I_L=\int_0^\infty r^{n+L-3}(1+r^2)^{2-n}\,dr .
\]
For the homogeneous degrees $2,4,6,8$, the identity
$I_{L+2}/I_L=(n+L-2)/(n-L-4)$ gives four positive leading principal
minors at $n=24$. At $n=25$ the first three remain positive and the
determinant is negative. This same family occurs in the analysis of
Khuri, Marques and Schoen. They prove in \cite[Proposition~A.4]{KMS}
that their Pohozaev quadratic form is coercive on the full space of
admissible metric Taylor terms through dimension $24$. For their
choice of a cubic $f$ giving a negative direction $f(|y|^2)H(y)$
in dimension $25$, see \cite[Proposition~A.8]{KMS}. In their
noncompactness construction, Brendle and Marques choose the coefficients
of $f$ by differentiating the centered energy and then check the
translation directions. We refer the reader to
\cite[Propositions~18 and~22 and Corollary~23]{BrMa} for these steps.

For the fourth- and sixth-order constant $Q$-curvature equations,
Gong, Kim and Wei denote the entries of the corresponding matrix by
$m_{ij}^{D,2}$, where $D,2$ specifies the divergence-free tensors
of angular degree $2$. Their sixth-order construction makes the
relation with the scale derivative explicit. With their quartic
coefficients $a_0,\ldots,a_4$ and bubble scale $\de$, the centered
energy correction is, up to a positive factor,
$-\sum_{i,j=0}^4 \frac{m_{ij}^{D,2}}{k_i+k_j}
 a_i a_j\de^{k_i+k_j}$.
Taking the negative scale derivative at $\de=1$ gives exactly the
quadratic form with entries $m_{ij}^{D,2}$, as in the proof of
\cite[Corollary~11.19]{GKW25}. For general even order $2N$, they
obtain the centered energy from the Hessian of total $Q_{2N}$-curvature
and differentiate it with respect to the bubble scale. For the resulting
matrix and its derivation, see \cite[formula~(1.5) and Proposition~5.3]{GKW26}.
At unit scale the
radial measure in $s=r^2$ is again proportional to
$s^{n/2+1}(1+s)^{-n}\,ds$, as in \eqref{rev-choice-moments}. The
integral pairs $f$ with a radial differential operator applied to $f$.
Its matrix records which homogeneous degrees can cancel. A spatial
balance of two nonnegative integrals as in
\eqref{threshold-inner-outer-balance} also requires a nonnegative
density.

For the fourth- and sixth-order equations the scale matrix on the
admissible coefficients is positive definite through dimensions $24$
and $26$, respectively. The same family of Weyl tensors multiplied
by radial polynomials already gives a negative direction at $25$ and
$27$. For the signs of these matrices, see
\cite[proofs of Lemmas~9.1 and~10.36]{GKW25}. The constructions of
Wei and Zhao in \cite{WZ13} and of Gong, Kim and Wei in \cite{GKW25}
realize a critical scale in the respective dimensions.
The choice of degree matters. As illustrated in \cite[Figure~1]{GKW26},
the determinant for the cubic family in the fourth-order problem first
changes sign at a real dimension of approximately $25.406$, whereas
the quartic used by Wei and Zhao already gives a
critical scale at $n=25$.

For the spinorial sequences of Theorem~\ref{centered-radial-bump-threshold},
the centering and exterior estimates make the local scale derivative
a necessary test even when all curvature jets at the limiting point
agree with those of the round metric. Theorem~\ref{rev-noncompact-sphere}
verifies the remaining center and spinor variations and constructs
solutions attaining the threshold. To treat centers that drift,
relative scales tending to zero or infinity, or exterior terms of the
same order, one would need further control of the solution away from
the bubble. In the scalar Yamabe problem, comparison of positive
solutions with the Green function supplies such a pointwise estimate
(see \cite[Proposition~2.5]{KMS}).

For the radial metric variations $f(|y|^2)H(y)$, we can thus obtain
the quadratic form in the coefficients of $f$ directly by differentiating
the centered reduced energy with respect to the bubble scale, without
using a local Pohozaev identity. The comparison with the classical
compactness proofs explains why this calculation detects their
dimension thresholds. Those proofs use coercivity of the Pohozaev
quadratic form on the full space of admissible metric Taylor terms to
deduce curvature vanishing, and then use mass estimates to exclude
concentration under their hypotheses (see, for instance,
\cite{KMS,GKW25,GKWambient}).
On our radial family the scale derivative agrees with the restriction
of that form up to sign and normalization. In each of the three scalar
problems this family already contains a negative direction in the
first dimension where the full form loses positive definiteness. For
the negative direction in the Yamabe problem, see
\cite[Proposition~A.8]{KMS}. The corresponding signs for all three
equations are given in \cite[Corollary~4.7]{GKWambient}.
This is why the centered calculation detects dimension $25$ for the
Yamabe and fourth-order equations and dimension $27$ for the
sixth-order equation.
\end{revision}

\begin{appendices}

\section{Calculations for the cubic radial profile}\label{app:cubic-calculations}
\subsection{The cut-off limit and the translation Hessian}\label{app:cubic-translations}

To establish the cut-off limit \eqref{Phi-infty}, substitute $h(x)=f(|x|^2)H(x)$ into
\eqref{widehat-Phi}. {We first consider the two integrals
defining $\Phi(\psi_{\lm,\xi,\ga})$ in Lemma~\ref{lemma: cl-expansion}, which have been specified in \eqref{Phi-on-M} and \eqref{Phi-on-M2}.
Write $\rho=|x|$ and let the parameters $\lm$ and $\xi$ vary in a fixed compact
set. Since $f$ is cubic and $H$ is homogeneous of degree two, each
$h_{ij}$ is $O(\rho^8)$ and each $\pa_i h_{jd}$ is $O(\rho^7)$ at
infinity. In the Clifford three-form term \eqref{Phi-on-M2}, the
coefficient $(\pa_i h_{jd})h_{dk}$ is therefore $O(\rho^{15})$,
and the spinor pairing is bounded by
$|\psi_{\lm,\xi,\ga}(x)|_{\ig_{\R^m}}^2=O(\rho^{2-2m})$ by
\eqref{bubble}. Including the polar volume factor $\rho^{m-1}$ gives that the integrand in \eqref{Phi-on-M2} behaves as
$O(\rho^{16-m})$ at infinity. For the term containing $h^2$, \eqref{Phi-on-M}
gives the sharper bound
\[
 \rho^{m-1}\frac{\lm^m\tr(h(x)^2)}{(\lm^2+|x-\xi|^2)^m}
 =O(\rho^{m-1+16-2m})=O(\rho^{15-m}).
\]
The powers $\rho^{16-m}$ and $\rho^{15-m}$ are integrable at infinity
for $m\ge18$ and $m\ge17$, respectively.}
The source in
\eqref{rev-source} is bounded by $C(1+|x|)^{8-m}$, which belongs to
$L^{2m/(m+1)}(\R^m)$ for $m>17$. The same bounds apply with the
cut-off $\eta(|x|/R)$, uniformly for $R\ge1$. Dominated convergence
therefore gives \eqref{Phi-infty} uniformly on compact sets of bubble parameters.

For the tensor in Example~\ref{rev-weyl-example}, we already know that
$\tfrac12|W|^2+P(W,\ga)=72$ for every unit spinor. We first show that
the Clifford three-form term in $\Phi(\psi_{\lm,\xi,\ga})$ vanishes
for every center $\xi$.
In the Euclidean frame $e_j=\partial/\partial x^j$, define the vector
fields $U$ and $V$ by
\[
 U(x)=U^j(x)e_j=-x^4e_1+x^1e_4,\qquad
 V(x)=V^j(x)e_j=x^3e_2-x^2e_3 .
\]
Thus $U^j$ and $V^j$ are their scalar components, and their squared
Euclidean lengths are $|U|^2=(x^1)^2+(x^4)^2$ and
$|V|^2=(x^2)^2+(x^3)^2$.
Direct substitution of the components in {Table~\ref{tab:weyl-vanishing-spinor-term}} gives
\begin{\equ}\label{rev-table1-structure}
\left\{
\aligned
&H=3(U\otimes V+V\otimes U), \\[0.3em]
&{\tr(H(x)^2)}
=18\big((x^1)^2+(x^4)^2\big)
\big((x^2)^2+(x^3)^2\big).
\endaligned 
\right.
\end{\equ}
For this tensor, all components of $H$ with an index greater than $4$ vanish.
To see why the Clifford three-form term in $\Phi$ vanishes for this
tensor, take $h=FH$ with $F$ an arbitrary smooth radial function.
Substituting $h=FH$ into the first term of $\Phi(\psi)$ in
Lemma~\ref{lemma: cl-expansion} gives
\begin{\equ}\label{three-form-F}
 -\frac{1}{32}
 \sum_{\substack{i,j,k \\ i\neq j\neq k\neq i}}
 \int_{\R^m}\left(\sum_d\partial_i(FH_{jd})\,F H_{dk}\right)
 \real\bigl(e_i\cdot_{\ig_{\R^m}}e_j\cdot_{\ig_{\R^m}}e_k
       \cdot_{\ig_{\R^m}}\psi,\psi\bigr)_{\ig_{\R^m}}\,d\vol_{\ig_{\R^m}} .
\end{\equ}
For pairwise distinct indices, interchanging any two of $i,j,k$
reverses the sign of the Clifford product. Thus, for each $i<j<k$,
the six permutations contribute the alternating sum of the scalar
coefficients multiplying the same ordered Clifford product.
When the derivative falls on $F$, the sum over $d$ gives
$F(\partial_iF)(H^2)_{jk}$. Since $H^2$ is symmetric, these terms
cancel in pairs when $j$ and $k$ are interchanged. When the derivative
falls on $H$, we use $H_{jd}=3(U^jV^d+V^jU^d)$ from
\eqref{rev-table1-structure}. After summing over $d$ and discarding
terms symmetric in $j,k$, the two remaining products, up to the
common factor $9F^2$, are
\[
 |V|^2(\partial_i U^j)U^k
 \quad\hbox{and}\quad
 |U|^2(\partial_i V^j)V^k.
\]
The first product can be nonzero only if $i,j,k\in\{1,4\}$, and the
second only if $i,j,k\in\{2,3\}$. Both vanish when the three indices
are pairwise distinct. Hence the three-form vanishes pointwise, and its
contribution to $\Phi(\psi_{\lm,\xi,\ga})$ is zero for every center $\xi$.
It follows that
\begin{\equ}\label{rev-positive-convolution}
 \Phi(\psi_{\lambda,\xi,\ga})
 =\frac{m^{m-1}\lambda^m}{16}
 \int_{\R^m}
 \frac{a^2\eta(r^{-1}|x|)^2f(b^{-2}|x|^2)^2\,{\tr(H(x)^2)}}
      {(\lambda^2+|x-\xi|^2)^m}\,dx
\end{\equ}
and this expression is independent of $\ga$. 

To obtain the full translation Hessian for Example~\ref{rev-weyl-example}, we also need to differentiate the term involving $K_{\lm,\xi,\ga}$ in \eqref{widehat-Phi}. After rescaling the bubble to unit scale, the linearized operator is $L_{\psi_{1,0,\ga}}$ as was considered in \eqref{equ-w0} and the Weyl tensor enters the right-hand side through $h$. We carry out this step for a general radial polynomial of degree at most three. To shorten the notation, we set
\[
 z(x)=u_0(x)(1-x)\cdot_{\ig_{\R^m}}\ga,
 \qquad u_0(x)=\frac{m^{(m-1)/2}}{(1+|x|^2)^{m/2}},
 \qquad s=|x|^2.
\]
For a translation in the $\xi^\ell$ direction, fix the vector
$e_\ell=\partial/\partial x^\ell$ in the Euclidean frame above.
Using the Euclidean metric to raise one index of $H$, we obtain the
vector field given by its $\ell$-th column,
\begin{\equ}\label{V-ell}
 V_\ell(x)=H(x)e_\ell=\sum_{j=1}^m H_{j\ell}(x)e_j .
\end{\equ}
By \eqref{rev-table1-structure}, only the first four columns can be
nonzero. For $1\le\ell\le4$, define the spinor
$C_\ell=D_{\ig_{\R^m}}(V_\ell\cdot_{\ig_{\R^m}}\ga)$, where
$D_{\ig_{\R^m}}$ is the Euclidean Dirac operator.
Then
$\diver V_\ell=0$,
$\langle x,V_\ell\rangle=0$ and $\De V_\ell=0$.
{Since $\ga$ is constant in the Euclidean spin frame,
the Dirac operator differentiates only the coefficients $H_{j\ell}$ of $V_\ell$.}
Pairing the terms with indices $(i,j)$ and $(j,i)$ gives
\[
C_\ell
 =\sum_{i,j}(\pa_i H_{j\ell})\,
 e_i\cdot_{\ig_{\R^m}}e_j\cdot_{\ig_{\R^m}}\ga
=-(\diver V_\ell)\ga
 +\sum_{i<j}(\pa_i H_{j\ell}-\pa_j H_{i\ell})\,
 e_i\cdot_{\ig_{\R^m}}e_j\cdot_{\ig_{\R^m}}\ga .
\]
{Since $\diver V_\ell=0$, only the terms with $i<j$ remain.}

We now take the radial factor in \eqref{three-form-F} to be a polynomial
in $s=|x|^2$. Thus $h(x)=F(s)H(x)$ with
$F(s)=\sum_{j=0}^N F_js^j$ of degree $N\le3$. Differentiating
$F_h(\psi_{1,\xi,\ga})$ from \eqref{rev-source} with respect to
$\xi^\ell$ gives
\begin{\equ}\label{rev-translation-source}
 \mff_\ell
 =\left.\partial_{\xi^\ell}F_h(\psi_{1,\xi,\ga})\right|_{\xi=0}
 =-\frac{m u_0F(s)}{2(1+s)}
   \big(V_\ell\cdot_{\ig_{\R^m}}\ga+x\cdot_{\ig_{\R^m}} V_\ell\cdot_{\ig_{\R^m}}\ga\big).
\end{\equ}
The right-hand side of \eqref{rev-translation-source} involves only
$V_\ell\cdot_{\ig_{\R^m}}\ga$ and
$x\cdot_{\ig_{\R^m}}V_\ell\cdot_{\ig_{\R^m}}\ga$.
Applying the Euclidean Dirac operator to these spinors also produces
$C_\ell$ and $x\cdot_{\ig_{\R^m}}C_\ell$. We therefore seek a
particular solution of $L_{\psi_{1,0,\ga}}w_\ell=\mff_\ell$ in the form
\begin{\equ}\label{rev-explicit-inverse}
 w_\ell(x)=u_0(x)\bigl[\varsigma(s)V_\ell(x)\cdot_{\ig_{\R^m}}\ga
 +\varkappa(s)x\cdot_{\ig_{\R^m}}V_\ell(x)\cdot_{\ig_{\R^m}}\ga\bigr]
 +u_0(x)\varpi(s)(1-x)\cdot_{\ig_{\R^m}}C_\ell(x),
\end{\equ}
where $s=|x|^2$ and $\varsigma,\varkappa,\varpi$ are real functions of $s$. To determine them,
we use the identities
\[
 D_{\ig_{\R^m}} C_\ell=0,\qquad
 D_{\ig_{\R^m}}(x\cdot_{\ig_{\R^m}} C_\ell)=-(m+2)C_\ell,
\]
\[
 D_{\ig_{\R^m}}(x\cdot_{\ig_{\R^m}} V_\ell\cdot_{\ig_{\R^m}}\ga)
 =-(m+4)V_\ell\cdot_{\ig_{\R^m}}\ga-x\cdot_{\ig_{\R^m}} C_\ell.
\]
The first two spinors in \eqref{rev-explicit-inverse} are pointwise
orthogonal to $z$ for the real inner product since
$\langle x,V_\ell\rangle_{\ig_{\R^m}}=0$. The factor $1-x$ in the
last term gives
\[
 \real\bigl(z,(1-x)\cdot_{\ig_{\R^m}}C_\ell\bigr)_{\ig_{\R^m}}
 =u_0(1+s)\real(\ga,C_\ell)_{\ig_{\R^m}}=0,
\]
{since $C_\ell$ is a real linear combination of
$e_i\cdot_{\ig_{\R^m}}e_j\cdot_{\ig_{\R^m}}\ga$ with $i<j$, and each
Clifford product $e_i\cdot_{\ig_{\R^m}}e_j$ is skew-adjoint.} Thus the real rank-one
term in the linearized operator vanishes on $w_\ell$ for every unit
spinor $\ga$.

We first apply $L_{\psi_{1,0,\ga}}$ to the full expression for
$w_\ell$ in \eqref{rev-explicit-inverse}.
We then divide both sides of the equation by $u_0(x)>0$ and compare
coefficients in
\[
 \frac{\bigl(L_{\psi_{1,0,\ga}}w_\ell\bigr)(x)}{u_0(x)}
 =\frac{\mff_\ell(x)}{u_0(x)}.
\]
On the left, the coefficients of $C_\ell$ and
$x\cdot_{\ig_{\R^m}}C_\ell$ are
$\varsigma+2s\varpi'+2\varpi$ and $2\varpi'-\varkappa$, respectively, where primes denote
derivatives in $s$. To cancel these terms, set
\begin{\equ}\label{rev-radial-coefficient-relations}
 \varsigma(s)=-2s\varpi'(s)-2\varpi(s),\qquad \varkappa(s)=2\varpi'(s).
\end{\equ}
With this choice, $V_\ell\cdot_{\ig_{\R^m}}\ga$ and
$x\cdot_{\ig_{\R^m}}V_\ell\cdot_{\ig_{\R^m}}\ga$ have the same
coefficient on the left. Comparing this common coefficient with
the one on the right, given by \eqref{rev-translation-source}
after division by $u_0(x)$, gives
\[
 2\varsigma'(s)-\frac{m}{1+s}\bigl(\varsigma(s)+\varkappa(s)\bigr)
 =-\frac{m}{2(1+s)}F(s).
\]
Substituting the expressions for $\varsigma$ and $\varkappa$ from
\eqref{rev-radial-coefficient-relations} and multiplying by $1+s$ gives the scalar equation
\[
 -4s(1+s)\varpi''+\big((2m-8)s-(2m+8)\big)\varpi'+2m\varpi
 =-\frac m2 F(s).
\]
Writing $\varpi(s)=\sum_{j=0}^N\varpi_j s^j$ and comparing powers of $s$
determines $\varpi_N,\ldots,\varpi_0$ successively. Set $\varpi_{N+1}=0$ and use
\begin{\equ}\label{rev-inverse-recurrence}
 \varpi_j=\frac{(m+2j+4)\varpi_{j+1}-mF_j/[4(j+1)]}{m-2j},
 \qquad j=N,N-1,\ldots,0.
\end{\equ}
Since $0\le j\le N\le3$ and $m\ge18$, each denominator satisfies
$m-2j\ge m-6>0$.
This gives the particular solution \eqref{rev-explicit-inverse}.
To compare this solution with the one obtained by applying
the normal inverse, we first check its decay. By \eqref{rev-radial-coefficient-relations},
the polynomial multiplying $u_0(x)$ in \eqref{rev-explicit-inverse}
has degree at most $2N+2\le8$ in $x$, since $V_\ell$ is quadratic
and $C_\ell$ is linear. As $|x|\to\infty$, it follows that
\[
 |w_\ell(x)|_{\ig_{\R^m}}=O(|x|^{8-m}),\qquad
 |\nabla w_\ell(x)|_{\ig_{\R^m}}=O(|x|^{7-m}).
\]
For $m\ge18$ both bounds are square integrable at infinity, and
$w_\ell$ is smooth near the origin. Hence
$w_\ell\in H^1(\R^m,\mbs(\R^m))\subset\msd^{1/2}(\R^m,\mbs(\R^m))$.
To see why we may use it in the pairing with $\mff_\ell$, take any
$Z\in T_{\psi_{1,0,\ga}}\cm$. Its translate $Z_\xi(x)=Z(x-\xi)$
belongs to $T_{\psi_{1,\xi,\ga}}\cm$. Since $\Ga$ vanishes on $\cm$
by Lemma~\ref{lemma: critical-manifold}~(4), its derivative in the
direction $Z_\xi$ vanishes. In terms of \eqref{rev-source}, this gives
\[
 \real\int_{\R^m}
 \bigl(F_h(\psi_{1,\xi,\ga}),Z_\xi\bigr)_{\ig_{\R^m}}\,dx=0.
\]
For the polynomial $h$ used here the same identity follows by removing
the cut-off, using the cut-off estimates at the beginning of Appendix~\ref{app:cubic-translations}.
Differentiate this identity with respect to $\xi^\ell$ at $\xi=0$.
The term obtained by differentiating $Z_\xi$ vanishes because
$F_h(\psi_{1,0,\ga})=0$. Thus \eqref{rev-translation-source} gives
\[
 \real\int_{\R^m}(\mff_\ell,Z)_{\ig_{\R^m}}\,dx
 =\inp{|D_{\ig_{\R^m}}|^{-1}\mff_\ell}{Z}_{\msd}=0.
\]
By the nondegeneracy in Lemma~\ref{lemma: critical-manifold}~(2),
$T_{\psi_{1,0,\ga}}\cm=\ker L_{\psi_{1,0,\ga}}$ in
$\msd^{1/2}(\R^m,\mbs(\R^m))$. Thus adding any element of this
kernel to $w_\ell$ leaves its pairing with $\mff_\ell$ unchanged.

For a polynomial $\cp(s)$, let $[\cp]_j$ denote the coefficient of $s^j$
in $\cp(s)$. Set
\begin{\equ}\label{rev-L-moments}
 L_j=M_j\frac{m-4-2j}{2}, \quad \text{for } 0\le j\le6.
\end{\equ}
The expression \eqref{rev-inverse-recurrence} applies to any radial
polynomial $F$ of degree at most three. We now use it to compute the
translation Hessian for the cubic polynomial $f$ chosen in
\eqref{rev-cubic-profile}. For a centered bubble of scale $\lm>0$,
the change of variables $x=\lm y$ turns the radial factor $f(|x|^2)$
into $f(\lm^2|y|^2)$. Thus, writing $s=|y|^2$ and $t=\lm^2$, we take
$F(s)=f(ts)$ in \eqref{rev-inverse-recurrence}. Its coefficients are $F_j=t^j[f]_j$
for $j=0,1,2,3$. We denote the resulting polynomial in $s$ by $\varpi(s,t)$.
Define
\begin{\equ}\label{rev-active-hessian-polynomial}
 \ch_m(t)=
 \sum_{j=0}^6q_jt^jL_j\left(2+\frac{8j}{m+4}\right)
 +8\sum_{j=0}^6[f(ts)\varpi(s,t)]_jL_j.
\end{\equ}
{Since $f(0)=1$ and \eqref{rev-inverse-recurrence} gives
$\varpi(s,0)=-1/4$, the two sums contribute $2L_0$ and $-2L_0$
to the constant term, giving $\ch_m(0)=0$.}
The second sum in \eqref{rev-active-hessian-polynomial} is the contribution of the negative
inverse-Hessian correction. Indeed, integrating the terms in
\eqref{rev-table1-structure} by the same coordinate symmetry as in
\eqref{I12} gives, for $1\le\ell,k\le4$,
\[
 \int_{S^{m-1}}\real
 (V_\ell(\theta)\cdot_{\ig_{\R^m}}\ga,
 V_k(\theta)\cdot_{\ig_{\R^m}}\ga)_{\ig_{\R^m}}\,d\theta
 =\frac{18\om_{m-1}}{m(m+2)}\de_{\ell k},
\]
\[
 \int_{S^{m-1}}\real
 (V_\ell(\theta)\cdot_{\ig_{\R^m}}\ga,
 \theta\cdot_{\ig_{\R^m}}C_k(\theta))_{\ig_{\R^m}}\,d\theta
 =-\frac{54\om_{m-1}}{m(m+2)}\de_{\ell k}.
\]
The four-form terms have zero angular average.  These identities hold for
all unit $\ga$.  In \eqref{rev-explicit-inverse}, the scalar combination
which remains after angular integration is
$\varsigma+s\varkappa+6\varpi=4\varpi$. A beta integral therefore gives the normalized pairing
$-8\sum_j[f(ts)\varpi(s,t)]_jL_j$, explaining the sign in
\eqref{rev-active-hessian-polynomial}.

For translations in directions $e_5,\ldots,e_m$, the differentiated
source vanishes since $H e_\ell=0$ for $\ell\ge5$. Differentiating
\eqref{rev-positive-convolution} with the cut-off removed gives $J_m$ in
\eqref{rev-transverse-translation-polynomial}. Together with the
coefficient $\ch_m$ in \eqref{rev-active-hessian-polynomial}, this
gives \eqref{rev-full-translation-hessian}. The first translation
derivative is zero at the center for every scale and every unit
spinor in this example.

\subsection{The six polynomial inequalities}\label{app:cubic-signs}

\begin{proof}[Proof of Lemma~\ref{rev-cubic-signs}]
Dimension $20$ first enters the noncompactness construction when we ask
at which scale $\lm>0$ the centered reduced energy
$\widehat\Phi_\infty(\lm,0,\ga)$ is stationary with $\ga$ fixed.
By \eqref{rev-central-profile}, the dependence on $W$ and $\ga$ is
through the factor $\tfrac12|W|^2+P(W,\ga)$, which is independent of
$\lm$. When this factor is nonzero, the scale equation reduces to
$P_m'(t)=tQ_m(t)=0$ for $t>0$,
where $t=\lm^2$. For the cubic polynomial in
\eqref{rev-cubic-profile}, evaluating \eqref{rev-radial-polynomials}
at $t=1$ gives $Q_{19}(1)>0$ and $Q_{20}(1)<0$. With this change of
sign in mind, we write $n=m-20$ and expand in powers of $n$ to show
that the required signs hold in dimension $20$ and persist in every
higher dimension.

To make the calculation easier to follow, we introduce two numbers
depending only on the dimension $m$ by setting
\[
 \Xi_0=\prod_{j=1}^6(m-2j-4) \quad \text{and} \quad
 \Xi_1=(m-2)(m-4)(m-6)^2(m-8)(m-10)(m-12)(m-14).
\]
The seven constants $M_0,\ldots,M_6$ in \eqref{rev-moments} can be
viewed as rational functions of $m$. Multiplying them by $\Xi_0$
makes each of them a polynomial in $m$. We use $\Xi_1$ for the same
purpose with the coefficients of $\ch_m(t)/t$. Indeed, for $0\le k\le6$, we have that
\[
 \Xi_0M_k=
 \prod_{j=1}^k(m+2j+2)\prod_{j=k+1}^6(m-2j-4).
\]
In \eqref{rev-inverse-recurrence}, $(m-2)(m-4)(m-6)$ clears the
denominators of the coefficients of $\varpi(s,t)$, since the factor
$m$ in the numerator of $\varpi_1$ cancels the denominator $m$ at
$j=0$. The remaining factors of $\Xi_1$ clear the denominators of
$L_0,\ldots,L_6$ in \eqref{rev-L-moments}. For $j\ge1$, the numerator
of $L_j$ also contains the factor $m+4$ that occurs in the denominator
of $8j/(m+4)$ in \eqref{rev-active-hessian-polynomial}.
Multiplication by $\Xi_0$ or $\Xi_1$ preserves signs in the range of the
lemma, since every factor in $\Xi_0$ is positive for $m>16$ and
every factor in $\Xi_1$ is positive for $m>14$.
{We work with these polynomials for $m\geq18$, the range in
which we obtained $\widehat\Phi_\infty$ by letting the cut-off radius
tend to infinity in \eqref{Phi-infty}. We then restrict to $m\geq20$
to prove the six sign inequalities.}

We first bound the degree of $-\Xi_1\ch_m(t)/t$ as a polynomial in
$t$ for fixed $m$. Since $\Xi_1$ depends only on $m$, it suffices to
consider $\ch_m(t)/t$. Applying \eqref{rev-inverse-recurrence} successively
for $j=3,2,1,0$ expresses each coefficient of $\varpi(s,t)$ as a linear
combination of the coefficients $F_j=t^j[f]_j$ of $f(ts)$. Thus $\varpi(s,t)$
has degree at most $3$ in $t$, and \eqref{rev-active-hessian-polynomial}
gives degree at most $6$ for $\ch_m(t)$. {Since $\ch_m(0)=0$
by \eqref{rev-active-hessian-polynomial}, the quotient $\ch_m(t)/t$
is a polynomial of degree at most $5$ in $t$.}

To bound the degree of $-\Xi_1\ch_m(t)/t$ in $n=m-20$, fix $t>0$
and let $m\to\infty$. By backward
substitution in \eqref{rev-inverse-recurrence}, the coefficients of
$\varpi(s,t)$ remain bounded. Setting $s=1$ in the scalar equation
preceding \eqref{rev-inverse-recurrence} then gives
$\varpi(1,t)=-f(t)/4+O(m^{-1})$. Since $L_j=m/2+O(1)$, the two sums
in \eqref{rev-active-hessian-polynomial} give
\[
 \ch_m(t)=m f(t)^2+4m f(t)\varpi(1,t)+O(1)=O(1).
\]
As $\Xi_1$ has degree $8$ in $m$, the polynomial $-\Xi_1\ch_m(t)/t$
has degree at most $8$ in $m$, and hence in $n=m-20$.

It remains to prove the six inequalities. With $m=n+20$, consider
the six polynomials
\[
\begin{gathered}
 \Xi_0Q_m(17/20),\qquad -\Xi_0Q_m(1),\qquad -\Xi_0Q_m'(t),\\[0.5em]
 -\Xi_0T_m(t),\qquad -\Xi_1\ch_m(t)/t,\qquad
 -(m+4)\Xi_0\bigl(J_m(t)+P_m(t)/t^2\bigr).
\end{gathered}
\]
Their positivity is equivalent to the six inequalities in the order
stated, since $\Xi_0$, $\Xi_1$, $m+4$, and $t$ are positive.
\begin{revision}
Let $R(n,t)$ be any one of these six polynomials. Using the degree
bounds $K$ in $n$ and $N$ in $t$ from the preceding calculations,
write
\[
 R(n,t)=\sum_{k=0}^K n^k\sum_{\ell=0}^N a_{k\ell}t^\ell.
\]
For the first two polynomials, $t$ has already been fixed at an
endpoint, so $N=0$ and it suffices to check that the coefficients
$a_{k0}$ are positive. For the remaining four, the coefficient of
each $n^k$ is itself a polynomial in $t$. We use the Bernstein basis
to give a positive lower bound for each such coefficient throughout
$17/20\le t\le1$.

First put $t=17/20+3\tau/20$, so that $0\le\tau\le1$. For each
fixed $k$, the binomial theorem gives
\[
 \sum_{\ell=0}^N a_{k\ell}
       \left(\frac{17}{20}+\frac{3}{20}\tau\right)^\ell
 =\sum_{j=0}^N b_{kj}\tau^j,\qquad
 b_{kj}=\left(\frac{3}{20}\right)^j
        \sum_{\ell=j}^N\binom{\ell}{j}
        \left(\frac{17}{20}\right)^{\ell-j}a_{k\ell}.
\]
We now express the same polynomial in $\tau$ in the Bernstein basis
of degree $N$. Denote its coefficients in this basis by
$r_{k0},\ldots,r_{kN}$, so that
\[
 \sum_{j=0}^N b_{kj}\tau^j
 =\sum_{i=0}^N r_{ki}\binom{N}{i}\tau^i(1-\tau)^{N-i}.
\]
To compute $r_{ki}$, write each $\tau^j$ on the left as
$\tau^j(\tau+(1-\tau))^{N-j}$ and expand. The coefficient of
$\tau^i(1-\tau)^{N-i}$ is then
$\sum_{j=0}^i\binom{N-j}{i-j}b_{kj}$. Dividing by $\binom{N}{i}$ gives
\begin{\equ}\label{rev-bernstein-coefficients}
 r_{ki}=\sum_{j=0}^i
 \frac{\binom{N-j}{i-j}}{\binom{N}{i}}\,b_{kj},
 \qquad 0\le k\le K,\quad 0\le i\le N.
\end{\equ}
Thus the original coefficients $a_{k\ell}$ determine the $b_{kj}$
by substitution, and these in turn determine the $r_{ki}$ by
\eqref{rev-bernstein-coefficients}. When $N=0$, this gives $r_{k0}=a_{k0}$.

For the cubic $f$ in \eqref{rev-cubic-profile}, we obtain the
$a_{k\ell}$ from the coefficients $q_j$ of $f^2$ in
\eqref{rev-cubic-square} and the formulas for $Q_m$, $\ch_m$, and
$T_m$ in \eqref{rev-radial-polynomials},
\eqref{rev-active-hessian-polynomial}, and \eqref{rev-transverse-translation-polynomial}.
For example, take $R(n,t)=-(m+4)\Xi_0\bigl(J_m(t)+P_m(t)/t^2\bigr)$.
Substituting $J_m(t)=2tT_m(t)/(m+4)$ and collecting the terms in
$t^j$ gives
\[
 \sum_{k=0}^K a_{kj}n^k
 =-q_j\,\Xi_0 M_j\bigl(m+4+2j(m-4-2j)\bigr),\qquad 0\le j\le6.
\]
Substituting $m=n+20$ on the right and collecting the powers of $n$
gives the individual coefficients $a_{kj}$.
\end{revision}
The following Table~\ref{tab:cubic-sign-coefficients} pairs each inequality with its
polynomial $R$. The last column gives a rational lower bound valid for
every $r_{ki}$ in that row, with $0\le k\le K$ and $0\le i\le N$.
\begin{center}
\begin{minipage}{\linewidth}
\centering
\captionof{table}{Coefficient bounds for the six inequalities in
Lemma~\ref{rev-cubic-signs}.}\label{tab:cubic-sign-coefficients}
\small
\renewcommand{\arraystretch}{1.5}
\begin{tabular}{|l|l|c|c|c|}
\hline
Desired inequality & $R(n,t)$ & $K$ & $N$ & Lower bound for $r_{ki}$ \\
\hline
$Q_m(17/20)>0$ & $\Xi_0Q_m(17/20)$ & $6$ & $0$ & $1/10$ \\
\hline
$Q_m(1)<0$ & $-\Xi_0Q_m(1)$ & $6$ & $0$ & $1/1000$ \\
\hline
$Q_m'(t)<0$ & $-\Xi_0Q_m'(t)$ & $6$ & $5$ & $1/2$ \\
\hline
$T_m(t)<0$ & $-\Xi_0T_m(t)$ & $7$ & $5$ & $2/5$ \\
\hline
$\ch_m(t)<0$ & $-\Xi_1\ch_m(t)/t$ & $8$ & $5$ & $1$ \\
\hline
$J_m(t)+P_m(t)/t^2<0$ & $-(m+4)\Xi_0(J_m(t)+P_m(t)/t^2)$ & $7$ & $6$ & $1/2$ \\
\hline
\end{tabular}
\end{minipage}
\end{center}
For example, in the second row $N=0$, so there are just seven
coefficients, one for each power of $n$ from $0$ to $6$. Each
$\Xi_0M_j$ is monic of degree $6$ in $n$. Thus the coefficient of $n^6$
in $-\Xi_0Q_m(1)$ is $-\sum_{j=0}^6(j+2)q_j=577/405769$.
This is the smallest of the seven coefficients and exceeds $1/1000$.

All six lower bounds are positive. Since the Bernstein basis functions
are nonnegative and sum to one on $[0,1]$, the coefficient of each
$n^k$ is at least the lower bound in its row for every $\tau\in[0,1]$.
In particular the constant coefficient is positive, and $n=m-20\ge0$.
Hence every $R(n,t)$ in the table is positive on its stated domain,
which proves all six inequalities.
In the two lower dimensions covered by the cut-off limit
\eqref{Phi-infty}, direct substitution in \eqref{rev-radial-polynomials}
gives $Q_{18}(1)>0$ and $Q_{19}(1)>0$. Thus the endpoint value $Q_m(1)$
changes sign from $m=19$ to $m=20$, and $20$ is the smallest dimension
$m\ge18$ for which the cubic in \eqref{rev-cubic-profile}
satisfies all six inequalities.
\end{proof}

For reference, when $m=20$ the critical polynomial is
\[
 Q_{20}(t)=2-\frac{324}{49}t+\frac{3114}{343}t^2
 -\frac{340}{49}t^3+\frac{14418}{4459}t^4
 -\frac{576}{637}t^5+\frac{4352}{31213}t^6,
\]
and its zero in $(17/20,1)$ satisfies
$0.9426915<t_{20}<0.9426917$.

\subsection{Perturbations of the cubic coefficients}\label{app:cubic-perturbations}

In \eqref{rev-cubic-profile}, the coefficients of $t$, $t^2$, and $t^3$,
in that order, form the vector $v=(v_1,v_2,v_3)=(-9/14,9/91,-2/637)$.
To verify the coefficient range in Remark~\ref{rem:cubic-coefficient-range},
let $\nu=(\nu_1,\nu_2,\nu_3)\in\R^3$ and replace each $v_j$ by
$v_j+\nu_j$, keeping the constant term $f(0)=1$. The resulting polynomial is
\[
 f_\nu(t)=1+(v_1+\nu_1)t+(v_2+\nu_2)t^2+(v_3+\nu_3)t^3.
\]
The degree estimates in the proof of Lemma~\ref{rev-cubic-signs} use
only that $f$ has degree at most $3$. Thus, with $f_\nu$ in place of $f$,
each of the six polynomials $R(n,t)$ in Table~\ref{tab:cubic-sign-coefficients}
has the same upper bounds $K$ and $N$ for its degrees in $n$ and $t$,
respectively.

Each coefficient $r=r_{ki}$ in \eqref{rev-bernstein-coefficients} is a
quadratic polynomial in the coefficients of $t$, $t^2$, and $t^3$ in $f$.
Write $r(v+\nu)=r(v)+L_r(\nu)+B_r(\nu)$, where $L_r$ is linear and
$B_r$ is homogeneous of degree~$2$. Let $\nu_{\max}=\max_{1\le j\le3}|\nu_j|$,
and let $\|L_r\|_1$ and $\|B_r\|_1$ denote the sums of the absolute
values of the respective monomial coefficients. Since $r(v)>0$ by
Table~\ref{tab:cubic-sign-coefficients}, the triangle inequality gives
\[
 \frac{|r(v+\nu)-r(v)|}{r(v)}
 \le \frac{\|L_r\|_1}{r(v)}\nu_{\max}
     +\frac{\|B_r\|_1}{r(v)}\nu_{\max}^2.
\]
To compute these ratios, we keep $\nu_1,\nu_2,\nu_3$ as variables.
For each $r_{ki}$ in the six rows, we use the coefficients of $f_\nu^2$
and the corresponding $\varpi$ in the formula for $R$, and apply
\eqref{rev-bernstein-coefficients}. Collecting the terms of degrees
$1$ and $2$ in $\nu$ gives $L_r$ and $B_r$, respectively. Taking the
largest of the resulting ratios in each case gives
\[
 \max_r\frac{\|L_r\|_1}{r(v)}=\frac{3304119}{577}<6000,\qquad
 \max_r\frac{\|B_r\|_1}{r(v)}=\frac{278995171}{57}<5\,000\,000 .
\]
Both maxima occur in the second row, which represents $-\Xi_0Q_m(1)$.
For the linear ratio the maximum occurs at the coefficient $r_{6,0}$ of
$n^6$. Since each $\Xi_0M_j$ is monic of degree $6$ in $n$,
\[
 r_{6,0}(v+\nu)=-2f_\nu(1)\bigl(f_\nu(1)+f_\nu'(1)\bigr).
\]
The two factors in this product are $577/1274+\sum_{j=1}^3\nu_j$
and $-1/637+\sum_{j=1}^3(j+1)\nu_j$, respectively. The constant term
is therefore $r_{6,0}(v)=577/405769$, and the sum of the absolute
values of the three linear coefficients is
$2\sum_{j=1}^3((j+1)577/1274-1/637)=57/7$. Their ratio is
$(57/7)/(577/405769)=3304119/577$.

For the quadratic ratio the maximum occurs at the constant coefficient
$r_{0,0}$ in $n$, obtained by setting $m=20$. By
\eqref{rev-radial-polynomials}, the quadratic part in $\nu$ is then
\[
 -\Xi_0\sum_{i,j=1}^3(i+j+2)M_{i+j}\nu_i\nu_j
 \qquad (m=20).
\]
All its monomial coefficients are negative. Summing their absolute
values gives $5766359040$, and $r_{0,0}(v)=5253120/4459$.
The ratio for this coefficient thus reduces to $278995171/57$.

If $\nu_{\max}\le 1/7000$, the relative change in every coefficient
is at most
\[
 \frac{6000}{7000}+\frac{5\,000\,000}{7000^2}=\frac{47}{49}.
\]
Hence $r(v+\nu)\ge (2/49)r(v)>0$ for every coefficient in every row.
As in the proof of Lemma~\ref{rev-cubic-signs}, the Bernstein basis
functions are nonnegative and sum to one on $[0,1]$. The coefficient
of each $n^k$ therefore stays positive throughout the interval.
Since $n=m-20\ge0$, each of the six polynomials remains positive
on its stated domain, proving the six strict inequalities for $f_\nu$.

The endpoint inequality $Q_m(1)<0$ limits how much we can vary the
three coefficients while retaining all six inequalities in every
dimension $m\ge20$. Since each $M_k$ in \eqref{rev-moments}
tends to one as $m\to\infty$, \eqref{rev-radial-polynomials} gives
\[
 \lim_{m\to\infty}Q_m(1)=2f(1)\bigl(f(1)+f'(1)\bigr).
\]
For $f$ in \eqref{rev-cubic-profile}, $f(1)=577/1274$ and
$f(1)+f'(1)=-1/637$. Taking $\nu_1=\nu_2=\nu_3=\nu_{\max}$
increases the latter quantity by $9\nu_{\max}$. Hence for
$\nu_{\max}>1/5733$ the perturbed polynomial fails this endpoint test in
all sufficiently large dimensions. This upper bound concerns the
test at the fixed relative scale $t=1$. Of course, a local maximum may persist
at a different scale.

\subsection{Stationarity at the centered extrema}\label{app:cubic-centered}

\begin{proof}[Proof of Proposition~\ref{rev-cubic-max}]
By Lemma~\ref{rev-cubic-signs}, $Q_m$ decreases strictly from a positive
value at $17/20$ to a negative value at $1$. The coefficient estimate in Appendix~\ref{app:cubic-perturbations} gives the same inequalities throughout the range
in Remark~\ref{rem:cubic-coefficient-range}. Since $P_m'(t)=tQ_m(t)$,
the unique zero $t_m$ of $Q_m$ in the interval $(17/20,1)$ satisfies
$P_m''(t_m)=t_mQ_m'(t_m)<0$. By \eqref{rev-central-profile0},
$P_m(t)>0$ for $t>0$. Thus $\lm\mapsto P_m(\lm^2)$ has a positive,
nondegenerate local maximum at $\sqrt{t_m}$.

{For fixed $W$, we regard $\mbs_m$
as a real Euclidean space with inner product $\real(\cdot,\cdot)_{\ig_{\R^m}}$.
On its unit sphere $S(\mbs_m)$ we integrate with respect to the induced
surface measure $d\ga$, normalized to have total mass $1$.} Each term in
$P(W,\ga)$ is a multiple of $\real(A\ga,\ga)_{\ig_{\R^m}}$, where
$A$ is the Clifford product of four distinct orthonormal basis vectors.
For any one of these vectors $e_i$, the isometry
$\ga\mapsto e_i\cdot_{\ig_{\R^m}}\ga$ preserves this measure and reverses
the sign of this term, since Clifford multiplication by $e_i$
anticommutes with $A$. Thus each term has zero integral, and
\[
 \frac12|W|^2+P(W,\ga_+)
 \ge \int_{S(\mbs_m)}\Big(\frac12|W|^2+P(W,\ga)\Big)\,d\ga
 =\frac12|W|^2.
\]
For every $W\ne0$ the maximum at $\ga_+$ is therefore positive.
Together with the positive local maximum of $P_m(\lm^2)$ at
$\lm=\sqrt{t_m}$, formula \eqref{rev-central-profile} gives a positive
local maximum of the centered energy at $(\sqrt{t_m},0,\ga_+)$.
The same formula gives a negative local minimum at $(\sqrt{t_m},0,\ga_-)$
provided that $\tfrac12|W|^2+P(W,\ga_-)<0$, and zero at every scale for any
spinor at which this factor vanishes.
For Example~\ref{Ex-2}, \eqref{nontrivial-Phi} and Fact~\ref{fact}
give the range $[-168,408]$. Fixing a spinor at which the value is
$-168$, we see by continuity in $W$ that the value remains negative
for all nearby algebraic Weyl tensors.

We now show that the derivative of $\widehat\Phi_\infty$ with respect
to the center $\xi$ vanishes at each of these extrema. The scalar term in
\eqref{Phi-on-M} is even in $\xi$. The correction in \eqref{widehat-Phi}
is quadratic in $\nabla\Ga(\psi_{\lm,\xi,\ga})$, which vanishes at
$\xi=0$ by Lemma~\ref{lemma: critical-manifold}~(4). It remains to
differentiate the Clifford term in \eqref{Phi-on-M2}.

We first collect the coefficient of
$e_i\cdot_{\ig_{\R^m}}e_j\cdot_{\ig_{\R^m}}e_k$ for $i<j<k$.
Recall that $h_{jq}(x)=f(|x|^2)H_{jq}(x)$ in the uncut tensor.
When the partial derivative $\pa_i$ acts on the radial factor in $\sum_q(\pa_i h_{jq})h_{kq}$,
the resulting coefficient is symmetric in $j,k$, since it contains
$\sum_q H_{jq}H_{kq}$. The corresponding Clifford products have opposite
signs when $j$ and $k$ are interchanged, so these contributions cancel.
For the terms in which $\pa_i$ acts on $H_{jq}$, grouping the six
permutations in pairs gives
$f(|x|^2)^2(\La_W(x))_{ijk}$, where the three-form $\La_W$ has components
\[
 (\La_W)_{ijk}
 =\sum_{\mathrm{cyc}(i,j,k)}\sum_{q=1}^m
   (\pa_iH_{jq}-\pa_jH_{iq})H_{kq}.
\]
Here the cyclic sum is taken over $(i,j,k)$, $(j,k,i)$, and $(k,i,j)$.

We claim that $x^i(\La_W(x))_{ijk}=0$ for all $j,k$.
Indeed, each $H_{jq}$ is homogeneous
of degree $2$, and $x^iH_{iq}=0$ since $H$ is symmetric and $Hx=0$.
Euler's identity and differentiation of $x^iH_{iq}=0$ with respect to $x^j$ give
\[
 x^i\pa_iH_{jq}=2H_{jq},\qquad
 x^i\pa_jH_{iq}=-H_{jq}.
\]
Contracting the three terms of the cyclic sum with $x$ therefore gives
\[
 x^i(\La_W(x))_{ijk}=\sum_q\bigl(3H_{jq}H_{kq}-3H_{kq}H_{jq}\bigr)=0.
\]

After substituting the bubble \eqref{bubble} into the Clifford term,
the numerator contains the product
\[
 (\lm+x-\xi)\cdot_{\ig_{\R^m}}\La_W(x)\cdot_{\ig_{\R^m}}(\lm-x+\xi).
\]
At $\xi=0$, the radial contraction just computed gives
$x\cdot_{\ig_{\R^m}}\La_W(x)=-\La_W(x)\cdot_{\ig_{\R^m}}x$,
so the product becomes
\[
 (\lm^2-|x|^2)\La_W(x)+2\lm\,x\cdot_{\ig_{\R^m}}\La_W(x).
\]
Consider first the contribution obtained by differentiating the weight
$(\lm^2+|x-\xi|^2)^{-m}$ with respect to $\xi^\ell$. At $\xi=0$ its
derivative is $2m x^\ell(\lm^2+|x|^2)^{-m-1}$. Since this is odd in
$x$, only the first term of the centered numerator contributes to the integral.

To separate the radial and spherical variables, write $x=\rho\theta$,
where $\rho=|x|$ and $\theta\in S^{m-1}$. Each component of $\La_W$ is
homogeneous of degree $3$, so $x^\ell\La_W(x)=\rho^4\theta^\ell\La_W(\theta)$.
After taking out the factors independent of $x$ and using
$dx=\rho^{m-1}\,d\rho\,d\theta$, we obtain the factorization
\[
\begin{aligned}
 &\int_{\R^m}\frac{(\lm^2-|x|^2)f(|x|^2)^2}{(\lm^2+|x|^2)^{m+1}}\,
        x^\ell(\La_W(x))_{ijk}\,dx\\[0.3em]
 &\qquad=
 \left(\int_0^\infty
 \frac{(\lm^2-\rho^2)f(\rho^2)^2\rho^{m+3}}
      {(\lm^2+\rho^2)^{m+1}}\,d\rho\right)
 \int_{S^{m-1}}\theta^\ell(\La_W(\theta))_{ijk}\,d\theta .
\end{aligned}
\]
The radial integral is the same for all indices and contains all the
dependence on $\lm$ and $f$. It remains to compute the integral over
$S^{m-1}$, which we denote by $\msi_{\ell ijk}$.
By definition, $\msi_{\ell ijk}$ is alternating in $i,j,k$. To see
why it is alternating in $\ell$ as well, use
$H_{kq}(x)=W_{k\al q\bt}x^\al x^\bt$ and the first Bianchi identity,
which gives $\pa_iH_{jq}-\pa_jH_{iq}=3W_{ijpq}x^p$.
Substituting these two identities into the component formula for $\La_W$
at $x=\theta$ gives
\[
 \begin{aligned}
 \msi_{\ell ijk}
 &=\int_{S^{m-1}}\theta^\ell
   \sum_{\mathrm{cyc}(i,j,k)}\sum_{q=1}^m
   (\pa_iH_{jq}-\pa_jH_{iq})(\theta)\,H_{kq}(\theta)\,d\theta\\
 &=3\sum_{\mathrm{cyc}(i,j,k)}\sum_{q=1}^m W_{ijpq}W_{k\al q\bt}
   \int_{S^{m-1}}\theta^\ell\theta^p\theta^\al\theta^\bt\,d\theta .
 \end{aligned}
\]
The Weyl components are independent of $\theta$, so it remains to
integrate the product of the four coordinates. By the same spherical
symmetry used to evaluate \eqref{I12}, we obtain
\begin{\equ}\label{sphere-four-coordinate-integral}
 \int_{S^{m-1}}\theta^\ell\theta^p\theta^\al\theta^\bt\,d\theta
 =\frac{\om_{m-1}}{m(m+2)}
  \bigl(\de^{\ell p}\de^{\al\bt}+\de^{\ell\al}\de^{p\bt}
        +\de^{\ell\bt}\de^{p\al}\bigr).
\end{\equ}
Replacing the integral in the second line of the expansion for
$\msi_{\ell ijk}$ by the right-hand side of \eqref{sphere-four-coordinate-integral}
contracts each of the three products of Kronecker deltas with
$W_{ijpq}W_{k\al q\bt}$. The term $\de^{\ell p}\de^{\al\bt}$
contributes zero because $\sum_{\al=1}^m W_{k\al q\al}=0$. The terms
$\de^{\ell\al}\de^{p\bt}$ and $\de^{\ell\bt}\de^{p\al}$ select the
components $W_{k\ell qp}$ and $W_{kpq\ell}$, respectively. Thus
\[
 \msi_{\ell ijk}=\frac{3\om_{m-1}}{m(m+2)}
 \sum_{\mathrm{cyc}(i,j,k)}\sum_{p,q}
 W_{ijpq}\bigl(W_{k\ell qp}+W_{kpq\ell}\bigr).
\]
Since $W_{ijpq}$ is antisymmetric in $p,q$, only the antisymmetric part
$(W_{kpq\ell}-W_{kqp\ell})/2=-W_{k\ell pq}/2$ contributes from
the second summand. Expanding the cyclic sum yields
\[
 \msi_{\ell ijk}=\frac{9\om_{m-1}}{2m(m+2)}
 \sum_{p,q}\bigl(
 W_{\ell i pq}W_{jk pq}-W_{\ell j pq}W_{ik pq}
 +W_{\ell k pq}W_{ij pq}\bigr),
\]
which is alternating in all four indices.

For a vector $v$, write $v^\flat$ for its Euclidean dual one-form and
$\iota_v$ for insertion of $v$ into the first argument of a differential form.
Clifford multiplication of the three-form $\La_W(\theta)$ by $\theta$ gives
\[
 \theta\cdot_{\ig_{\R^m}}\La_W(\theta)
 =\theta^\flat\wedge\La_W(\theta)-\iota_\theta\La_W(\theta).
\]
The first term is a four-form. The second is a two-form whose $(j,k)$
component is $\theta^i(\La_W(\theta))_{ijk}$, which vanishes by the
identity $x^i(\La_W(x))_{ijk}=0$ at $x=\theta$. Integrating the remaining
four-form over the unit sphere, we obtain
\[
 \msf_W=\int_{S^{m-1}}\theta^\flat\wedge\La_W(\theta)\,d\theta .
\]
Its components are
\[
 (\msf_W)_{\ell ijk}
 =\msi_{\ell ijk}-\msi_{i\ell jk}+\msi_{j\ell ik}-\msi_{k\ell ij}
 =4\msi_{\ell ijk}.
\]
Thus for any translation direction $v$ we obtain
\begin{\equ}\label{rev-weyl-angular-contraction}
 \int_{S^{m-1}}\langle\theta,v\rangle\La_W(\theta)\,d\theta
 =\frac14\iota_v\msf_W .
\end{\equ}

In the centered energy itself, the weight is even in $x$, so only the
term $2\lm\,x\cdot_{\ig_{\R^m}}\La_W(x)$ of the numerator contributes. Its integral
over the sphere at radius $\rho$ gives $\msf_W$ times a radial factor.
Thus the spinor-dependent part of the centered energy is a nonzero
scalar multiple of $\real(\msf_W\cdot_{\ig_{\R^m}}\ga,\ga)_{\ig_{\R^m}}$.
In particular, $\ga_+$ and $\ga_-$ are eigenspinors of the
self-adjoint Clifford action of $\msf_W$.

We can now differentiate in any center direction $v$ at $\xi=0$.
In the differentiated numerator, the part even in $x$ is
\[
 v\cdot_{\ig_{\R^m}}\La_W(x)\cdot_{\ig_{\R^m}}x
 +x\cdot_{\ig_{\R^m}}\La_W(x)\cdot_{\ig_{\R^m}}v .
\]
Using $\La_W(x)\cdot_{\ig_{\R^m}}x=-x\cdot_{\ig_{\R^m}}\La_W(x)$
in the first term and integrating over the unit sphere gives
$\msf_W\cdot_{\ig_{\R^m}}v-v\cdot_{\ig_{\R^m}}\msf_W$.
For the derivative of the weight, \eqref{rev-weyl-angular-contraction}
gives the spherical factor $\iota_v\msf_W/4$. Since $\msf_W$ is a four-form,
\[
 \msf_W\cdot_{\ig_{\R^m}}v-v\cdot_{\ig_{\R^m}}\msf_W=2\iota_v\msf_W.
\]
Thus both contributions are multiples of the difference
\[
 \real\bigl(\msf_W\cdot_{\ig_{\R^m}}(v\cdot_{\ig_{\R^m}}\ga),\ga\bigr)_{\ig_{\R^m}}
 -
 \real\bigl(v\cdot_{\ig_{\R^m}}(\msf_W\cdot_{\ig_{\R^m}}\ga),\ga\bigr)_{\ig_{\R^m}}.
\]
If $\msf_W\cdot_{\ig_{\R^m}}\ga=\beta\ga$, self-adjointness of the action
of $\msf_W$ allows us to move it to the second slot in the first pairing.
Both pairings then equal $\beta\real(v\cdot_{\ig_{\R^m}}\ga,\ga)_{\ig_{\R^m}}$.
Hence the translation derivative vanishes at both centered spinor extrema.

Finally, in the coordinates of \eqref{ab-a1} the radial cut-off
multiplies $f(\rho^2)$ by $\eta(\rho/R)$, where $R=r/b$.
In the centered energy \eqref{rev-central-profile0} this inserts
$\eta(\rho/R)^2$ in the radial integral. The dependence on $\ga$
and the preceding calculation of the translation derivatives are
therefore the same for every $R$. In particular, each extremizing
spinor still gives zero spinor and translation derivatives at the center.

On a compact scale interval around $\sqrt{t_m}$, the radial integrand
in \eqref{rev-central-profile0} and its first two scale derivatives are
$O(\rho^{15-m})$ at infinity, uniformly in the scale. Since $m\ge20$,
this gives $C^2$ convergence of the cut-off radial function to
$\kappa_m P_m(\lm^2)$ on that interval. The limit has a nondegenerate
maximum at $\sqrt{t_m}$. On a sufficiently small interval about this
scale its second derivative is negative, and its first derivative
has opposite signs at the two endpoints. These signs persist for
large $R$, giving a unique critical scale $\ell_R$ in that interval
with $\ell_R\to\sqrt{t_m}$. At this scale all first derivatives vanish
for the extremizing spinors.
Formula \eqref{ab-a1} takes the normalized scale $\ell_R$ to the scale
$b\ell_R$ of $\widehat\Phi_{a,b,r}$.
\end{proof}

\subsection{Joint variations of the center and the spinor}\label{app:cubic-joint}
We first derive the criterion \eqref{rev-weyl-sign-condition} for a
general algebraic Weyl tensor $W$. We then verify it on both
extremal eigenspaces for Example~\ref{Ex-2}.

To compute the first term from \eqref{widehat-Phi}, differentiate
$\Phi(\psi_{\ell_m,\xi,\ga})$ twice in direction $v$ at $\xi=0$ and subtract
$\inp{K_{\ell_m,0,\ga}Y_v}{Y_v}_{\msd}$, where $Y_v$ is the
translation derivative of the Hilbert representative of the source
\eqref{rev-source} at $\xi=0$ for the uncut tensor. All other terms
from differentiating the normal correction vanish because the centered
source is zero. For any algebraic $W$, the identities
$\diver(Hv)=0$, $\De(Hv)=0$ and $\langle x,Hv\rangle=0$
allow one to use \eqref{rev-explicit-inverse}--\eqref{rev-inverse-recurrence}
with $V=Hv$ after rescaling to unit bubble scale.

The second term in \eqref{rev-joint-translation-test} accounts for the
variation of the spinor. At scale $\ell_m$ the translation derivative
as a quadratic form on unit spinors is represented by the endomorphism
\[
 \ga\ \mapsto\ \frac{G_m}{2\ell_m}\Bigl(
 \msa_W(v\cdot_{\ig_{\R^m}}\ga)
 -v\cdot_{\ig_{\R^m}}(\msa_W\ga)\Bigr).
\]
Indeed, in the calculation following \eqref{rev-weyl-angular-contraction}
the centered numerator contributes $2\ell_m\,x\cdot_{\ig_{\R^m}}\La_W(x)$,
and the even part of its translation derivative is
\[
 (x\cdot_{\ig_{\R^m}}\La_W(x))\cdot_{\ig_{\R^m}}v
 -v\cdot_{\ig_{\R^m}}(x\cdot_{\ig_{\R^m}}\La_W(x)).
\]
The term from the denominator is proportional to
\[
 \int_0^\infty
 \frac{(\ell_m^2-r^2)r^{m+3} f(r^2)^2}{(\ell_m^2+r^2)^{m+1}}\,dr=0
\]
by the scale equation $P_m'(t_m)=0$. Since $\msa_W\ga=\bt_*\ga$,
the mixed derivative in directions $v$ and $\eta$ is
$\dfrac{G_m}{\ell_m}\real\bigl((\msa_W-\bt_*\operatorname{Id})
(v\cdot_{\ig_{\R^m}}\ga),\eta\bigr)_{\ig_{\R^m}}$
for every tangent spinor $\eta$ to the unit spinor sphere at $\ga$.
Completing the square in $\eta$ gives the joint Hessian
\[
 \ch_{W,\ga}(v)
 +2G_m\real\left((\msa_W-\bt_*\operatorname{Id})
   \left(\eta+\frac{v\cdot_{\ig_{\R^m}}\ga}{2\ell_m}\right),
   \eta+\frac{v\cdot_{\ig_{\R^m}}\ga}{2\ell_m}\right)_{\ig_{\R^m}}.
\]
The scale has no mixed term with either variable. This follows from
\eqref{rev-central-profile} and from translation stationarity at every
scale for $\ga\in S(E_{\bt_*})$.

The square in the spinor variation has the sign of the spinor Hessian,
which is strict normal to $S(E_{\bt_*})$ in the unit spinor sphere.
Under \eqref{rev-weyl-sign-condition}, the remaining translation term
has the same strict sign for every $v\ne0$. Together with the scale
Hessian, this gives \eqref{rev-uniform-local-gap}.

Let $E=\operatorname{span}\{e_1,e_2,e_3,e_4\}$ and let $A$ denote Clifford multiplication by the oriented volume form of $E$. Fix a unit spinor
$\ga$ satisfying $A\ga=\sigma\ga$, with $\sigma\in\{-1,1\}$, and set
$p=P(W,\ga)$. For this Weyl tensor, \eqref{nontrivial-Phi} gives
$p=-288\sigma$ and $\tfrac12|W|^2=120$.

With $t=\lm^2$, the claimed translation Hessian is
\eqref{table-two-translation-hessian}, where $v_E$ and $v_{E^\perp}$
are the orthogonal components of $v$.
For comparison, in Example~\ref{rev-weyl-example} we have $p=0$ and
$\tfrac12|W|^2=72$. The corresponding coefficients of $|v_E|^2$ and
$|v_{E^\perp}|^2$ are $72\ch_m(t)$ and $72J_m(t)$, respectively,
by \eqref{rev-full-translation-hessian}.

Assuming \eqref{table-two-translation-hessian} for the moment, we can
read off its sign at $t=t_m$. Lemma~\ref{rev-cubic-signs} gives $\ch_m(t_m)<0$ and $J_m(t_m)<0$,
since $J_m(t_m)=2t_mT_m(t_m)/(m+4)$. On the positive branch $p=288$,
the coefficients of $|v_E|^2$ and $|v_{E^\perp}|^2$ in
\eqref{table-two-translation-hessian} are
$264\ch_m(t_m)+288J_m(t_m)$ and $408J_m(t_m)$, respectively, so both
are negative. On the negative branch $p=-288$, the same coefficients
are $-24\ch_m(t_m)-288J_m(t_m)$ and $-168J_m(t_m)$, so both are positive.
The mixed scale and translation
derivatives vanish since the translation derivative vanishes at the
center for every scale at these spinors.

We now prove \eqref{table-two-translation-hessian}. In \eqref{widehat-Phi}
we must differentiate both $\Phi$ and the correction involving the inverse
of the linearized operator. The former involves $\tr(H^2)$ and the
three-form $\La_W$. For the latter we use \eqref{rev-translation-source}
and \eqref{rev-explicit-inverse} with the columns $V_i=He_i$ and spinors
$C_i=D_{\ig_{\R^m}}(V_i\cdot_{\ig_{\R^m}}\ga)$ of Example~\ref{Ex-2}.

After integrating in the radial variable, it remains to integrate
polynomials of degrees $4$ and $6$ in $\theta=(\theta^1,\ldots,\theta^m)$
over $S^{m-1}$. Formula \eqref{sphere-four-coordinate-integral} gives
the integrals of degree $4$. For degree $6$, reflection in a coordinate
hyperplane gives zero whenever that coordinate has odd exponent.
For distinct $i,j,k$, the remaining possibilities are
\[
 \frac{1}{\om_{m-1}}\int_{S^{m-1}}
 \Big((\theta^i)^6,\,(\theta^i)^4(\theta^j)^2,\,
                    (\theta^i)^2(\theta^j)^2(\theta^k)^2\Big)\,d\theta
 =\frac{(15,3,1)}{m(m+2)(m+4)}.
\]
Indeed, for a homogeneous polynomial $\cp$ of degree $5$ and a fixed
$i$, the divergence theorem on the unit ball and integration in polar
coordinates give
\[
 \int_{S^{m-1}}\theta^i\cp(\theta)\,d\theta
 =\int_{B_1(0)}\pa_i\cp(x)\,dx
 =\frac{1}{m+4}\int_{S^{m-1}}(\pa_i\cp)(\theta)\,d\theta .
\]
Removing one factor $\theta^i$ from each of the three monomials and
applying this identity reduces all three integrals to
\eqref{sphere-four-coordinate-integral}.

We apply these formulas to the entries of Table~\ref{tab:weyl-indefinite-spinor-term}.
Since $\ga$ is a unit spinor, the real inner product of
$V_i\cdot_{\ig_{\R^m}}\ga$ and $V_j\cdot_{\ig_{\R^m}}\ga$ equals
$\langle V_i,V_j\rangle$. Substituting
$H_{ri}(\theta)=W_{r\al i\bt}\theta^\al\theta^\bt$ into this inner
product and using \eqref{sphere-four-coordinate-integral} gives, for
$1\le i,j\le4$,
\[
 \frac{m(m+2)}{\om_{m-1}}\int_{S^{m-1}}
             \langle V_i(\theta),V_j(\theta)\rangle\,d\theta
 =\sum_{r,\al,\bt=1}^4
    W_{r\al i\bt}\bigl(W_{r\al j\bt}+W_{r\bt j\al}\bigr)
 =30\de_{ij}.
\]
Here the term involving the trace of $W$ vanishes. Summing over $i=j$
and using $\tr(H^2)=\sum_{i=1}^4|V_i|^2$ gives
\[
 \int_{S^{m-1}}\tr(H(\theta)^2)\,d\theta
 =\frac{120\om_{m-1}}{m(m+2)} .
\]
Inserting $\theta^i\theta^j$ in the same integral and using the
degree-$6$ formula gives
\[
 \frac{m(m+2)(m+4)}{\om_{m-1}}\int_{S^{m-1}}
 \theta^i\theta^j\tr(H(\theta)^2)\,d\theta
 =\begin{cases}
    240 & i=j\le4,\\
    120 & i=j>4,\\
    0   & i\ne j .
   \end{cases}
\]
For instance, when $i=j=1$, the coefficients of the even monomials
of the three types displayed above sum to $4$, $52$ and $24$,
respectively, giving $15\cdot4+3\cdot52+24=240$. When $i=j>4$,
$\tr(H^2)$ is independent of $\theta^i$, so the extra square leaves
the numerator $120$ unchanged.

It remains to compute the integral involving $C_j$. Since
$\langle\theta,V_i(\theta)\rangle=0$, moving Clifford multiplication by
$V_i(\theta)$ to the second slot gives the Clifford product of $\theta$,
$V_i(\theta)$ and the two-form acting on $\ga$ in $C_j(\theta)$,
in that order. Its scalar part is $-3\langle V_i(\theta),V_j(\theta)\rangle$
because $\theta^k(\pa_kH_{\ell j}-\pa_\ell H_{kj})(\theta)=3H_{\ell j}(\theta)$.
Its two-form part has zero real inner product with $\ga$ after acting
on $\ga$. For the four-form part, any coordinate outside $E$ integrates
to zero by reflection. Write $\epsilon_{rsk\ell}$ for the sign of the
permutation $(r,s,k,\ell)$ of $(1,2,3,4)$, and set it equal to zero for
repeated indices. Using
$(\pa_kH_{\ell j}-\pa_\ell H_{kj})(\theta)=3W_{k\ell\al j}\theta^\al$
and \eqref{sphere-four-coordinate-integral}, the coefficient of
$A$ in the integrated Clifford action, after multiplication by
$m(m+2)/\om_{m-1}$, is
\[
 3\sum_{r,s,\al=1}^4\ \sum_{1\le k<\ell\le4}
 \epsilon_{rsk\ell}\bigl(W_{sr i\al}+W_{s\al ir}\bigr)W_{k\ell\al j}
 =72\de_{ij}.
\]
Since $A\ga=\sigma\ga$ and $p=-288\sigma$, the scalar and four-form
parts together give, for $1\le i,j\le4$,
\[
 \frac{m(m+2)}{\om_{m-1}}\int_{S^{m-1}}
 \real\bigl(V_i(\theta)\cdot_{\ig_{\R^m}}\ga,\,
           \theta\cdot_{\ig_{\R^m}}C_j(\theta)\bigr)_{\ig_{\R^m}}\,d\theta
 =(-3\cdot30+72\sigma)\de_{ij}
 =\Big(-90-\frac{p}{4}\Big)\de_{ij}.
\]
Both spinor integrals vanish if either index is greater than $4$,
since $V_i=C_i=0$ for $i>4$.

We first evaluate the contribution of $\Phi$. Substituting the same Weyl
components into the formula for $\La_W$ gives the pointwise identity
\[
 \La_W(x)=12|x_E|^2\iota_{x_E}\operatorname{vol}_E,
\]
where $\operatorname{vol}_E$ is the oriented volume form of $E$.
It follows that the part of $\Phi(\psi_{\lm,\xi,\ga})$ even in $\xi$
has the same integral as in \eqref{rev-positive-convolution}, with
$a=b=1$, the cut-off removed, and the factor ${\tr(H(x)^2)}$ replaced by
\[
 \tr(H(x)^2)+\frac p{24}|x_E|^2
       \bigl(|x_E|^2-\langle x_E,\xi_E\rangle\bigr).
\]
Indeed, in the Clifford product for the three-form term only
$2\lm(x-\xi)^\flat\wedge\La_W$ contributes to the even part, and the
four-forms with exactly three indices in $E$ interchange the two
eigenspaces of $A$. Let $U_m(t)=\sum_{j=0}^6 q_jt^j L_j$.
Differentiating the kernel twice gives the contribution
$(120+p)(2U_m+2J_m)$ along each unit direction in $E$ in the
normalization of \eqref{table-two-translation-hessian}. Differentiating
the term linear in $\xi_E$ and the kernel once each contributes
$-pU_m$. Thus $\Phi$ contributes
\[
 (120+p)(2U_m+2J_m)-pU_m
 =\Big(120+\frac{p}{2}\Big)(2U_m+2J_m)+pJ_m .
\]
The additional factor $(1+s)^{-1}$ in the rescaled radial
integral replaces $M_j$ by $\frac{m-4-2j}{2m}M_j$.
For the inverse-Hessian term, the sphere integrals involving $V_i$ and $C_j$ change
the coefficient of $8\sum_j[f(ts)\varpi(s,t)]_jL_j$ in
\eqref{rev-active-hessian-polynomial} from $72$ to $120+p/2$, since
the combination in \eqref{rev-explicit-inverse} becomes
$30(\varsigma+s\varkappa+6\varpi)+(p/2)\varpi=(120+p/2)\varpi$.
By \eqref{rev-transverse-translation-polynomial},
$2U_m+2J_m$ is the first sum in \eqref{rev-active-hessian-polynomial}.
Adding the inverse-Hessian correction therefore gives the coefficient
$\bigl(120+p/2\bigr)\ch_m+pJ_m$ of $|v_E|^2$.
In $E^\perp$ the columns of $H$ vanish, so only the kernel derivatives
contribute, giving $(120+p)J_m$ as the coefficient of $|v_{E^\perp}|^2$.
The mixed entries between $E$ and $E^\perp$ vanish by parity.
This proves \eqref{table-two-translation-hessian}.

It remains to substitute the two eigenvalues into
\eqref{rev-joint-translation-test}. On $A\ga=\sigma\ga$ one has
$9\bt_*=120+p$ and $p=-288\sigma$. Clifford multiplication by $v_E$
interchanges the eigenspaces of $A$, while multiplication by
$v_{E^\perp}$ preserves them. Hence
\[
 \real\bigl((\msa_W-\bt_*\operatorname{Id})(v\cdot_{\ig_{\R^m}}\ga),
             v\cdot_{\ig_{\R^m}}\ga\bigr)_{\ig_{\R^m}}
 =-\frac{2p}{9}|v_E|^2 .
\]
By \eqref{table-two-translation-hessian}, the coefficients of
$|v_E|^2$ and $|v_{E^\perp}|^2$ in
$\ch_{W,\ga}(v)/(\kappa_m t_m)$ are therefore
\[
 (120+p/2)\ch_m(t_m)+p\big(J_m(t_m)+P_m(t_m)/t_m^2\big) \quad \text{and} \quad
 \qquad (120+p)J_m(t_m),
\]
respectively. In the coefficient of $|v_E|^2$, allowing the spinor to
vary adds $pP_m(t_m)/t_m^2$ to the pure translation coefficient in
\eqref{table-two-translation-hessian}. This is where we use the
stronger bound $J_m(t_m)+P_m(t_m)/t_m^2<0$ from
Lemma~\ref{rev-cubic-signs}. Together with $\ch_m(t_m)<0$ and
$J_m(t_m)<0$, it makes both coefficients negative for $p=288$ and
both positive for $p=-288$. Thus
\eqref{rev-uniform-local-gap} holds on disjoint neighborhoods of the
positive maximum and the negative minimum.

For each fixed $m$ and $f$ these conclusions persist for all algebraic
Weyl tensors sufficiently close to the respective example. Indeed,
$D^2_{\xi\xi}\widehat\Phi_\infty$ depends continuously on the unit
spinor and quadratically on $W$ by \eqref{widehat-Phi}. If $W_j\to W$,
the largest and smallest eigenvalues of $\msa_{W_j}$ converge to
the corresponding eigenvalues of $\msa_W$. Any sequence of
associated unit eigenspinors has a subsequence converging to a unit
eigenspinor for $W$. Applying \eqref{rev-joint-translation-test} along
this subsequence and using compactness of the unit translation sphere
shows that the strict sign in \eqref{rev-weyl-sign-condition} persists.
This argument also covers splitting of an extremal eigenvalue.

\end{appendices}

\section*{AI Declaration Statement}
We began this work in 2025. The ODE analysis on product manifolds in Section 4 and the Weyl tensor constructions in Examples 5.9 and 5.10 are solely due to the authors. We had completed this work by July 2025, using \textit{Mathematica} as our only computational tool to obtain the Weyl tensor data and verify the algebraic identities in the two examples. On August 2, 2025, the second author presented these examples at the 1st Yunnan Symposium on Nonlinear Analysis and Nonlinear Partial Differential Equations, an international conference held in Kunming, China.

In the subsequent development and revision of the manuscript, we used OpenAI’s ChatGPT to assist with further derivations, carry out and check symbolic calculations more efficiently, and improve the consistency of the notation and the clarity of the exposition. We take sole responsibility for all mathematical ideas, arguments, and results presented in this manuscript.

\bigskip

\noindent
\authoraddress{Juncheng Wei\\
    Department of Mathematics,\\
    The Chinese University of Hong Kong,\\
    Shatin, New Territories, Hong Kong}
    {wei@math.cuhk.edu.hk}%
\hfill
\authoraddress{Tian Xu\\
    Department of Mathematics,\\
    Zhejiang Normal University,\\
    Jinhua, Zhejiang 321004, China}
    {xutian@amss.ac.cn}\par
	
	\end{document}